**<u>Новоженова О.Г.</u>**

# БИОГРАФИЯ И НАУЧНЫЕ ТРУДЫ
# АЛЕКСЕЯ НИКИФОРОВИЧА ГЕРАСИМОВА

## О ЛИНЕЙНЫХ ОПЕРАТОРАХ,
## УПРУГО-ВЯЗКОСТИ, ЭЛЕВТЕРОЗЕ
## И
## ДРОБНЫХ ПРОИЗВОДНЫХ

По единственной статье, известной в мировой научной литературе (Прикладная Математика и Механика,1948, № 3), расшифрованы инициалы, восстановлена биография, найдены еще 11 сохранившихся работ советского механика Алексея Никифоровича Герасимова, на 20 лет раньше Капуто предложившего использование дробной производной для задач вязко-упругости (производная Герасимова-Капуто).

## БИОГРАФИЯ

Алексей Никифорович Герасимов родился 24 марта 1897 г. в станице Кочетовская Донской области, в семье московского врача Никифора Илларионовича и домохозяйки Веры Сергеевны Герасимовых.

В 1906 г. пошел в приготовительный класс московской гимназии № 10, с 7-го класса начал помогать отцу (родились еще два брата и сестра), давая частные уроки по математике и физике.

В 1916 г. поступил на математическое отделение физмата МГУ, где и проучился до мобилизации в 1919 г. по профсоюзной линии в РККа в качестве рядового.

Одновременно в 1917–1918 гг. работал преподавателем физики и математики и завучем школы I ступени № 21 в г. Москве.

В 1919–1922 гг. служил в рядах Красной Армии. Участвовал в гражданской войне на Юго-Восточном фронте в составе 56 дивизии. После войны состоял на воинском учёте в качестве запасного минёра второй очереди Военно-морского флота.

С 1922 г. по май 1930 г. работал учителем физики и математики в школах II ступени Замоскворецкого района г. Москвы. С 1930 г. по 1938 г. работал преподавателем физики и термодинамики на кафедре физики Московского текстильного института, возглавляемой И.И. Васильевым. В 1932 г. опубликовал статью «Теория рычажных весов с постоянной чувствительностью» в журнале «За реконструкцию текстильной промышленности».

Параллельно с работой экстерном учился в I Московском Государственном Университете на физико-математическом факультете, который окончил в 1936 г. по специальности «Прикладная математика», защитив диплом I степени на тему «Принцип соответствия в теории линейных операторов»(научный руководитель И.И. Привалов).

С 1937 г. по 1940 г. состоял заочником в аспирантуре НИИМ МГУ, сдавал экзамены из программы кандидатского минимума, работал над диссертацией под руководством И.М.Привалова. Опубликовал четыре работы в журнале «Прикладная математика и механика» АН СССР: «Проблема упругого последействия и внутреннее трение»./ ПММ,1938, I,493-536;Письмо в редакцию/ ПММ, ОТН,1938,т.II, в. 1, с.137;Основания теории деформаций упруго-вязких тел/ ПММ, ОТН, 1938, т. II, в.3, С.379-388;К вопросу о малых колебаниях упруго-вязких мембран/ ПММ,1939,т.II,в.4, 467-486.

По обстоятельствам квартирного характера с беременной женой Еленой Алексеевной Калиной и ее дочерью от первого брака Вероникой с 1938 по 1940 г был вынужден уехать в г. Йошкар-Ола. Здесь он работал в Марийском пединституте им. Н.К. Крупской зав. кафедрой математики. Читал курсы теоретической механики, астрономии, теорию аналитических функций, спец. курсы анализа, теорию вероятностей.

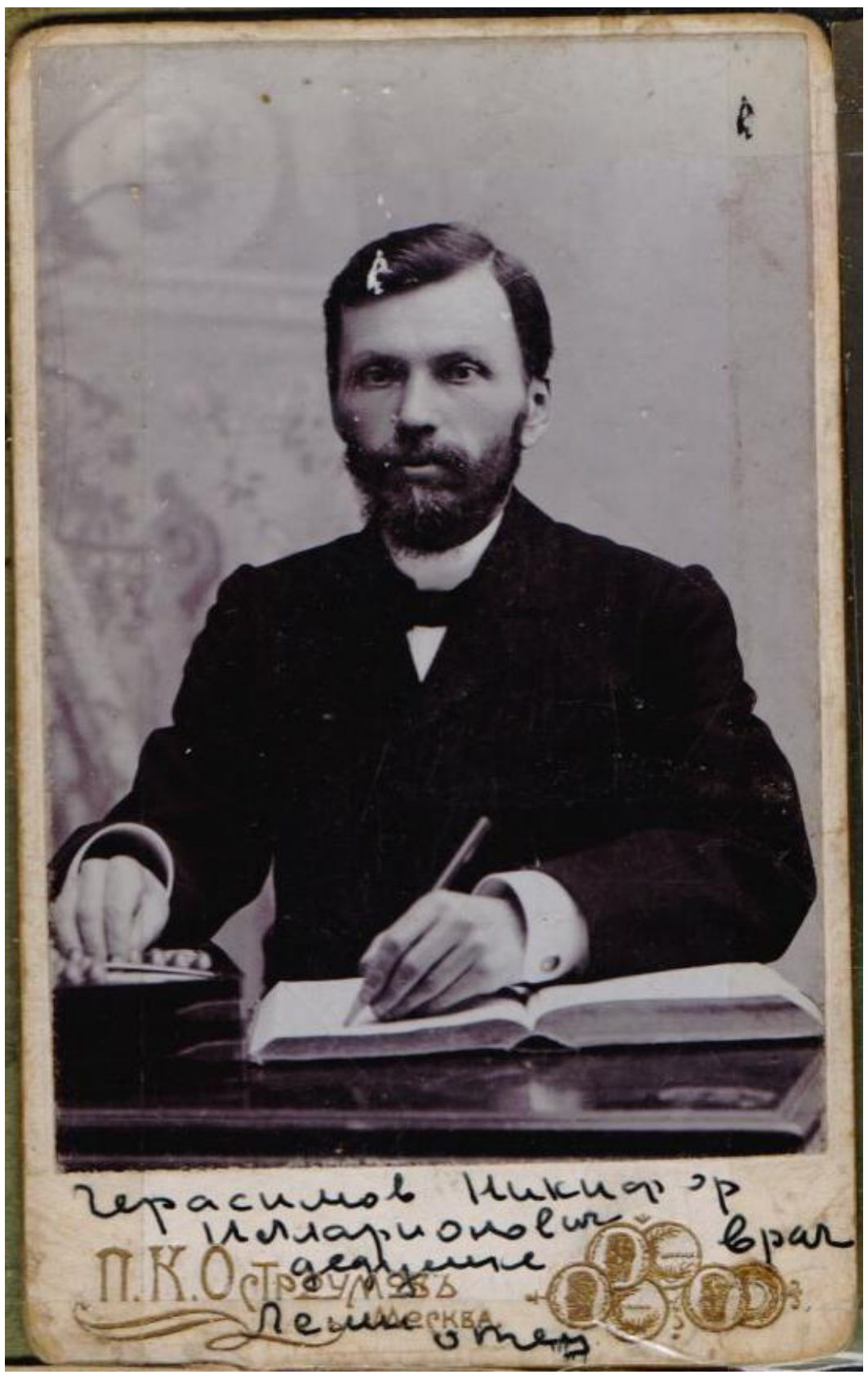

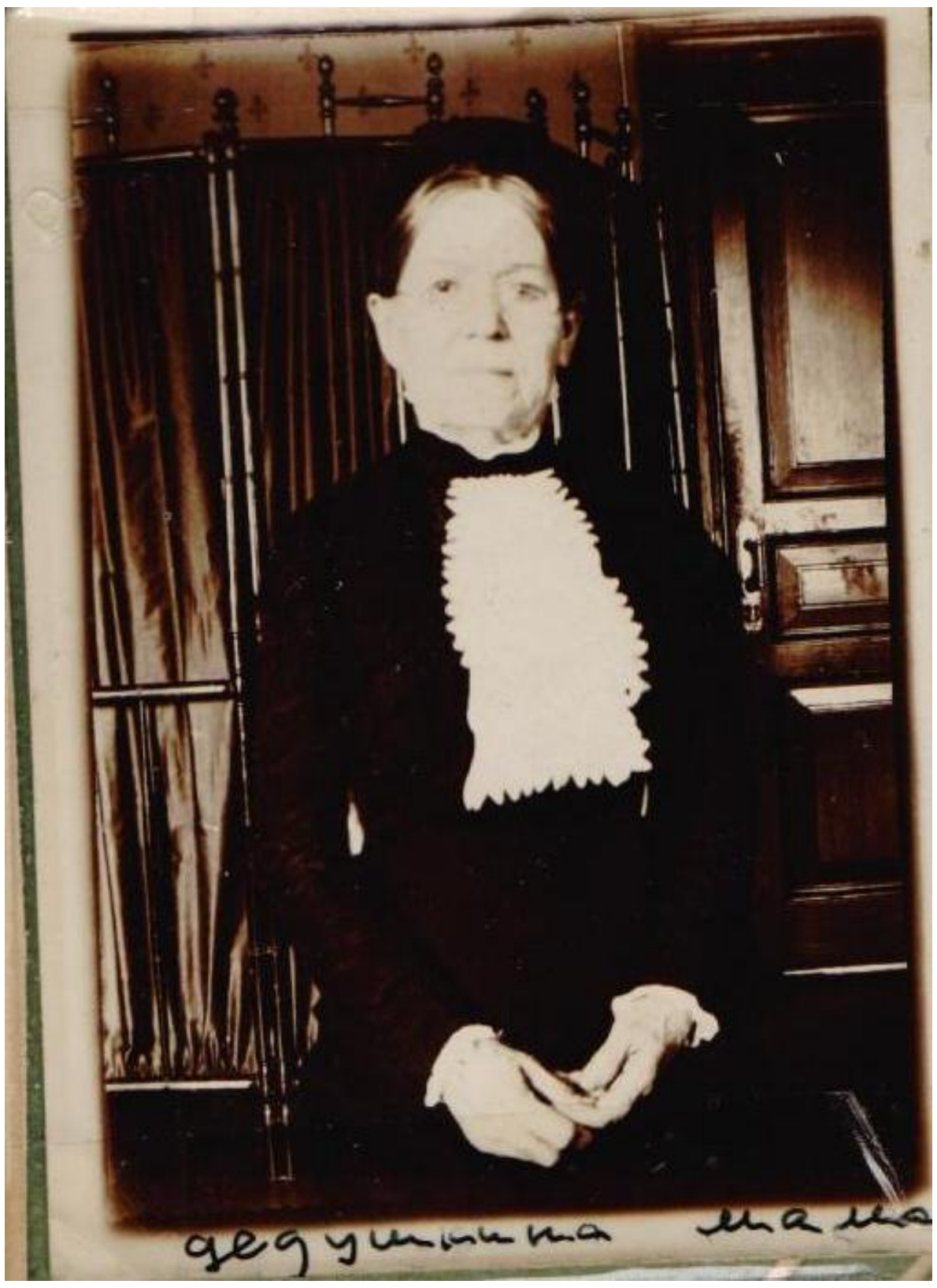

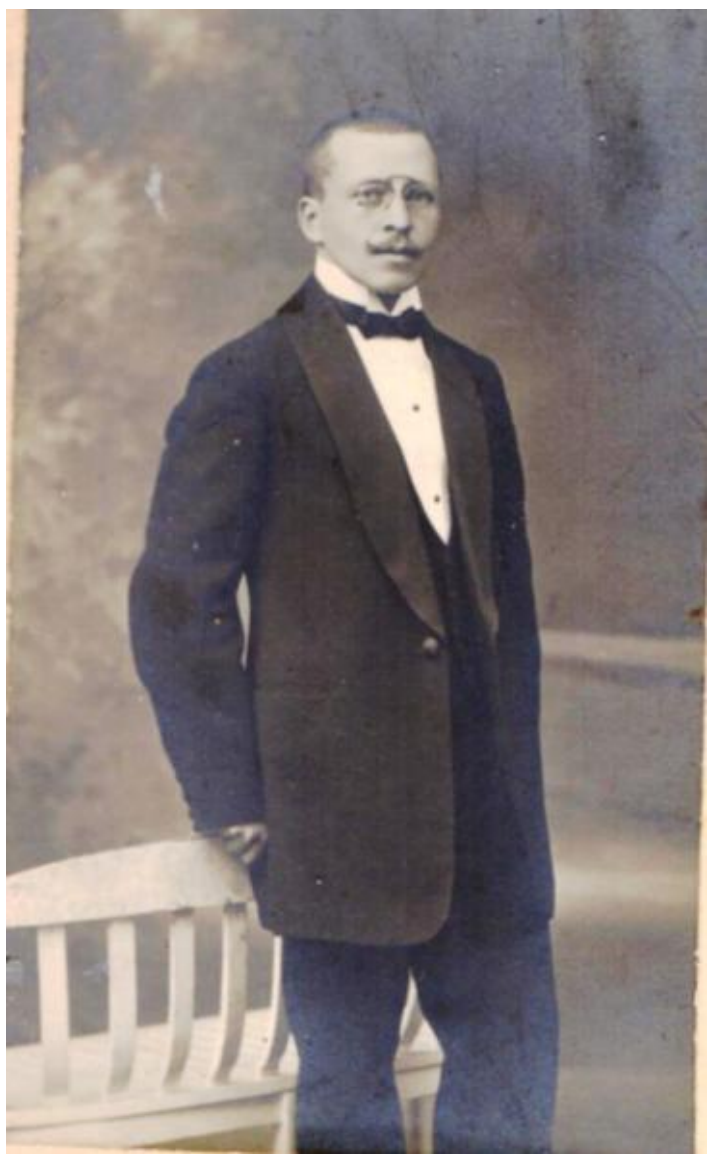

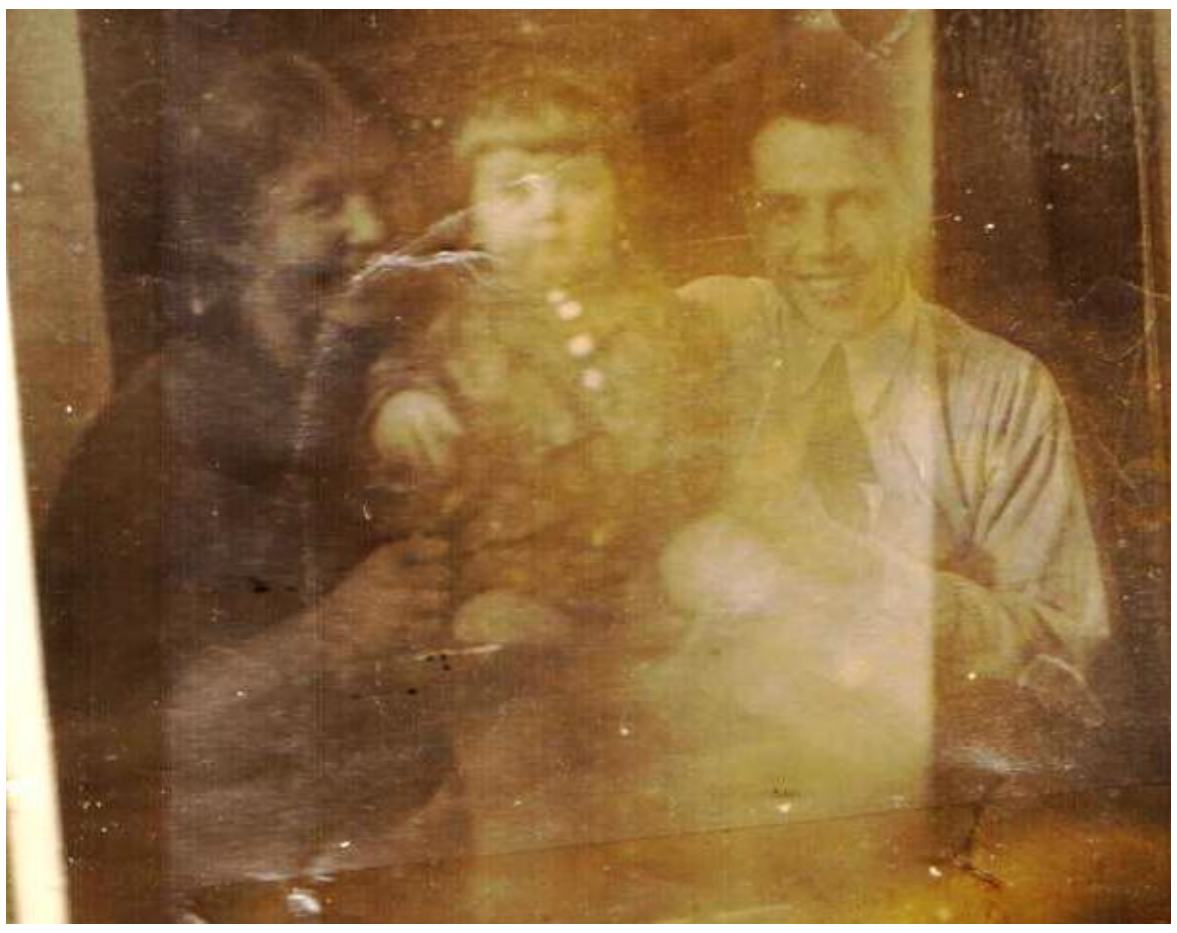

Отец жены - дворянин

г. Йошкар-Ола,1939 г.

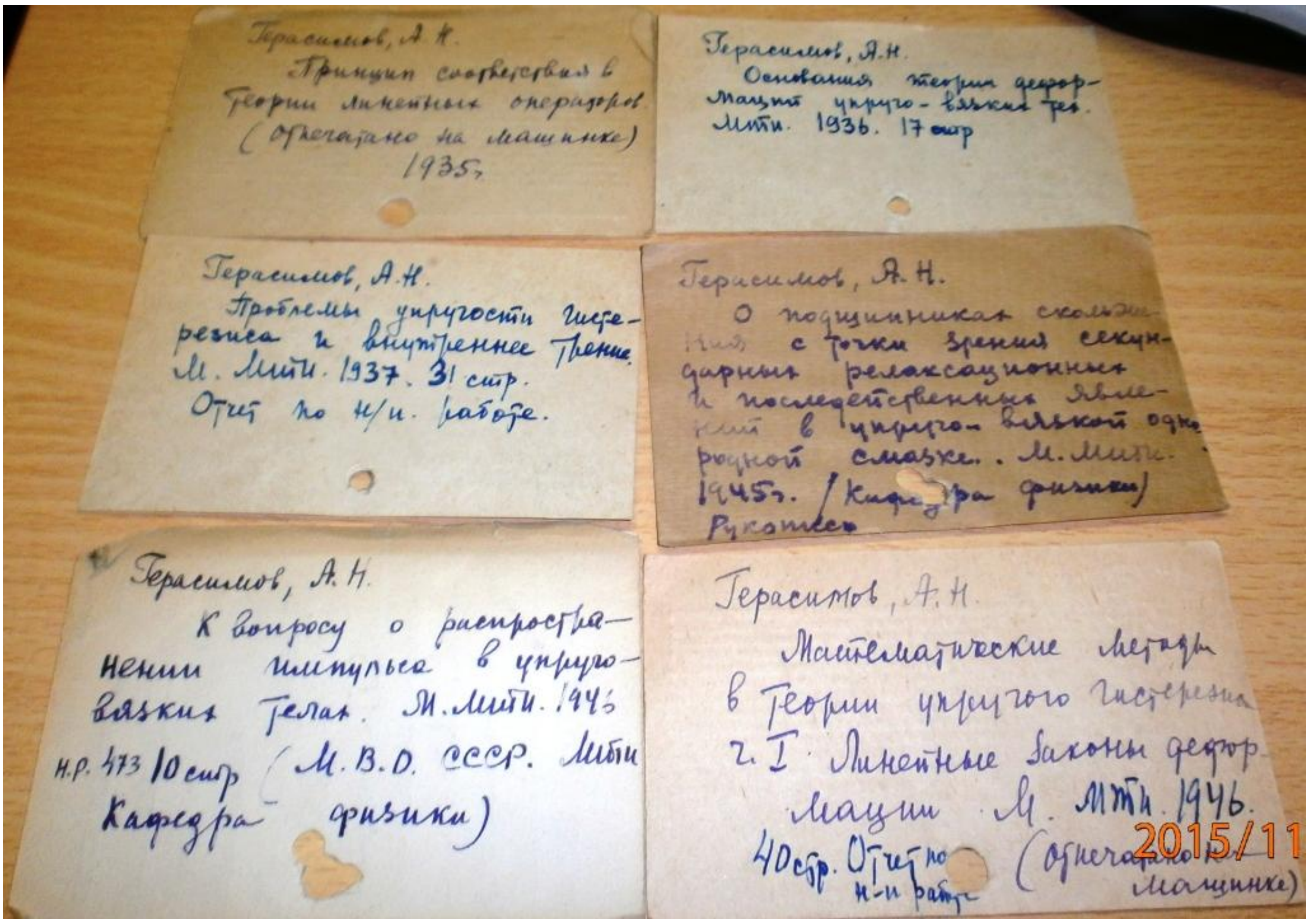
Герасимов, А.Н.
Принцип соответствия в теории линейных операторов. (Отпечатано на машинке) 1935г.

Герасимов, А.Н.
Основания теории деформации упруго-вязких тел. Мти. 1936. 17 стр

Герасимов, А.Н.
Проблемы упругости гистерезиса и внутреннее трение. М. Мти. 1937. 31 стр.
Отчет по н/и. работе.

Герасимов, А.Н.
О подшипниках скольжения с точки зрения секундарных релаксационных и последейственных явлений в упруго-вязкой однородной смазке. М. Мти. 1945г. (Кафедра физики) Рукопись

Герасимов, А.Н.
К вопросу о распространении импульса в упруго-вязких телах. М. Мти. 1946 н.р. 473 10 стр (М.В.О. СССР. Мти Кафедра физики)

Герасимов, А.Н.
Математические методы в теории упругого гистерезиса ч. I. Линейные законы деформации М. МТИ. 1946. 40 стр. Отчет по н-и работе (отпечатано на машинке)

Карточки из библиотеки МТИ

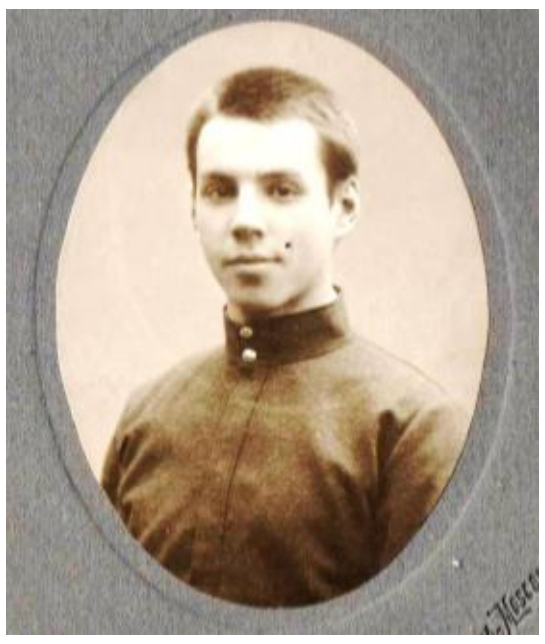
МГУ,1918

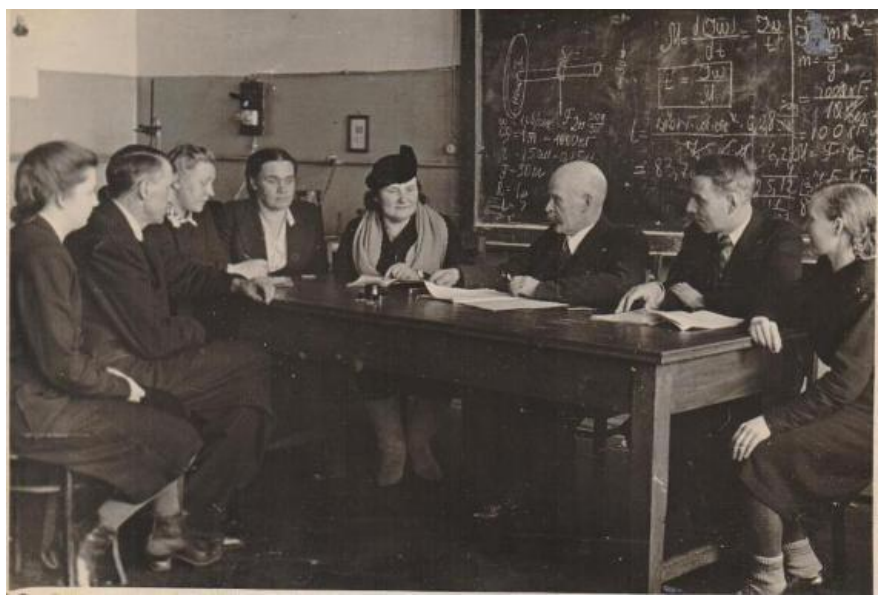
МТИ,1935

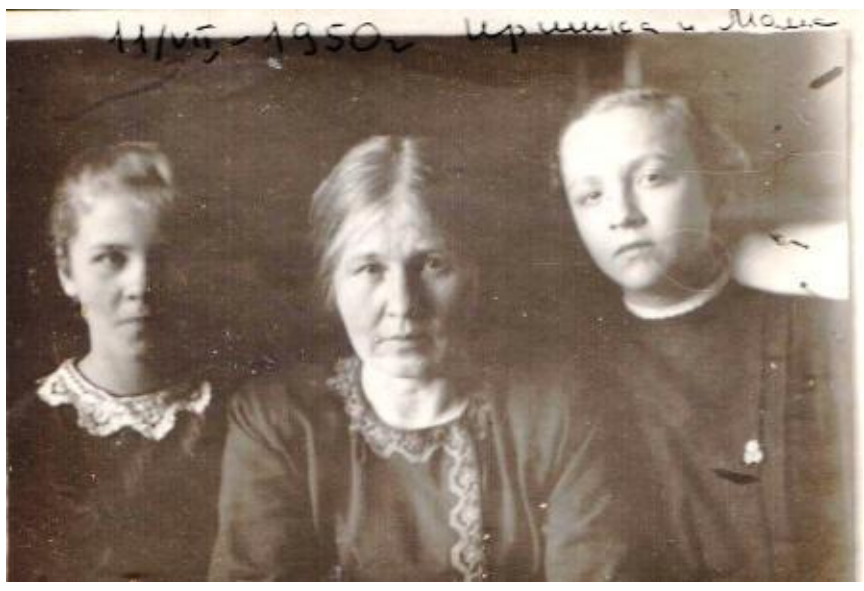
Семья 11 июля 1950

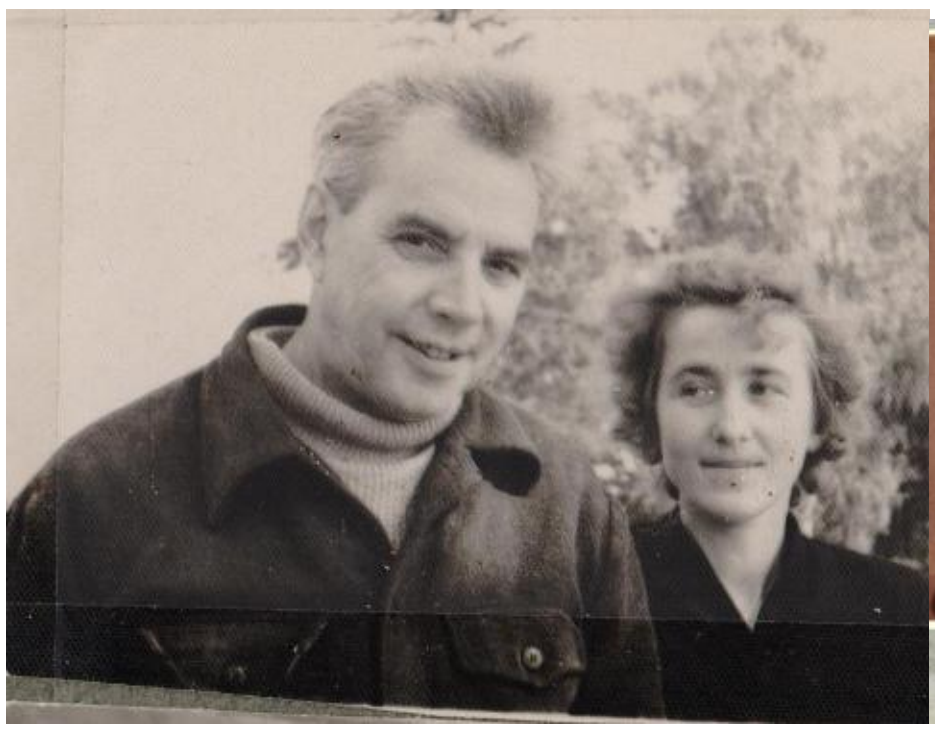
С дочерью Роной 1956

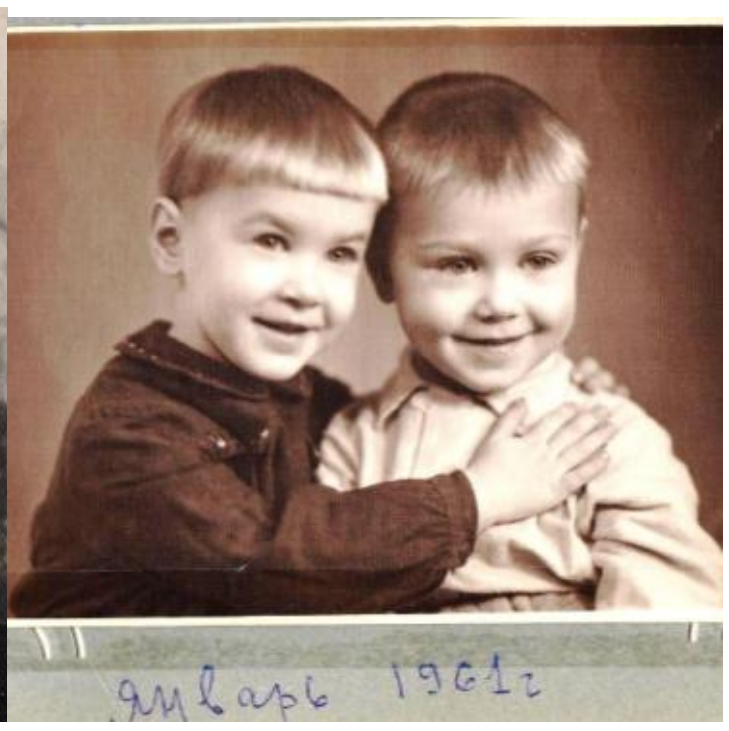

Внуки 1961

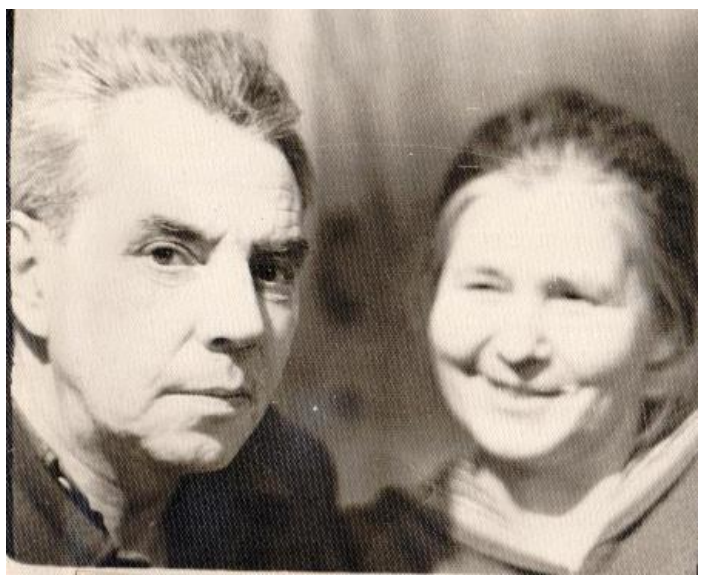
Супруги 1960

После рождения дочери Ирины в 1940 г. вернулись в Москву и восстановили прописку, обменяв маленькую квартиру на Большой Якиманке на нормативное жилье в Угловом переулке. С 1940 г. по 1943 г. А.Н. работал в МАИ на кафедре физики, с октября 1943г.- снова в МТИ, доцентом.

Надо отметить, что А.Н. тяжело переносил голод в войну. Старшая дочь Вероника рассказывала, как она собирала крапиву по московским пустырям, как давала через час младшенькой хлебушка с ноготок, как ругались родители, хотели разводиться, да потом пожалели отца.

В 1942 г., после гибели в 1940 г. научного руководителя И.И.Привалова, А.Н. Герасимов закончил заочно аспирантуру и представил уже под руководством В.В. Степанова к защите диссертацию на тему: «Некоторые задачи теории упругости с учетом последействия и релаксации по линейному закону», а 8 марта 1943 г. в НИИ механики МГУ ее защитил.

5 июня 1943 г. А.Н. был утвержден в ученом звании доцента МАИ по кафедре «Высшая математика».

В 1945 -1947гг. А.Н. делает ряд докладов на конференциях МТИ, на семинарах механики в МГУ и в ОТН АН СССР по готовящейся докторской диссертации «*Производные дробных порядков и их применение к динамике волокнистых материалов»*.

В характеристике от 24.04.1947 г: «В настоящее время он занят оформлением своей докторской диссертации на тему: «*Математические методы теории деформации упруго-вязких тел»*». 29 мая 1947 г. А.Н. делает доклад в Институте механики МГУ, 5 июня 1947 г. сдает в журнал ПММ статью «Обобщение линейных законов деформирования и его применение к задачам внутреннего трения» (ПММ, ОТН, 1948,т.XII, №3, с.251-260). Эта статья оказалась первой публикацией в мировой литературе, в которой дробная производная использовалась для описания вязкоупругих материалов.

В 1949 А.Н. Герасимов на конференции МТИ делает доклад извлечений из докторской диссертации «*Применение производных дробного порядка к динамике вытяжного поля*».

В характеристике от 30.06.1950 г.: «…А.Н. Герасимов готовит диссертацию на степень доктора физ.-мат. наук «*Применение производных дробного порядка к динамике волокнистых материалов*», которая должна быть закончена в декабре 1950 г. Одновременно с этим состоит (*беспартийный*!) слушателем Ун-та марксизма-ленинизма (2-ой год)».

К сожалению, 22 июля 1950 г после продолжительной болезни умирает и второй научный руководитель - Вячеслав Васильевич Степанов. В МТИ, МГУ и ОТН АН ССР критикуют докторскую диссертацию Герасимова, не принимают использования аппарата дробного дифференцирования: «Какая от этого практическая польза?»

А.Н. Герасимов находит иной объект для применения своей идеи. Магнитная добротность, величина, обратная тангенсу угла магнитных потерь, характеризует соотношение упругих и вязких сил, которые существуют в подвижной системе динамика вблизи частоты резонанса. Но и опять его не понимают. Статья «К расчету магнитных полей многослойной катушки», сданная в журнал «Электричество» в 1950 г., пролежала в редакции 3 года и так и не была опубликована.

Минимум четырежды перекраиваемая диссертация теряет свою «дробность». В характеристике от 10 апреля 1954 г. читаем: «… готовит *докторскую диссертацию по вопросам вытягивания*, которую должен закончить в декабре 1954 г.».

В характеристике от 18 декабря 1954 г.: « Доцент Герасимов А.Н. –хороший лектор и педагог. Он умеет преподносить студентам материал курса в живой и увлекательной форме, что подтверждается отзывами на стенограммы его лекций. Для его лекций особенно характерны точность формулировок и высокий научный уровень». Здесь уже не упоминается о докторской диссертации, зачеркнута последняя строка «~~Характеристика дана в связи с защитой им докторской диссертации»~~, но сделана надпись «для физфака МГУ».

В 1954 г. А.Н. Герасимов устраивается на полставки на физфак МГУ. К этому времени старшая дочь Вероника (Рона) уже замужем за Александром Афанасьевым, работает после пединститута учителем физики, младшая восьмиклассница Ирина закончила музыкальную школу. В 56 г. Ирина, золотая медалистка, почему-то идет не на мехмат в МГУ, а тоже в педагогический ин-т.

При первом нашем телефонном разговоре Рона сказала: «Я знаю Работнова, отец много с ним работал». Возможно, это было связано с последними публикациями Герасимова в 1956-58 гг. в журнале « Известия АН СССР, ОТН», где Ю.Н. был в это время главным редактором. Герасимов публикует четыре статьи (*Кинетика процесса вытягивания.1.Стационарный процесс; Кинетика процесса вытягивания. II. Нестационарный процесс; О скоростях волокон в процессе вытягивания; Квазистационарный процесс работы вытяжного прибора*), в которых уже нет ни слова о дробной производной, лишь проскальзывают фразы об отсутствии последействия в технологическом процессе вытягивания волокон.

В 1958 г. 23 марта у Роны рождаются близнецы – Елена и Андрей, которые-таки по стопам деда пошли на мехмат. Елена защитила кандидатскую диссертацию, а Андрей после аспирантуры, кончив курсы бухгалтеров, сначала стал банкиром, а после дорос до помощника министра финансов.

Алексей Никифорович так и не защитил докторскую диссертацию. Он умер в 1968 г. 14 марта, а в декабре того же года от запущенного рака умерла Алена, два последние года ухаживавшая за парализованным мужем. Прах их покоится в 34-й секции 22-го колумбария Донского монастыря.

# I. ПРЕДИСЛОВИЕ

Представляем 12 уцелевших из найденных наименований работ Алексея Никифоровича Герасимова – советского механика, получившего по единственной статье 1948 г. [1] известность в качестве «пионера применения дробного дифференцирования», до настоящего времени не имевшего в научной литературе даже собственного имени [2-4]. В этой статье было всего 3 ссылки - на Л. Больцмана, А.И. Лурье и А.Ю.Ишлинского. Только благодаря В.С. Постникову [5], одновременно упоминавшему статьи Герасимова [6] и Ишлинского [7], нам удалось найти более ранние публикации А.Н.Г. в ПММ, а также два его домашних адреса -московский (с местом работы в МТИ) и в г. Йошкар-Ола. Поначалу длительный перерыв в публикациях послужил тщетному поиску А.Н.Г. среди репрессированных. По ответу на запрос из Марийского государственного краеведческого музея I-й московский государственный университет им. Покровского был идентифицирован с МГУ имени М.В. Ломоносова, в архиве которого был найден диплом А.Н.Г., а в архиве физфака – личное дело совместителя с упоминанием о кандидатской работе, защищенной в МГУ. По запросу в МТИ обнаружились лишь библиотечные карточки на утраченные отчеты и публикация 56 г. в журнале «Известия АН ССР» [8] (потом здесь же найдены еще 3 публикации). И уже во время симпозиума 2016 г., посвященного 105-летию А.А. Ильюшина, в библиотеке МГУ обнаружилась (вот удача!) не списанная с 1943 г. рукопись кандидатской диссертации А.Н.Г. По смене московских прописок нашлась старшая дочь – Калина В. И. с внуками, чьи любезно предоставленные фото и воспоминания пополнили биографию А.Н.Г. К сожалению, архив А.Н. Герасимова оказался утрачен после смерти в 1996 г. его младшей дочери И. А. Томилиной и последующей продажи квартиры. Но и сохранившиеся работы дают достаточно полный научный портрет А.Н. Герасимова - дважды мехматянина, ученика профессоров И.И. Привалова[9] и В.В. Степанова, к краткому изложению которых мы и переходим.

**1.***В статье «ТЕОРИЯ РЫЧАЖНЫХ ВЕСОВ С ПОСТОЯННОЙ ЧУВСТВИТЕЛЬНОСТЬЮ»* получена общая математическая теория рычажных весов с постоянной чувствительностью и *равномерной шкалой*. Рассмотрены наиболее распространенные виды рычажных весов, так называемые «квадранты», используемые для измерения веса данного мотка пряжи или номера последней, а также «динамометры», используемые для определения разрывного усилия в весовом счете. Найдена такая *форма края лекала* коромысла квадранта, основанного на принципе лекала и ленты, при которой шкалу квадранта можно равномерно градуировать по измеряемой величине (весу). В предположении невесомости ленты, на которую крепится вес, приближенное решение общей задачи о форме края лекала получено путем решения уравнения Клеро. Изучен вопрос о величине ошибки вследствие пренебрежения весом ленты.

Получены уравнения, которые определяют *форму края лекала для квадранта с равномерной шкалой, измеряющей вес пряжи:*

$$x = a\frac{\sin^2\varphi}{\varphi^2}\ , \quad y = a\frac{\varphi - \sin\varphi\cos\varphi}{\varphi^2}\ ,$$

где угол поворота коромысла $\varphi$ меняется в пределах от $\varphi_1$ до $\varphi_2$, постоянная квадранта $a = ClG$ имеет размерность длины и пропорциональна расстоянию $l$ от оси вращения коромысла до его центра тяжести и весу коромысла $G$ .

Для квадранта с *равномерной шкалой,* измеряющего *номер N пряжи,* $N = \frac{\lambda}{P}$, и $b = \frac{lG}{C\lambda}$ ($b$ , $\lambda$ –постоянные, имеющие размерность длины), получены параметрические уравнения, которые определяют *форму края лекала*:

$$x = b\sin^2\varphi, \quad y = b[(\varphi_0 - \varphi) - \sin\varphi\cos\varphi],$$

где $\varphi$ меняется согласно $b$ способностью прибора в пределах от $\varphi_1$ до $\varphi_2$.

Было отмечено, что этот же принцип может быть применен в текстильном деле и к динамометрам, *в том числе на основе электромагнита (соленоида и катушки),* определяющим разрывную длину нити, поскольку разрывная длина выражается эмпирически в функции веса нити.

**2**. *Диплом «ПРИНЦИП СООТВЕТСТВИЯ В ТЕОРИИ ЛИНЕЙНЫХ ОПЕРАТОРОВ»* по специальности «Прикладная математика», написанный на мехмате МГУ под научным руководством И. И. Привалова, тесно связан с работами фон Дж. Неймана [10] и Ф. Рисса [11-13]. Во Введении изложены основные свойства комплексного гильбертова пространства Ꜧ, конкретно функционального пространства из функций, обладающих интегрируемым квадратом модуля в некоторой основной области Ω, изоморфного с Ꜧ. Дано понятие линейного оператора и внутреннего произведения, рассмотрены область определения и образы. Изучены возможности продолжения того принципа, который составляет сущность метода изучения ограниченных операторов Ф. Рисса [12]. Согласно этому принципу, если R ограниченный самосопряженный оператор и (Ru,u) при $|u| = 1$изменяется в пределах от m до M, то всякому полиному f(μ) с действительными коэффициентами для значений μ, заключенных между m и M, соответствует ограниченный самосопряженный оператор f(R), причем это соответствие положительного типа, мультипликативно и дистрибутивно. Исходя из этого положения, Ф. Рисс распространяет его на все ограниченные функции f(μ), являющиеся пределами возрастающих (или убывающих) последовательностей непрерывных функций и получает для оператора f(R) выражение в виде интеграла Стилтьеса

$$f(R) = \int_m^M f(\mu)dE_\mu,$$

где $E_\mu$ –особым образом сконструированный оператор, зависящий от непрерывно меняющегося параметра μ, обладающий свойствами

$E_\mu=0,\ \mu<m;\ E_\mu=E,\ \mu>M;\ \dot{E}_\mu E_\lambda u = E_\mu u, \mu \le \lambda.$

А.Н.Герасимовым показано, что принцип соответствия Рисса, с сохранением своих свойств в основном и лишь с небольшими изменениями в деталях, распространяется и на класс операторов *неограниченных.* Проделана *реставрация* того пути, по которому Ф.Рисс шел первоначально и который был впоследствии заменен другим.

Изучены свойства некоторых ограниченных операторов, связанных с данным самосопряженным оператором, (ограниченным или неограниченным)/§ 1/. При этом исходной точкой взято уравнение, рассмотренное Ф.Риссом [11].
В § 2 показано, каким образом можно установить соответствие между рациональными дробями (с комплексными корнями и полюсами) и операторами. Описаны характерные свойства этого соответствия (его положительность, дистрибутивность, мультипликативность и пр.)
В §3 доказано, что соответствие с сохранением всех основных черт может быть продолжено на класс непрерывных функций, а также и на такие разрывные функции, которые являются пределами последовательностей непрерывных.
В §4 соответствие распространено на неограниченные функции, показана применимость всего предыдущего построения выражения оператора в виде интеграла Стилтьеса, условия сходимости последнего и способ получения обратного оператора для R-μE (где E –оператор тождественного преобразования). Показано, что в точках спектра оператор $R-\mu$ не имеет себе обратного.

**3.***В статье «ПРОБЛЕМА УПРУГОГО ПОСЛЕДЕЙСТВИЯ И ВНУТРЕННЕЕ ТРЕНИЕ»* для случая одномерного тела (нити) построена система уравнений, описывающая закон движения частиц упруго-вязкого тела, являясь обобщением уравнений Ламе́ и Навье-Стокса.
Работа состоит из трех частей. *Первая часть* посвящена малым поперечным колебаниям нити, закрепленной на концах. В § 1 дается вывод уравнения движения. В § 2 показано, что уравнения (N-S) гидродинамики приводят к тому же закону движения нити при условии надлежащего истолкования гидродинамического давления. В § 3 приводится доказательство единственности решения полученного уравнения с частными производными третьего порядка. В § 4 проведено решение этого уравнения по способу Д. Бернулли. Решение получается в виде ряда, и в § 5 доказывается возможность почленного дифференцирования этого ряда нужное количество раз. Далее, в § 6 устанавливается, что полученный в § 4 ряд действительно удовлетворяет всем условиям задачи. В § 7 выясняется физический смысл найденного решения. В § 8 разбирается вопрос об использовании этого решения для опытного определения коэффициентов упругости и вязкости вещества нити путем наблюдения периода колебания и декремента затухания. Далее исследуется неоднородное уравнение: в § 9 приводится решение для тяжелой нити, в § 10 рассматривается общий случай вынужденного движения нити. В § 11 исследуется вопрос о сходимости рядов и устанавливается соблюдение всех условий задачи.

*Вторая часть* этой работы представляет собой довольно сжатый набросок решения задачи о крутильных колебаниях нити, один из концов которой закреплен, а другой несет некоторый дополнительный момент инерции, причем разбирается лишь случай колебаний в пустоте и с небольшой амплитудой. В § 12 показано, что задача о крутильных колебаниях приводит к тому же уравнению движения, что и задача о поперечных колебаниях(§ 1). В § 13 исследуются условия на свободных концах нити, в § 14 дается решение уравнения, в § 15 устанавливается, что это решение удовлетворяет всем условиям задачи.

*Третья часть посвящена упругому гистерезису*. В § 17 приведен упрощенный вывод уравнения Больцмана; § 18 посвящен способу определения ядра уравнения для того случая, когда рассматриваемая частица нити не лежит в особой точке. С использованием дискретизации времени в уравнении В. Вольтерра

$$u(x,t) = \frac{1}{M}[f(x,t) - \rho\ddot{u}(x,t)] + \int_0^t K(x,t-\theta)[f(x,\theta 0 - \rho\ddot{u}(x,\theta)]d\theta$$

(где x– абсцисса по оси нити, u –линейное перемещение или угол закручивания, M–масса бесконечно тонкого слоя нити, f(x,t) – приложенная сила, ρ– плотность нити) и метода Прони дано приближенное уравнение для K($x$,$t$-θ) в виде суммы показательных функций

$$K(x,t-\theta) = \sum_{p=1}^{N} A_p(x)e^{\alpha_p(x)[t-\theta]},$$

где $A_p(x)$ , $\alpha_p(x)$ − надлежащим способом подобранные для данного x числа. Были описаны вынужденные колебания *анизотропной* нити, доказана единственность решения.

В § 19 дается общий способ определения ядра.

**4.** *ПИСЬМО В РЕДАКЦИЮ (ПОПРАВКА К СТАТЬЕ « ПРОБЛЕМА УПРУГОГО ПОСЛЕДЕЙСТВИЯ И ВНУТРЕННЕЕ ТРЕНИЕ»*. Предложено уточненное доказательство единственности решения данного там в § 3 уравнения колебаний упругой вязкой нити. Кстати, свой вариант доказательства единственности напечатал А. Ю. Ишлинский [14].

**5**. *Статья «ОСНОВАНИЯ ТЕОРИИ ДЕФОРМАЦИЙ УПРУГО-ВЯЗКИХ ТЕЛ»*. Начав с одномерной теории нитей, обладающих одновременно свойствами упругости и вязкости, А.Н. Герасимов вводит *понятие упруго-вязкости* (ныне общепринято «вязко-упругость»). Им получен тензор упруго-вязких напряжений, выраженный в функции тензоров деформаций смещений и скоростей и их первых инвариантов, показана его симметричность.

6. *В статье « К ВОПРОСУ О МАЛЫХ КОЛЕБАНИЯХ УПРУГО-ВЯЗКИХ МЕМБРАН»* рассмотрена задача о малых колебаниях достаточно тонкой пластинки, вещество которой обладает одновременно и упругими и вязкими свойствами. Решение получено путем разложения в ряд Фурье по фундаментальным функциям нормированной и ортогональной системы. А.Н.Герасимовым доказан и определен *элевтероз* - постепенное освобождение основного тона от сопутствующих обертонов вследствие внутреннего трения, полученный для малых колебаний упруго-вязких нитей и мембран. Эффект элевтероза позволяет в найденном решении ограничиваться лишь конечным числом членов. На основании этого эффекта, как указывал А.Н. Герасимов, «можно было бы применить акустический метод к практическому определению модуля сдвига и коэффициента внутреннего трения (кстати, позволяющего вычислять время релаксации по модели Максвелла)».

**7.** *Кандидатская диссертация «НЕКОТОРЫЕ ВОПРОСЫ ТЕОРИИ УПРУГОСТИ С УЧЕТОМ РЕЛАКСАЦИИ И ПОСЛЕДЕЙСТВИЯ ПО ЛИНЕЙНОМУ ЗАКОНУ».* Перечисленные выше статьи, опубликованные А.Н.Герасимовым в 1939-1940 году, посвященные колебаниям упруго-вязких тел, были представлены в качестве диссертации на соискание степени кандидата физико-математических наук. Он снабдил их общим предисловием, обзором литературы по упруго-вязким и пластическим деформациям, критикой сообщаемым в них результатов и дополнением, содержащим некоторые новые результаты, полученные автором в 40-42 годах.

Приведен обзор литературы по данному вопросу и по смежным с ним задачам математической теории упругости. В этом обзоре одна часть касается теории пластических деформаций, другая посвящена специально линейным законам упругого последействия и релаксаций.

Вкратце изложено содержание собственных статей и дана критика сообщаемых там выводов.

Последняя часть содержит изложение задачи о колебаниях упруго-вязкой плоской мембраны, релаксирующей и последействующей по линейному закону. В этой же части дается обобщение разработанной теории, полученное А.Н.Герасимовым уже после опубликования приложенных работ, на трехмерный случай.

**8.** *Известная статья «ОБОБЩЕНИЕ ЛИНЕЙНЫХ ЗАКОНОВ ДЕФОРМАЦИИ И ИХ ПРИЛОЖЕНИЕ К ЗАДАЧАМ ВНУТРЕННЕГО ТРЕНИЯ»,* сделанная по докладу в Институте механики АН СССР 29 мая 1947 г., за 20 лет до публикации работы с дробной производной Капуто [15]. Однако правильнее ее называть дробной *производной Герасимова−Капуто* [3], так как в этой статье А.Н. Герасимов ввел аналогичную производную при рассмотрении дифференциальных уравнений с частными дробными производными.

А.Н. Герасимов упругое последействие по теории наследственности описывает линейным интегральным уравнением Больцмана

$$\sigma(t) = E\varepsilon(t) + \int_0^{\infty} G(\tau)\varepsilon(t-\tau)d\tau ,$$

где  E – упругая постоянная и $G(\tau)$–наследственная функция, определяемая из опыта, наследственная часть напряжения $\sigma(t)=\int_0^\infty G(\tau)\varepsilon(t-\tau)d\tau$.

В вязком случае, когда напряжение зависит от скорости деформации, наследственная часть напряжения равна $\sigma(t)=\int_0^\infty K(\tau)\dot{\varepsilon}(t-\tau)d\tau$. Для материалов волокнистой структуры наследственная функция имеет вид $K(\tau)=\frac{A}{\tau^\alpha}$, где постоянная A>0 и 0<α<1. Положив $A=\frac{\chi}{\Gamma(1-\alpha)}$, где χ>0 зависит от свойств вещества и Γ есть Эйлеров интеграл второго рода (гамма–функция), получает

$$\sigma(t)=\chi\frac{1}{\Gamma(1-\alpha)}\int_0^\infty\frac{\dot{\varepsilon}(t-\tau)d\tau}{\tau^\alpha}=\chi\frac{\partial^\alpha\varepsilon(t)}{\partial t^\alpha}\quad(0<\alpha<1)\qquad(1)$$

так как производная от ε(t) по t порядка α будет $\frac{1}{\Gamma(1-\alpha)}\int_0^\infty\frac{\dot{\varepsilon}(t-\tau)d\tau}{\tau^\alpha}$.

Это линейное соотношение между ε и σ при α=0 обращается в закон Гука, при α=1 – в закон Ньютона для внутреннего трения (вязкости).

С использованием зависимости (1) рассмотрена задача о движении среды между двумя параллельными плоскостями, из которых одна неподвижна, а другая движется прямолинейно параллельно первой по закону $y(x,t)|_{x=l}=\varphi(t)$, ($\varphi(0)=0$, $\varphi'(0)=0$), где φ(t)–данная функция времени. По методу Хевисайда, с учетом теоремы Бореля, для предельных случаев α=0 и α=1 получены соответственно уравнение вынужденных поперечных волн в упругой среде и распределение смещений в жидкости, данное А.И. Лурье. Для случая $\alpha=\frac{1}{2}$ снова по методу Хевисайда получено

$$\sigma|_{x=l}=\chi cv\left\{\frac{1}{\Gamma(m)t^{1-m}}+2\sum_{k=0}^{\infty}\sum_{j=0}^{\infty}(-1)^j\frac{[(2k+2)cl]^j}{j!\Gamma(-jm+m)t^{1+(j-1)m}}\right\},\ m=1-\frac{\alpha}{2}.\qquad(2)$$

Во втором примере проанализировано движение смазки в подшипнике скольжения при рассмотрении движения жидкости между коаксильными цилиндрическими поверхностями, вращающимися по заданному закону. С помощью интегралов Стилтьеса показана возможность описания уравнением (1) деформирования не только по линейному закону $\sigma+\frac{1}{q}\dot{\sigma}=E(\varepsilon+\frac{1}{n}\dot{\varepsilon})$, но и для вторичных релаксационно-последейственных процессов, описываемых соотношением

$$\omega+\frac{1}{q}\dot{\sigma}+\lambda\ddot{\sigma}=E(\varepsilon+\frac{1}{n}\dot{\varepsilon}+\mu\ddot{\varepsilon}),$$

при определенных допущениях.

**9.** *Статья «КИНЕТИКА ПРОЦЕССА ВЫТЯГИВАНИЯ. 1. СТАЦИОНАРНЫЙ ПРОЦЕСС»* была опубликована после большого перерыва в 1956 г., после четырех неудачных защит докторской диссертации по использованию аппарата дробного дифференцирования. В ней и в последующих трех работах, А.Н.Герасимов отказывается даже об упоминании не только дробной производной, но и явления последействия. Тем не менее, судя по архивным данным, с большой долей уверенности можно предположить, что в последних четырех статьях изложены те же результаты, но полученные другими методами.

Рассмотрен простейший вытяжной аппарат с двумя парами валиков. Начало координат О поместим в зажиме питающей (задней) пары, ось $Ox$ направим по движению ленты. Растягиваемый волокнистый объект мы называем лентой. При стационарном процессе индивидуальные скорости отдельных волокон в границах вытяжного поля представлены почти исключительно только двумя значениями, причем каждое волокно меняет свою скорость мгновенно (и не больше одного раза) с меньшего значения $v_0$ на большее $v_1$. Это положение мы будем называть в дальнейшем *принципом отсутствия промежуточных скоростей.* Предполагалось, что в заднем зажиме имеются волокна, уже получившие большую из двух возможных скоростей, тогда как в переднем зажиме (по движению ленты) имеются волокна, еще сохранившие меньшую скорость. Считалось, что процесс вытягивания сопровождается проскальзыванием волокон в обоих зажимах Кроме того, учтена способность волокон растягиваться. Показано, как из этих положений, без привлечения каких-либо иных гипотез, методами классической механики подвижных сред строится вся кинетика процесса.

<u>Обозначено:</u>

$q$ – секундный расход массы (масса волокнистого материала, проходящего через сечение поля в единицу времени – секунду). В стационарном процессе $q$ – одно и то же для всех сечений поля. Вместе с тем $q$ является линейной плотностью количества движения.

$\Lambda$ – средняя линейная плотность массы волокна, обратно пропорциональная номеру волокна;

$\lambda(x)$ – линейная плотность массы ленты в сечении поля; величина, обратно пропорциональная номеру ленты в этом сечении.

Показано, что *средняя арифметическая скорость* $v(x)$ волокон, пронизывающих сечение $x$, есть не что иное, как *скорость центра массы этих волокон.*

Получено соотношением для связи фактической вытяжки $b$ с предельной $B$ через меры *проскальзывания* волокон в заднем зажиме $\alpha$ и проскальзывания волокон в переднем зажиме $\beta$. При этом максимальное отношение средней квадратичной скорости $u(x)$ к квадрату средней арифметической скорости по совокупности всех волокон, пронизывающих сечение $x$, $v(x)$, равно

$$\sigma_m = \frac{(v_0 + v_1)^2}{4v_0v_1} = \frac{(B+1)^2}{4B}$$

Через разность *линейной плотности кинетической энергии* в сечении $x$ поля

$$E(x) = \frac{\lambda(x)[u(x)]^2}{2} = \frac{q}{2v}\left[(v_0 + v_1)v - v_0v_1\right] = \frac{qv_0}{2} + \frac{qv_1}{2} - \frac{qv_0v_1}{2v(x)}$$

и линейной плотности кинетической энергии для случая ,если бы все волокна, пронизывающие это сечение, имели единую скорость $v(x)$, равной

$$E'(x) = \frac{\lambda(x)[v(x)]^2}{2} = \frac{qv(x)}{2},$$

получено выражение для *диссипативной функции*

$$D(x) = E(x) - E'(x)$$

или
$$D(v) = \frac{q}{2}\frac{(v_1 - v)(v - v_0)}{v}$$

В наиболее общем случае плотность $f$ растягивающей силы представляется в виде степенного ряда по диссипативной функции

$$f = a_0 + a_1D + a_2D^2 + ... \quad ,$$

в линейном приближении

$$f = r + \delta D$$

или
$$q\frac{dv}{dx} = r + \frac{q\delta}{2}\frac{(v_1 - v)(v - v_0)}{v},$$

решая которое, получено *уравнение распределения средних арифметических скоростей по сечениям поля:*

$$\frac{v - w_0}{(w_1 - v)^{\varepsilon}} = PQ^x .$$

Здесь положительные постоянные $P$ и $Q$, имеют значения:

$$P = \frac{v(0) - w_0}{[w_1 - v(0)]^{\varepsilon}}, \quad Q = \frac{v(1) - w_0}{v(0) - w_0}\left(\frac{w_1 - v(0)}{w_1 - v(1)}\right)^{\varepsilon} \quad ,$$

$$\varepsilon = \frac{w_1}{w_0} > 1, \quad w_0 = \frac{r}{\delta q} + \frac{v_0 + v_1}{2} - \sqrt{\left(\frac{r}{\delta q} + \frac{v_0 + v_1}{2}\right)^2 - v_0v_1} \quad ,$$

$$w_1 = \frac{r}{\delta q} + \frac{v_0 + v_1}{2} + \sqrt{\left(\frac{r}{\delta q} + \frac{v_0 + v_1}{2}\right)^2 - v_0v_1} .$$

Величина, обратная $\delta$, имеет размерность длины. Легко видеть, что $\delta$ зависит только от вытяжки (предельной) $B$ и от обоих коэффициентов проскальзывания $\alpha$ и $\beta$; а именно:

$$\frac{1}{\delta} = S = \frac{B - 1}{2\ln\left[(1-\alpha)^B(1-\beta)/(\alpha\beta^B)\right]}$$

Рассмотрены различные примеры, построены кривые изменения скоростей и утонения.

**10.** *В статье «КИНЕТИКА ПРОЦЕССА ВЫТЯГИВАНИЯ. II. НЕСТАЦИОНАРНЫЙ ПРОЦЕСС»* приведены уравнение неразрывности

$$\frac{\partial(\lambda v)}{\partial x} = -\frac{\partial \lambda}{\partial t} \qquad (1)$$

и уравнение движения на основе общих соображений, разработанных И.В. Мещерским. Если в момент $t$ масса тела есть $M$, скорость его в этот момент есть $v$ и если за время $dt$ тело присоединяет массу $dM$, которая до присоединения имела скорость $u$, то, называя буквой $F$ движущую силу в момент $t$, имеем

$$F = \frac{d(Mv)}{dt} - u\frac{dM}{dt}$$

Показано, что движущая сила $F(x,t)$, действующая в момент $t$ на отрезок $(0,x)$ ленты, определяется равенством

$$F(x,t) = q(x,t)v(x,t) - q(0,t)v(0,t) + \int_0^x \frac{\partial q(\xi,t)}{\partial t} d\xi \quad , \qquad (2)$$

или, так как

$$q = \lambda v \quad , \qquad (3)$$

$$\frac{\partial[F - \lambda v^2]}{\partial x} = \frac{\partial(\lambda v)}{\partial t}$$

Система трех независимых уравнений (1) - (3) связывает четыре функции: $F, q, \lambda, v$ переменных $x$ и $t$. Эта система должна быть дополнена четвертым уравнением, не зависящим от трех предыдущих и связывающим, вообще говоря, все четыре функции. Это уравнение состояния ленты должно отражать *связь между структурными свойствами объекта и параметрами процесса*. Рассмотрены случаи определения $q$, $v$ и $F$ по известной $\lambda(x,t)$ и секундному расходу массы $q(0,t)$ на входе в поле, а также определения $\lambda$, $v$ и $F$ по известному секундному расходу массы $q(x,t)$ и распределению линейной плотности массы $\lambda(x,0)$ в начальный момент времени, отмечена единственность решения в обоих случаях.

Для доказательства единственности решения при определении линейная плотности массы $\lambda(x,t)$ по распределению скоростей по сечениям поля в различные моменты времени $v(x,t)$ вместо уравнения неразрывности

$$\frac{\partial q}{\partial x} + \frac{\partial \lambda}{\partial t} = 0 \qquad (N)$$

рассмотрено несколько более общее уравнение

$$F_\alpha(\lambda) \equiv \frac{\partial(v\lambda)}{\partial x} + \alpha\frac{\partial \lambda}{\partial t} = 0 \qquad (NO)$$

где $\alpha$ – безразмерный параметр. При $\alpha = 1$ это уравнение совпадает с (N). Решение для задачи

$$F_{\alpha}(\lambda) = 0\,,\ \lambda(x,t)\big|_{x=0} = \lambda(0,t) = \varphi(t)$$

ищется в виде ряда по степеням параметра $\alpha$

$$\lambda(x,t) = \lambda_0(x,t) + \alpha\lambda_1(x,t) + \alpha^2\lambda_2(x,t) + \ldots$$

Составляя формально ряды для производных, подставляя в (NO) и сравнивая множители при одинаковых степенях параметра $\alpha$, получена последовательность рекуррентных уравнений

$$\frac{\partial(v\lambda_0)}{\partial x} = 0\,,\ \frac{\partial(v\lambda_k)}{\partial x} = -\frac{\partial\lambda_{k-1}}{\partial t} \qquad (k = 1,2,\ldots)$$

Из этой системы определяется сначала $\lambda_0(x,t)$, которая должна чтобы она обращаться в заданную функцию $\varphi(t)$ при $x = 0$. Затем шаг за шагом находятся все $\lambda_k(x,t)$, $k = 1,2,\ldots$, причем каждая из них должна обращаться в нуль при $x = 0$. Тогда ряд разложения формально дает решение задачи. Для того, чтобы решение было не только формальным, но и фактическим, нужно иметь уверенность в равномерной и абсолютной сходимости его во всей полуполосе $0 \le x \le l\,, t \ge 0$. Если это действительно так при $\alpha = 1$, то ряд

$$\lambda(x,t) = \sum_{k=0}^{\infty} \frac{(-1)^k}{v(x,t)} \int_0 \frac{\partial}{\partial t} \ldots \frac{1}{v(x,t)} \int_0 \frac{\partial}{\partial t}\left[\varphi(t)\frac{v(0,t)}{v(x,t)}\right] \underbrace{dx\ldots dx}_{k\cdot раз}$$

является решением задачи

$$\frac{\partial\lambda}{\partial t} + v\frac{\partial\lambda}{\partial x} + \lambda\frac{\partial v}{\partial x} = 0\,,\ \lambda(x,t)\big|_{x=0} = \varphi(t)$$

Показано, что это решение будет единственно возможным. Отмечено, непосредственно из физических соображений, что знание желаемого распределения скоростей и условия на входе еще недостаточно для того, чтобы процесс мог фактически осуществляться. Этому отвечает то обстоятельство, что многократно применяемый оператор

$$\frac{-1}{v}\int_0^x \frac{\partial}{\partial t} \ldots dx$$

в ряде, дающем $\lambda(x,t)$, не ограничен, и надлежащая сходимость ряда, вообще говоря, ничем не обеспечена. В каждом конкретном случае сходимость ряда должна проверяться применительно к этому случаю.

**11.** В *статье «О СКОРОСТЯХ ВОЛОКОН В ПРОЦЕССЕ ВЫТЯГИВАНИЯ»* приведено недостающее четвертое уравнение, полученное экспериментально, путем промеров секундного расхода массы и линейной плотности ленты в различных сечениях поля и в различные моменты времени

Предполагалось *во-первых*, что процесс, вообще нестационарный, протекает *без последействия*, т.е. распределение скоростей волокон зависит лишь от мгновенных значений, определяющих процесс обстоятельств, а не от тех, которые эти обстоятельства имели в предшествующие моменты времени. *Во-вторых*, принято, что два любых волокна с одинаковой длиной, занимающие в данный момент времени один и тот же участок поля, в точках с одинаковыми абсциссами имеют одинаковые мгновенные скорости. *В-третьих*, допускалось наличие переходов от одной скорости к другой по ряду промежуточных значений. Наконец, *в четвертых*, считалось, что проскальзывания волокон нет ни на входе, ни на выходе.

Полученные из тщательно обставленного опыта сведения о распределении скоростей волокон в вытяжном поле предложено стараться подбором констант $\mu$ и $\nu$ уложить в уравнение

$$\frac{[v(l,t)-w]^{\mu}}{w-v(0,t)}=[v(l,t)-v(0,t)]^{\mu-1}\left[\frac{l-x}{x}\right]^{\nu}$$

Отсюда можно, в принципе, найти мгновенную скорость $w$ для всех тех частиц ленты, которые в момент $t$ находятся в сечении поля с абсциссой $x$.

Значение $w$ отождествлялось с искомой средней скоростью ленты в момент $t$ в сечении $x$.

*Было сделано обобщение*, что $\mu>1$ *не постоянная, а функция времени*, равная вытяжке $B(t)$, аналогично $\nu=\kappa B(t)$,где $\kappa>1$ не зависит ни от $x$, ни от $t$, причем $[\kappa]=1$. Тогда уравнение перепишется в виде

$$\frac{[v(l,t)-w]^{B(t)}}{w-v(0,t)}=[v(l,t)-v(0,t)]^{B(t)-1}\left[\frac{l-x}{x}\right]^{\kappa B(t)},$$

откуда

$$\frac{\partial w}{\partial x}=\frac{\kappa B(t)l}{[B(t)-1]x(l-x)}\frac{[w-v(0,t)][v(l,t)-w]}{w}.$$

Обнаружен и смысл введенной в конце статьи **9** «динамической длины» волокна

$$S=\frac{B-1}{2\kappa B}l.$$

Из четырех независимых уравнений:

$$q=\lambda v,\quad \frac{\partial q}{\partial x}=-\frac{\partial\lambda}{\partial t},\quad \frac{\partial(F-qv)}{\partial x}=-\frac{\partial q}{\partial t}.$$

$$\frac{[v(l,t)-v]^{B(t)}}{v-v(0,t)}=[v(l,t)-v(0,t)]^{B(t)-1}\left[\frac{l-x}{x}\right]^{\kappa B(t)}$$

возможно (правда, лишь в принципе) определить все четыре неизвестные функции $v, \lambda, q$ и $F$ независимых переменных $x$ и $t$ с надлежащими краевыми и начальными условиями.

**12.** Статья «КВАЗИСТАЦИОНАРНЫЙ ПРОЦЕСС РАБОТЫ ВЫТЯЖНОГО ПРИБОРА».
В системе уравнений:

$$q = \lambda v\,, \quad \frac{\partial q}{\partial x} = -\frac{\partial \lambda}{\partial t}\,, \quad \frac{\partial (F - qv)}{\partial x} = -\frac{\partial q}{\partial t}\,.$$

$$\frac{[v(l,t) - v]^{B(t)}}{v - v(0,t)} = [v(l,t) - v(0,t)]^{B(t)-1}\left[\frac{l - x}{x}\right]^{\kappa B(t)} \qquad (1)$$

где $q-$ секундный расход массы ленты, $\lambda-$линейная плотность, $v-$скорость, $F$ – сила, $B(t)$ – полная вытяжка в момент $t$, равная

$$B(t) = \frac{v(l,t)}{v(0,t)} \qquad (2)$$

положительная постоянная $\kappa$ подбирается из опыта и является динамической характеристикой вещества.

Система (1) определяет все четыре неизвестные функции $v, \lambda, q$ и $F$ независимых переменных, как бы не менялись с течением времени скорости $v(0,t)$ и $v(l,t)$ на входе и на выходе. Но конструкция прибора такова, что за исключением случайных малых отклонений полная вытяжка остается постоянной для любого $t$, потому что скорость $v(l,t)$ на выходе во всякий момент пропорциональна скорости $v(0,t)$ на входе. С учетом сказанного, предположено, что

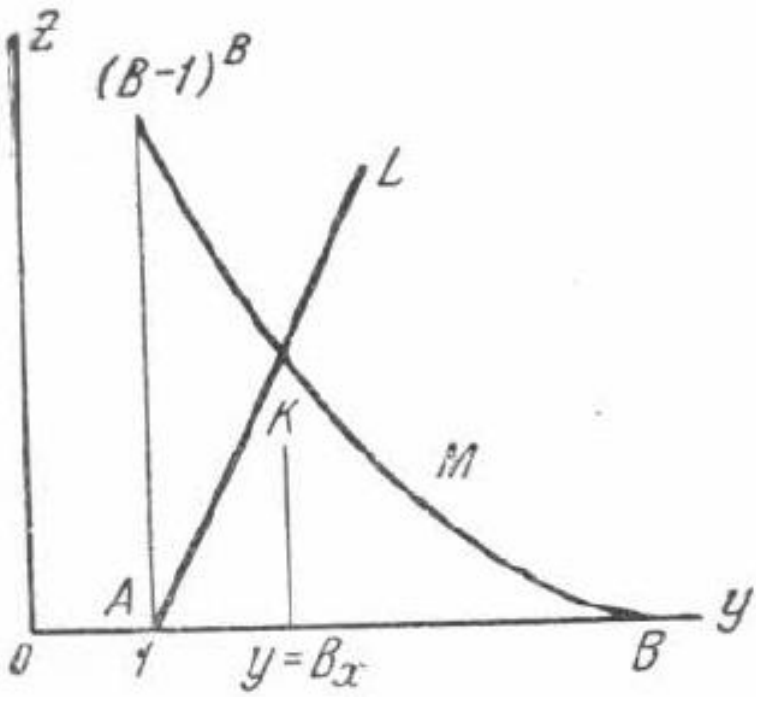

$$v(0,t) = \sigma(t)v(0)\ ,\ v(l,t) = \sigma(t)v(l)$$

где $\sigma(t)-$ данная функция $t$, причем $\sigma(t) > 0$ и $\sigma(0) = 1$. С другой стороны, *отсутствие последействия в волокнах* позволяет считать, что и вообще

$$v(x,t) = \sigma(t)v(x)$$

так что переменные $x,t$ в выражении для скорости разделяются естественными условиями процесса.

Вследствие однородности последнее уравнение (1) переписано в виде

$$\frac{(B - B_x)^B}{B_x - 1} = (B - 1)^{B-1}\left(\frac{l - x}{x}\right)^{\kappa B}$$

Для любого $B \geq 1$ и $x$ $(0 \leq x \leq l)$ показано, как графическим приемом можно определить $B_x$. Для движущей силы получено соотношение

$$F(x,t) - F(0,t) = q(t)[v(x,t) - v(0,t)] + \frac{\partial q}{\partial t}x$$

и тем самым решение задачи доведено до конца.

Сделано такое *наблюдение*. Первые три уравнения написанной выше системы (1) верны для процесса вытягивания любого линейно-протяженного тела. Вводя обозначение

$$qv - F = \Phi$$

и дифференцируя второе из уравнений (1) по $t$, а третье – по $x$, получим

$$\frac{\partial^2 \Phi}{\partial x^2} = \frac{\partial^2 \lambda}{\partial t^2}.$$

Это указывает на специфический характер процесса. Например, если бы каждая из двух функций $\Phi$ и $\lambda$ зависела бы от одной и той же величины $u$, в свою очередь являющейся функцией только двух независимых переменных $x$ и $t$, так что

$$\Phi = \Phi(u), \quad \lambda = \lambda(u)$$

то получилось бы хорошо изученное волновое одномерное уравнение

$$\frac{\partial}{\partial x}\left[\Phi'(u)\frac{\partial u}{\partial x}\right] = \frac{\partial}{\partial t}\left[\lambda'(u)\frac{\partial u}{\partial t}\right],$$

решаемое известными методами.

## ЛИТЕРАТУРА

---

В заключение составитель считает своим долгом выразить глубокую благодарность дочери А.Н.Герасимова *Веронике Ивановне Калиной*, директору Марийского мемориального народного музея истории ГУЛАГа *Николаю Александровичу Аракчееву,* директору Музея истории города Йошкар-Олы *Надежде Эриковне Логиновой*, документоведу Архива МГУ *Наталье Петровне Каргиной*, а также заведующей архива Московского Текстильного Института *Ольге Викторовне Кулаковой* .

*Москва, 13.01.2018 info-viniti@yandex.ru*

*Новоженова Ольга Георгиевна*

# II. ТРУДЫ А. Н. ГЕРАСИМОВА

## 1.ТЕОРИЯ РЫЧАЖНЫХ ВЕСОВ С ПОСТОЯННОЙ ЧУВСТВИТЕЛЬНОСТЬЮ

Одним из наиболее употребительных и распространенных видов рычажных весов является так наз. «квадранты» и «динамометры», применяемые в текстильной, бумажной, резиновой и кожевенной промышленности. Квадранты служат для измерения некоторых определенных технических величин, функционально связанных с весом исследуемого тела. В текстильном деле, например, квадранты употребляются для измерения веса данного мотка пряжи, номера последней и пр. Что же касается динамометров, то они применяются для определения разрывного усилия в весовом счете, дающего представление о крепости испытуемого материала.

Эти два прибора, имеющие один и тот же принцип рычажных весов, весьма разнообразные по особенностям применения и конструкции, обладают одними и теми же механическими свойствами и могут быть обобщены в настоящем исследовании. Приборы того и другого типа я буду называть в дальнейшем «квадрантами», хотя все излагаемое ниже с таким же успехом может быть отнесено и к динамометрам.

В наиболее распространенной системе квадрантов коромысло прибора представляет собой трехплечий рычаг. На одно плечо подвешивают тот груз, для которого желают найти числовое значение измеряемой величины, другое плечо играет роль противовеса, а третье имеет вид стрелки. Конец стрелки может перемещаться по круговой шкале, градуированной по измеряемой величине и имеющей центр на оси вращения коромысла. Последнее повертывается под действием подвешиваемого груза, и конец стрелки устанавливается против некоторого деления на шкале, указывая непосредственно искомое числовое значение.

В силу механических свойств подобной системы шкала прибора, вообще говоря, должна быть градуирована неравномерно, так как длина деления на шкале, а также и чувствительность прибора величины меняются в зависимости от подвешиваемого груза. Кроме того, процесс градуировки шкал таких квадрантов крайне кропотлив, сложен и требует от исполнителя специальной теоретической подготовки. Особенно большие затруднения встречаются в том случае, когда на шкале прибора приходится градуировать несколько (до восьми) различных поясов с делениями, как это делается, например, при изготовлении так наз. «универсальных» квадрантов системы Шоппер.

Изготовление квадрантов у нас в СССР стало возможным лишь после того, как заведующий кафедрой физики в Московском текстильном институте И.И. Васильев разработал теорию этого дела. Результаты, полученные им по этому вопросу, опубликованы в « Известиях московского текстильного института» т.II, 1929 г. в статье под названием «Теория построения шкал рычажных весов». В этой же статье описывается «шкально-микроскопический» метод, разработанный И.И. Васильевым и старшим механиком лаборатории А.П.Матвеевым. Метод этот был впоследствии усовершенствован доцентом А.Н.Марковым.

Благодаря работам И.И Васильева, А.П. Матвеева и А.Н. Маркова, физическая лаборатория Московского текстильного института и в настоящее время является, по-видимому, единственным в СССР местом, где производятся дальнейшие изменения в области усовершенствования квадрантов. Так, например, осенью прошлого года А.П. Матвеев совместно с механиком Г.К. Петровым разработал и предложил четыре конструкции квадрантов, названных «револьверными». При помощи остроумных приспособлений в квадрантах Матвеева-Петрова семь-восемь параллельных специальных поясов с делениями заменяются лишь одним безличным поясом. Конструкции Матвеева-Петрова уже вытесняют «универсальные» шопперовские приборы, имея перед последними некоторые преимущества.

Но шкала прибора и в системе Матвеева-Петрова должна быть градуирована неравномерно.

В конце прошлого года проф. Васильев порекомендовал нам заняться исследованием возможностей такой системы квадранта, которая имела бы равномерно градуированную шкалу. Лучше всего этому требованию может удовлетворить система, построенная на принципе «лекала и ленты».

Принцип этот не нов, квадранты такой конструкции встречаются в технике, но теоретические данные, относящиеся к ним, очень скудны и носят эмпирический характер. Например, в книжке J.Zingler “Theorie der Zusammensetzten Waagen” (Verl. von J. Springar, Berlin,1928) имеется лишь несколько страниц, посвященных теории этого принципа, причем автор ограничивается рассмотрением малоинтересных частностей. Развиваемые там теоретические соображения малоубедительны и лишены общности, а рассуждения о величине допускаемой ошибки заставляют желать лучшего. Кроме того, Zingler совсем не касается интересующего нас вопроса о возможности применения принципа лекала и ленты для достижения постоянства чувствительности прибора.

Предлагаемая статья имеет в виду пополнить существенный пробел в технической литературе, посвященной теории весов. В этой статье мы сначала даем общую математическую теорию рычажных весов с постоянной чувствительностью и равномерной шкалой, затем исследуем вопрос о верхнем пределе допускаемой ошибки и в заключение рассматриваем в качестве примеров два частных случая, имеющих большое значение для текстильной промышленности.

Коромысло квадрантов, основанное на принципе лекала и ленты, состоит (рис.1) из неизменно связанных между собой лекала ODBCO, край OBC которого вырезывается по некоторой определенной кривой, и стрелки указателя OA, конец которой при повороте коромысла вокруг оси O пробегает по круговой шкале, имеющей центр на оси O и градуированной по измеряемой величине. На краю лекала в точке D закрепляется один конец тонкой и гибкой стальной ленты DBK с крючком K на другом конце. На этот крючок подвешивается исследуемое тело. Под действием силы тяжести часть ленты DB охватывает дугу края лекала, а часть BK, свешиваясь свободно, заставляет коромысло повернуться на некоторый угол. Конец стрелки останавливается против определенного, соответствующего исследуемому грузу, деления на шкале и указывает непосредственно числовое значение измеряемой величины.

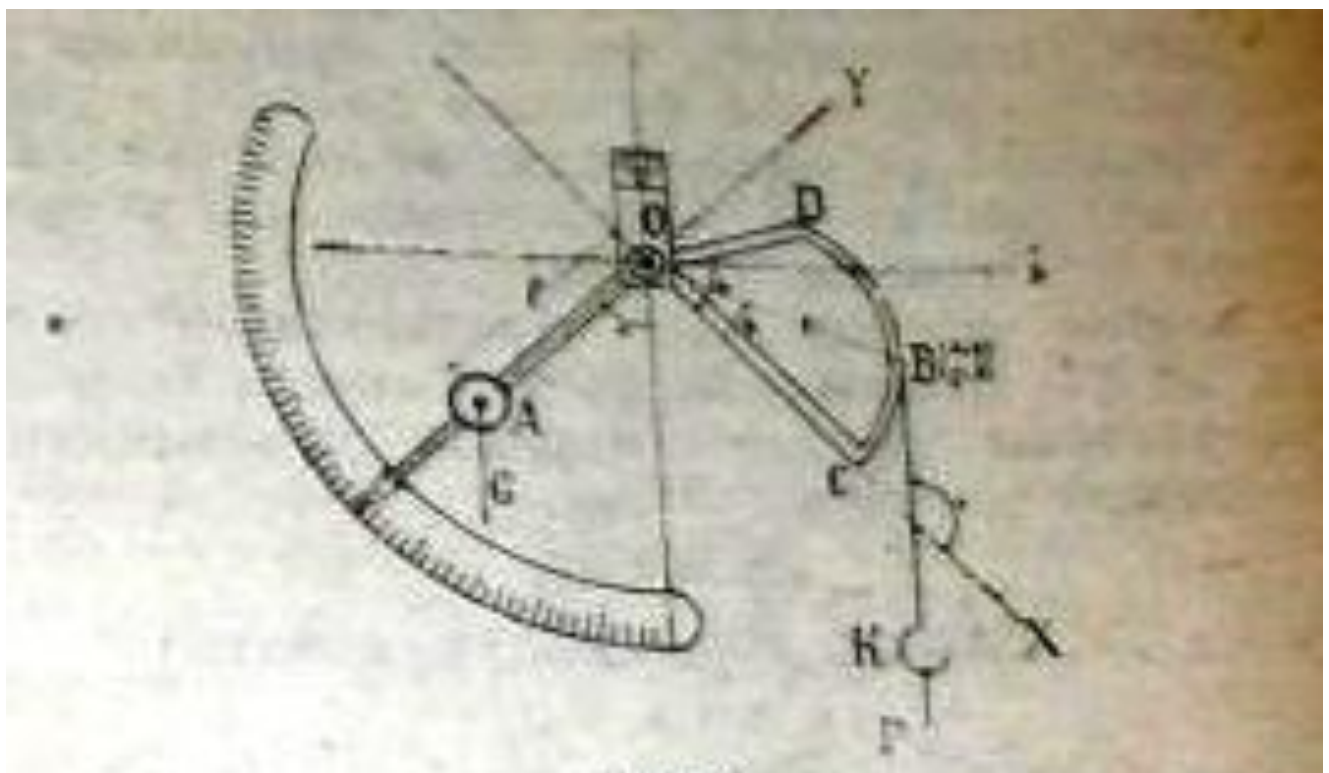

Рис. 1

Общая постановка задачи заключается в том, чтобы *найти такую форму края* DBC *лекала, при которой шкалу прибора нужно будет равномерно градуировать по измеряемой величине.*

Для решения задачи заметим, во-первых, что коромысло весов всегда можно сбалансировать так, что *центр тяжести* его будет лежать *на оси стрелки*. Не нарушая общности рассуждений, мы это и будем предполагать. Во-вторых, мы покажем дальше, что *весом ленты можно пренебречь*, не выходя за пределы допустимой ошибки. Поэтому ограничимся лишь приближенным решением задачи, предполагая, что лента невесома.

Введем следующие обозначения (рис.1):

G – вес коромысла;

P – вес исследуемого тела;

$\xi 0 \eta$ –неподвижная система координат, ось $0\xi$ которой горизонтальна;

x0y – система координат, связана с плоскостью лекала, при отсутствии груза Р совпадающая с системой $\xi 0 \eta$;

x,y – координаты в плоскости лекала той точки В его края, в которой лента сходит с лекала при нагруженном коромысле;

$\rho$, $\theta$ – полярные координаты той же точки В в плоскости лекала, причем полярной осью является ось ох;

$\varphi$ –угол поворота коромысла, соответствующий грузу Р;

$\tau$ –угол между касательной к краю лекала в точке В и осью ох;

l –расстояние от оси вращения О коромысла до его центра тяжести А.

Математическая постановка задачи сводится к следующим условиям:

1.Условию равновесия системы, заключающемуся в том, чтобы момент силы G с плечом $l\sin\varphi$ (рис.1) относительно точки О был бы равен моменту силы Р с плечом $\rho\cos(\theta-\varphi)$ относительно той же точки. Таким образом, это условие напишется в виде

$$lG\sin\varphi = \rho P\cos(\theta-\varphi). \tag{1}$$

2. Условие постоянства чувствительности прибора. Пусть

$$P = F(z) \quad , \tag{2}$$

где $z$ - подлежащая измерению величина, а

$F$ - знак функции, дифференцируемой в пределах между $z_1$ и $z_2$, зависящими от назначения и способности прибора. Относительно $F$ будем предполагать, что *логарифмическая производная* ее *обратной функции* $z = f(P)$ в указанных пределах не превосходит некоторого конечного количества М. Мы покажем, что в приложениях это условие соблюдается.

Условие постоянства чувствительности требует, чтобы отклонения $d\rho$ коромысла были пропорциональны соответствующим изменениям $dz$ измеряемой величины. Итак,

$$d\rho = Cdz,$$

где С – постоянная. Интегрируя, найдем

$$\rho = Cz + \rho_0, \tag{3}$$

где $\rho_0$ - другая постоянная. Так как $z$ по предложению может меняться лишь в пределах от $z_1$ до $z_2$, то, согласно (3) угол $\varphi$ должен меняться в пределах $Cz_1 + \varphi_0$ до $Cz_2 + \varphi_0$.

Пределы $z_1, z_2$ и $\varphi_1, \varphi_2$ оказываются связанными соотношениями:

$$z_1 = \frac{\varphi_1 - \varphi_0}{C} \ , \ z_2 = \frac{\varphi_2 - \varphi_0}{C} \ .$$

Определив из (3) и подставив в (2), найдем условие постоянства чувствительности в виде

$$P = F\left(\frac{\varphi - \varphi_0}{C}\right), \tag{4}$$

где $\varphi$ меняется в пределах между $\varphi_1$ и $\varphi_2$, а постоянная С есть мера чувствительности прибора. Она положительная, если $z$ возрастает вместе с $\varphi$ и отрицательна, если $z$ убывает при возрастании $\varphi$, как видно из уравнения (3). На чувствительность прибора мы здесь смотрим следовательно в несколько обобщенном смысле, приписывая ей знак.

3. Третье и последнее условие задачи мы назовем условием сохранения отвесного положения свешивающейся части ВК ленты. Дело в том, что вследствие неизменяемости коромысла при отклонении последнего на угол $d\varphi$– на такой же угол повертывается и каждая касательная к краю лекала. Пусть в процессе такого поворота точка разлуки ленты с краем лекала переходит из положения В в бесконечно близкое положение $B'$. Если в начале этого процесса отвесное положение занимала касательная к краю ленты в точке $B$, то в конце его отвесное положение займет касательная к $B'$. Следовательно, угол смежности $d\tau$ между касательными к $B$ и $B'$ непременно должен быть равен углу $d\varphi$ поворота коромысла. Итак:

$$d\varphi = d\tau \ .$$

Интегрируя, находим

$$\varphi = \tau + A,$$

где $A$–постоянная. Но при $\varphi = 0$ мы имели бы $\tau = \frac{\pi}{2}$. Следовательно $A = -\frac{\pi}{2}$, и условие сохранения отвесного положения конца ленты напишется в окончательном виде так:

$$\varphi = \tau - \frac{\pi}{2} \ . \tag{5}$$

Впрочем, этот результат можно получить и непосредственно, из геометрических соображений.

Итак, предполагая ленту невесомой, мы свели задачу к системе трех уравнений (1),(4)-(5). Исключим из (1) $P$ и $\varphi$ при помощи (4) и (5). Тогда (1) сможем написать в виде

$$lG\sin(\tau-\frac{\pi}{2})=F\left(\frac{\tau-\frac{\pi}{2}-\varphi_0}{C}\right)\rho\cos(\theta-\tau+\frac{\pi}{2}).$$

Но по известным формулам:

$$\tau-\frac{\pi}{2}=-(\frac{\pi}{2}-arctg\frac{dy}{dx})=-arcctg\frac{dy}{dx}\;;$$

$$\sin(\tau-\frac{\pi}{2})=-\sin(\frac{\pi}{2}-\tau)=-\cos\tau=-\frac{dx}{ds},$$

$$\cos(\theta-\tau+\frac{\pi}{2})=-\sin(\theta-\tau)=-(\sin\theta\cos\tau-\cos\theta\sin\tau)=-(\frac{y}{\rho}\cdot\frac{dx}{ds}-\frac{x}{\rho}\frac{dy}{ds}).$$

Здесь $ds$–элемент дуги края лекала.

Подставляя все это в вышеописанное уравнение, получим

$$lG\frac{dx}{ds}=F\left(\frac{-\varphi_0-arcctg\frac{dy}{dx}}{C}\right)\cdot\rho\left(\frac{y}{\rho}\frac{dx}{ds}-\frac{x}{\rho}\frac{dy}{ds}\right),$$

откуда, сокращая на $\rho$ и на $\frac{dx}{ds}$, найдем

$$y=px+\frac{lG}{F\left(\frac{-\varphi_0-arcctg\,p}{C}\right)}\,, \qquad (6)$$

если обозначить, как это делают обычно, $\frac{dy}{dx}$ через $p$.

Дифференциальное уравнение 1-го порядка (6) принадлежит к легко интегрируемым уравнениям класса Clairot. Дифференцируем его по $x$. После сокращения по $p$ в левой и правой части, получим:

$$\frac{dp}{dx}\left(x-\frac{lGF'}{C(1+p^2)F^2}\right)=0\;. \qquad (7)$$

Этому же удовлетворим, во-первых, положив $\frac{dp}{dx}=0$. Но тогда $p = const$, а при этом условии (6) определяет в плоскости лекала семейство прямых, зависящее от одного параметра $p$. В таком предположении уравнение (6) и есть общий интеграл рассматриваемого уравнения Clairot.

Во-вторых, уравнению (7) удовлетворим, положив

$$x = \frac{lGF'}{C(1+p^2)F^2} \quad . \tag{8}$$

Тогда, согласно (6)

$$y = \frac{plGF'}{C(1+p^2)F^2} + \frac{lG}{F} \quad . \tag{9}$$

Исключая из (8) и (9) параметр $p$, получим такое соотношение между $x$ и $y$, которое не может быть получено из общего интеграла (6) ни при каком специальном значении параметра $p$ и которое является особым интегралом уравнения Clairot. В плоскости лекала *XOY* оно определит огибающую семейства прямых (6).

Уравнениям (8) и (9) можно придать более удобный для исследования вид, если вместо $p$ ввести в качестве параметра угол $\varphi$ отклонения коромысла, воспользовавшись формулами:

$$p = -ctg\varphi \ , \ \frac{1}{1+p^2} = \sin^2\varphi \ , \ \frac{p}{1+p^2} = -\sin\varphi\cos\varphi$$

Тогда

$$x = \frac{lGF'\left(\frac{\varphi-\varphi_0}{C}\right)\sin^2\varphi}{C\left[F\left(\frac{\varphi-\varphi_0}{C}\right)\right]^2}, \tag{A$_1$}$$

$$y = -\frac{lGF(\frac{\varphi-\varphi_0}{C})\sin\varphi\cos\varphi}{C\left[F(\frac{\varphi-\varphi_0}{C})\right]^2} + \frac{lG}{F\left(\frac{\varphi-\varphi_0}{C}\right)} \tag{A$_2$}$$

В уравнениях (A) $\varphi$ может меняться в пределах от $\varphi_1$ до $\varphi_2$, смысл которых разъяснен выше.

Формулы (А) и дают приближенное решение общей задачи о форме края лекала в предположении невесомости ленты. Именно, лекало должно быть вырезано по кривой (А), чтобы шкалу квадранта нужно было бы градуировать равномерно по измеряемой величине $z$, связанной с весом исследуемого тела зависимостью $P = F(z)$.

Решим наконец *вопрос о величине той ошибки*, которую мы делаем, пренебрегая весом ленты. Мы утверждаем, что при наложенных на функцию $F$ ограничениях ширину ленты всегда можно выбрать настолько малой, что весом ленты можно будет пренебречь, не выходя за пределы допустимой ошибки.

Пусть Q–вес ленты, $\sigma$–отношение плеча силы Q к плечу приложенного к точке В веса исследуемого тела Р. Перенесем силу Q в точку В. Тогда система сил, действующих на коромысло квадранта, сведется к следующим двум:

а) весу Р коромысла, приложенного к точке А

и

б) весу $P + \sigma Q$ исследуемого тела и ленты вместе, приложенного к точке В.

Заметим, что $\sigma$ есть число, не превосходящее некоторого конечного числа $\sigma^*$.

Если $z = f(P)$, то вместо действительного числового значения величины $z = f(P)$ из-за влияния веса ленты мы прочитаем на шкале прибора значение $\varsigma = f(P + \sigma Q)$. Абсолютная ошибка в отсчете окажется равной

$$\varsigma - z = f(P + \sigma Q) - f(P),$$

а относительная ошибка $\varepsilon$ в процентах будет тогда $\varepsilon = \dfrac{100[f(P + \sigma Q) - f(P)]}{f(P)}$.

Считая $\sigma Q$ величиной малой, разлагая разность, написанную в квадратных скобках по формуле Тейлора и ограничиваясь только первой степенью $\sigma Q$, получим

$$\varepsilon = \frac{100 f'(P)}{f(P)} \sigma Q.$$

По предположению, сделанному относительно функции $f(P)$, ее логарифмическая производная $f'(P) : f(P)$ не превосходит М, где М – конечное количество. Поэтому

$$\varepsilon \le 100 M \sigma Q.$$

Это неравенство усилится, если вместо $\sigma$ ввести то наибольшее числовое значение $\sigma^*$, которое достигает отношение $\sigma$ в пределах изменения угла $\varphi$. Следовательно, à fortiori

$$\varepsilon \le 100 M \sigma^* Q.$$

Пусть Е % есть наибольшая, допустимая на данном квадранте, относительная ошибка. Потребуем, чтобы

$$100M\sigma^{*}Q < E\,.$$

Тогда $$Q < \frac{E}{100M\sigma^{*}}\quad . \qquad (10)$$

Ширина ленты должна быть выбрана таким образом, чтобы неравенство (10) удовлетворялось. В таком случае ошибка в отсчете величины $z$ на данном квадранте не превзойдет $E\%$ и весом ленты свободно можно будет пренебречь.

Применим эти общие теоретические соображения к ***двум частным случаям***, практически важным для текстильного дела.

Вес Р мотка пряжи .Посмотрим, каким образом результаты общей теории, приведшей нас к формулам (А), могут быть приложены к тому случаю, когда измеряемая на квадранте величина $z$ есть вес Р мотка пряжи. Тогда

$$P = F(z) = F\left(\frac{\varphi-\varphi_0}{C}\right) = z = \frac{\varphi-\varphi_0}{C}$$

$$F\left(\frac{\varphi-\varphi_0}{C}\right) = 1$$

Пусть $\varphi_0 = 0$. В таком случае

$$P = \frac{\varphi}{C} \qquad (2)$$

Если $\varphi$ меняется в пределах от $\varphi_1$ от $\varphi_2$, то $P$ меняется от $P_1 = \frac{\varphi_1}{C}$ до $P_2 = \frac{\varphi_2}{C}$. В этих пределах должны соблюдаться условия, налагаемые общей теорией на функцию $F$ для того, чтобы весом ленты можно было пренебречь. А так как логарифмическая производная $f'(P) : f(P)$ в рассматриваемом случае принимает вид $\frac{1}{P}$ и достигает максимума $M = \frac{1}{P_1}$ при $P = P_1$, то *необходимо, чтобы $P_1$ было больше нуля*. Это условие будем предполагать соблюденным, причем заметим, что если $P_1$ мало, то стальная лента должна быть заменена тонкой, практически невесомой шелковой нитью.

Итак, имеем

$$F\left(\frac{\varphi-\varphi_0}{C}\right) = \frac{\varphi}{C}, \qquad F'\left(\frac{\varphi-\varphi_0}{C}\right) = 1 \quad .$$

Подставляя это в (А) и вводя величину

$$a = ClG \tag{11}$$

которая имеет размерность длины и может быть названа постоянной квадранта, получим

$$x = a\frac{\sin^2\varphi}{\varphi^2} \;,$$

$$y = a\frac{\varphi - \sin\varphi\cos\varphi}{\varphi^2} \;, \tag{B}$$

где $\varphi$ меняется в пределах от $\varphi_1$ до $\varphi_2$.

Уравнения (В) и определяют форму края лекала для квадранта с равномерной шкалой, измеряющей вес пряжи.

Из этих уравнений видно, что все лекала этого типа для различных значений постоянной **а** в одном и том же интервале $(\varphi_1, \varphi_2)$ между собой подобны и подобно расположены, имея центр подобия в начале координат. Поэтому по данному лекалу, построенному для какого-нибудь специального значения постоянной **а**, пантографическим методом сможем построить лекало и для любого иного **а** в тех же пределах углов отклонения.

Впрочем, можно пользоваться для измерения веса пряжи и каким-нибудь одним определенным коромыслом в данном интервале $(\varphi_1, \varphi_2)$, употребляя револьверное приспособление Матвеева-Петрова.

В самом деле, величину С можно менять, оставляя **l** и **а** постоянными (т.е. не меняя коромысла) и заставляя меняться лишь вес G коромысла. Для этого в центре тяжести А (рис.1) коромысла устраивается шпилька s, на которую можно надевать надлежащим образом подобранные цилиндрические грузики m в форме шайбочек. Эти грузики можно, например, выбрать так, что по безличной шкале прибора сможем прочитывать вес пряжи в любой системе или достигать любой чувствительности прибора.*

**Примечание редакции. Утверждение автора о грузиках верно только теоретически, при отсутствии трения.*

Ниже приведена таблица значений x и y, вычисленных согласно (В) для a=1 и углов $\varphi$ в пределах от 0 до $360^0$ через каждые $10^0$. Эта таблица может быть, конечно, продолжена как угодно далеко. Вид кривой, определяемой уравнениями (В), показан на рис.2; числа, поставленные около кривой, указывают соответствующие значения углов поворота коромысла.

| $\varphi^0$ | x | y | $\varphi^0$ | x | y | $\varphi^0$ | x | y |
|---|---|---|---|---|---|---|---|---|
| 0 | 1,000 | 0,000 | 130 | 0,114 | 0,536 | 250 | 0,045 | 0,213 |
| 10 | 0,990 | 0,115 | 140 | 0,069 | 0,492 | 260 | 0,046 | 0,212 |
| 20 | 0,960 | 0,227 | 150 | 0,036 | 0,444 | 270 | 0,047 | 0,212 |
| 30 | 0,912 | 0,331 | 160 | 0,015 | 0,389 | 280 | 0,041 | 0,211 |
| 40 | 0,848 | 0,422 | 170 | 0,003 | 0,357 | 290 | 0,035 | 0,210 |
| 50 | 0,771 | 0,499 | 180 | 0,000 | 0,318 | 300 | 0,027 | 0,207 |
| 60 | 0,684 | 0,568 | 190 | 0,009 | 0,286 | 310 | 0,020 | 0,202 |
| 70 | 0,591 | 0,603 | 200 | 0,009 | 0,260 | 320 | 0,013 | 0,195 |
| 80 | 0,498 | 0,629 | 210 | 0,019 | 0,241 | 330 | 0,008 | 0,187 |
| 90 | 0,405 | 0,637 | 220 | 0,028 | 0,227 | 340 | 0,003 | 0,178 |
| 100 | 0,319 | 0,629 | 230 | 0,036 | 0,219 | 350 | 0,001 | 0,169 |
| 110 | 0,240 | 0,608 | 240 | 0,043 | 0,214 | 360 | 0,000 | 0,159 |
| 120 | 0,171 | 0,576 | | | | | | |

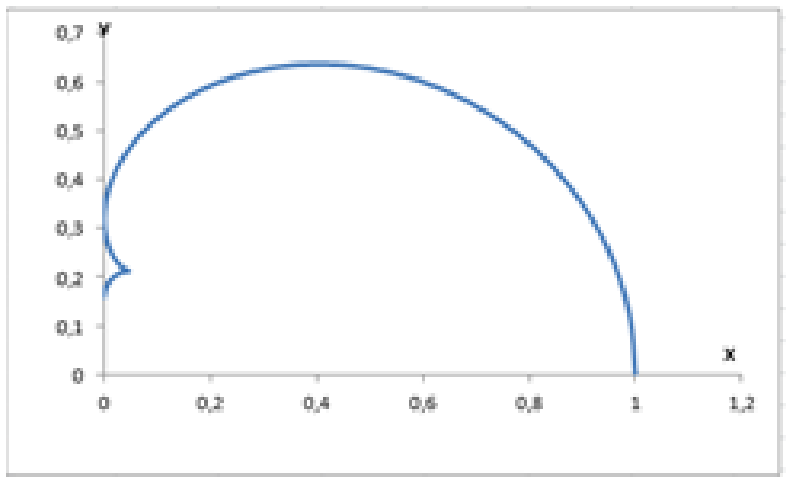

Рис.2

Расчет квадранта с равномерной шкалой для измерения веса пряжи можно вести по следующей схеме. Задаваясь пределами $P_1$ и $P_2$ измеряемого веса пряжи и чувствительностью С прибора, определяют предельные углы $\varphi_1$ и $\varphi_2$ отклонения коромысла по ф-лам

$$\varphi_1 = CP_1 \ , \ \varphi_2 = CP_2 .$$

Таким образом, устанавливается диапазон шкалы прибора и выделяется тот участок кривой (В), по которой должен быть вырезан край лекала. Подбирают затем **a**, l и G и строят коромысло. Шкалу прибора равномерно градуируют по весу в пределах между $P_1$ и $P_2$, а ширину ленты устанавливают на основании общих соображений при помощи ф-лы (10).

*Номер N мотка пряжи.* Применим формулы (А) к тому случаю, когда измеряемой на квадранте величиной $z$ является *номер пряжи*. Как известно, номер пряжи есть длина единицы ее веса. Поэтому между номером $N$ и весом $P$ данного мотка пряжи существует зависимость

$$N = \frac{\lambda}{P} \quad , \qquad (12)$$

где $\lambda$ – постоянная, имеющая размерность длины.

В таком случае соотношение $P = F(x) = F\left(\frac{\varphi - \varphi_0}{C}\right)$ общей теории принимает вид

$$P = F(z) = F\left(\frac{\varphi - \varphi_0}{C}\right) = \frac{\lambda}{z} = \frac{\lambda}{N} = \frac{C\lambda}{\varphi - \varphi_0} \ .$$

Так как с возрастанием веса $P$, а следовательно и угла $\varphi$ отклонения коромысла номер $N$ пряжи убывает, то С нужно считать отрицательным.

Если же разуметь под С его абсолютную величину, то соотношение между $P$ и $\varphi$ мы должны будем писать в виде

$$P = \frac{C\lambda}{\varphi_0 - \varphi} \ . \qquad (13)$$

Так мы и будем поступать в дальнейшем.

Для того, чтобы весом ленты можно было пренебречь согласно с требованиями общей теории необходимо, чтобы логарифмическая производная $f'(P) : f(P)$ не превосходила некоторого конечного числа $M$. В нашем случае

$$f(P) = N = \frac{\lambda}{P}; \ f'(P) = -\frac{\lambda}{P^2} \ ,$$

следовательно

$$\frac{f'(P)}{f(P)} = -\frac{\lambda}{P^2} : \frac{\lambda}{P} = -\frac{1}{P} \ .$$

Называя буквой $P_1$ наименьший измеряемый груз, мы видим, что условие конечности $M$ равноценно условию

$$P_1 > 0 \qquad (14)$$

Значит, желая измерять на квадранте очень высокие номера пряжи; соответствующие очень малым весам $P$ мотка пряжи данной длины $\lambda$, мы должны будем заменить стальную ленту практически невесомой нитью.

Условие (14) будем предполагать выполненным.

Из (13) получается, что при изменении веса мотка пряжи от $P_1$ до $P_2$ угол отклонения $\varphi$ коромысла будет изменяться в пределах от $\varphi_1$ до $\varphi_2$, причем

$$\varphi_1 = \varphi_0 - \frac{C\lambda}{P_1} \ , \ \varphi_2 = \varphi_0 - \frac{C\lambda}{P_2} \ .$$

Соответственно с этим номер $N$ пряжи, определяемый формулой (12), будет претерпевать изменение от наибольшего значения $N_1 = \frac{\lambda}{P_1}$ до наименьшего $N_2 = \frac{\lambda}{P_2}$ .

Обратимся теперь к вопросу о форме лекала. Считая С положительным, мы имеем

$$F\left(\frac{\varphi - \varphi_0}{C}\right) = \frac{C\lambda}{\varphi_0 - \varphi} \ ; \qquad F'\left(\frac{\varphi - \varphi_0}{C}\right) = -\frac{C^2\lambda}{(\varphi_0 - \varphi)^2} \ .$$

Подставляя это в формулы (А) общей теории, получим

$$x=\frac{lG}{C\lambda}\sin^2\sigma\ ,\ y=\frac{lG}{C\lambda}\left[(\varphi_0-\varphi)-\sin\varphi\cos\varphi\right]\ ,$$

или, вводя в качестве постоянной квадранта величину $b=\frac{lG}{C\lambda}$,имеющую размерность длины, сможем написанным выше уравнениям придать вид:

$$x=b\sin^2\varphi,\quad y=b\left[(\varphi_0-\varphi)-\sin\varphi\cos\varphi\right], \qquad (C)$$

где $\varphi$ меняется в установленных назначением $b$ способностью прибора в пределах от $\varphi_1$ до $\varphi_2$.

*Параметрические уравнения (С) и определяют форму края лекала,* предназначенного для *квадранта с равномерной шкалой, измеряющего номера N пряжи* в пределах от $N_1$до $N_2$, причем

$$N_2=\lambda:\frac{C\lambda}{\varphi_0-\varphi_2}=\frac{\varphi_0-\varphi_2}{C},$$

$$N_1=\lambda:\frac{C\lambda}{\varphi_0-\varphi_1}=\frac{\varphi_0-\varphi_1}{C}.$$

Величина $\varphi_0$ вполне произвольна, но выбор ее должен быть целесообразен. Положим, например, $\varphi_0=\pi$. Тогда формулам (С) можно придать вид, более приспособленный для вычислений, если в качестве параметра ввести другой угол $\psi$, связанный с $\varphi$ соотношением $\pi-\varphi=\psi$.

Тогда

$$x=b\sin^2\psi,\quad y=b\left[\psi+\sin\psi\cos\psi\right] \qquad (C')$$

причем угол $\psi$ меняется в пределах между $\psi_1$ и $\psi_2$, где

$$\psi_1=\pi-\varphi_1\quad,\ \psi_2=\pi-\varphi_2$$

Далее приведена таблица значений, вычисленных по ф-лам $(C')$ для b=1 и углов $\psi$, меняющихся через каждые $10^0$ в пределах от $0^0$ до $180^0$.

| $\psi$ | x | Y | $\psi$ | x | y | $\psi$ | x | y |
|---|---|---|---|---|---|---|---|---|
| $0^0$ | 0,000 | 0,000 | 60 | 0,750 | 1,480 | 130 | 0,587 | 1,777 |
| 10 | 0,030 | 0,346 | 70 | 0,883 | 1,543 | 135 | 0,500 | 1,856 |
| 20 | 0,017 | 0,670 | 80 | 0,970 | 1,567 | 140 | 0,413 | 1,951 |
| 30 | 0,250 | 0,957 | 90 | 1,000 | 1,571 | 150 | 0,250 | 2,185 |
| 40 | 0,413 | 1,191 | 100 | 0,970 | 1,574 | 160 | 0,117 | 2,471 |
| 45 | 0,50 | 1,285 | 110 | 0,883 | 1,559 | 170 | 0,0 0 | 2,796 |
| 50 | 0,587 | 1,365 | 120 | 0,750 | 1,661 | 180 | 0,000 | 3,142 |

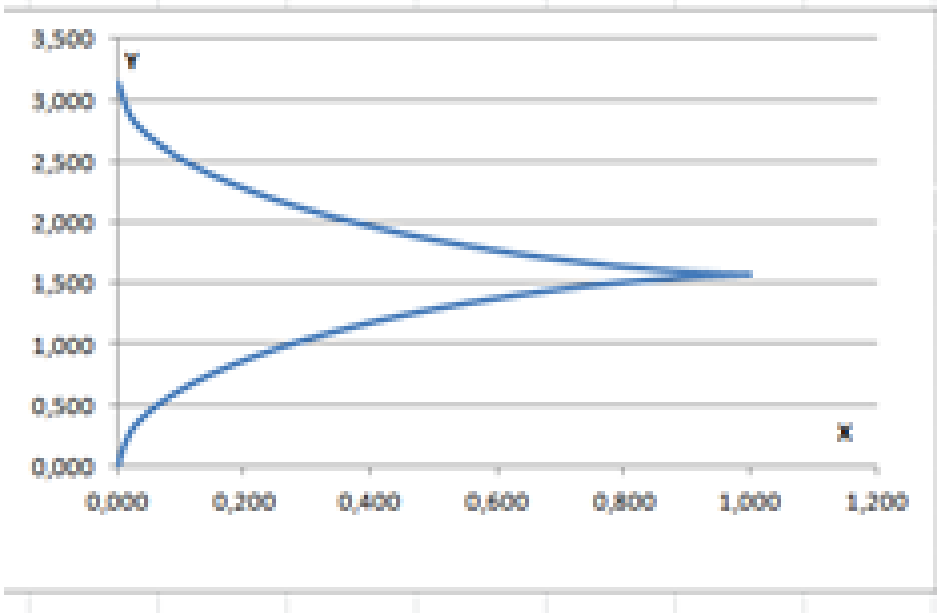

На рис. 3 показан вид кривой, определяемой уравнениями $(C')$, причем числа на рисунке указывают значения угла отклонения $\varphi$ коромысла для соответствующих точек кривой.

Из уравнений (С) края лекала видно, то и в этом случае можно прибегнуть к помощи пантографа для того, чтобы по данному лекалу для данного $\varphi_0$ построить в тех пределах углов отклонения другое, для любого иного значения постоянное b. В самом деле, из формул (С) вытекает, что все лекала для различных значений постоянной b и данного $\varphi_0$ между собой подобны, подобно расположены и имеют центр подобия в начале координат.

Изменение же угла $\varphi_0$ равноценно сдвигу всей кривой вдоль оси $O\gamma$.

Так же как это мы имели в первом примере, можно ограничиться, впрочем, одним типом коромысла квадранта для данного интервала $(\varphi_1, \varphi_2)$ отклонений, прибегнув к револьверному приспособлению Матвеева-Петрова. Для этого и в квадрантах этого типа в центре тяжести А (рис.1) коромысла устраивается шпилька s, на которую можно надевать по желанию, каждый из специально подобранных цилиндрических грузиков m. От этого будет меняться лишь вес коромысла G и чувствительность С прибора при постоянных b, $\lambda$ и l.

Ширина ленты прибора должна быть установлена таким образом, чтобы удовлетворялось требование (10) общей теории. Расчет остальных частей прибора можно вести так: задаваясь определенными С,$N_1$ и $N_2$, $\lambda$ и подбирая наиболее удобным образом другие постоянные ($\varphi_0$,b,C), определяют диапазон шкалы, которую и градуируют по номеру в установленных пределах. Из этих же данных находят и тот участок кривой (С), по которому должен быть вырезан край лекала.

Чтобы коромысло квадранта не опрокидывалось, находясь в нерабочем состоянии (а это иногда может случиться), прибегают к помощи специальных упоров, на которые и ложится стрелка – указатель, когда коромысло не нагружено.

Ограничиваясь рассмотрением двух приведенных примеров, замечу *в заключение*, что это не единственные случаи в технике, когда общая теория квадранта с постоянной чувствительностью, построенного на принципе лекала и ленты, могла бы быть полезной. Этот же принцип может быть применен в текстильном деле к динамометрам, определяющим разрывную длину нити [1], поскольку разрывная длина выражается эмпирически в функции веса нити. Далее в ювелирном деле, например, считают, что стоимость драгоценного камня пропорциональна квадрату его веса. Пользуясь общими формулами (А), можно было бы построить квадрант с равномерной шкалой, при помощи которого мы могли бы непосредственно измерять стоимость камня, вместо того, чтобы ее вычислять. И по всей вероятности существует в технике немало аналогичных случаев, когда принцип лекала и ленты, обычно игнорируемый, мог бы оказаться плодотворным.

Считаю долгом выразить здесь глубокую благодарность проф. И.И. Васильеву, помогавшему мне ценными замечаниями, и А.П. Матвееву, который, несмотря на болезнь, взял на себя составление чертежей и работу по просмотру рукописи.

1) *Применение электромагнита (соленоида и катушки) для динамометров впервые указано проф. И.И. Васильевым*

## От редакции

Работа т. Герасимова разбирает очень интересный и важный для нашего приборостроения вопрос. Тема проработана достаточно детально теоретически в применении к собственно весам. В отношении применения данного принципа к динамометрам вопрос не достаточно полно проработан, так как не разобраны вопросы о законе изменения действующей силы по времени условии сравнимости показаний различных динамометров.

Автором дано неверно утверждение, что динамометр определяет разрывную длину.

Динамометр определяет только силу, деформирующую образец. Разрывная длина находится вычислением на основе определенной динамометром силы и веса единицы длины продукта.

## **2.**ПРИНЦИП СООТВЕТСТВИЯ В ТЕОРИИ ЛИНЕЙНЫХ ОПЕРАТОРОВ

*(диплом по специальности «Прикладная математика»,*
*Научный руководитель И. И. Привалов)*

ь

*СОДЕРЖАНИЕ*

Введение

Основным понятием в теории линейных преобразований является понятие абстрактного гильбертова пространства, которое мы будем обозначать Ꞩ. Одной из возможных интерпретаций этого понятия является многообразие всех элементов «u» бесконечно большого числа измерений, причем сумма $\sum|u_k|^2$ квадратов модулей компонентов каждого из элементов **u** предполагается сходящейся.

Другой возможной конкретизацией Ꞩ, и притом очень существенной, является функциональное пространство из функций, обладающих интегрируемым квадратом модуля в некоторой основной области $\Omega$ .Для такого рода конкретных пространств можно доказать, как это делает J.von Neumann (Allgemeine Eigenwerttheorie Hermitsescher Functionaloperatoren, Math.Ann.102(1929), 49-131), что они изоморфны с Ꞩ. Это значит, что всякому элементу из рассматриваемого конкретного пространства может быть поставлен во взаимно-однозначное соответствие элемент из Ӈ, причем соотношениям между элементами из Ӈ отвечают соответственные соотношения между соответствующими элементами из рассматриваемого конкретного пространства.

Абстрактное гильбертово пространство Ӈ характеризуется вполне посредством приведенных ниже пяти свойств, которые обычно постулируются. Вследствие изоморфизма, теми же пятью постулируемыми положениями характеризуется и функциональное пространство из функций с интегрируемым квадратом. Таким образом, основная в функциональном анализе проблема линейных преобразований такого рода функций адекватна проблеме линейных преобразований комплексного гильбертова пространства.

*Абстрактное гильбертово пространство* Ӈ мы определяем наличием следующих пяти положений:

1.Для элементов $u_1, u_2, \ldots$ пространства Ӈ определены:

*сложение* $u_1 + u_2$; *вычитание* $u_1 - u_2$; *умножение на скаляр* a, au ,по тем принципам, которые установлены в обычной векторной алгебре. Это свойство Ӈ называют свойством линейности.

Из Ӈ можно выделить часть, имеющую любое (меньшее, чем у Ӈ) число измерений. Если эта часть такова, что, содержа в себе элементы $u_1, u_2, \ldots, u_n$, она содержит в себе и любую их линейную комбинацию $\sum_1^n c_k u_k$ ,($c_k$ – постоянные скаляры), то эту часть мы называема *линейным многообразием*. Линейные многообразия $\mathfrak{W}_1, \mathfrak{W}_2, \ldots, \mathfrak{W}_n$ *линейно независимы* между собой, при $u_1 \subset \mathfrak{W}_1$, $u_2 \subset \mathfrak{W}_2$, $u_n \subset \mathfrak{W}_n$, из $u_1+u_2+\ldots+u_n=0$ следует, что $u_1=u_2=\ldots=u_n=0$. Если $\mathfrak{W}'$ – *какая-нибудь часть* Ӈ , то множество всех $a_1u_1 + a_2u_2+\ldots + a_nu_n$ , где $u_1, u_2, \ldots, u_n \subset \mathfrak{W}'$, $a_1, a_2, \ldots, a_n$ – любые числа, представляет собой наименьшее линейное многообразие $\mathfrak{W}$, охватывающее $\mathfrak{W}'$.

2. Для элементов u,v… пространства Ӈ по правилам, применяемым в комплексной векторной алгебре, устанавливается понятие внутреннего произведения (u,v) со следующими свойствами:

1) (au, v)=a(u, v); 2)$(u_1+u_2,v) = (u_1,v) + (u_2,v)$; 3)$(u,v)=\overline{(v,u)}$ , где $\bar{x}$ обозначает комплексную сопряженную с $x$;

4) (u,u) $=|u|^2 \geq 0$ (причем знак равенства возможен тогда и только тогда, когда u=0.

Из этих постулируемых свойств внутреннего произведения вытекает далее:

5) (u, av)= $\bar{a}\overline{(v,u)} = \bar{a}(u,v) = (\bar{a}u, v)$;

6) (u,$v_1$+$v_2$)=$\overline{(v_1+v_2,u)} = \overline{(v_1,u)} + \overline{(v_2,u)} = (u,v_1)+(u,v_2)$;

Далее

($u_1$+$u_2$, $v_1$+$v_2$)=($u_1$,$v_1$+$v_2$)+ ($u_2$,$v_1$+$v_2$)= ($u_1$,$v_1$)+ ($u_1$,$v_2$) +($u_2$,$v_1$) + ($u_2$,$v_2$)= (u1,v1)+(u1,v2)+$\overline{(v_1,u_2)}$ +($u_2$,$v_2$).

Необходимо заметить, что внутреннее произведение (u,v) *непрерывно* относительно каждого из своих аргументов. В этом можно убедиться, положив по 2) $u_1$=u+Δu, $u_2$=-u и проследив переход к пределу при Δu ⟶0.

Затем имеет место *неравенство Cauchy:*

$$|(u,v)| \le |u|\cdot|v|.$$

Для доказательства возьмем:

$$\left(\frac{u+v}{2},\frac{u+v}{2}\right) = \frac{1}{4}\left[|u|^2+|v|^2+2\operatorname{Re}(u,v)\right], \qquad *)$$

$$\left(\frac{u-v}{2},\frac{u-v}{2}\right) = \frac{1}{4}\left[|u|^2+|v|^2-2\operatorname{Re}(u,v)\right],$$

откуда

$$\operatorname{Re}\,(u,v)=\left(\frac{u+v}{2},\frac{u+v}{2}\right)-\left(\frac{u-v}{2},\frac{u-v}{2}\right).$$

Так как разность двух неотрицательных количеств не больше, по абсолютной величине, их суммы, которая в данном случае равна

$$\left(\frac{u+v}{2},\frac{u+v}{2}\right)+\left(\frac{u-v}{2},\frac{u-v}{2}\right)=\frac{|u|^2+|v|^2}{2},$$

то имеем $\operatorname{Re}(u,v) \le \frac{1}{2}\left[|u|^2+|v|^2\right]$.

Пусть a>0. Заменим u,v соответственно на au, v. Левая часть последнего неравенства не изменится при изменении a; что касается правой, то она, очень большая при малых a, с возрастанием a уменьшается. Легко видеть, что $\frac{1}{2}\left[|au|^2+\left|\frac{v}{a}\right|^2\right]$ достигает минимума при $a^2=\frac{|v|}{|u|}$ и что этот минимум равен

*) *Здесь Re x обозначает действительную часть комплексного числа x.*

$|u|\cdot|v|$. Таким образом, Re(u,v) не больше этого минимума.

Заменив же u,v соответственно на ϑu, v ($|\vartheta|=1$), мы оставим правую часть постоянной при изменении ϑ. Левая часть станет равной Re ϑ(u,v). Пусть $(u,v)=|(u,v)|e^{i\psi}, \vartheta=e^{i\varphi}$. Тогда $\vartheta(u,v)=|(u,v)|e^{i(\varphi+\psi)}$ и максимум Re ϑ(u,v) есть $|(u,v)|$.

Сопоставление обоих выводов приводит к неравенству Cauchy:

Так как Re(u,v) $\mathrm{Re}(u,v)\le|u|\cdot|v|$, то

$$|u+v|^2=|u|^2+|v|^2+2\,\mathrm{Re}(u,v)\le(|u|+|v|)^2,$$

откуда имеем $|u+v|^2\le|u|+|v|$.

Два линейных многообразия $\mathfrak{M}_1$ ,$\mathfrak{M}_2$ взаимно ортогональны, если для любых их элементов $u_1\subset\mathfrak{M}_1$, $u_2\subset\mathfrak{M}_2$ имеем

$(u_1,u_2)=0$ ,

то есть, если два любых элемента этих многообразий взаимно ортогональны.

3. Ӈ есть пространство бесконечного множества измерений, Ӈ обладает любым конечным числом независимых элементов.

4. Пространство Ӈ полно, то есть, если последовательность элементов в Ӈ: $u_1,u_2,\ldots,u_n$ удовлетворяет условиям сходимости Cauchy (для данного $\varepsilon>0$ можно найти такое N=N(ε), что из ε>0 можно найти такое N=N(ε), что из m, n>N следует $|u_m-u_n|<\varepsilon$), то эта последовательность сходится к элементу u, который принадлежит к Ӈ.

Заметим, что когда $u^{(\nu)}\to u$, то $(u^{(\nu)},v)\to(u,v)$; $(v,u^{(\nu)})\to(v,u)$.

Если $u^{(\nu)}\to u, v^{(\mu)}\to v$, то $(u^{(\nu)},v^{(\mu)})\to(u,v)$.

5.Наконец, последнее свойство пространства Ӈ из числа тех постулируемых, которыми Ӈ вполне характеризуется, мы установим следующим образом: Пусть $\mathfrak{M}$-линейное многообразие, выбранное на Ӈ. Пусть $\mathfrak{M}$ не всюду плотно в Ӈ, то-есть, в Ӈ существует по крайней мере один такой элемент w, для которого нельзя найти последовательности элементов $u_1$, $u_2$, …,$u_n$, выбранной на $\mathfrak{M}$ и сходящейся к w. Тогда существует в Ӈ такой элемент, который, будучи отличен от нуля, ортогонален всем элементам $u\subset\mathfrak{M}$.

Из перечисленных пяти основных свойств пространства Ӈ *вытекает* следующее *важное обстоятельство*.

«***Пусть*** l(u)- функция от элемента u из некоторого линейного пространства 𝔚⊂ Ꞗ. Предположим, что l(u) ограничена в 𝔚 и удовлетворяет условию дистрибутивности:

$$l(c_1u_1 + c_2u_2) = c_1l(u_1) + c_2l(u_2),$$

где $u_1,u_2 \subset$ 𝔚, $c_1,c_2$-любые постоянные. Пусть, наконец, 𝔚 всюду плотно в Ꞗ.

***Тогда*** существует один и только один элемент v* из Ꞗ такой, что

l(u)=(u,v*).»-

_Действительно_, для того, чтобы найти v*, достаточно выбрать из Ꞗ такой элемент (v*), компоненты которого были бы количества, комплексно сопряженные с коэффициентами при компонентах u в линейной функции l(u). Очевидно, такой элемент v* всегда может быть найден. Сомнение могло бы возникнуть лишь относительно *единственности v**. Однако, если бы существовало два подобных элемента $v_1^*,v_2^*$ , то мы имели бы $(u, v_1^*)=( u, v_2^*)$, так что для всех элементов u из 𝔚 всюду плотного в Ꞗ ⊃ 𝔚 было бы тождественно

$$(u, v_1^* - v_2^*) \equiv 0,$$

откуда $v_1^* = v_2^*$, _что и надо было доказать._

Пусть теперь установлены правила, посредством которых элемент u преобразуется в элемент $v$, так что можно написать $v = Ru,$ где $R$ - символ того закона, при помощи которого осуществляется преобразование элемента u в элемент v. Этот символ мы называем _оператором_, часть Ꞗ, содержащую в себе те и только те u, для которых $R$ имеет смысл, - _областью определения_ (или областью существования) $R$, а элементы $v = Ru$ - _образами_ элементов u.

Множество всех образов $v$ мы называем _областью значений_ оператора $R$. Если 𝔚 есть область определения $R$ ,($u \subset$ 𝔚), а 𝔑-область его значений($v \subset$ 𝔑), то можно сказать, что $R$ _отображает_ 𝔚 на 𝔑.

Оператор $R$ называется _линейным_, если он удовлетворяет условию дистрибутивности:

$$R\sum c_k u_k = \sum c_k R u_k \ ,$$

где $c_k$ – любые постоянные, а R имеет смысл как для каждого из элементов $u_k$ ,так и для их линейной комбинации $\sum c_k u_k$ .

Внутреннее произведение $(Ru, v)$ /билинейная форма/, имеет смысл для всех u из области определения R, которую обозначим 𝔚 для всех $v \subset$ Ꞗ. Для данного $v$ это произведение представляет собой дистрибутивную функцию от u , если R – линейный оператор. При надлежащим образом выбранной $v$, $(Ru,v)$ есть *ограниченная* функция u для всех $u \subset \mathfrak{W}$. Мы будем рассматривать только такие операторы, которые определены на всюду плотных линейных многообразиях $\mathfrak{W}$. В таком случае, для данного $v$, для которого $(Ru,v)$ ограничена, можно найти один (и только один) элемент $v^* \subset \mathfrak{H}$ такой, что $(Ru,v)=(u,v^*)$.[x)]

Все те $v$, для которых это имеет место, образуют линейное многообразие $\mathfrak{W}$*. Если и $\mathfrak{W}$* всюду плотно в $\mathfrak{H}$, то можно найти единственный линейный оператор R* такой, что $v^*$ = R* $v$, следовательно, тогда $(Ru,v)=(u,R^*v)$.

Удовлетворяющий этому условию оператор R* мы называем, следуя Ф.Риссу, *сопряженным* с R. Аналогично, для R* можно определить сопряженный ему оператор R**. *Сопряженный к R* оператор R*** имеет смысл по крайней мере для всех тех u, для которых определен смысл R, вообще же R** является *продолжением* R.

Если R= R*, то R= R**. В этом случае оператор R мы называем *самосопряженным*. Если R – самосопряженен, то $(Ru,v)=(u,Rv)$ для всех $u,v \subset \mathfrak{W}$.В таком случае квадратичная форма $(Ru,u)=(u,Ru)$ может принимать только действительные значения, как это видно из того свойства внутреннего произведения, по которому $(a,b)=\overline{(b,a)}$.

Самосопряженный линейный оператор R мы будем называть *определенным положительно*, если положительно определена квадратичная форма $(Ru,u)$, ему соответствующая.

Пусть R – линейный оператор (не обязательно самосопряженный). Мы скажем, что R *ограничен на линейном многообразии* $\mathfrak{W} \subset \mathfrak{H}$, если существует такое не зависящее от $u \subset \mathfrak{W}$ число M>0, для которого $|Ru| \le M|u|$, $u \subset \mathfrak{W}$.

Наименьшее число $M_R$ из определяемых таким образом чисел M мы называем *гранью* оператора R.

Если R ограничен на всюду плотном в $\mathfrak{H}$ линейном многообразии $\mathfrak{W}$, то его продолжение R** определено и ограничено для всего $\mathfrak{H}$ (см. статью Riesz F. Sur la décomposition des opérations fonctionelles linéaires/Acta Scientiarum Mathematicarum,1928) . Если самосопряженный (ограниченный или неограниченный) оператор R таков, что для всех u, для которых R имеет смысл (Ru,u)=0, то Ru=0 и мы считаем

[x)] *Это условие может быть заменено и другим, обеспечивающим однозначность* $v^*$.

тогда R=0 . При этом существенно иметь ввиду, что область 𝔚, определения R предполагается всюду плотной в Ḩ.

Оператор R мы будем называть ограниченным, если он ограничен в Ḩ. Пусть $R_1, R_2, \ldots$ *–последовательность самосопряженных ограниченных операторов*. Мы скажем, что она *сходится к оператору R*, когда для всех $u \subset$ Ḩ: $R_n u \to Ru$ и когда все грани $M_{R_1}, M_{R_2}, \ldots$ ограничены в совокупности.

***«Если*** $R_n \to R$, то для всякого ограниченного самосопряженного оператора S: $R_n S \to RS; SR_n \to SR$. Если $R_n \to R, S_n \to S$ то $R_n S_n \to RS$.»–

*В самом деле*, так как $M_{R_n} M_{S_n} \geq M_{R_n S_n}$, то совокупность граней $M_{R_n S_n}$ ограничена. Кроме того $R_n S_n u = R_n S u - R_n (S - S_n) u$. Но $R_n S \to RS$, а относительно $R_n (S - S_n)$ легко доказать, что это выражение при $n \to \infty$ стремится к нулю. Наконец, заметим, что для самосопряженных операторов ограниченность адекватна наличию таких двух чисел m и M, для которых операторы ME-R и R-mE были бы положительно определены. *)

Предлагаемая работа посвящается возможностям продолжения того *принципа,* который составляет сущность метода изучения ограниченных операторов Ф. Рисса. Согласно этому принципу, *если R ограниченный самосопряженный оператор и (Ru,u) при* $|u| = 1$ *изменяется в пределах от m до M, то всякому полиному f(μ) с действительными коэффициентами для значений μ, заключенных между m и M, соответствует ограниченный самосопряженный оператор f(R), причем это соответствие обладает следующими свойствами:*

a/ оно положительного типа, т.е. если f(μ)>0 при $m \leq \mu \leq M$ , то и $(f(R)u, u) > 0$ при $|u| = 1$.

b/ оно мультипликативно: это значит, что произведению двух полиномов $f_1(\mu)$, $f_2(\mu)$ в рассматриваемом интервале соответствует произведение соответствующих операторов $f_1(R) f_2(R) = f_1 f_2(R) = f_2 f_1(R) = f_2(R) f_1(R)$.

c/ это соответствие дистрибутивно, т.е. если $f_1(\mu)$, $f_2(\mu)$-полиномы с действительными коэффициентами и $c_1, c_2$ –действительные числа, то полиному

*) *Здесь и в дальнейшем E обозначает оператор тождественного преобразования*

*от непрерывно меняющегося параметра μ, обладающий свойствами*

$$E_\mu = 0,\ \mu < m;\ E_\mu = E,\ \mu > M;\ \dot{E}_\mu E_\lambda u = E_\mu u, \mu \leq \lambda.$$

$c_1 f_1(\mu) + c_2 f_2(\mu) = f(\mu)$ соответствует самосопряженный ограниченный оператор f(R), причем $f(R) = c_1 f_1(R) + c_2 f_2(R)$.

Исходя из этого положения, Ф. Рисс распространяет его на все ограниченные функции f(μ), являющиеся пределами возрастающих (или убывающих) последовательностей непрерывных функций и получает для оператора f(R) выражение в виде интеграла Стилтьеса

$$f(R) = \int_{m}^{M} f(\mu) dE_{\mu},$$

где $E_{\mu}$ –особым образом сконструированный оператор, зависящий

Что касается неограниченных операторов, то для изучения их свойств Ф.Рисс предпочитает избрать иной метод.*) В статье, напечатанной в Acta Scient.(см. сноску), Ф.Рисс, путем рассмотрения некоторых специальных ограниченных операторов, доказывает возможность разложения оператора E тождественного преобразования на два: $E_-$ и $E_+$, $E = E_- + E_+$, каждый из которых коммутирует **) как с данным самосопряженным оператором R , так и со всяким коммутирующим с R , оператором. При этом $E_-R + E_+R = R$ и операторы $E_-R$, $E_+R$ соответственно отрицательно и положительно определены. Далее, рассматривается разложение единицы $E_{\lambda}$ для самосопряженного оператора $R - \lambda E$ ($\lambda$- действительно) и выводятся свойства операторов $E_{\lambda}$, которые позволяют построить выражение для оператора в виде интеграла Стилтьеса и установить необходимые и достаточные признаки его сходимости.

***Предлагаемая работа имеет целью*** показать, что принцип соответствия, с сохранением своих свойств в основном и лишь с небольшими изменениями в деталях, распространяется и на класс операторов неограниченных.

Она является, в сущности, *реставрацией того пути*, по которому Ф.Рисс шел первоначально и который был впоследствии заменен другим.

Мы начнем с изучения свойств некоторых ограниченных операторов, связанных с данным самосопряженным оператором, (ограниченным или неограниченным)/**§ 1**/. При этом исходной точкой будет уравнение, рассмотренное Ф.Риссом в цитированной статье из Acta Szeped.

*) *Acta Scient. ac litt.,1931, $\bar{V}$. Über die linearen Transform.d. kompl.Hilbertschen Raumen*

**) *Пусть S имеет смысл для всех u⊂ Ӈ; при u⊂ 𝔚⊂ Ӈ, где 𝔚 область определения R, SRu также имеет смысл. Если, при этом RSu имеет смысл, и если RSu=SRu,то R и S коммутируют между собой. Очевидно, тогда SR есть продолжение RS в том смысле, что SR определен по крайней мере для всех тех элементов u, для которых RS имеет смысл.*

Затем **в § 2** будет показано, каким образом можно установить соответствие между рациональными дробями (с комплексными корнями и полюсами) и операторами. Будут констатированы характерные свойства этого соответствия (его положительность, дистрибутивность, мультипликативность и пр.)

В §**3** мы покажем, что соответствие с сохранением всех основных черт может быть продолжено на класс непрерывных функций, а также и на такие разрывные функции, которые являются пределами последовательностей непрерывных.

Наконец, §**4** содержит в себе распространение соответствия на неограниченные функции и применение всего предыдущего построения выражения оператора в виде интеграла Стилтьеса , условия сходимости последнего и способ получения обратного оператора для R-μE.

**§ 1.***Некоторые специальные ограниченные операторы, связанные с данным самосопряженным оператором*

***Пусть*** R- линейный самосопряженный (ограниченный или неограниченный) оператор, $\lambda = \lambda' + i\lambda''$ – комплексное число. Рассмотрим уравнение

$$(R - \lambda E)u = v \tag{1}$$

В однородном уравнении, когда $v = 0$, имеется одно и только одно решение u=0.–

<u>*В самом деле*</u>, если u– решение уравнения $Ru - \lambda u = 0$, то должно быть

$$(Ru,u) - (\lambda u,u) = ((R - \lambda E)u,u) = (0,u) = 0\,,$$

откуда $(Ru,u) = \lambda(u,u)$.

Так как $(Ru,u),(u,u)$– действительны, то при комплексном λ последнее равенство приводит к заключению

$$(u,u) = |u|^2 = 0, \therefore u = 0,$$

<u>*что и утверждалось*</u>.

***Отсюда следует***, что ее неоднородное уравнение (1) при данном v имеет не больше одного решения.

<u>*Действительно,*</u> если бы (1) имело два различных решения $u_1$,$u_2$, так что $Ru_1 - \lambda u_1 = v, Ru_2 - \lambda u_2 = v,$

то вычитание привело бы нас к однородному уравнению

$R(u_1 - u_2) - \lambda(u_1 - u_2) = 0,$

которое имело бы ненулевое решение, что противоречит доказанному выше.

***Пусть*** теперь элемент u в уравнении (1) пробегает линейное многообразие 𝔚, где 𝔚 – область определения R. Тогда $v$ пробегает линейное многообразие, которое обозначим 𝔑, 𝔑⊂Ꞩ. Утверждаем, что <u>*𝔑 – всюду плотно в Ꞩ*</u>.–

<u>*В самом деле,*</u> если бы это было не так, то в Ꞩ существовал бы такой элемент w, для которого было бы:

$0 = (v, w) = ((R - \lambda E)u, w) = (Ru, w) - \lambda(u, w), \therefore (Ru, w) = (u, \bar{\lambda} w).$

Так как $\lambda(u, w)$ конечно, то $(Ru, w)$ также конечно, следовательно, w– элемент продолжения R.

По определению <u>*самосопряженного*</u> оператора R, w должен принадлежать к области определения $\mathfrak{W}$ и R должно иметь смысл для w, причем $(Ru, w) = (u, Rw)$. Итак, для всех $u \subset \mathfrak{W}$ имеем: $(u, Rw) = (u, \bar{\lambda} w)$.

Отсюда следует, что $Rw = \bar{\lambda} w.$ Тогда $(Rw, w) = \bar{\lambda}(w, w)$ - *чего не может быть* при комплексном λ вследствие действительности форм $(Rw, w)$ и $(w, w)$.

Затем *доказем, что* $\mathfrak{N}$ <u>*замкнуто.*</u>

*Действительно,* пусть имеем последовательность элементов $v^{(\nu)} \subset \mathfrak{N}$ такую, что $v^{(\nu)}$ сходится к элементу $v$ (может быть, $v \not\subset \mathfrak{N}$) (по доказанному, такая последовательность существует), причем $v^{(\nu)} = (R - \lambda E)u^{(\nu)}$.

*Надо доказать*:

1) $u^{(\nu)} \to u$;

2) $(R - \lambda E)u = v$ , так что $v \subset \mathfrak{N}$.–

<u>*Мы имеем*</u>

$$\left|v^{(\mu)} - v^{(\nu)}\right|^2 = \left|(R - \lambda' E)(u^{(\mu)} - u^{(\nu)}) - i\lambda''(u^{(\mu)} - u^{(\nu)})\right|^2 =$$

$$\left|(R - \lambda' E)(u^{(\mu)} - u^{(\nu)})\right|^2 + \lambda''^2\left|u^{(\mu)} - u^{(\nu)}\right|^2 \geq \lambda''^2\left|u^{(\mu)} - u^{(\nu)}\right|^2.$$

Значит, когда последовательность $v^{(\nu)}$ сходится, то и последовательность $u^{(\nu)}$ также сходится. Назовем предел последней $u$ , $u^{(\nu)} \to u$ .

Далее, для любого элемента $w$ мы должны иметь, вследствие ограниченности

$$(w, (R - \lambda' E)u^{(\nu)}) = (w, v^{(\nu)} + i\lambda'' u^{(\nu)}):$$

С другой стороны

$$(w, (R - \lambda' E)u^{(\nu)}) = ((R - \lambda' E)w, u^{(\nu)}) \to ((R - \lambda' E)w, u) = (w, (R - \lambda' E)u).$$

$$(w, (R - \lambda' E)u^{(\nu)}) = (w, v^{(\nu)} + i\lambda'' u^{(\nu)}) \to (w, v + i\lambda'' u).$$

На основании равенства левых частей для любого $v$ , должны быть равны и пределы справа; поэтому

$$(w, (R - \lambda' E)u) = (w, v + i\lambda'' u); \; (w, (R - \lambda E)u - v) = 0 \; ,$$

причем последнее должно соблюдаться тождественно для всех w. Но тогда

$$(R - \lambda E)u = v ,$$

*что и требовалось доказать.*

Но замкнутое линейное многообразие в Ӈ есть само Ӈ, поэтому, когда u пробегает область 𝔚 определения R, v в уравнении (1) пробегает все Ӈ.

Таким образом, уравнение (1) для каждого данного $v \subset$ Ӈ имеет решение u⊂ 𝔚 и мы видели, что это решение единственно. Обозначим это решение символом

$$u = P_\lambda v \qquad (2)$$

Обратный относительно $R - \lambda E$ оператор $P_\lambda$ определен, таким образом, для всех элементов гильбертова пространства и преобразует последнее в область 𝔚 определения R.

Кроме того, из неравенства

$$\left|v^{(\mu)} - v^{(\nu)}\right|^2 \geq \lambda''^2 \left|u^{(\mu)} - u^{(\nu)}\right|^2 \quad ,$$

которое мы уже имели раньше, следует, что оператор $P_\lambda$ ограничен.

Действительно, положим в написанном неравенстве

$$u^{(\nu)} = 0; u^{(\mu)} = u; v^{(\nu)} = 0; v^{(\mu)} = v \,,$$

оно примет вид:

$$|v|^2 \geq \lambda''^2 |u|^2 \,,$$

или, на основании (2):

$$|P_\lambda v| \leq \frac{1}{|\lambda''|} \cdot |v| ,$$

*откуда и следует утверждаемое*.

Наряду с $P_\lambda$ *рассмотрим* ограниченный, определенный для всего Ӈ, оператор $P_{\bar{\lambda}}$, обратный относительно $R - \bar{\lambda}E$. Докажем, что $P_\lambda$ и $P_{\bar{\lambda}}$ *взаимно сопряжены.*

*Действительно,* так как $(R - \lambda E)P_\lambda = E,\ (R - \bar{\lambda}E)P_{\bar{\lambda}} = E,$

то $\qquad (P_\lambda u, v) = (P_\lambda u, (R - \bar{\lambda}E)P_{\bar{\lambda}} v) = ((R - \lambda E)P_\lambda u, P_{\bar{\lambda}} v) = (u, P_{\bar{\lambda}} v),$

*что и требовалось доказать*.

Далее, рассмотрим два оператора $P_{\lambda_1}, P_{\lambda_2}$, соответствующие двум различным значениям комплексного параметра λ. Очевидно,

$$P_{\lambda_1} = P_{\lambda_2}(R - \lambda_2 E)P_{\lambda_1} ,$$

$$P_{\lambda_2} = P_{\lambda_1}(R - \lambda_1 E)P_{\lambda_2}.$$

Почленное вычитание дает:

$$P_{\lambda_1} - P_{\lambda_2} = (\lambda_1 - \lambda_2)P_{\lambda_2}P_{\lambda_1} \qquad (3)$$

А переставив значки 1 и 2, найдем

$$P_{\lambda_1} - P_{\lambda_2} = (\lambda_1 - \lambda_2)P_{\lambda_1}P_{\lambda_2}.$$

Сравнивая (3) с только что полученным результатом, видим, что ограниченные операторы $P_{\lambda_1}, P_{\lambda_2}$ для различных значений комплексного λ коммутируют между собой. В частности, $P_\lambda$ коммутирует с $P_{\bar{\lambda}}$.

Рассмотрим теперь оператор

$$Q_\lambda = P_\lambda P_{\bar{\lambda}} = P_{\bar{\lambda}} P_\lambda \; . \qquad (4)$$

Он имеет смысл для всех элементов гильбертова пространства и, как произведение ограниченных операторов, ограничен. Кроме того, $Q_\lambda$ самосопряжен, так как

$$(Q_\lambda u, v) = (P_{\bar{\lambda}} P_\lambda u, v) = (P_\lambda u, P_\lambda v) = (u, P_{\bar{\lambda}} P_\lambda v) = (u, Q_\lambda v).$$

Отсюда, между прочим, видно, что $Q_\lambda$ определен положительно, так как при $v = u$:

$$(Q_\lambda u, u) = (P_\lambda u, P_\lambda u) = |P_\lambda u|^2 \geq 0.$$

Наконец, легко видеть, что при $|u| = 1$ совокупность значений формы $(Q_\lambda u, u)$ заключена между 0 и $\lambda''^{-2}$.

Затем, пусть $\sigma = \sigma' + i\sigma''$, λ – два комплексных числа.

Оператор

$$S_{\sigma\lambda} = [R - \sigma E]P_\lambda \qquad (5)$$

имеет смысл для всех элементов u гильбертова пространства, так как посредством $P_\lambda$ последнее преобразуется в $\mathfrak{W}$, для элементов которого определены R и $R - \sigma E$.

Сопоставляя (5) с тождеством $E = (R - \lambda E)P_\lambda$,

получаем

$$S_{\sigma\lambda} = E + (\lambda - \sigma)P_\lambda \; . \qquad \mathbf{(6)}$$

А так как $P_\lambda$ – ограничен, то и $S_{\sigma\lambda}$ тоже ограничен.

Для двух пар комплексных чисел $\sigma_1, \lambda_1$; $\sigma_2, \lambda_2$ мы имеем

$$S_{\sigma_1\lambda_1} = E + (\lambda_1 - \sigma_1)P_{\lambda_1}$$

$$S_{\sigma_2\lambda_2} = E + (\lambda_2 - \sigma_2)P_{\lambda_2}$$

Произведение $S_{\sigma_1\lambda_1}S_{\sigma_2\lambda_2}$ есть оператор, имеющий смысл для всех элементов пространства Ꞑ и, очевидно, ограниченный. При этом, так как $P_{\lambda_1}, P_{\lambda_2}$ взаимно коммутируют, то и $S_{\sigma_1\lambda_1}, S_{\sigma_2\lambda_2}$ взаимно коммутируют. В частности,

$$S_{\sigma\lambda}S_{\bar{\sigma}\bar{\lambda}} = S_{\bar{\sigma}\bar{\lambda}}S_{\sigma\lambda}.$$

Далее, легко видеть, что *операторы* $S_{\sigma\lambda}, S_{\bar{\sigma}\bar{\lambda}}$ *взаимно сопряжены*;

для доказательства достаточно воспользоваться их выражениями через $P_\lambda, P_{\bar{\lambda}}$;

затем мы имеем:

$$(S_{\bar{\sigma}\bar{\lambda}}S_{\sigma\lambda}u, v) = ([E + (\bar{\lambda} - \bar{\sigma})P_{\bar{\lambda}}][E + (\lambda - \sigma)P_\lambda]u, v) =$$

$$= ([E + (\lambda - \sigma)P_\lambda]u, [E + (\lambda - \sigma)P_\lambda]v) =$$

$$= (u, [E + (\bar{\lambda} - \bar{\sigma})P_{\bar{\lambda}}][E + (\lambda - \sigma)P_\lambda]v) = (u, S_{\bar{\sigma}\bar{\lambda}}S_{\sigma\lambda}v).$$

Из этого рассуждения видно, что ограниченный оператор

$$T_{\sigma\lambda} = S_{\sigma\lambda}S_{\bar{\sigma}\bar{\lambda}} = S_{\bar{\sigma}\bar{\lambda}}S_{\sigma\lambda} \qquad \textbf{(7)}$$

*самосопряжен*. Если в предыдущей цепи равенств положить $v = u$, то станет очевидным, что $T_{\sigma\lambda}$ определен положительно, так как

$$(T_{\sigma\lambda}u, u) = (S_{\sigma\lambda}u, S_{\bar{\sigma}\bar{\lambda}}u) = |S_{\sigma\lambda}|^2 \geq 0$$

Пусть теперь имеем два ряда комплексных чисел $\sigma_1, \sigma_2, ..., \sigma_n; \lambda_1, \lambda_2, ..., \lambda_n$.

*Рассмотрим оператор*

$$T = T_{\sigma_1\lambda_1}T_{\sigma_2\lambda_2}...T_{\sigma_n\lambda_n} \qquad (8)$$

Как произведение ограниченных, самосопряженных и взаимно коммутирующих операторов, он ограничен и самосопряжен.

Далее, при помощи соотношения (6) легко показать, что

$$S_{\sigma_i\lambda_i}S_{\sigma_k\lambda_k} = (R - \sigma_i E)P_{\lambda_i}(R - \sigma_k E)P_{\lambda_k} = (R - \sigma_k E)P_{\lambda_k}(R - \sigma_i E)P_{\lambda_i} =$$

$$= (R - \sigma_k E)P_{\lambda_i}(R - \sigma_i E)P_{\lambda_k} = S_{\sigma_k\lambda_i}S_{\sigma_i\lambda_k}.$$

Отсюда следует, что в операторе T, представленном в виде

$$T=(R-\sigma_1 E)P_{\lambda_1}(R-\overline{\sigma_1}E)P_{\overline{\lambda}_1}....(R-\sigma_n E)P_{\lambda_n}(R-\overline{\sigma_n}E)P_{\overline{\lambda_n}} ,$$

можно менять местами:

1) все множители вида $R-\sigma E$,

2) все множители вида $P_\lambda$ между собой.

Наконец, оператор вида

$$H=T_{\sigma_1\lambda_1}T_{\sigma_2\lambda_2}....T_{\sigma_n\lambda_n}Q_{\lambda'}Q_{\lambda''}....Q_{\lambda^{(\nu)}}$$

где $Q_{\lambda^{(\nu)}}, T_{\sigma_j\lambda_j}$ определены посредством (4) и (7),) имеет смысл для всех элементов гильбертова пространства, ограничен и самосопряжен, причем в его развернутом выражении в виде произведения множителей типа $R-\sigma E, P_\lambda$ можно менять местами как множители первого, так и второго типа.

**§ 2.** *Соответствие между рациональными дробями и операторами*

Мы будем исходить из уравнений

$$P_{\lambda_1}P_{\lambda_2}=\frac{1}{\lambda_1-\lambda_2}P_{\lambda_1}+\frac{1}{\lambda_2-\lambda_1}P_{\lambda_2}, \quad S_{\sigma\lambda}=(R-\sigma E)P_\lambda=E+(\lambda-\sigma)P_\lambda ,$$

установленных в § 1 (уравнения (3) и (6)). Здесь $\lambda_1,\lambda_2,\sigma,\lambda-$любые комплексные числа.

Пусть $\sigma_1,\sigma_2,....\sigma_n$ и $\lambda_1,\lambda_2,....\lambda_n-$ две системы любых комплексных чисел. Оператор

$$S_{\sigma_1\lambda_1}S_{\sigma_2\lambda_2}....S_{\sigma_n\lambda_n}=(R-\sigma_1 E)P_{\lambda_1}(R-\sigma_2 E)P_{\lambda_2}....(R-\sigma_n E)P_{\lambda_n} , \qquad (1)$$

хотя и не самосопряжен, но, как произведение ограниченных операторов типа $(R-\sigma E)P_\lambda=S_{\sigma\lambda}$, ограничен и мы видели, что он допускает перестановку как множителей вида $R-\sigma_j E$ между собой, так и множителей вида $P_{\lambda_k}$ между собой.

Мы **<u>*утверждаем*</u>:**

*Оператор (1) может быть представлен в виде линейного агрегата из операторов* $P_{\lambda_k}$*:*

$$S_{\sigma_1\lambda_1}....S_{\sigma_n\lambda_n}=E+\sum_{k=1}^{n}c_k P_{\lambda_k} ,$$

*причем комплексные коэффициенты $c_k$ равны соответствующим коэффициентам в разложении рациональной дроби*

$$\frac{(\mu-\sigma_1)(\mu-\sigma_2)....(\mu-\sigma_n)}{(\mu-\lambda_1)(\mu-\lambda_2)....(\mu-\lambda_n)}=1+\sum_{k=1}^{n}\frac{c_k}{\mu-\lambda_k}.$$

*Действительно*, рассмотрим сначала тот случай, когда среди n комплексных чисел $\lambda_k$ нет равных между собой. Доказательство проведем по методу полной индукции.

Случай n=1 очевиден, так как с одной стороны

$$(R-\sigma_1 E)P_{\lambda_1}=1+(\lambda_1-\sigma_1)P_{\lambda_1},$$

а с другой

$$\frac{\mu-\sigma_1}{\mu-\lambda_1}=\frac{\mu-\lambda_1+(\lambda_1-\sigma_1)}{\mu-\lambda_1}=1+\frac{\lambda_1-\sigma_1}{\mu-\lambda_1},$$

так что в этом случае $c_1^1=\lambda_1-\sigma_1$.

Предположим теперь, что доказываемое положение верно для случая n множителей, так что имеем:

$$S_{\sigma_1\lambda_1}....S_{\sigma_n\lambda_n}=E+\sum_{k=1}^{n}c_k^{(n)}P_{\lambda_k}, \qquad \frac{(\mu-\sigma_1)....(\mu-\sigma_n)}{(\mu-\lambda_1)....(\mu-\lambda_n)}=1+\sum_{k=1}^{n}\frac{c_k^{(n)}}{\mu-\lambda_k}$$

с соответственно равными коэффициентами $c_k^{(n)}$. Докажем, что в таком случае это положение верно и для n+1 множителей, то есть, что в разложениях для

$$S_{\sigma_1\lambda_1}....S_{\sigma_n\lambda_n}S_{\sigma_{n+1}\lambda_{n+1}}=[E+\sum_{k=1}^{n}c_k^{(n)}P_{\lambda_k}][E+(\lambda_{n+1}-\sigma_{n+1})P_{\lambda_{n+1}}]$$

и для

$$\frac{(\mu-\sigma_1)....(\mu-\sigma_n)(\mu-\sigma_{n+1})}{(\mu-\lambda_1)....(\mu-\lambda_n)(\mu-\lambda_{n+1})}=[1+\sum_{k=1}^{n}\frac{c_k^{(n)}}{\mu-\lambda_k}][1+\frac{\lambda_{n+1}-\sigma_{n+1}}{\mu-\lambda_{n+1}}]$$

мы имеем соответственно равные коэффициенты $c_k^{(n+1)}$. Тогда теорема будет доказана.

Но:

$$S_{\sigma_1\lambda_1}....S_{\sigma_{n+1}\lambda_{n+1}}=E+\sum_{k=1}^{n}c_k^{(n)}P_{\lambda_k}+(\lambda_{n+1}-\sigma_{n+1})P_{\lambda_{n+1}}+\sum_{k=1}^{n}c_k^{(n)}(\lambda_{n+1}-\sigma_{n+1})P_{\lambda_k}P_{\lambda_{n+1}}.$$

Кроме того $P_{\lambda_k} P_{\lambda_{n+1}} = \frac{1}{\lambda_k - \lambda_{n+1}} P_{\lambda_k} + \frac{1}{\lambda_{n+1} - \lambda_k} P_{\lambda_{n+1}}$ .

Подставляя, получим

$$S_{\sigma_1\lambda_1} \ldots S_{\sigma_{n+1}\lambda_{n+1}} = E + \sum_{k=1}^{n} c_k^{(n)} \frac{\lambda_k - \sigma_{n+1}}{\lambda_k - \lambda_{n+1}} P_{\lambda_k} + [(\lambda_{n+1} - \sigma_{n+1}) + \sum_{k=1}^{n} c_k^{(n)} \frac{\lambda_{n+1} - \sigma_{n+1}}{\lambda_{n+1} - \lambda_k}] P_{\lambda_{n+1}} \quad ,$$

или:

$$S_{\sigma_1\lambda_1} \ldots S_{\sigma_{n+1}\lambda_{n+1}} = E + \sum_{k=1}^{n+1} c_k^{(n+1)} P_{\lambda_k} \quad ,$$

где

$$c_k^{(n+1)} = c_k^{(n)} \cdot \frac{\lambda_k - \sigma_{n+1}}{\lambda_k - \lambda_{k+1}} \text{ для } k = 1,2,\ldots,n \;, \quad c_{n+1}^{(n+1)} = (\lambda_{n+1} - \sigma_{n+1}) + \sum_{k=1}^{n} c_k^{(n)} \frac{\lambda_{n+1} - \sigma_{n+1}}{\lambda_{n+1} - \lambda_k} \;. \quad (3)$$

С другой стороны легко проверить, что

$$\frac{(\mu - \sigma_1)\ldots(\mu - \sigma_{n+1})}{(\mu - \lambda_1)\ldots(\mu - \lambda_{n+1})} = 1 + \sum_{k=1}^{n} c_k^{(n)} \frac{\lambda_k - \sigma_{n+1}}{\lambda_k - \lambda_{n+1}} \frac{1}{\mu - \lambda_k} + [(\lambda_{n+1} - \sigma_{n+1}) + \sum_{k=1}^{n} c_k^{(n)} \frac{\lambda_{n+1} - \sigma_{n+1}}{\lambda_{n+1} - \lambda_k}] \cdot \frac{1}{\mu - \lambda_{n+1}} \;.$$

Таким образом, коэффициенты разложения дроби и операторы совпадают для n+1 множителей, если они совпадают для случая n множителей. Так как доказано раньше совпадение коэффициентов для случая одного множителя, то оно доказано тем самым для любого их числа. Из рекуррентных формул (3) легко найти, что

$$c_k^{(n)} = \frac{(\lambda_k - \sigma_1)\ldots(\lambda_k - \sigma_n)}{(\lambda_k - \lambda_1)\ldots(\lambda_k - \lambda_n)} \quad , \quad k = 1,2,\ldots,n, \quad n \text{-любое.}$$

В *знаменателе отсутствует тот множитель* $\lambda_k - \lambda_i$ *, для которого* $i = k$ .

***Рассмотрим*** *теперь оператор*

$$S_{\sigma_1\lambda} S_{\sigma_2\lambda} \ldots S_{\sigma_k\lambda} = (R - \sigma_1 E) P_\lambda \ldots (R - \sigma_k E) P_\lambda$$

*и рациональную дробь*

$$\frac{(\mu - \sigma_1)\ldots(\mu - \sigma_k)}{(\mu - \lambda)^k} \;.$$

*Докажем,* что и в этом случае рассматриваемый *оператор разлагается единственным образом в сумму элементарных операторов*

$$S_{\sigma_1\lambda} \ldots S_{\sigma_k\lambda} = E + \sum_{j=1}^{k} c_j^{(k)} P_\lambda^j \quad ,$$

*причем коэффициенты этого разложения таковы же, что и коэффициенты разложения рассматриваемой дроби*

$$\frac{(\mu-\sigma_1)...(\mu-\sigma_k)}{(\mu-\lambda)^k}=1+\sum_{j=1}^{k}\frac{c_j^{(k)}}{(\mu-\lambda)^j} \ .$$

*Действительно, при* $k=1$ это очевидно, так как

$$S_{\sigma_1\lambda}=(R-\sigma_1 E)P_\lambda=E+(\lambda-\sigma_1)P_\lambda \ ,$$

в то время, как

$$\frac{\mu-\sigma_1}{\mu-\lambda}=\frac{\mu-\lambda+(\lambda-\sigma_1)}{\mu-\lambda}=1+(\lambda-\sigma_1)\frac{1}{\mu-\lambda} \ .$$

Далее, *при* $k=2$ находим:

$$S_{\sigma_1\lambda}S_{\sigma_2\lambda}=[E+(\lambda-\sigma_1)P_{\lambda_1}][E+(\lambda-\sigma_2)P_\lambda]=E+[(\lambda-\sigma_1)+(\lambda-\sigma_1)]P_\lambda+(\lambda-\sigma_1)(\lambda-\sigma_2)P_\lambda^2 \ ,$$

а с другой стороны

$$\frac{(\mu-\sigma_1)(\mu-\sigma_2)}{(\mu-\lambda)^2}=1+\frac{(\lambda-\sigma_1)+(\lambda-\sigma_2)}{\mu-\lambda}+\frac{(\lambda-\sigma_1)(\lambda-\sigma_2)}{(\mu-\lambda)^2} \ ,$$

следовательно,

$$c_1^{(2)}=(\lambda-\sigma_1)+(\lambda-\sigma_2) \ , \ c_2^{(2)}=(\lambda-\sigma_1)(\lambda-\sigma_2) \ .$$

Идя так далее, легко доказать утверждаемое положение для любого количества $k$ сомножителей.

Таким образом, теорема о разложении операторов вида $S_{\sigma_1\lambda_1}...S_{\sigma_n\lambda_n}$ в сумму элементарных сохраняет свою силу и для того случая, *когда в ряде комплексных чисел* $\lambda_1,\lambda_2,...,\lambda_n$ *имеются равные.* Этот результат можно было бы получить, перемножив операторы, соответствующие функциям лишь с простыми корнями знаменателя и применив уравнения (3) и (6).

Далее, заметим, что *оператор* $P_{\lambda_1}P_{\lambda_2}...P_{\lambda_m}$ , где все $\lambda_j$ -различны, также разлагается в линейный агрегат из $P_{\lambda_j}$ по формуле

$$P_{\lambda_1}P_{\lambda_2}...P_{\lambda_m}=\sum_{j=1}^{m}c_j^{(m)}P_{\lambda_j}$$

с теми же коэффициентами $c_j^{(m)}$, которые имеются в разложении

$$\frac{1}{(\mu-\lambda_1)...(\mu-\lambda_m)}=\sum_{j=1}^{m}c_j^m\frac{1}{\mu-\lambda_j} \ .$$

*Действительно,* сначала имеем

$$P_{\lambda_1}P_{\lambda_2}=\frac{1}{\lambda_1-\lambda_2}P_{\lambda_1}+\frac{1}{\lambda_2-\lambda_1}P_{\lambda_2}\ .$$

Затем:

$$P_{\lambda_1}P_{\lambda_2}P_{\lambda_3}=\frac{1}{\lambda_1-\lambda_2}P_{\lambda_1}P_{\lambda_3}+\frac{1}{\lambda_2-\lambda_1}P_{\lambda_2}P_{\lambda_3}=\frac{1}{\lambda_1-\lambda_2}[\frac{1}{\lambda_1-\lambda_3}P_{\lambda_1}+\frac{1}{\lambda_3-\lambda_1}P_{\lambda_3}]+\frac{1}{\lambda_2-\lambda_1}[\frac{1}{\lambda_2-\lambda_3}P_{\lambda_2}+\frac{1}{\lambda_3-\lambda_2}P_{\lambda_3}]=$$

$$=\frac{1}{(\lambda_1-\lambda_2)(\lambda_1-\lambda_3)}P_{\lambda_1}+\frac{1}{(\lambda_2-\lambda_1)(\lambda_2-\lambda_3)}P_{\lambda_2}+\frac{1}{(\lambda_3-\lambda_1)(\lambda_3-\lambda_2)}P_{\lambda_3}\ .$$

С другой стороны, имеются легко проверяемые тождества

$$\frac{1}{(\mu-\lambda_1)(\mu-\lambda_2)}=\frac{1}{\lambda_1-\lambda_2}\frac{1}{\mu-\lambda_1}+\frac{1}{\lambda_2-\lambda_1}\frac{1}{\mu-\lambda_2}\ ,$$

$$\frac{1}{(\mu-\lambda_1)(\mu-\lambda_2)(\mu-\lambda_3)}=\frac{1}{(\lambda_1-\lambda_2)(\lambda_1-\lambda_3)}\frac{1}{\mu-\lambda_1}+\frac{1}{(\lambda_2-\lambda_1)(\lambda_2-\lambda_3)}\frac{1}{\mu-\lambda_2}+\frac{1}{(\lambda_3-\lambda_1)(\lambda_3-\lambda_2)}\frac{1}{\mu-\lambda_3}$$

и так далее, откуда и видна справедливость утверждаемого положения для любого количества m множителей.
Теперь легко убедиться, что любой оператор вида

$$\begin{aligned}&(R-\sigma_1'E)P_{\lambda'}(R-\sigma_2'E)P_{\lambda'}.....(R-\sigma_{k_1}'E)P_{\lambda'}\\&(R-\sigma_1''E)P_{\lambda''}(R-\sigma_2''E)P_{\lambda''}.....(R-\sigma_{k_2}''E)P_{\lambda''}\\&\ldots\ldots\ldots\ldots\ldots\ldots\ldots\ldots\ldots\ldots\ldots\ldots\\&(R-\sigma_1^{(m)}E)P_{\lambda^{(m)}}(R-\sigma_2^{(m)}E)P_{\lambda^{(m)}}.....(R-\sigma_{k_m}^{(m)}E)P_{\lambda^{(m)}}P_{\lambda_1}^{m_1}P_{\lambda_2}^{m_2}...P_{\lambda_k}^{m_k}\end{aligned}\tag{1}$$

может быть представлен в виде разложения

$$\begin{aligned}&(c_1^{'}P_{\lambda'}+....+c_{k_1}^{'}P_{\lambda'}^{k_1})+(c_1^{''}P_{\lambda''}+...+c_{k_2}^{''}P_{\lambda''}^{k_2})+...+(c_1^{(m)}P_{\lambda^{(m)}}+...+c_{k_m}^{(m)}P_{\lambda_m}^{k_m})+\\&(a_1^{'}P_{\lambda_1}+...+a_{m_1}^{'}P_{\lambda_1}^{m_1})+...+(a_1^{(k)}P_{\lambda_k}+...a_{m_k}^{(k)}P_{\lambda_k}^{m_k})\end{aligned}\ ,\tag{2}$$

коэффициенты которого идентичны с соответственными коэффициентами разложения рациональной дроби

$$\frac{(\mu-\sigma_1^{'})...(\mu-\sigma_{k_1}^{'})}{(\mu-\lambda^{'})^{k_1}}...\frac{(\mu-\sigma_1^{(m)})...(\mu-\sigma_{k_m}^{(m)})}{(\mu-\lambda^{(m)})^{k_m}}\frac{1}{(\mu-\lambda_1)^{m_1}}...\frac{1}{(\mu-\lambda_k)^{m_k}}$$

на элементарные. Здесь $R-$ любой самосопряженный оператор, $\lambda,\sigma-$любые комплексные числа.

Рассмотренные до сих пор операторы типа (1), вообще говоря, не были самосопряженными. Перейдем теперь к самосопряженным операторам вида (1), для которых доказанная выше теорема о разложении остается в силе.

Всякой рациональной функции

$$H(\mu)=\frac{f(\mu)}{F(\mu)}=\frac{\prod_{j=1}^{m}(\mu-\sigma_j)(\mu-\overline{\sigma_j})}{\prod_{k=1}^{n}(\mu-\lambda_k)(\mu-\overline{\lambda_k})},\qquad m\le n \tag{3}$$

с комплексными нулями и полюсами мы ставим в соответствие оператор $H(R)$:

$$H(R)=(R-\sigma_1E)P_{\lambda_1}(R-\overline{\sigma_1}E)P_{\overline{\lambda_1}}...(R-\sigma_mE)P_{\lambda_m}(R-\overline{\sigma_m}E)P_{\overline{\lambda_m}}P_{\lambda_{m+1}}P_{\overline{\lambda_{m+1}}}....P_{\lambda_n}P_{\overline{\lambda_n}}$$

Так как в выражении справа можно менять как порядок множителей вида $R-\sigma E$, так и порядок множителей вида $P_\lambda$, что соответствует изменению порядка множителей в числителе и знаменателе дроби $H(\mu)$, не меняющей при этом своего значения, то *способ установления соответствия* *однозначен.*

Затем, *легко убедиться*, что таким образом установленное соответствие мультипликативно: произведение двух рациональных функций типа (3) есть функция того же типа и ей соответствует оператор, равный произведению соответствующих операторов.

Далее, мы *покажем,* что это соответствие положительного типа, то есть, каждой не отрицательной для всех действительных $\mu$ функции $H(\mu)$ соответствует (один и только один) самосопряженный ограниченный оператор $H(R)$, и притом определенный положительно.

*Для доказательства* этого положения мы предварительно докажем следующую *теорему*, заимствованную нами из статьи Ф.Рисса.

*Произведение самосопряженных положительно определенных и взаимно-коммутирующих ограниченных операторов есть оператор самосопряженный и определенный положительно.–*

*Действительно,* пусть имеем два ограниченных, самосопряженных и положительно определенных оператора $M, N$, причем $MN=NM$.

Не нарушая общности рассуждений, можно считать, что они имеют одну и ту же грань, равную единице. Определим операторы $M_k$ рекуррентной формулой

$$M_k=M;\ M_{k+1}=M_k-M_k^2.$$

Очевидно, все $M_k, E-M_k$ определены положительно, как это следует путем перехода от k к k+1 из соотношений

$$M_{k+1} = M_k^2(E-M_k) + M_k(E-M_k)^2 ,$$

$$E - M_{k+1} = (E-M_k) + M_k^2 .$$

Далее, имеем

$$M = M_1^2 + M_2^2 + ... + M_n^2 + M_{n+1} ,$$

следовательно

$$(M_1^2 u,u) + (M_2^2 u,u) + ... + (M_n^2 u,u) \le (Mu,u) .$$

Значит, ряд $\sum_k (M_k^2 u,u) = \sum_k |M_k u|^2$ сходится, а тогда $|M_k u| \to 0$ при $k \to \infty$.

Поэтому можно утверждать, что $M = \sum M_k^2$ .

Аналогично $N = \sum N_i^2$ . Тогда оператор $MN$ может быть представлен как сумма членов вида $M_k^2 N_i^2$ .Каждый из этих членов определен положительно, так как

$$(M_k^2 N_i^2 u,u) = (M_k N_i u, M_k N_i u) \ge 0 .$$

Тогда и сумма таких операторов также определена положительна, <u>*чем и доказывается вторая часть теоремы. Первая ее часть очевидна.*</u>

Методом полной индукции *легко распространить* доказанное положение *на случай любого числа операторов.*

В частности, все операторы вида $S_{\sigma_j \lambda_k} S_{\overline{\sigma_j}\,\overline{\lambda_k}}$, $P_{\lambda_k} P_{\overline{\lambda_k}}$ ограничены, самосопряжены, взаимно коммутируют и определены положительно. Тогда, по доказанному, определено положительно и их произведение, которое есть не что иное, как $H(R)$, соответствующий неотрицательной функции $H(\mu)$. Этим самым <u>установлен положительный характер соответствия.</u>

***Из*** *положительности соответствия* вытекает, что *значения квадратичной формы* $(H(R)u,u) = (u,H(R)u)$ *при* $|u| = 1$*заключены между наименьшим* $m$ *и наибольшим М значением функции* $H(\mu)$ *при всех действительных* $\mu$ .–

*В самом деле,* положительным рациональным функциям $M - H(\mu)$ и $H(\mu) - m$ ответствуют положительно определенные операторы

$$ME - H(R) \ , \ H(R) - mE \ ,$$

т.е. такие, для которых

$$((ME - H(R))u,u) \geq 0, \quad ((H(R) - mE)u,u) \geq 0,$$

откуда при $|u| = 1$:

$$M \geq (H(R)u,u) \geq m,$$

*что и требовалось доказать.*

Наконец, *покажем,* что установленное соответствие между рациональными дробями и операторами дистрибутивно, то есть, что линейной комбинации (с действительными коэффициентами) из функций типа (3) соответствует оператор, равный такой же линейной комбинации из соответствующих операторов.

*Это положение становится очевидным*, если заметить, что в то время, как функция $H(\mu)$ может быть разложена в сумму элементарных дробей:

$$H(\mu) = 1 + \sum_{j,k} \left[ \frac{c_j^{(k)}}{(\mu - \lambda_j)^k} + \frac{\overline{c_j^{(k)}}}{(\mu - \overline{\lambda_j})^k} \right],$$

оператор $H(R)$, соответствующий этой функции, по доказанному выше, представляется в виде суммы:

$$H(R) = E + \sum_{j,k} \left[ c_j^{(k)} P_{\lambda_j}^k + \overline{c_j^{(k)}} P_{\lambda_j}^k \right].$$

**§ 3.** *Распространение соответствия на класс функций, более широкий*

В предыдущем параграфе было показано, каким образом можно установить соответствие между рациональными дробями и операторами. Там же было показано, то такое соответствие, будучи однозначно определенным, есть мультипликативное и дистрибутивное соответствие положительного типа. Можно пойти дальше и показать, что это соответствие с сохранением своих свойств распространяется, во-первых, на класс всех ограниченных функций, равномерно непрерывных для всех действительных значений аргумента (существенно заметить при этом, речь идет лишь о таких функциях, которые непрерывны в «замкнутом» интервале $(-\infty, +\infty)$ и, следовательно, равномерно непрерывны для всех значений аргумента).

В самом деле, такие функции могут быть с любой степенью точности равномерно аппроксимированы последовательностью ограниченных рациональных дробей того вида, который был изучен в §2.

Пусть имеем последовательность $H^{(\nu)}$ ограниченных рациональных дробей, равномерно сходящуюся к ограниченной и непрерывной для всех (считая и бесконечно удаленные точки) действительных значений $u$ функции $f(\mu)$. Это значит, что для всякого, как угодно малого $\varepsilon > 0$ существует такое целое $n > 0$, не зависящее от $\mu$,
что

$$\left|H^{(\nu_2)}(\mu) - H^{(\nu_1)}(\mu)\right| \leq \varepsilon, \text{ если } \nu_1, \nu_2 \geq \nu.$$

Но, согласно §2, функциям $H^{(\nu_2)}, H^{(\nu_1)}$ соответствуют операторы $H^{(\nu_2)}(R), H^{(\nu_1)}(R)$ и квадратичные формы $(H^{(\nu_2)}(R)u,u), (H^{(\nu_1)}(R)u,u)$,
причем для всех $u$ должно быть:

$$\left|(H^{(\nu_2)}(R)u,u) - (H^{(\nu_1)}(R)u,u)\right| \leq \varepsilon E(u,u).$$

Это значит, что последовательность квадратичных форм равномерно сходится. В таком случае, равномерно сходится и последовательность операторов. Предел этой последовательности мы считаем оператором $f(R)$, соответствующим функции $f(\mu)$ и этот предел не зависит от выбора аппроксимирующих функций.

Таким образом продолженное соответствие сохраняет свою мультипликативность и дистрибутивность.

В самом деле, если

$$H_1^{(\nu_1)}(\mu) \to f_1(\mu), \ H_2^{(\nu_2)}(\mu) \to f_2(\mu),$$

то $f_1(\mu) = H_1^{(\nu_1)}(\mu) + \varepsilon_1, \ \ f_2(\mu) = H_2^{(\nu_2)}(\mu) + \varepsilon_2,$

где $\varepsilon_1, \varepsilon_2 > 0$ и не зависят от $\mu$, $\nu_1, \nu_2$–достаточно большие целые числа.

Функциям $f_1, f_2; H_1^{(\nu_1)}, H_2^{(\nu_2)}$ соответствуют квадратичные формы $(f_1(R)u,u); (f_2(R)u,u);$ $(H_1^{(\nu_1)}(R)u,u); (H_2^{(\nu_2)}(R)u,u)$ , причем для форм последних двух типов установлена в §2 мультипликативность и дистрибутивность соответствия. Тогда из написанных равенств путем надлежащих действий и переходом к пределу при $\nu_1, \nu_2 \to \infty$, $\varepsilon_1, \varepsilon_2 \to 0$, мы констатируем наличие тех же свойств и у операторов $f_1(R), \ f_2(R)$.

Наконец, из соотношения $f(R) = H^{(\nu)}(R) + \varepsilon E$ и из положительности соответствия для операторов вида $H^{(\nu)}(R)$, вследствие произвола в выборе $\varepsilon$ при надлежаще подобранном $\nu$ вытекает, что и для всех ограниченных, непрерывных равномерно для всех действительных $\mu$ функций вида $f(\mu)$ установленное соответствие сохраняет положительный характер.

Как следствие, отсюда, в свою очередь получается, что значения квадратичной формы $(f(R)u,u)$ при $|u|=1$остаются заключенными между наименьшим и наибольшим из значений $f(\mu)$ при всех действительных $\mu$.

Установив соответствие для класса ограниченных, равномерно непрерывных функций и констатировав его свойства для этих функций, *можно распространить соответствие с сохранением тех же свойств на класс функций, еще более широкий.* Именно, речь идет о *всех ограниченных функциях (может быть, разрывных), которые являются пределами монотонных последовательностей равномерно непрерывных функций.*

*Пусть* $f(\mu)$ – такая функция. Будем считать $f(\mu)$ положительной.

Рассмотрим последовательность $f^{(\nu)}(\mu)$непрерывных неотрицательных функций, сходящуюся к $f(\mu)$. Каждой функции $f^{(\nu)}(\mu)$из этой последовательности соответствует ограниченный самосопряженный оператор (и форма) $f^{(\nu)}(R)$, определенные положительно. При этом все формы $f^{(\nu)}(R)$ ограничены в совокупности и образуют не убывающую последовательность; то же относится и к формам вида $f^{(\nu)^2}(R)$.

Пусть $M>0$есть верхняя точная граница функции $f(\mu)$. Тогда $\left|f^{(\nu)}(R)u\right|\le M$ при $|u|=1$для всех $\nu$ и $\left|f^{(\nu)^2}(R)u\right|\le M^2$. Отсюда вытекает, что

$$((f^{(\nu)}(R)-f^{(\mu)}(R))^2u,u)=\left|f^{(\nu)}(R)u-f^{(\mu)}(R)u\right|^2 \text{ при } |u|=1$$

и при достаточно больших $\nu,\mu$ может быть сделано по модулю меньше любой, как угодно малой, положительной величины.

В таком случае последовательность $f^{(\nu)}(R)$, будучи ограниченной, сходится к некоторому пределу (сильная сходимость). Этот предел мы ставим в соответствие с функцией $f(\mu)$ и обозначаем его $f(R)$.

Мы *докажем независимость этого предела от выбора последовательности аппроксимирующих функций.* Следуя Ф.Риссу(см. ссылку выше, часть V), мы *докажем сначала более общее положение.*

*Пусть* $H_1^{(\nu)}(\mu),H_2^{(\nu)}(\mu)$ – две возрастающие последовательности ограниченных рациональных функций, равномерно сходящихся к функциям, ограниченным и непрерывным для всех $\mu$ (в смысле, о котором говорилось выше). Пусть $f_1(\mu)\le f_2(\mu)$ для всех $\mu$. Докажем, что *тогда* $(f_1(R)u,u)\le(f_2(R)u,u)$ *для всех* $u$.

*Действительно,* можно зафиксировать такое $m$ и такое $\varepsilon > 0$, что для всех достаточно больших $n$ имеем $f_2^{(n)}(\mu) > f_1^{(n)}(\mu) - \varepsilon$, потому что если бы это было не так, то те значения $\mu$, для которых $f_1^{(n)}(\mu) - \varepsilon \geq f_2^{(1)}(\mu), \ldots, f_1^{(m)}(\mu) - \varepsilon \geq f_2^{(2)}(\mu), \ldots$ образовывали бы последовательность замкнутых множеств, каждое из которых заключало бы в себе все следующие.

Тогда непременно нашлась бы точка (по крайней мере, одна) $\mu^*$, такая, которая принадлежала бы всем этим множествам. В этой точке мы имели бы

$$f_1\,(\mu) > f_1^{(m)}(\mu^*) - \varepsilon \geq f_2\,(\mu^*)$$

вопреки предположению $f_1 \leq f_2$ для всех $\mu$. Итак, для достаточно большого $n$ при выбранных $m, \varepsilon$:

$$f_1^{(m)}(\mu) - \varepsilon < f_2^{(n)}(\mu).$$

Но тогда, на основании соответствия для рациональных дробей, имеем

$$(f_1^{(m)}(R)u,u) - \varepsilon E(u,u) \leq (f_2^{(n)}(R)u,u).$$

Вследствие произвола в выборе $m, \varepsilon$, имеем тогда

$$(f_1\,(R)u,u) \leq (f_2^{(n)}(R)u,u) \leq (f_2(R)u,u),$$

*что и хотели доказать.*

В частности, если $f_1(\mu) = f_2(\mu)$,

то и $$f_1(R) = f_2(R),$$

чем и доказывается *независимость предельного оператора* $f(R)$ *от выбора приближающих дробей*.

Далее, очевиден положительный характер соответствия. А из положительности соответствия, в свою очередь, вытекает, что *значения формы* $(f(R)u,u)$ *при* $|u| = 1$ *заключены между точными, верхней и нижней, границами* $f(\mu)$.

Затем, *можно показать*, что *соответствие распространяется и на разности* положительных функций только что рассмотренного вида.

*Действительно,* если $f_1(\mu), f_2(\mu); f_1^{'}(\mu), f_2^{'}(\mu) -$ функции, принадлежащие к классу рассмотренных, то и их суммы

$$f_1(\mu) + f_1^{'}(\mu),\ f_2(\mu) + f_2^{'}(\mu)$$

также принадлежат к классу рассмотренных.

Вследствие очевидной аддитивности соответствия, этим функциям соответствуют операторы

$$f_1(R)+f_1^{'}(R),\ f_2(R)+f_2^{'}(R),$$

так что если имеет место равенство

$$f_1(\mu)+f_1^{'}(\mu)=f_2(\mu)+f_2^{'}(\mu),$$

то имеет место также

$$f_1(R)+f_1^{'}(R)=f_2(R)+f_2^{'}(R).$$

Но, в таком случае функции

$$f_1(\mu)-f_2(\mu)=f_2^{'}(\mu)-f_1^{'}(\mu)$$

должен быть поставлен в соответствие оператор

$$f_1(R)-f_2(R)=f_2^{'}(R)-f_1^{'}(R),$$

*что и требовалось доказать*.

Продолженное таким образом *соответствие*, очевидно, *дистрибутивно.*

*Докажем,* что оно сохраняет *мультипликативность* для этого *расширенного класса функций.*

*В самом деле,* пусть имеем две возрастающие [x)] ограниченные последовательности

В силу мультипликативности соответствия для класса непрерывных функций равномерно непрерывных функций $f^{(\nu)}(\mu), g^{(\nu)}(\mu),$ соответственно сходящиеся к ограниченным функциям $f(\mu), g(\mu)$. Этим функциям соответствуют операторы, причем для последних мы имеем

$$f^{(\nu)}(R)\to f(R), g^{(\nu)}(R)\to g(R).$$

Рассмотрим последовательность произведений из функций

$$f^{(\nu)}(\mu)g^{(\nu)}(\mu)=f^{(\nu)}g^{(\nu)}(\mu).$$

Эта последовательность не убывает и ограничена, следовательно, сходится к пределу, который обозначаем через $fg(\mu)$.

[x)] *Очевидно, это требование не нарушает общности рассуждений*

Для соответствующих операторов имеем:

$$f^{(\nu)}g^{(\nu)}(R) \to fg(R).$$

$$f^{(\nu)}g^{(\nu)}(R) = f^{(\nu)}(R)g^{(\nu)}(R),$$

так что $f^{(\nu)}(R)g^{(\nu)}(R) \to fg(R)$ .

С другой стороны, последовательность $f^{(\nu)}(R)g^{(\nu)}(R)$ произведений ограниченных в совокупности операторов удовлетворяет условию сильной сходимости, о которой была речь во «Введении». По доказанной там теореме имеем:

$$f^{(\nu)}(R)g^{(\nu)}(R) \to f(R)g(R).$$

Сопоставление обоих соотношений приводит к выводу $fg(R) = f(R)g(R)$,
*что и требовалось доказать*.
Из мультипликативности же соответствия следует коммутативность всех операторов для класса ограниченных функций, являющихся пределами последовательностей функций, равномерно непрерывных для всех, как угодно больших по абсолютной величине, значений $\mu$.

*В заключение* этого параграфа *рассмотрим примеры*, играющие значительную роль в дальнейшем.

*Пример 1.* Пусть $\mu$–действительное число. Рассмотрим функцию $\Psi_\mu(\lambda)$, определенную для всех действительных $\lambda$ следующим образом:

$$\Psi_\mu(\lambda) = 1 \quad \text{для } \lambda < \mu;$$

$$\Psi_\mu(\lambda) = 0 \quad \text{для } \lambda \ge \mu.$$

Эта формула является пределом неубывающей последовательности неотрицательных непрерывных функций (например, функций $\Psi_\mu^{(\nu)}(\lambda)$, определенных таким образом

$$\Psi_\mu^{(\nu)}(\lambda) = 1 \qquad \text{при } \lambda \le \mu - \varepsilon^{(\nu)};$$

$$\Psi_\mu^{(\nu)}(\lambda) = \frac{\mu - \lambda}{\mu - \varepsilon^{(\nu)}} \qquad \text{при } \mu - \varepsilon^{(\nu)} \le \lambda \le \mu;$$

$$\Psi_\mu^{(\nu)}(\lambda) = 0 \qquad \text{при } \lambda \ge \mu$$

для всех $\nu$, причем $\varepsilon^{(\nu)}(\nu = 1,2,...)$ есть последовательность положительных убывающих чисел, поэтому ей соответствует положительно определенный, ограниченный и самосопряженный оператор $\Psi_\mu(R)$, который обозначим $E_\mu$. При этом квадратичная положительно определенная форма $(E_\mu u, u)$ при $|u| = 1$ имеет совокупность значений, заключенную между 0 и 1.Далее, так как $[\Psi_\mu(\lambda)]^2 = \Psi_\mu(\lambda)$, то $E_\mu^2 = E_\mu$.

Пусть теперь $\mu_j$ возрастая, пробегает множество действительных чисел $\mu_1,\ \mu_2, \ldots$ Построим последовательность функций $\Psi_{\mu_k}(\lambda)$ и соответствующих операторов $E_{\mu_k}$. Так как последовательность функций $\Psi_{\mu_k}(\lambda)$ не убывает, то и последовательность квадратичных форм $(E_{\mu_k}u,u)$ также не убывает.

При $\mu_k \to -\infty$ имеем $\Psi_{\mu_k}(\lambda) = 0$. Значит, $\lim\limits_{\mu\to-\infty}(E_\mu u,u) = 0$.

Наоборот, легко видеть, что

$$\lim_{\mu\to+\infty}(E_\mu u,u) = (Eu,u) = |u|^2 .$$

Так как, далее $\Psi_{\mu_k}(\lambda)\Psi_{\mu_l}(\lambda) = \Psi_{\mu_k}(\lambda)$, если $\mu_k \le \mu_l$

то, по принципу соответствия $E_{\mu_k}E_{\mu_l} = E_{\mu_k}$ при этом условии.

Наконец, в силу мультипликативности соответствия и коммутативности перемножаемых операторов, все $E_{\mu_k}$; $k = 1.2,...,$, коммутируют между собой.

Далее, рассмотрим разности функций

$$\Delta\Psi_{\mu_k}(\lambda) = \Psi_{\mu_{k+1}}(\lambda) - \Psi_{\mu_k}(\lambda).$$

Все они ограничены и положительны, причем

$$\Delta\Psi_{\mu_k}(\lambda) = 0 \text{ при } \lambda < \mu_k \text{ и } \lambda \ge \mu_{k+1}$$

$$\Delta\Psi_{\mu_k}(\lambda) = 1 \text{ при } \mu_k \le \lambda \le \mu_{k+1} .$$

Каждой из них отвечает оператор (и квадратичная форма), который мы обозначим $\Delta E_{\mu_k}$. Значения $(\Delta E_{\mu_k}u,u)$ при $|u| = 1$лежат между 0 и 1. Рассмотрим два интервала $(\mu_k, \mu_{k+1}) = \delta_k$, $(\mu_j, \mu_{j+1}) = \delta_j$ вещественной оси, не имеющие общих точек. Так как

$$\Delta\Psi_{\mu_k}(\lambda)\Delta\Psi_{\mu_j}(\lambda) = 0,$$

то $$\Delta E_{\mu_k}\Delta E_{\mu_j} = 0.$$

Это значит, что операторы $\Delta E_{\mu_k}, \Delta E_{\mu_j}$ для интервалов $\delta_k$, $\delta_j$ без общих точек взаимно ортогональны.

Наконец, легко видеть, что $\sum\limits_k \Delta E_{\mu_k} = E$.

Откуда следует сходимость ряда слева. Возводя в квадрат и принимая во внимание ортогональность всех $\Delta E_{\mu_k}$, получим также

$$\sum_k |\Delta E_{\mu_k} u|^2 = |u|^2 \text{ для всех } u \subset Ӈ.$$

Последнее равенство выражает <u>*условие полноты для систем*</u> $\Delta E_{\mu_k}$.

Положив для краткости $\Delta E_{\mu_k} u = u_k = \Delta_k u$,

последние два соотношения можно представить в виде

$$\sum_k u_k = u_j, \qquad \sum_k |u_k|^2 = |u|^2.$$

<u>*Пример 2.*</u> Рассмотрим функцию $\Pi_\mu(\lambda)$, определенную для всех вещественных $\lambda$ при вещественном $\mu$, таким образом

$$\Pi_\mu(\lambda) = 0 \text{ для } \lambda \geq \mu; \quad \Pi_\mu(\lambda) = 1 \text{ для } \lambda = \mu.$$

*Покажем*, что она является пределом последовательности непрерывных функций, ей отвечают оператор $P_\mu$ [*)] и квадратичная форма $(P_\mu u, u)$, определенные положительно, причем значения формы при $|u| = 1$ лежат между 0 и 1.

<u>*Действительно,*</u> пусть $\mu$ пробегает, возрастая, множество точек $\mu_1, \mu_2, ...$ Построим последовательности функций $\Pi_{\mu_k}(\lambda)$ и операторов $P_{\mu_k}$.

Так как $|\Pi_\mu(\lambda)|^2 = \Pi_\mu(\lambda)$ ограничена, то $P^2{}_{\mu_k} = P_{\mu_k}$.

С другой стороны,

$$\Pi_{\mu_k}(\lambda)\Pi_{\mu_j}(\lambda) = 0 \text{ при } \mu_k \neq \mu_j.$$

Значит, $P_{\mu_k} \cdot P_{\mu_j} = 0$.

Таким образом, видим, что операторы $P_{\mu_k}$ взаимно ортогональны и нормированы. Легко убедиться, что все они взаимно коммутируют и что для них удовлетворяется неравенство Бесселя $\sum P^2_{\mu_k} \leq E$, откуда следует сходимость ряда, стоящего в левой части неравенства.

*) *Оператор $P_\mu$ не имеет ничего общего с теми операторами $P_\lambda$, которые определены были раньше при помощи соотношений $(R - \lambda E)P_\lambda = E$.*

## § 4. *Неограниченные операторы. Интегральные представления. Резольвента для $R - \mu E$*

В предыдущих двух параграфах было установлено, что принцип соответствия имеет силу для всех ограниченных функций, которые являются пределами последовательностей функций равномерно непрерывных. Мы видели, что соответствие между функциями и операторами положительно, дистрибутивно и мультипликативно. Исходной точкой наших рассуждений была при этом ограниченная дробь вида

$$\frac{(\mu-\sigma_1)(\mu-\sigma_2)...(\mu-\sigma_{m)}}{(\mu-\lambda_1)(\mu-\lambda_2)...(\mu-\lambda_n)} \;\; ; \;\; m \le n ,$$

где $\sigma_1, \sigma_2, ..., \sigma_m$, $\lambda_1, \lambda_2, ..., \lambda_n$ –комплексные числа (§ 2).

Как мы видели, функция

$$f_n(\mu) = \frac{(\mu-\sigma_1)...(\mu-\sigma_n)}{(\mu-\lambda_1)...(\mu-\lambda_n)} = 1 + \sum_k \frac{c_k}{\mu-\lambda_k}$$

и оператор

$$f_n(R) = (R-\sigma_1 E)P_{\lambda_1} ...(R-\sigma_n E)P_{\lambda_n} = E + \sum_k c_k P_{\lambda_k}$$

стоят в соответствии друг с другом.

Пусть теперь $\sigma$ – какое-нибудь число. Рассмотрим функцию

$$\mu - \sigma + \sum_k c_k + \sum_k \frac{c_k(\lambda_k - \sigma)}{\mu - \lambda_k} = (\mu-\sigma) f_n(\mu) .$$

Оператор $R - \sigma E + \sum_k c_k E + \sum_k c_k (\lambda_k - \sigma) P_{\lambda_k}$ ,

представляющий сумму неограниченного оператора $R - \sigma E$ с ограниченным $\sum_k c_k E + \sum_k c_k (\lambda_k - \sigma) P_{\lambda_k}$ , очевидно, определен для всех элементов $u$ из области $\mathfrak{W}$ определения $R$ .

Этот оператор может быть представлен в виде

$$R - \sigma E + \sum_k c_k [(R - \lambda_k E)P_{\lambda_k} + \lambda_k P_{\lambda_k}] - \sum_k c_k \sigma P_{\lambda_k} =$$

$$= R - \sigma E + \sum_k c_k R P_{\lambda_k} - \sum_k c_k \sigma P_{\lambda_k} = R - \sigma E + (R - \sigma E)\sum_k c_k P_{\lambda_k}$$

или в виде

$$(R - \sigma E)(E + \sum_k c_k P_{\lambda_k}) = (R - \sigma E) f_n(R) .$$

С другой стороны, вследствие коммутативности $R$ и $P_\lambda$ , мы имеем также

$$R - \sigma E + \sum_k c_k P_{\lambda_k} R - \sum_k c_k \sigma P_{\lambda_k} = R - \sigma E + \sum_k c_k P_{\lambda_k} (R - \sigma E) =$$

$$= (E + \sum_k c_k P_{\lambda_k})(R - \sigma E) = f_n(R)(R - \sigma E).$$

Следовательно, оператор

$$(R-\sigma E)f_n(R)=f_n(R)(R-\sigma E),$$

определенный для всех элементов $u \subset \mathfrak{W}$, должен быть поставлен в соответствие с *функцией*

$$(\mu-\sigma)f_n(\mu),$$

*которая, очевидно, не ограничена.* Из приведенных выше соображений видно, что это *соответствие дистрибутивно*: линейной комбинации из функций вида $(\mu-\sigma)f_n(\mu)$ соответствует такая же линейная комбинация из соответствующих операторов.

Пусть теперь последовательность функций типа $f_n(\mu)$ равномерно сходится к ограниченной непрерывной функции $f(\mu)$. В таком случае последовательность функций вида

$$\varphi_n(\mu)=\sum_k c_k+\sum_k \frac{c_k(\lambda_k-\sigma)}{(\mu-\lambda_k)}=(\mu-\sigma)[f_n(\mu)-1],$$

каждая из которых ограничена, так как соответствует ограниченному оператору вида

$$\varphi_n(R)=\sum_k c_k E+\sum_k c_k(\lambda_k-\sigma)P_{\lambda_k} \ ,$$

сходится к предельной ограниченной функции

$$(\mu-\sigma)[f(\mu)-1].$$

Вместе с тем последовательность операторов $\varphi_n(R)$ сходится к предельному оператору $(R-\sigma E)[f(R)-E]$, который непрерывно ограничен. Но это значит, что при любой равномерно непрерывной функции $f(\mu)$ устанавливается соответствие *между неограниченным оператором* $(R-\sigma E)f(R)$ и *неограниченной же функцией* $(\mu-\sigma)[f(\mu)-1]$. При этом для всех $u \subset \mathfrak{W}$ имеем

$$(R-\sigma E)f(R)=f(R)(R-\sigma E).$$

Этот вывод остается *в силе также и для всех разрывных ограниченных функций* $f(\mu)$, которые являются пределами последовательностей непрерывных функций рассмотренного вида.

*Для доказательства этого положения* нужно рассмотреть последовательности ограниченных непрерывных функций $(\mu-\sigma)[f_n(\mu)-1]$ и соответствующих им ограниченных операторов $(R-\sigma E)[f_n(R)-E]$, причем предполагается, что $f_n \to f$, где $f$ – ограниченная (разрывная) функция. В таком случае, на основании свойств соответствия, изложенных в § 3, мы имеем:

$$(\mu-\sigma)[f_n(\mu)-1] \to (\mu-\sigma)[f(\mu)-1],$$

$$(R-\sigma E)[f_n(R)-E] \to (R-\sigma E)[f(R)-E].$$

Поэтому функции $(\mu-\sigma)[f(\mu)-1]$ соответствует оператор $(R-\sigma E)[f(R)-E]$ и, следовательно, функции $(\mu-\sigma)f(\mu)$ отвечает ограниченный оператор $(R-\sigma E)f(R)$. Далее, *соответствие сохраняет свою дистрибутивность* и в этом случае.

Наконец, для всех элементов из области определения $R$ мы имеем

$$(R-\sigma E)f(R) = f(R)(R-\sigma E).$$

Наконец, заметим, что все высказанные выше положения a fortiori имеют силу в том случае, когда не только $f(\mu)$, но и $(\mu-\sigma)f(\mu)$ есть функция ограниченная.

*Пример, важный для дальнейшего.* В конце § 3 мы видели, что *соответствие распространяется на* случай функции вида

$$\Psi_\mu(\lambda) = 1 \text{ при } \lambda < \mu,$$

$$\Psi_\mu(\lambda) = 0 \text{ при } \lambda \ge \mu,$$

где $\mu$ – действительное число,

а *также на разности таких функций*:

$$\Delta\Psi_\mu(\lambda) = \Psi_{\mu''}(\lambda) - \Psi_{\mu'}(\lambda),$$

причем $\Delta\Psi_\mu(\lambda)$, очевидно, равна нулю при $\lambda < \mu'$ и при $\lambda \ge \mu''$, а при $\mu' \le \lambda < \mu''$ равняется единице.

Далее, *соответствие распространяется*, с сохранением его свойств, *и на ограниченную функцию*

$$(\lambda-\sigma)\Delta\Psi_\mu(\lambda), \ \sigma \text{ –действительное},$$

которая равна нулю при $\lambda < \mu'$, $\lambda \ge \mu''$, а при $\mu' \le \lambda < \mu''$ становится равной $\lambda - \sigma$.

Этой функции соответствует ограниченный самосопряженный оператор

$$(R-\sigma E)\Delta E_\mu,$$

где $\Delta E_\mu$ – оператор, рассмотренный в конце § 3 и соответствующий функции $\Delta\Psi_\mu(\lambda)$.

На основании того, что было сказано выше, *утверждаем,* что *для всех элементов и из области $\mathfrak{W}$ определения $R$ имеет место равенство*

$$(R-\sigma E)\Delta E_\mu = \Delta E_\mu (R-\sigma E),$$

выражающее *коммутативность* $\Delta E_\mu$ и $(R-\sigma E)$.

При этом, вследствие положительности соответствия, значения квадратичной формы $((R-\sigma E)\Delta E_\mu u, u)$ при $|u|=1$ заключены между $\mu'-\sigma$ и $\mu''-\sigma$.

В частности, при $\sigma=0, (R\Delta E_\mu u, u)$ при $|u|=1$ имеет совокупность значений между $\mu'$ и $\mu''$, а при $\mu'<\sigma<\mu''$ имеем:

$$\left|\left((R-\sigma E)\Delta E_\mu u, u\right)\right| \le (\mu''-\mu')\cdot |u|^2,$$

так что $\quad |R\Delta E_\mu u - \sigma\Delta E_\mu u| \le (\mu''-\mu')\cdot |u|$ .

В силу ограниченности операторов $R\Delta E_\mu, \Delta E_\mu$, они имеют смысл для любого $u \subset$ Ꜧ. Подставляя вместо $u$ элемент $\Delta E_\mu u$ и помня, что $\quad (\Delta E_\mu)^2 u = \Delta E_\mu u$ ,

получаем основное неравенство

$$(\mu'<\sigma<\mu'') \quad |R\Delta E_\mu u - \sigma\Delta E_\mu u| \le |\mu''-\mu'|\cdot|\Delta E_\mu u|.$$

Вернемся к функции $(\lambda-\sigma)\Delta\Psi_\mu(\lambda)$,

$$(\lambda-\sigma)\Delta\Psi_\mu(\lambda)=0 \qquad \text{при} \quad \lambda<\mu', \ge \mu'',$$

$$(\lambda-\sigma)\Delta\Psi_\mu(\lambda)=\lambda-\sigma \quad \text{при} \quad \lambda \ge \mu', < \mu''$$

и соответствующему ей ограниченному самосопряженному оператору $(R-\sigma E)\Delta E_\mu$.

Пусть интервал $(\mu',\mu'')$ расширяется, так что $\mu'\to-\infty, \mu''\to+\infty$.

Тогда, как было доказано в конце § 3 относительно сходимости оператора $E_\mu$, видим, что $\quad E_{\mu'}\to 0, E_{\mu''}\to E$ и, следовательно, $\Delta E_\mu \to E$.

Поэтому

$$(R-\sigma E)\Delta E_\mu \to R-\sigma E.$$

С другой стороны, функция $(\lambda-\sigma)\Delta\Psi_\mu(\lambda)$, определенная выше, при $\mu'\to-\infty, \mu''\to+\infty$, сходится к функции $\lambda-\sigma$ ,определенной для всех действительных $\lambda$.

Таким образом, *неограниченной функции* $\lambda-\sigma$ *соответствует неограниченный оператор* $R-\sigma E$.

Пусть действительная ось $\lambda$ разделена на счетное число не налегающих друг на друга интервалов $(\mu_k', \mu_k'')$. Каждому из этих интервалов соответствует ограниченная функция

$$\Delta\Psi_\mu(\lambda) = \begin{cases} 0 \\ 1 \end{cases} \text{при } \lambda \begin{cases} <\mu_k', \geq \mu_k'' \\ \geq \mu_k', <\mu_k'' \end{cases},$$

а каждой из таких функций, в свою очередь, соответствует оператор $\Delta E_\mu$.

При этом $(\Delta E_{\mu_k})^2 = \Delta E_{\mu_k}$, $\Delta E_{\mu_k} \Delta E_{\mu_j} = 0$ при $k \neq j$

и, наконец,

$$\sum \Delta E_{\mu_k} = E.$$

Рассмотрим элемент $u$ такой, для которого оператор $R$ имеет смысл и пусть $Ru = v$.

Положим для краткости

$$\Delta E_{\mu_k} u = u_k, \quad \Delta E_{\mu_k} v = v_k. \tag{1}$$

Тогда $(u_k, u_j) = 0, k \neq j;\ (u_k, u_k) = \left|u_k^2\right|$,

$$\sum |u_k|^2 = (u,u) = |u|^2, \sum u_k = u,$$

и аналогично $(v_k, v_j) = 0, k \neq j; (v_k, v_k) = \left|v_k^2\right|,\ \sum |v_k|^2 = |v|^2$

или $$\sum_k \left|R\Delta E_{\mu_k} u\right|^2 = |Ru|^2, \tag{2}$$

откуда видно, что для того, чтобы «$u$» принадлежало к области определения оператора $R$, <u>необходима сходимость</u> ряда $\sum_k \left|R\Delta E_{\mu_k} u\right|^2$.

Кроме того, так как $v_k = \Delta E_{\mu_k} v = \Delta E_{\mu_k} Ru = R\Delta E_{\mu_k} u = Ru_k$,

значит $$Ru = v = \sum v_k = \sum Ru_k. \tag{3}$$

Предположим, что ряд в левой части (2) сходится для данного $u$.

Тогда

$$\left|\sum_k^n Ru_k - \sum_k^m Ru_k\right|^2 = \left|\sum_{k=m+1}^n Ru_k\right|^2 = \sum_{k=m+1}^n |Ru_k|^2$$ при любом $n > m$ и при $m \to \infty$

стремится к нулю, а это, в свою очередь, означает, что ряд $\sum Ru_k$ сходится.

*Пусть $u^*$ есть соответствующий предельный элемент.* Если $R$ имеет смысл для какого-нибудь элемента $w$, то для всякого такого элемента

$$(Rw, \sum_1^n u_k) = (w, \sum_1^n Ru_k),$$

при $n \to \infty$ левая часть стремится к $(Rw, u)$, а правая – к $(w, u^*)$; поэтому в пределе

$$(Rw, u) = (w, u^*),$$

откуда, согласно самосопряженности $R$, $Ru = u^*$, следовательно, $Ru$ имеет смысл.

*Итак, сходимость ряда* $\sum_k |Ru_k|^2$ *есть не только необходимое*, но и достаточное условие того, чтобы для данного $u$ *выражение* $Ru$ *имело смысл.*

Предположим теперь, что минимум длины интервалов $(\mu'_k, \mu''_k)$ для всех $k$ не превосходит $\delta > 0$. Выберем внутри каждого из них $\sigma_k$, причем $\mu' \le \sigma_k \le \mu''$. Тогда на основании (1), очевидно, будем иметь

$$|Ru_k - \sigma_k u_k| \le (\mu''_k - \mu'_k)|u_k| \le \delta |u_k|.$$

Возводя эти неравенства в квадрат и складывая, получим

$$\sum_k |Ru_k - \sigma_k u_k|^2 \le \delta_2 \sum_k |u_k|^2 = \delta^2 |u|^2 \quad , \qquad (4)$$

так что ряд слева непременно сходится.

Отсюда далее следует, что при всяком $u$:

$$\left| (\sum_{k=1}^{n} |Ru_k|^2)^{1/2} - (\sum_{k=1}^{n} |\sigma_k|^2 |u_k|^2)^{1/2} \right| \le (\sum_{k=1}^{n} \left|Ru_k - \sigma_k u_k\right|^2)^{1/2} \le \delta |u|.$$

Это значит, что сходимость одного из рядов $\sum_{k=1}^{n} |Ru_k|^2$, $\sum_{k=1}^{n} |\sigma_k|^2 |u_k|^2$ влечет за собой сходимость другого.

Таким образом, для того, чтобы оператор $R$ имел смысл для элемента $u$, *необходимо,* чтобы при всяком разбиении оси $\lambda$ на интервалы $(\mu'_k, \mu''_k)$ ограниченной длины и при любом выборе $\sigma_k$ на этих интервалах ряд $\sum_{k=1}^{n} |\sigma_k|^2 |u_k|^2$ сходился.

Но, *достаточна* сходимость этого ряда в каком-нибудь одном частном случае.

Далее, из сходимости ряда в левой части (4) вытекает сходимость ряда

$$\sum_k (Ru_k - \sigma_k u_k).$$

При этом, если $Ru$ имеет смысл, то $\sum_k Ru_k = Ru$ сходится, а тогда сходится и $\sum_k \sigma_k u_k$. Наконец,

$$\left| Ru - \sum_k \sigma_k u_k \right| \le \delta |u|. \qquad (5)$$

Пусть теперь $\delta \to 0$. Тогда вместо рядов получим в пределе соответствующие интегралы Стилтьеса. При этом увидим, *что сходимость интеграла*

$$\int_{-\infty}^{+\infty} |\lambda|^2 d(E_\lambda u, u)$$

*есть необходимое и достаточное условие того, чтобы для данного* $u$ $Ru$ *имел бы смысл.* Если это условие выполнено, то

$$Ru = \int_{-\infty}^{+\infty} \lambda d(E_\lambda u).$$

Пусть теперь $F(\lambda)-$ ограниченная и равномерно непрерывная для всех действительных $\lambda$ функция. *Разделим ось* $\lambda$ на *две системы интервалов*: $(\mu'_k, \mu''_k), (\nu'_j, \nu''_j)$ и построим системы соответствующих функций $\Delta\Psi_{\mu_k}(\lambda)$, $\Delta\Psi_{\nu_j}(\lambda)$, которым отвечают операторы $\Delta E_{\mu_k}$, $\Delta E_{\nu_j}$.

Пусть колебание функций $F(\lambda)$ на каждом из интервалов не превосходит $\varpi > 0$. Имея в виду положить в пределе $\varpi = 0$, можем принять, что значение функции $F(\lambda)$ в точках интервалов $(\mu'_k, \mu''_k)$ равно значению ее в точке $\mu'_k$.

Тогда $\left|\sum F(\mu'_k)\Delta\Psi_{\mu_k}(u) - \sum F(\nu_j)\Delta\Psi_{\nu_j}(u)\right|$ не превзойдет $2\omega$.

Переходя, подобно тому, как это было сделано выше, от функций к соответствующим операторам, увидим, что

$$\left|\sum F(\mu'_k)\Delta E_{\mu_k} u - \sum F(\nu'_j)\Delta E_{\nu_j} u\right| \le 2\omega \cdot |u|.$$

При $\omega \to 0$ это неравенство убеждает нас в том, что $\sum F(\mu'_k)\Delta E_{\mu_k} u$ имеет предел, не зависящий от способа подразделения оси $\lambda$ на частные интервалы, и что предел этот есть интеграл Стилтьеса $\int_{-\infty}^{+\infty} F(\lambda) dE_\lambda u$, который сходится для всех $u \subset$ Ꜧ.

Этот интеграл мы обозначим

$$F(R) = \int_{-\infty}^{+\infty} F(\lambda) dE_\lambda u.$$

Из самого способа его получения явствует, что он коммутирует с $R$, так что $f(R)R$ есть продолжение $Rf(R)$.

Пусть $R$ имеет смысл для «$u$». Тогда

$$F(R)Ru = RF(R)u = \int_{-\infty}^{+\infty} F(\lambda) R dE_\lambda u.$$

Так как разность $\int_{-\infty}^{+\infty} F(\lambda) R dE_\lambda u - \int_{-\infty}^{+\infty} F(\lambda) \lambda dE_\lambda u$

можно рассматривать как предел суммы $\sum F(\sigma_k)(Ru_k - \sigma_k u_k)$;

и так как

$$\sum |Ru_k - \sigma_k u_k|^2 \le \delta^2 |u|^2,$$

то $$\left|\sum F(\sigma_k)(Ru_k - \sigma_k u_k)\right|^2 \le (\max|F(\lambda)|)^2 \delta^2 |u|^2.$$

Отсюда, вследствие произвола в выборе $\delta$, разность обеих сумм в левой части в пределе стремится к нулю и мы имеем

$$F(R)Ru = RF(R)u = \int_{-\infty}^{+\infty} F(\lambda)RdE_\lambda u = \int_{-\infty}^{+\infty} F(\lambda)\lambda dE_\lambda u .$$

Как применение полученных результатов, <u>рассмотрим функцию</u>

$$f(\lambda) = \frac{1}{\lambda - \mu},$$

где $\lambda$ – действительная переменная, $\mu$ – комплексный параметр.

Так как $f(\lambda)$ – ограничена и непрерывна, то ей соответствует оператор

$$f(R) = \int_{-\infty}^{+\infty} \frac{dE_\lambda}{\lambda - \mu}.$$

Пусть $R$ – любой (самосопряженный) оператор, $\mathfrak{W}$ – его область определения, $u \subset \mathfrak{W}$. Тогда

$$f(R)Ru = \int_{-\infty}^{+\infty} \frac{\lambda dE_\lambda u}{\lambda - \mu},$$

или $$Rf(R)u = f(R)Ru = \int_{-\infty}^{+\infty} \frac{\lambda - \mu + \mu}{\lambda - \mu} dE_\lambda u = Eu + \mu f(R)u;$$

откуда

$$f(R)[R - \mu E]u = [R - \mu E]f(R)u = Eu .$$

Значит, оператор

$$f(R) = \int_{-\infty}^{+\infty} \frac{dE_\lambda}{\lambda - \mu}$$

обратен по отношению к $R - \mu E$ для всех комплексных $\mu$. Но $\int_{-\infty}^{+\infty} \frac{dE_\lambda}{\lambda - \mu}$ имеет смысл и *для таких действительных значений* $\mu$ в достаточно малой окрестности которых $E_\lambda$ *остается постоянной, несмотря на то, что при* $\lambda = \mu$ *дробь* $\frac{1}{\lambda - \mu}$ *обращается в бесконечность*.

<u>*В самом деле,*</u> рассмотрим интервал $(\mu - h, \mu + h)$, где $h$ достаточно мало. Пусть для всех $\lambda$ внутри этого интервала $dE_\lambda \equiv 0$. Истинное значение отношения двух непрерывных функций $dE_\lambda$ и $\lambda$ - $\mu$ при $\lambda = \mu$ есть $\left[\frac{d(dE_\lambda)}{d(\lambda - \mu)}\right]_{\lambda=\mu} = 0$, откуда и следует высказанное положение.

<u>*Исключение*</u> в этом смысле составляют лишь такие точки действительной оси, в достаточно малых окрестностях, для которых $dE_\lambda \neq 0$. Эти исключительные значения $\mu$ образуют замкнутое множество, называемое ***спектром формы*** $(Ru, u)$.

*Докажем*, что в <u>*точках спектра оператор*</u> $R-\mu$ <u>*не имеет себе обратного*</u>.

<u>*Действительно,*</u> пусть $\mu-$ точка действительной оси, в достаточно малой окрестности $(\mu-h,\mu+h)$которой неубывающая функция $E_\lambda$ не остается постоянной. Определим функцию $\varphi(\lambda)$ так:

$\varphi(\lambda)=1$ внутри $(\mu-h,\mu+h)$ и в точке $\mu-h$,

$\varphi(\lambda)=0$ в прочих точках.

$\varphi(\lambda)$ есть предел последовательности непрерывных функций, следовательно, ей соответствует оператор (квадратичная форма), которую обозначим $\varphi(R)$.

Далее, так как $[\varphi(\lambda)]^2=\varphi(\lambda)$,

то $[\varphi(\lambda)]^2$соответствует оператор (и форма) $[\varphi(R)]^2=\varphi(R)$.

Кроме того:

$$\varphi(R)=\int_{-\infty}^{+\infty}\varphi(\lambda)dE_\lambda=\int_{-\mu+h}^{\mu+h}1dE_\lambda=E_{\mu+h}-E_{\mu-h},$$

значит форма $(\varphi(R)u,u)$ положительна по крайней мере для одного u.

Положим $v=\varphi(R)u$. Тогда $|v|^2=[\varphi(R)u]^2>0$, следовательно, $v$, соответствующий $u\neq 0$, также $\neq 0$.

Далее $\varphi(R)v=\varphi(R)\varphi(R)u=\varphi(R)u=v$.Поэтому

$$[R-\mu E]^2\varphi(R)v=[R-\mu E]^2v.$$

Но оператору $[R-\mu E]^2\varphi(R)$ соответствует функция $(\lambda-\mu)^2\varphi(\lambda)$,которая, как легко видеть, на интервале $(\mu-h\leq\lambda<\mu+h)$ принимает значения от 0 до $h^2$. Следовательно, в этих же пределах заключены и значения формы $([R-\mu E]^2v,v)$, $|v|=1$.

Положим

$$[R-\mu E]v=v'.$$

Тогда

$$|v'|^2=|([R-\mu E]^2v,v)|\leq h^2|v|^2.$$

Так как $h-$ как угодно мало, то левая часть может быть только нулем.

Итак, в уравнении

$[R-\mu E]v=v'$ имеем $v'=0$ при $v\neq 0$.

Для того, чтобы уравнение $[R-\mu E]v=v'$ могло быть обращено, необходимо (хотя и недостаточно), чтобы

$$|([R-\mu E]^2v,v)|=|v'|^2\geq m^2|v|^2,$$

где $m\neq 0-$ действительное число.

Так как в рассматриваемом случае $v'=0$, $v\neq 0$, то это условие не соблюдается, откуда следует, что для выбранного значения $\mu$ оператор $R-\mu E$ не имеет себе обратного, то есть, $\mu$ *принадлежит спектру* $R$, *<u>что и требовалось доказать.</u>*

Москва, май, 1936 г. А Герасимов. Мехмат/заочное отделение

*Прикладная математика и механика, 1938, I,493-536.*

## 3.ПРОБЛЕМА УПРУГОГО ПОСЛЕДЕЙСТВИЯ И ВНУТРЕННЕЕ ТРЕНИЕ

Теория упругости дает закон малых деформаций в виде системы уравнений:[1)]

$$(\lambda+\mu))\frac{\partial\theta}{\partial x}+\mu\Delta u+\rho(X-\frac{\partial^2 u}{\partial t^2})=0,$$

$$(\lambda+\mu))\frac{\partial\theta}{\partial y}+\mu\Delta v+\rho(Y-\frac{\partial^2 v}{\partial t^2})=0, \qquad (L)$$

$$(\lambda+\mu))\frac{\partial\theta}{\partial z}+\mu\Delta w+\rho(Z-\frac{\partial^2 w}{\partial t^2})=0,$$

Здесь:

u,v,w – компоненты смещения,

X,Y,Z – компоненты внешних сил, рассчитанных на единицу массы,

$\lambda,\mu$ – постоянные Ламé,

$$\theta=\frac{\partial u}{\partial x}+\frac{\partial v}{\partial y}+\frac{\partial w}{\partial z}$$

и $\Delta$ обозначает оператор Лапласа

$$\frac{\partial^2}{\partial x^2}+\frac{\partial^2}{\partial y^2}+\frac{\partial^2}{\partial z^2}.$$

Уравнения (L) содержат, таким образом, частные производные второго порядка относительно компонент смещения. Физическая основа их заключается в законе Гука, по которому во всяком упругом теле при деформации последнего и направленные в сторону ее уменьшения. Такого рода упругие силы принадлежат к классу консервативных сил, зависящих от взаимного расположения частиц деформируемого тела.

С другой стороны, гидродинамика описывает движение вязких жидкостей при помощи уравнений Навье – Стокса:[2)]

$$-\frac{\partial p}{\partial x}+\eta\left\{\frac{1}{3}\frac{\partial\Theta}{\partial x}+\Delta v_x\right\}+\rho\left\{X-\frac{dv_x}{dt}\right\}=0,$$

[1)] *Ляв, Матем. теория упругости, стр. 144, ОНТИ; 1935; М.М. Филоненко-Бородич. Основы теории упругости, ГСИ, стр. 43, 1933.*

[2)]*Вебстер А.Г., Механика материальных точек, твердых, упругих и жидких тел, ГТТИ, стр.592.1933.*

$$-\frac{\partial p}{\partial y}+\eta\left\{\frac{1}{3}\frac{\partial\Theta}{\partial y}+\Delta v_y\right\}+\rho\left\{Y-\frac{dv_y}{dt}\right\}=0, \qquad \text{(N-S)}$$

$$-\frac{\partial p}{\partial z}+\eta\left\{\frac{1}{3}\frac{\partial\Theta}{\partial z}+\Delta v_z\right\}+\rho\left\{Z-\frac{dv_z}{dt}\right\}=0,$$

где $p$ есть давление, $v_x, v_y, v_z$ – компоненты скорости,

$$\Theta=\frac{\partial v_x}{\partial x}+\frac{\partial v_y}{\partial y}+\frac{\partial v_z}{\partial z} \quad \text{и}$$

$\eta$ – коэффициент вязкости.

Относительно компонент смещения уравнения (N-S) уже третьего порядка. В их основу положен закон внутреннего трения, установленный Ньютоном. По этому закону сила внутреннего трения, или сила вязкости, возникающая между двумя скользящими друг по другу жидкими слоями, пропорциональна градиенту скорости и поверхности раздела и направлена в сторону уменьшения этого градиента.

Таким образом, если упругие силы зависят от распределения смещений по объему деформированного тела, то силы вязкости, возникающие при деформации вязких тел, зависят от распределения скоростей.

Но в природе широко распространен класс веществ, по всем признакам обладающих одновременно как вязкими, так и упругими свойствами. Достаточно назвать смолы, хлебное тесто и пр. К этой категории веществ относятся, все коллоиды. Такого рода вещества я называю упруго-вязкими.

Казалось бы очень важным построить систему уравнений, описывающую закон движения частиц упруго-вязкого тела. Такая система уравнений должна содержать в себе как уравнения Ламé, так и уравнения Навье-Стокса, являясь их обобщением.

В предлагаемой работе сделан первый шаг на пути такого обобщения, ограничиваясь случаем одномерных тел. Именно, объектом настоящего исследования является тело, один из линейных размеров которого преимущественно развит по сравнению с линейными размерами в прочих направлениях. Такого рода тело условимся называть нитью, не вкладывая, впрочем, в это слово физического содержания. Исключительно большой размер нити естественно называть длиной нити, а тот отрезок кривой, около которого сосредоточено все вещество нити - ее осью. Всякое сечение нити плоскостью, перпендикулярной к касательной к оси, будет называться в дальнейшем просто сечением. Размеры любого сечения нити очень малы по сравнению с ее длиной.

Будем считать, что нить состоит из вещества однородного, но не обязательно изотропного, т.е. обладая одинаковыми физическими свойствами во всех своих точках, нить может обладать различными свойствами в различных направлениях. Площади сечений нити будем считать равными, не предполагая, однако, равенства самих сечений. Переход от одной формы сечения к другой предполагаемы непрерывным, так что на любом участке длины нити можно взять два настолько близких сечения, что выделяемый ими объем может быть принят за прямой цилиндр с очень маленькой высотой.

Наконец, мы предполагаем, что для изучаемого тела опытно установлены плотность $\rho$ вещества, коэффициент вязкости $\eta$ и коэффициенты Ламé.

Работа состоит из трех частей. Первая часть посвящена малым поперечным колебаниям нити, закрепленной на концах. В § 1 дается вывод уравнения движения. В § 2 показано, что уравнения (N-S) гидродинамики приводят к тому же закону движения нити при условии надлежащего истолкования гидродинамического давления. В § 3 приводится доказательство единственности решения полученного уравнения с частными производными третьего порядка. В § 4 проведено решение этого уравнения по способу Д. Бернулли. Решение получается в виде ряда, и в § 5 доказывается возможность почленного дифференцирования этого ряда нужное количество раз. Далее, в § 6 устанавливается, что полученный в § 4 ряд действительно удовлетворяет всем условиям задачи. В § 7 выясняется физический смысл найденного решения. В § 8 разбирается вопрос об использовании этого решения для опытного определения коэффициентов упругости и вязкости вещества нити путем наблюдения периода колебания и декремента затухания. Далее исследуется неоднородное уравнение: в § 9 приводится решение для тяжелой нити, в § 10 рассматривается общий случай вынужденного движения нити. В § 11 исследуется вопрос о сходимости рядов и устанавливается соблюдение всех условий задачи.

Вторая часть этой работы представляет собой довольно сжатый набросок решения задачи о крутильных колебаниях нити, один из концов которой закреплен, а другой несет некоторый дополнительный момент инерции, причем разбирается лишь случай колебаний в пустоте и с небольшой амплитудой. В § 12 показано, что задача о крутильных колебаниях приводит к тому же уравнению движения, что и задача о поперечных колебаниях(§ 1). В § 13 исследуются условия на свободных концах нити, в § 14 дается решение уравнения, в § 15 устанавливается, что это решение удовлетворяет всем условиям задачи.

Третья часть посвящена упругому гистерезису. В § 17 приведен упрощенный вывод уравнения Больцмана; § 18 посвящен способу определения ядра уравнения для того случая, когда рассматриваемая частица нити не лежит в особой точке; в § 19 дается общий способ определения ядра.

**§ 1**. Пусть нить с длиной $l$ и с площадью сечения $s$ натянута так, что ее ось располагается на сегменте $[0,l]$ прямой $OX$, причем концы нити закреплены. Тогда все сечения последней расположатся перпендикулярно к $OX$.

Пусть в момент времени $t=0$ оси нити придана форма, определяемая уравнением:

$$u|_{t=0} = \Phi(x) \tag{1}$$

так что $\Phi(0)=\Phi(l)=0$.

Эту начальную форму нити мы предполагаем плоской и обладающей настолько малой кривизной, что сечения нити можно считать попрежнему перпендикулярными к $OX$.

Кроме того, допустим, что в тот же момент времени $t=0$ всем элементам нити сообщены начальные скорости, перпендикулярные к $OX$ и направленные в плоскости кривой (1). Таким образом дано начальное распределение скоростей:

$$\frac{\partial u}{\partial t}\Big|_{t=0} = \varphi(x). \tag{2}$$

Составим уравнение движения нити. Рассмотрим силы, действующие на бесконечно тонкий слой B нити (фиг. 1),

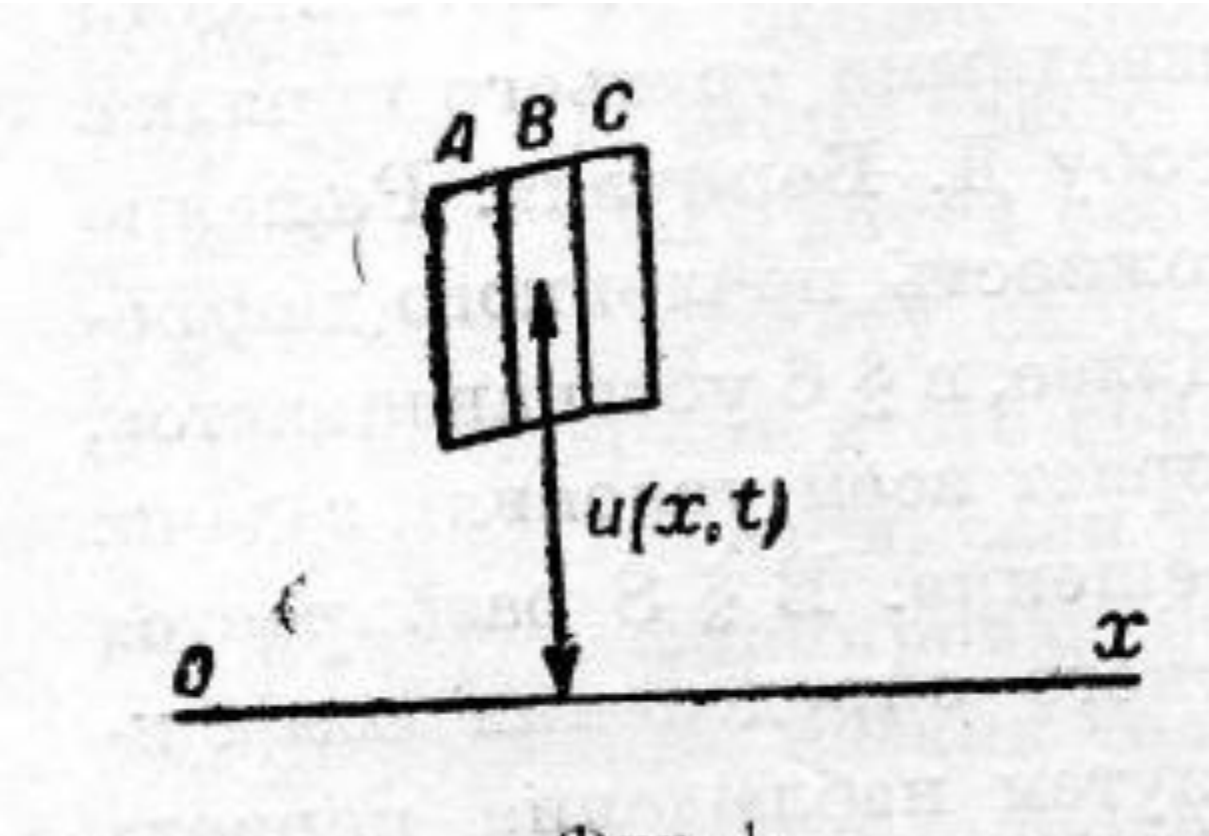

Фиг. 1.

выделяемой из ее тела двумя бесконечно близкими сечениями с абсциссами:

$$x-\frac{dx}{2} \text{ и } x+\frac{dx}{2}.$$

Тогда сила инерции элемента B есть

$$-\rho s dx \frac{\partial^2 u}{\partial t^2} \tag{3}$$

Если считать, что направление смещений $u$ вертикально и обращено кверху, то вес рассматриваемого элемента оказывается равным

$$-\rho g s dx \tag{4}$$

Пусть, далее, нить натянута так, что на единицу площади ее сечения действует продольная сила $T$. Тогда на всю площадь $s$ сечения действует сила натяжения, равная $Ts$. Ее вертикальная составляющая в точке с абсциссой $x-\frac{dx}{2}$ есть (фиг.2):

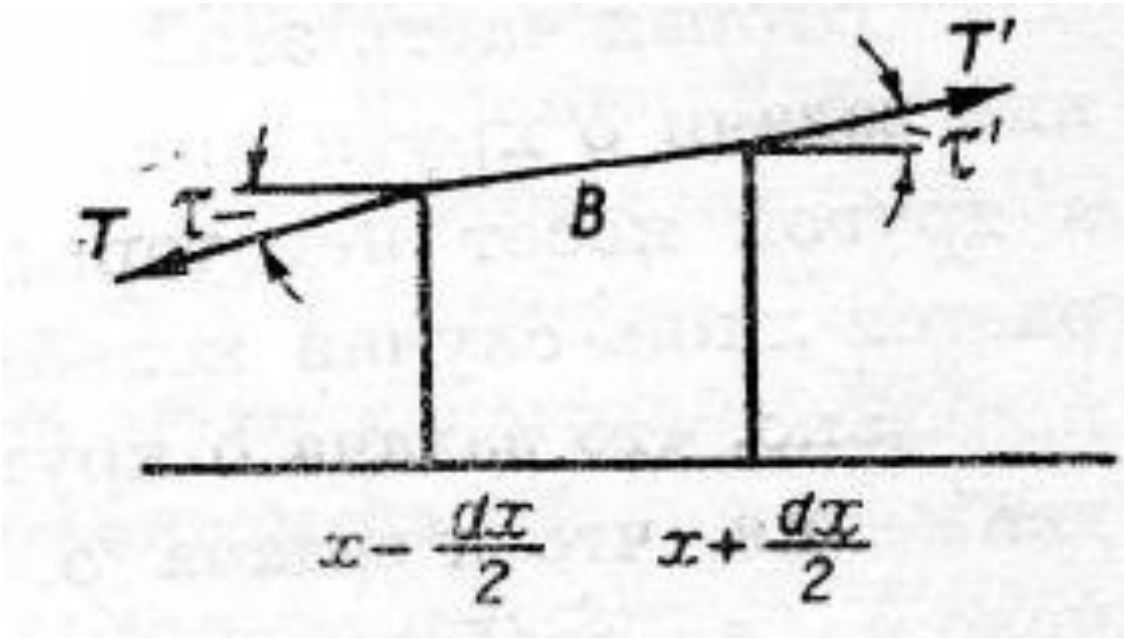

Фиг. 2.

$$-Ts\sin\tau=-Tstg\,\tau=-Ts\left(\frac{\partial u}{\partial x}\right)_{x-\frac{dx}{2}},$$

где $\tau$ обозначает угол наклона касательной к $OX$. В точке с абсциссой $x+\frac{dx}{2}$ натяжение дает вертикальную составляющую

$$+Ts\left(\frac{\partial u}{\partial x}\right)_{x+\frac{dx}{2}}.$$

Таким образом, вследствие натяжения нити, ее элемент В испытывает силу, направленную вниз и равную

$$-Ts\left[\left(\frac{\partial u}{\partial x}\right)_{x+\frac{dx}{2}}-\left(\frac{\partial u}{\partial x}\right)_{x-\frac{dx}{2}}\right]=-Tsdx\frac{\partial^2 u}{\partial x^2}, \qquad (5')$$

а сила, направленная в сторону смещения, получается равной

$$+Tsdx\frac{\partial^2 u}{\partial x^2}. \qquad (5)$$

Силу (5’) , пропорциональную кривизне средней линии $\frac{\partial^2 u}{\partial x^2}$ и направленную в сторону уменьшения этой кривизны, можно было бы назвать, следуя Планку, силой инерции формы. Однако в таком случае постоянный коэффициент $-T$ должен быть истолкован несколько иначе. Я хочу сказать следующее: если бы натяжение $T$ было равно нулю, то и в таком режиме движения нить испытывала бы местные деформации сдвига. Мы рассмотрим этот вопрос подробнее. Вследствие того, что в сечении с абсциссой $x-\frac{dx}{2}$ появляется деформация сдвига, измеряемая величиной угла $\tau$ , или, что то же, значением производной

$$\left(\frac{\partial u}{\partial x}\right)_{x-\frac{dx}{2}},$$

по закону Гука в данном сечении нити возникает сила упругости, равная

$$-\mu\left(\frac{\partial u}{\partial x}\right)_{x-\frac{dx}{2}}\cdot s,$$

где $\mu$ – модуль сдвига.

Аналогично в сечении с абсциссой $x+\dfrac{dx}{2}$ появляется сила упругости

$$+\mu\left(\frac{\partial u}{\partial x}\right)_{x+\frac{dx}{2}}\cdot s.$$

Поэтому силой инерции формы нужно считать не (5’), а величину

$$-Msdx\frac{\partial^2 u}{\partial x^2},$$

где положено

$$M=T+\mu. \tag{6}$$

При отсутствии натяжения $T$ коэффициент $M$ делается тождественным модулю сдвига $\mu$.

Итак, результирующая сила упругости и натяжения, направленная в сторону смещения $u$ и действующая на элемент В нити, имеет величину:

$$+Msdx\frac{\partial^2 u}{\partial x^2},\quad M=T+\mu. \tag{7}$$

Наконец, благодаря различию в скоростях в каждый момент времени слои А,В,С,…, выделяемые из тела нити рядом бесконечно близких сечений, вынуждены скользить друг по другу. Вследствие вязкости вещества они взаимодействуют между собой. В основу этого взаимодействия я кладу закон Ньютона, по которому сила действия одного слоя на другой, соседний с ним, должна быть пропорциональна площади s сечения нити и градиенту скорости $\dfrac{\partial\dot{u}}{\partial x}=\dfrac{\partial^2 u}{\partial x\partial t}$; кроме того, она должна быть направлена в сторону уменьшения этого градиента.

Поэтому слой А (фиг.1) действует на слой В с силой

$$-\eta s\left(\frac{\partial\dot{u}}{\partial x}\right)_{x-\frac{dx}{2}};$$

по той же причине действие В на С есть

$$-\eta s\left(\frac{\partial \dot{u}}{\partial x}\right)_{x+\frac{dx}{2}};$$

Следовательно, действие C на B равно:

$$+\eta s\left(\frac{\partial \dot{u}}{\partial x}\right)_{x+\frac{dx}{2}};$$

Следовательно, действие слоев A и C на лежащий между ними слой B получается равным

$$+\eta s\left[\left(\frac{\partial \dot{u}}{\partial x}\right)_{x+\frac{dx}{2}}-\left(\frac{\partial \dot{u}}{\partial x}\right)_{x-\frac{dx}{2}}\right]=\eta s\left(\frac{\partial^2 \dot{u}}{\partial x^2}\right). \tag{8}$$

По принципу Даламбера, сумма сил (3),(4),(7) и (8) должна тождественно равняться нулю, и мы приходим к закону малых поперечных колебаний упруго-вязкой нити:

$$\eta\frac{\partial^3 u}{\partial x^2 \partial t}+M\left(\frac{\partial^2 u}{\partial x^2}\right)=\rho\left[g+\left(\frac{\partial^2 u}{\partial t^2}\right)\right] \tag{A}$$

при условии (1) и (2) для начального момента $t=0$ и граничных условиях:

$$u|_{x=0}=u|_{x=t}=0. \tag{9}$$

При $\eta=0$ уравнение (A) превращается в хорошо изученное со времени Даламбера уравнение колебаний упругой тяжелой струны, натянутой вдоль оси $u=0$ и закрепленной на концах.

**§ 2.** Так как в рассуждениях § 1 имеются пункты, которые можно было бы оспаривать, я даю здесь другой вывод уравнения (A), исходя из гидромеханических соображений.

Рассмотренное в § 1 движение частиц A,B,C,… нити происходит так, как если бы нить находилась в поле гидродинамического давления, линии градиента которого параллельны направлению смещений.

Направим ось $OY$ в плоскости движения в сторону положительного отсчета смещений, так что $y\equiv u$. Предположим сначала, что нить имеет прямоугольное сечение с площадью $s$ и длиной ребра, параллельного $OY$, равной $\Delta y$. Следовательно, длина другой стороны сечения нити равна $\frac{s}{\Delta y}$,исходящая вследствие натяжения и упругости нити и действу и так как длина слоя B по-прежнему считается равной $dx$, то площадь $s'$ грани B, перпендикулярной к $OY$, есть

$$s' = \frac{s}{\Delta y} dx. \qquad (10).$$

Существенно отметить, что $\Delta y$ не есть прирост смещения $u$. В § 1 мы видели (7), что сила, происходящая вследствие натяжения и упругости нити и действующая на слой B в сторону положительных ординат $OY$, равна:

$$+ Msdx \frac{\partial^2 u}{\partial x^2},$$

Следовательно, на единицу площади $s' \perp OY$ действует сила

$$+ M \frac{s}{s'} dx \frac{\partial^2 u}{\partial x^2} \;,$$

Которую можно отождествлять с убылью фиктивного гидродинамического давления в направлении $OY$ на отрезке $\Delta y$. Этот недостаток равен

$$-\frac{\partial p}{\partial y} \Delta y.$$

Сопоставление двух последних выражени1й на основании (10) дает:

$$-\frac{\partial p}{\partial y} = M \frac{\partial^2 u}{\partial x^2}. \qquad (11)$$

Мы предполагали сечение нити прямоугольным. Однако, каково бы оно ни было, его всегда можно разбить на достаточно узкие прямоугольники посредством прямых, параллельных $OY$. Повторение приваеденных рассуждений даст тот же результат (11), который носит поэтому общий характер.

Далее, мы имеем систему уравнений Навье-Стокса, которую напишем в виде:

жидкостей при помощи уравнений Навье – Стокса:[1]

$$-\frac{\partial p}{\partial x} + \eta \left\{ \frac{1}{3} \frac{\partial \theta}{\partial x} + \Delta v_x \right\} + \rho \left\{ g_x - \frac{dv_x}{dt} \right\} = 0,$$

$$-\frac{\partial p}{\partial y} + \eta \left\{ \frac{1}{3} \frac{\partial \theta}{\partial y} + \Delta v_y \right\} + \rho \left\{ g_y - \frac{dv_y}{dt} \right\} = 0, \qquad \text{(N-S)}$$

$$-\frac{\partial p}{\partial z} + \eta \left\{ \frac{1}{3} \frac{\partial \theta}{\partial z} + \Delta v_z \right\} + \rho \left\{ g_z - \frac{dv_z}{dt} \right\} = 0,$$

здесь

$$\theta = \frac{\partial v_x}{\partial x} + \frac{\partial v_y}{\partial y} + \frac{\partial v_z}{\partial z} \quad ,$$

$v_x, v_y, v_z$ обозначают проекции скорости жидкой частицы на оси координат, $g_x = 0,\ \ g_y = -g,\ \ g_z = 0$ – проекция ускорения $g$ силы тяжести, если ось $OY$ направлена, как в нашем случае, вертикально вверх; символ $\frac{d}{dt}$ обозначает полную производную по времени, так что, например:

$$\frac{dv_x}{dt} = \frac{\partial v_x}{\partial t} + \frac{\partial v_x}{\partial x} v_x + \frac{\partial v_x}{\partial y} v_y + \frac{\partial v_x}{\partial z} v_z \ .$$

В рассматриваемом случае

$$\frac{\partial p}{\partial x} = 0 \ , \ -\frac{\partial p}{\partial y} = M \frac{\partial^2 u}{\partial x^2} ,\ \frac{\partial p}{\partial z} = 0 \ .$$

Далее,

$$v_x = 0,\ \ v_y = \frac{\partial u}{\partial t} = \dot{u} \ ,\ v_z = 0,$$

так что

$$\theta = \frac{\partial}{\partial y}\left(\frac{\partial u}{\partial t}\right) = \frac{\partial}{\partial t}\left(\frac{\partial u}{\partial y}\right) = \frac{\partial}{\partial t}(1) = 0 \, .$$

Затем находим:

$$\Delta v_x = 0,\ \ \Delta v_y = \Delta \frac{\partial u}{\partial t} = \frac{\partial}{\partial t}\left(\frac{\partial^2 u}{\partial x^2}\right),\ \ \Delta v_z = 0 \, ;$$

наконец,

$$\frac{dv_x}{dt} = 0 \, ,\ \frac{dv_y}{dt} = \frac{\partial}{\partial t}\left(\frac{\partial u}{\partial t}\right) = \frac{\partial^2 u}{\partial t^2} ,\ \frac{dv_z}{dt} = 0 \, .$$

Подставляя все это в уравнения Навье-Стокса, увидим, чт о остается только второе уравнение, принимающее вид:

$$\eta \frac{\partial^3 u}{\partial x^2 \partial t} + M \frac{\partial^2 u}{\partial x^2} = \rho \left\{ g + \frac{\partial^2 u}{\partial} t^2 \right\} ,$$

Вполне совпадающий с полученным в § 1.

§ **3.**Вместо выведенного в предыдущих параграфах уравнения (A) мы рассмотрим теперь более общее уравнение:

$$Ru \equiv \eta \frac{\partial^3 u}{\partial x^2 \partial t} + M \frac{\partial^2 u}{\partial x^2} - \rho \frac{\partial^2 u}{\partial t^2} - H \frac{\partial u}{\partial t} = -F(x,t), \qquad \text{(B)}$$

описывающее вынужденные колебания упруго-вязкой нити в среде, сопротивляющейся с силой, пропорциональной первой степени скорости. Здесь $\eta, M, \rho, H$ – положительные коэффициенты.

В случае отсутствия вынуждающей силы

$$f(x,t) \equiv 0$$

имеем однородное уравнение:

$$Ru \equiv \eta \frac{\partial^3 u}{\partial x^2 \partial t} + M \frac{\partial^2 u}{\partial x^2} - \rho \frac{\partial^2 u}{\partial t^2} - H \frac{\partial u}{\partial t} = 0, \qquad \text{(C)}$$

которое может быть представлено в виде

$$\frac{\partial^2 v}{\partial x^2} = \frac{\partial w}{\partial t}, \qquad (12)$$

если положить

$$v = \eta \frac{\partial u}{\partial t} + Mu \quad , \; w = \rho \frac{\partial u}{\partial t} + Hu \qquad (13)$$

Рассмотрим тот случай, когда

$$\eta H - \rho M \neq 0. \qquad (14)$$

Очевидно, что посредством (13) $v$ и $w$ однозначно определяются во всех тех точках плоскости независимых переменных $x$ и $t$, в которых известны значения $u$ и $t$, в которых известны значения $u$ и $\frac{\partial u}{\partial t}$.

С другой стороны, пусть $P$ и $Q$ – две какие-угодно функции переменных $x, t$ пусть $P$ и $Q$ имеют для всех значений $x$ и $t$, непрерывные производные: $P$ – до второго, по $x$, $Q$ – до первого, по $t$, порядка включительно.

Рассмотрим область $S$:

$$0 \le t \le +\infty, \; 0 \le x \le l \qquad (15)$$

и контур $C$, ее ограничивающий. Будем считать, что положительная ось $t$ направлена влево, если смотреть вдоль положительной оси $x$.

Применение формулы Гаусса к области $S$, а также к области $S_+$ являющейся дополнением к $S$ относительно всей плоскости $x,t$ дает следующие два соотношения:

$$\iint\limits_{S}\left\{\frac{\partial^2 P}{\partial x^2}-\frac{\partial Q}{\partial t}\right\}dxdt=\int\limits_{\underline{C}}\frac{\partial P}{\partial x}dt+Qdx\ ,$$

$$\iint\limits_{S_+}\left\{\frac{\partial^2 P}{\partial x^2}-\frac{\partial Q}{\partial t}\right\}dxdt=\int\limits_{\underline{C}}\frac{\partial P}{\partial x}dt+Qdx\,. \tag{16}$$

Положим

$$P=\frac{1}{2}v^2,\ Q=\int\limits_0^t v\frac{\partial w}{\partial t}dt\ . \tag{17}$$

Тогда

$$\frac{\partial P}{\partial x}=v\frac{\partial v}{\partial x}\ ,\quad \frac{\partial^2 P}{\partial x^2}=\left(\frac{\partial v}{\partial x}\right)^2+v\frac{\partial^2 v}{\partial x^2}=\left(\frac{\partial v}{\partial x}\right)^2+v\frac{\partial w}{\partial t},\ \ \frac{\partial Q}{\partial x}=v\frac{\partial w}{\partial t}\ .$$

Чтобы иметь право применить формулу (16) как к $S$, так и к $S_+$, достаточно предположить, что $\frac{\partial v}{\partial x}$ непрерывна для всех значений $x$ и $t$, и что $v$ и $\frac{\partial w}{\partial t}$ непрерывны при $t=0, t=+\infty$, $x=0$ и $x=l$. Подстановка (17) в (16) дает тогда:

$$\iint\limits_{S}\left(\frac{\partial v}{\partial x}\right)^2dxdt=\int\limits_0^{+\infty}\left[v\frac{\partial v}{\partial x}\right]_{x=l}dt-\int\limits_0^{+\infty}\left[v\frac{\partial v}{\partial x}\right]_{x=0}dt+\int\limits_0^{l}dx\int\limits_0^{0}v\frac{\partial w}{\partial t}dt-\int\limits_0^{l}dx\int\limits_0^{+\infty}v\frac{\partial w}{\partial t}dt\ .$$

Теперь можно показать, что если $v$ и $w$ исчезают на контуре С, то они исчезают тождественно для всех значений $x$ и $t$.

Действительно, рассмотрим правую часть первой из написанных формул. Так как $v|_{x=l}=v|_{x=0}=0$, то первые два интеграла исчезают. Очевидно исчезновение и третьего интеграла. Остальное дает:

$$\iint\limits_{S}\left(\frac{\partial v}{\partial x}\right)^2dxdt=-\int\limits_0^{l}dx\int\limits_0^{+\infty}v\frac{\partial w}{\partial t}dt\quad .$$

Но

$$v\frac{\partial w}{\partial t} = \left[\eta\frac{\partial u}{\partial t} + Mu\right]\left[\rho\frac{\partial^2 u}{\partial t^2} + H\frac{\partial u}{\partial t}\right] = \frac{\partial}{\partial t}\left[\frac{\eta\rho}{2}\left(\frac{\partial u}{\partial t}\right)^2 + \frac{HM}{2}u^2 + \rho Mu\frac{\partial u}{\partial t}\right] + ,$$

$$+(\eta H - \rho M)\left(\frac{\partial u}{\partial t}\right)^2,$$

следовательно,

$$\int_0^{+\infty} v\frac{\partial w}{\partial t}dt = \left[\frac{\eta\rho}{2}\left(\frac{\partial u}{\partial t}\right)^2 + \frac{HM}{2}u^2 + \rho Mu\frac{\partial u}{\partial t}\right]_{t=0}^{t=+\infty} + (\eta H - \rho M)\int_0^{+\infty}\left(\frac{\partial u}{\partial t}\right)^2 dt .$$

Но выражение, стоящее в квадратных скобках, обращается в нуль при $t = 0$ и $t = +\infty$, так как, вследствие исчезновения $u$ и $w$ на контуре $C$, имее

$$u\big|_{t=0} = u\big|_{t=+\infty} = 0, \qquad \frac{\partial u}{\partial t}\Big|_{t=0} = \frac{\partial u}{\partial t}\Big|_{t=+\infty} = 0.$$

Поэтому

$$-\int_0^l dx\int_0^{+\infty} v\frac{\partial w}{\partial t}dt = -(\eta H - \rho M)\int_0^l dx\int_0^{+\infty}\left(\frac{\partial u}{\partial t}\right)^2 dt .$$

Формулы (16) принимают вид:

$$\iint_S \left(\frac{\partial v}{\partial x}\right)^2 dxdt = -(\eta H - \rho M)\int_0^l dx\int_0^{+\infty}\left(\frac{\partial u}{\partial t}\right)^2 dt ,$$

$$\iint_{S_+} \left(\frac{\partial v}{\partial x}\right)^2 dxdt = +(\eta H - \rho M)\int_0^l dx\int_0^{+\infty}\left(\frac{\partial u}{\partial t}\right)^2 dt.$$

Почленное сложение их дает:

$$\int\int_{S+S_+} \left(\frac{\partial v}{\partial x}\right)^2 dxdt = 0 ,$$

где $S + S_+$ обозначает всю плоскость переменных $x, t$. Так как все члены интеграла слева неотрицательны, то отсюда следует:

$$\frac{\partial v}{\partial x} \equiv 0 .$$

Значит, $v$ может зависеть только от $t$. Однако при $x = 0$ мы имеем $v = 0$, следовательно,

$$v \equiv 0 .$$

Поэтому $\frac{\partial^2 v}{\partial x^2} \equiv 0$. Но $\frac{\partial^2 v}{\partial x^2} = \frac{\partial w}{\partial t}$, а это значит, что $w$ может быть функцией только $x$; при $t = 0$ мы имеем $w = 0$. Отсюда следует, что

$$w \equiv 0,$$

*что и утверждалось.*

На основании замечания, сделанного по поводу уравнений (13) при условии (14), доказанное положение может быть сформулировано так:

Если $u$ есть решение уравнения (C), непрерывное и обладающее непрерывными производными $\frac{\partial u}{\partial t}, \frac{\partial^2 u}{\partial t^2}, \frac{\partial u}{\partial x}$ и $\frac{\partial^2 u}{\partial x \partial t}$ для всех значений $x$ и $t$ и обращающееся в нуль при $x = 0$ и при $x = l$ для всех $t \geq 0$, а также при $t = 0$ и при $t = +\infty$ для всех $x$ из интервала $0 \leq x \leq l$, то $u \equiv 0$.

Пусть теперь заданы какие-нибудь значения $u$ и $\frac{\partial u}{\partial t}$ на контуре С. Это значит, что на этом контуре заданы значения $v$ и $w$. В таком случае существует единственное решение $u$ уравнения (C), непрерывное вместе со своими производными до второго порядка включительно во всей плоскости переменных $x$ и $t$ и принимающее на С вместе с $\frac{\partial u}{\partial t}$ заданные последовательности значений.,

В самом деле, если бы таких решений было два, например, $u_1$ и $u_2$, то в силу дистрибутивных свойств оператора $R$ разность $u = u_1 - u_2$ также была бы решением для (C).Она удовлетворяла бы поставленным условиям непрерывности и обращались бы вместе со своей производной $\frac{\partial u}{\partial t}$ в нуль для точек контура С, не будучи тождественна равна нулю, а это находится в противоречии с доказанным положением.

Обращаясь теперь к неоднородному уравнению, покажем, что *уравнение (B) при условии (14) имеет единственное решение непрерывное со своими производными первого и второго порядка и удовлетворяющее условиям:*

$$u\big|_{t=0} = \Phi(x), \quad \frac{\partial u}{\partial t}\Big|_{t=0} = \varphi(x),$$

$$0 \leq x \leq l;$$

$$u\big|_{t=+\infty} = \Psi(x), \quad \frac{\partial u}{\partial t}\Big|_{t=+\infty} = \psi(x)$$

$$u\big|_{x=0} = \mathrm{X}(t), \quad \frac{\partial u}{\partial t}\Big|_{x=0} = \chi(x), \qquad t \geq 0,$$

$$u\big|_{x=l} = \Omega(t), \quad \frac{\partial u}{\partial t}\Big|_{x=l} = \omega(x).$$

*причем, конечно,*

$$\Phi(0)=\mathrm{X}(0)\ ,\quad \varphi(0)=\chi(0);\quad \Phi(l)=\Omega(0),\quad \varphi(l)=\omega(l);$$

$$\Psi(0)=\mathrm{X}(+\infty)\,,\ \psi(0)=\chi(+\infty);\quad \Psi(l)=\Omega(+\infty)\,,\ \psi(l)=\omega(+\infty)\,.$$

*Действительно,* если бы существовало два различных решения $u_1$ и $u_2$ Уравнения (B), то однородное уравнение (C) имело бы неисчезающее тождественно решение $u=u_1-u_2$, которое удовлетворяло бы условиям непрерывности и обращалось бы вместе с $\frac{\partial u}{\partial t}$ в нуль на границе области

$$0\le t,\ \ 0\le x\le l\ ,$$

чего, как мы видели, не может быть.

До сих пор речь шла о том случае, когда $\eta H-\rho M\neq 0$.

Пусть теперь $\eta H-\rho M=0$. Тогда, положив

$$\rho:\eta=H:M=\lambda,\ \lambda>0,$$

мы имели бы:

$$w=\lambda v\,,$$

и уравнение (C) приняло бы вид:

$$\frac{\partial^2 v}{\partial x^2}=\lambda\frac{\partial v}{\partial t}\,. \tag{18}$$

Относительно уравнения (18) доказано,[1)] что если его интеграл $v$, непрерывный вместе со своими первыми производными, исчезает для любого $x$ из интервала $0\le x\le l$ при $t\ge 0$ как при $x=0$, так и при $x=l$, то $v\equiv 0$ внутри области $S:0\le t,\ \ 0\le x\le l$. Поэтому тогда

$$\eta\frac{\partial u}{\partial t}+Mu\equiv 0\,,$$

откуда следует, что

$$u=\alpha e^{-\frac{M}{\eta}t}\,,$$

где $\alpha$ может зависеть только от $x$. Так как при $t=0$ должно быть $u=0$, то $\alpha=0$. В таком случае $u\equiv 0$. Остальные рассуждения по поводу ненулевых граничных условий и по поводу неоднородного уравнения могут быть буквально повторены и в этом случае, поэтому выводы этого параграфа нужно считать общими.

Связь излагаемой в этой работе теории с уравнением (18) теплопроводности не может быть случайной; во время вязкого трения слоев нити друг о друга вырабатывается тепло и процесс передачи последнего вдоль нити не может быть не связанным с процессом движения слоев нити.

Назовем интеграл уравнения (B) *правильным* в некоторой области изменения переменных $x,t$, если этот интеграл вместе со своими производными до второго порядка включительно непрерывен в данной области.

*[1)] Э.Гурса , Курс математического анализа, т.III, ч. I, ГТТИ, стр.260-262, 1933.*

В таком случае можно считать доказанным, что существует[1] единственный правильный во всей плоскости переменных $x, t$ интеграл $u$ уравнения (B), принимающей на контуре C вместе с $\frac{\partial u}{\partial t}$ заданные последовательности значений. Более того, на основании формулы

$$\iint\limits_{S}\left(\frac{\partial v}{\partial x}\right)^2 dxdt = -(\eta H - \rho M)\int\limits_0^l dx\int\limits_0^{+\infty}\left(\frac{\partial u}{\partial t}\right)^2 dt = -(\eta H - \rho M)\iint\limits_{S}\left(\frac{\partial u}{\partial t}\right)^2 dxdt$$

при $(\eta H - \rho M) > 0$ имеем

$$\frac{\partial u}{\partial t} \equiv 0 \quad , \ 0 \le t \ , \ 0 \le x \le l \, .$$

Как и выше, отсюда мы заключаем, что правильное в $S$ решение $u$ однородного уравнения (C) исчезает тождественно внутри S, если оно равно нулю на границе $S$ вместе с $\frac{\partial u}{\partial t}$ и если $(\eta H - \rho M) > 0$. Значит, при $(\eta H - \rho M) > 0$ неоднoгродное уравнение (B) имеет единственное правильное в $S$ решение, принимающее на границе $S$ вместе с $\frac{\partial u}{\partial t}$ заданные значения. Непосредственное доказательство единственности правильного в $S$ интеграла уравнения (B), принимающего вместе с $\frac{\partial u}{\partial t}$ заданные значения на C, для случая $(\eta H - \rho M) < 0$ было бы довольно затруднительным. Однако из формулы

$$\iint\limits_{S_+}\left(\frac{\partial v}{\partial x}\right)^2 dxdt = +(\eta H - \rho M)\int\limits_0^l dx\int\limits_0^{+\infty}\left(\frac{\partial u}{\partial t}\right)^2 dt = +(\eta H - \rho M)\iint\limits_{S}\left(\frac{\partial u}{\partial t}\right)^2 dxdt,$$

$$(\eta H - \rho M) < 0$$

вследствие положительности обоих интегралов и подинтегральных выражений выводим опять

$$\frac{\partial u}{\partial t} \equiv 0 \text{ при } 0 \le t \ , \ 0 \le x \le l \, ,$$

где под $u$ разумеется интеграл однородного уравнения (C). Отсюда, как и раньше, получаем, что общее уравнение (B) и в случае $(\eta H - \rho M) < 0$ имеет единственное, правильное в области $0 \le t$, $0 \le x \le l$ решение касается случая , которое принимает на границе этой области вместе с $\frac{\partial u}{\partial t}$ заданную последовательность значений. Что $(\eta H - \rho M) = 0$, то, как мы видели, он не заставляет сомневаться в единственности решения; нужно заметить, что в этом случае задание $u$ и $\frac{\partial u}{\partial t}$ при $t = +\infty$ является излишним.

[1] *Факт существования решения будет установлен в следующих параграфах. Изложенное позволяет лишь утверждать, что может быть не больше одного решения, удовлетворяющего условиям текста.*

<u>*Примечание.*</u> Из соображений, приведенных выше, следует, что интеграл уравнения

$$Ru \equiv \eta \frac{\partial^3 u}{\partial x^2 \partial t} + M \frac{\partial^2 u}{\partial x^2} - \rho \frac{\partial^2 u}{\partial t^2} - H \frac{\partial u}{\partial t} = -f(x,t) \qquad \text{(C)}$$

однозначно определен во всей плоскости переменных $x, t,$ если известны значения $u, \frac{\partial u}{\partial t}$ на контуре C:

$$u\big|_{t=0} = \Phi(x), \quad \frac{\partial u}{\partial t}\Big|_{t=0} = \varphi(x),$$

$$u\big|_{t=+\infty} = \Psi(x), \quad \frac{\partial u}{\partial t}\Big|_{t=+\infty} = \psi(x),$$

$$u\big|_{x=0} = \mathrm{X}(t), \quad \frac{\partial u}{\partial t}\Big|_{x=0} = \chi(t),$$

$$u\big|_{x=l} = \Omega(t), \quad \frac{\partial u}{\partial t}\Big|_{x=l} = \omega(t).$$

Заметим, во-первых, что из

$$u\big|_{x=0} = \mathrm{X}(t), \qquad u\big|_{x=l} = \Omega(t)$$

непосредственным дифференцированием по $t$ получается

$$\frac{\partial u}{\partial t}\Big|_{x=0} = \mathrm{X}'(t) \qquad \frac{\partial u}{\partial t}\Big|_{x=l} = \Omega'(t).$$

Следовательно, должно быть

$$\chi(t) = \mathrm{X}'(t), \quad \omega(t) = \Omega'(t),$$

и $\chi(t)$ и $\omega(t)$ не могут быть заданы произвольно, если даны $\mathrm{X}(t)$ и $\Omega(t)$. Если даны две последние функции, то уже тем самым даны и значения их производных по $t$ при $x = 0$ и при $x = l$.

Во-вторых, условия при $t = +\infty$:

$$u\big|_{t=+\infty} = \Psi(x), \qquad \frac{\partial u}{\partial t}\Big|_{t=+\infty} = \psi(x)$$

производят несколько странное впечатление. Действительно, может ли быть речь об упругом последействии, если значения смещений $u$ заранее известны при $t = +\infty$? Как связать этот пункт рассуждений с физической постановкой задачи об упругом гистерезисе? Я утверждаю, что *условие*

$$u\big|_{t=+\infty} = \Psi(x)$$

*является лишним в приведенном выше доказательстве единственности решения.* Достаточно считать данным распределение скоростей при $t = +\infty$:

$$\frac{\partial u}{\partial t}\Big|_{t=+\infty} = \psi(x)\ .$$

Чтобы доказать это, достаточно предположить, что для однородного уравнения

$$Ru = 0$$

при $t = +\infty$ лишь $\dfrac{\partial u}{\partial t}$ исчезает для всех значений $x$ из интервала $(0, l)$. Тогда равенство

$$\int_0^l \int_0^{+\infty} \left(\frac{\partial v}{\partial x}\right)^2 dxdt = -\int_0^l dx \int_0^{+\infty} v \frac{\partial w}{\partial t} dt = -(\eta H - \rho M)\iint_S \left(\frac{\partial u}{\partial t}\right)^2 dxdt -$$

$$-\int_0^l \left[\frac{\eta\rho}{2}\left(\frac{\partial u}{\partial t}\right)^2 + \frac{HM}{2}u^2 + \rho M u \frac{\partial u}{\partial t}\right]_{t=0}^{t++\infty} dx$$

принимает вид

$$\int_0^l \int_0^{+\infty} \left(\frac{\partial v}{\partial x}\right)^2 dxdt = -(\eta H - \rho M)\int_0^l \int_0^{+\infty} \left(\frac{\partial u}{\partial t}\right)^2 dxdt = -\frac{HM}{2}\int_0^l \left(u^2\right)_{t=+\infty} dx\ .$$

При $(\eta H - \rho M) > 0$ оба члена справа отрицательны, в то время как слева имеем положительное качество. Равенство возможно при исчезновении всех членов. Далее, можно повторить приведенные в тексте соображения о единственности решения неоднородного уравнения.

Итак, можно считать доказанным, что уравнение

$$Ru = f(x, t)$$

имеет единственное решение $u$, которое на сторонах прямоугольника С удовлетворяет условиям:

$$u\big|_{t=0} = \Phi(x)\,, \qquad \frac{\partial u}{\partial t}\Big|_{t=0} = \varphi(x),$$

$$u\big|_{x=0} = \mathrm{X}(t)\,, \qquad u\big|_{x=l} = \Omega(t) \qquad\qquad (X)$$

$$\frac{\partial u}{\partial t}\Big|_{t=+\infty} = \psi(x)\,.$$

Последнее условие предполагает, что известно лишь распределение скоростей вдоль нити при $t = +\infty$.

С физической точки зрения дело обстоит дело именно так. Действительно, скорости движения слоев нити при $t = +\infty$ должны непременно исчезнуть, так что

$$\psi(x) \equiv 0$$

для всех $x$ из интервала $(0,l)$. Это следует хотя бы из уравнения для свободной энергии $E$ в изолированной системе, полученного Н.А.Наседкиным:

$$E = E_0 e^{-\alpha^2 t}\,,$$

из которого видно, что свободная энергия $E$ такой системы стремится к нулю при $t \to +\infty$, переходя в теплоту. А со свободной энергией стремятся к нулю и скорости движения слоев $\frac{\partial u}{\partial t}$.

Поэтому во всех дальнейших рассуждениях мы можем считать данным условие

$$\frac{\partial u}{\partial t}\Big|_{t=+\infty} \equiv 0\,,$$

обеспечивающее вместе с четырьмя прочими условиями (X) единственность решения задачи.

**§ 4.** Вернемся к однородному уравнению (C):

$$\eta\frac{\partial^3 u}{\partial x^2 \partial t} + M\frac{\partial^2 u}{\partial x^2} - \rho\frac{\partial^2 u}{\partial t^2} - H\frac{\partial u}{\partial t} = 0 \quad , \qquad \text{(C)}$$

где $u$ обозначает поперечное смещение частиц упруго-вязкой нити, закрепленной в точках $x = 0$ и $x = l$ и движущейся в сопротивляющейся среде. Так как вынуждающая внешняя сила отсутствует, то движение, описываемое уравнением (C), мы будем называть свободным колебанием, не предполагая периодичности последнего. Мы будем искать решение $u(x,t)$ уравнения (C), удовлетворяющее условиям:

$$u\big|_{t=0} = \Phi(x)\,, \quad \frac{\partial u}{\partial t}\Big|_{t=0} = \varphi(x),$$

$$u\big|_{x=0} = 0\,, \quad u\big|_{x=l} = 0 \quad , \qquad (19)$$

$$\frac{\partial u}{\partial t}\Big|_{t=+\infty} = 0 .$$

Как мы видели, может существовать не больше одного такого решения, правильного при $t \geq 0$ и при $0 \leq x \leq l$.

Пользуясь тем, что переменные $x$ и $t$ в (С) разделяются, я применяю способ Д.Бернулли.

Пусть

$$u(x,t) = XT , \qquad (20)$$

где $X, T$ зависят соответственно только от $x$ и от $t$. Уравнение (С) принимает вид:

$$\frac{X''}{X} = \frac{\rho T'' + HT'}{\eta T' + MT} .$$

Благодаря независимости переменных $x$ и $t$ каждая из частей написанного равенства может быть только постоянной. Обозначая последнюю $-\kappa$, получаем два обыкновенных уравнения:

$$X'' + \kappa X = 0 , \qquad (21)$$

$$\rho T'' + (\kappa\eta + H)T' + \kappa MT = 0 . \qquad (22)$$

Интеграл первого из них, исчезающий при $x = 0$, есть

$$X = a \sin \sqrt{x} x.$$

Так как, по условию задачи,X должно обращаться в нуль при $x = l$, то

$$\sqrt{\kappa} = \frac{m\pi}{l}, \; m \text{ – целое.} \qquad (23)$$

Не нарушая общности рассуждений, $m$ можно считать положительным, так как выбор знака определяется лишь выбором направления отсчета смещений. Тривиальное значение $m = 0$ считается исключенным.

Характеристическое уравнение для (22) есть

$$\rho q^2 + (\kappa\eta + H)Q + \kappa M = 0 ,$$

так что

$$q', q'' = \frac{-(\kappa\eta + H) \pm \sqrt{(\kappa\eta + H)^2 - 4\kappa\rho M}}{2\rho} , \qquad (24)$$

$$\kappa = \frac{m^2\pi^2}{l^2},\ m = 1,2,,...$$

$$T = C'e^{q't} + C''e^{q''t}$$

с произвольными постоянными $C'$ и $C''$ есть общий интеграл (22). Отсюда следует, что общий интеграл уравнения (С), обращающийся в нуль при $x = 0$ и при $x = l$, имеет вид:

$$u(x,t) = \sum_{m=1}^{\infty}\left(C'_m e^{q'_m t} + C''_m e^{q''_m t}\right)\sin\frac{m\pi x}{l}, \qquad (25)$$

где

$$q'_m, q''_m = -\frac{\frac{m^2\pi^2\eta}{b^2} + H}{2\rho} \pm \left[\frac{\left(\frac{m^2\pi^2\eta}{b^2} + H\right)^2}{4\rho^2} - \frac{m^2\pi^2 M}{\rho l^2}\right]^{\frac{1}{2}},$$

а постоянные $C', C''$ остаются пока произвольными. Мы их выберем, требуя соблюдения условий:

$$u\big|_{t=0} = \Phi(x), \quad \frac{\partial u}{\partial t}\Big|_{t=0} = \varphi(x).$$

Дифференцируя (25) почленно по $t$ (законность этого действия будет установлена в следующем параграфе), находим

$$\frac{\partial u}{\partial t} = \sum_{m=1}^{\infty}\left(q'_m C'_m e^{q'_m t} + q''_m C''_m e^{q''_m t}\right)\sin\frac{m\pi x}{l}, \qquad (26)$$

так что , положив в (25) и (26) $t = 0$, имеем:

$$u\big|_{t=0} = \sum_{m=1}^{\infty}(C'_m + C''_m)\sin\frac{m\pi x}{l},$$

$$\frac{\partial u}{\partial t}\Big|_{t=0} = \sum_{m=1}^{\infty}(q'_m C'_m + q''_m C'')\sin\frac{m\pi x}{l}. \qquad (27)$$

Предположим, с другой стороны, что

$$\Phi(x) = \sum_{m=1}^{\infty} A_m \sin\frac{m\pi x}{l}, \qquad (28)$$

$$\varphi(x) = \sum_{m=1}^{\infty} B_m \sin\frac{m\pi x}{l},$$

так что при $x = 0$ и при $x = l$

$$\Phi(0) = \varphi(0) = \Phi(l) = \varphi(l).$$

Тогда постоянные $A_m$ и $B_m$ определяется формулами:

$$A_m = \frac{2}{l}\int_0^l \Phi(x)\sin\frac{m\pi x}{l},$$

$$(m = 1,2,...) \qquad (29)$$

$$B_m = \frac{2}{l}\int_0^l \varphi(x)\sin\frac{m\pi x}{l} \quad .$$

Для того, чтобы соблюдались условия

$$u\big|_{t=0} = \Phi(x), \ \frac{\partial u}{\partial t}\Big|_{t=0} = \varphi(x), \qquad (30)$$

достаточно положить

$$C'_m + C''_m = A_m,$$

$$(m = 1,2,...) \qquad (31)$$

$$q'_m C'_m + q''_m C''_m = B_m .$$

При всех $m$, для которых $q'_m \neq q''_m$ (формул а 24) из (31) однозначно определяются коэффициенты $C'_m$ и $C''_m$ по данным $A_m$ и $B_m$. Но если существует такое целое $m = m_0 \geq 1$, для которого $q'_{m_0} = q''_{m_0}$, то найти $C'_{m_0}$ и $C''_{m_0}$ из (31) нельзя, да в этом и нет надобности, так как тогда достаточно знать сумму

$C'_{m_0} + C''_{m_0}$. Последняя может быть определена посредством любого из двух уравнений (31) при условии:

$$B_{m_0} = q_{m_0} A_{m_0} \quad (q_{m_0} = q'_{m_0} = q''_{m_0}).$$

Итак, поперечные колебания упруго-вязкой нити, закрепленной на концах и движущейся в сопротивляющейся среде, в однородном случае свободных колебаний определяется уравнением:

$$u(x,t) = \sum_{m=1}^{\infty}\left(C'_m e^{q'_m t} + C''_m e^{q''_m t}\right)\sin\frac{m\pi x}{l} \ , \qquad (32)$$

где

$$C'_m + C''_m = \int_0^l \Phi(x)\sin\frac{m\pi x}{l} \quad ,$$

$$q'_m C'_m + q''_m C''_m = \int_0^l \varphi(x)\sin\frac{m\pi x}{l}$$

и где $q'_m, q''_m$, будучи корнями уравнений

$$\rho q^2 + (\kappa_m \eta + H)q + \kappa_m M = 0, \ \kappa_m = \frac{m^2\pi^2}{l^2} \ , \ (m = 1,2,...),$$

имеют значения

$$q'_m, q''_m = -\frac{\dfrac{m^2\pi^2\eta}{l^2}+H}{2\rho} \pm \left[\frac{\left(\dfrac{m^2\pi^2\eta}{l^2}+H\right)^2}{4\rho^2} - \frac{m^2\pi^2 M}{\rho l^2}\right]^{\frac{1}{2}}.$$

Таким образом общее колебание складывается из частных вида:

$$\left(C'_m e^{q'_m t} + C''_m e^{q''_m t}\right)\sin\frac{m\pi x}{l},$$

Каждое из которых, не смотря на возможную апериодичность, мы будем называть обертоном соответствующего. Обертон первого порядка $(m=1)$ будет считаться основным тоном.

Из выражения для $q'_m, q''_m$ видно, что при возрастании $m$ значения $\frac{1}{|q'_m|}, \frac{1}{|q''_m|}$ представляют собой величины того же порядка малости, как и $\frac{1}{m^2}$.

**§ 5.** Теперь перейдем к вопросу о сходимости полученных рядов. Пусть введенные в предыдущем параграфе функции $\Phi(x), \varphi(x)$ непрерывны вместе с их производными $\Phi'(x), \varphi'(x)$ во всех точках интервала $0 \le x \le l$ и пусть их вторые производные $\Phi''(x), \varphi''(x)$ на том же интервале удовлетворяют условиям Дирихле, т.е. остаются конечными, имеют конечное число максимумов и минимумов и конечное число точек разрыва первого рода.

В таком случае ряды

$$\sum_{m=1}^{\infty}\left(C'_m + C''_m\right)\frac{m\pi}{l}\cos\frac{m\pi x}{l},$$

$$\sum_{m=1}^{\infty}\left(q'_m C'_m + q''_m C''_m\right)\frac{m\pi}{l}\cos\frac{m\pi x}{l}$$

представляют собой производные $\Phi'(x)$ и $\varphi'(x)$ во всех точках интервала $0 \le x \le l$. Так как последние сами непрерывны и имеют производные $\Phi''(x)$ и $\varphi''(x)$, удовлетворяющие условиям Дирихле в том же интервале, то коэффициенты их разложений

$$\left(C'_m + C''_m\right)\frac{m\pi}{l}, \left(C'_m q'_m + C''_m q''_m\right)\frac{m\pi}{l}$$

имеют порядок малости [1)] не ниже, чем $\frac{1}{m^2}$. Поэтому ряды, получающиеся из двух вышенаписанных путем почленного дифференцирования, равномерно сходятся и представляют собой соответственно разложения для $\Phi''(x)$ и $\varphi''(x)$.

[1)] *Крылов А.Н. , Лекции о приближенных вычислениях, стр. 172, 194, 2-е изд., 1933*

Итак, ряды

$$\sum_{m=1}^{\infty}(C'_m e^{q'_m t} + C''_m e^{q''_m t})\sin\frac{m\pi x}{l},$$

$$\sum_{m=1}^{\infty}\left(C'_m e^{q'_m t} + C''_m e^{q''_m t}\right)\left(-\frac{m^2\pi^2}{l^2}\right)\sin\frac{m\pi x}{l},$$

$$\sum_{m=1}^{\infty}(C'_m q'_m e^{q'_m t} + C''_m q''_m e^{q''_m t})\sin\frac{m\pi x}{l},$$

$$\sum_{m=1}^{\infty}\left(C'_m q'_m e^{q'_m t} + C''_m q''_m e^{q''_m t}\right)\left(-\frac{m^2\pi^2}{l^2}\right)\sin\frac{m\pi x}{l}$$

при $t = 0$ соответственно равномерно сходятся к функциям

$$\Phi(x) = u\big|_{t=0},\ \Phi''(x) = \frac{\partial^2 u}{\partial x^2}\Big|_{t=0},\ \varphi(x) = \frac{\partial u}{\partial t}\Big|_{t=0},\ \varphi''(x) = \frac{\partial^3 u}{\partial x^2 \partial t}\Big|_{t=0}.$$

Характер корней $q'_m, q''_m$ позволяет утверждать, что при любом $t > 0$ и при каком угодно $A$

$$\left|Ae^{q'_m t}\right| < |A|,\ \left|Ae^{q''_m t}\right| < |A|.$$

Следовательно, выписанные выше четыре ряда, сходясь равномерно при $t = 0$, тем более сходятся равномерно при любом $t > 0$. Кроме того, ряд

$$\sum_{m=1}^{\infty}\left(C'_m q'^2_m e^{q'_m t} + C''_m q''^2_m e^{q''_m t}\right)\sin\frac{m\pi x}{l}$$

имеет при всяком $t \geq 0$ коэффициенты при $\sin\frac{m\pi x}{l}$ того же порядка малости, что и ряд

$$\sum_{m=1}^{\infty}\left(C'_m q'_m e^{q'_m t} + C''_m q''_m e^{q''_m t}\right)\left(-\frac{m^2\pi^2}{l^2}\right)\sin\frac{m\pi x}{l}.$$

А так как последний, как мы видели, сходится равномерно при всяком $t \geq 0$, то и первый равномерно сходится для всех $t \geq 0$. Но этот ряд получается из

$$\sum_{m=1}^{\infty}\left(C'_m q'_m e^{q'_m t} + C''_m q''_m e^{q''_m t}\right)\sin\frac{m\pi x}{l}$$

путем почленного дифференцирования по $t$; следовательно, такое дифференцирование законно и ряд

$$\sum_{m=1}^{\infty}\left(C'_m q'^2_m e^{q'_m t} + C''_m q''^2_m e^{q''_m t}\right)\sin\frac{m\pi x}{l}$$

при всяком $t \geq 0$ равномерно сходится к $\frac{\partial^2 u}{\partial t^2}$.

Мы доказали, таким образом, что при сделанных предположениях относительно функций $\Phi(x)$ и $\varphi(ч)$, найденное в § 4 *выражение для смещения* $u(x,t)$, которое мы получили в виде бесконечного ряда, *допускает почленное дифференцирование по x и t столько раз, сколько нужно для составления всех производных, входящих в уравнение (С).*

Приведенные здесь рассуждения тем более имеют полную силу в том случае, когда $\Phi''(x)$ непрерывна.

**§ 6.** В § 4 было найдено, что уравнение

$$\eta\frac{\partial^3 u}{\partial x^2\partial t}+M\frac{\partial^2 u}{\partial x^2}-\rho\frac{\partial^2 u}{\partial t^2}-H\frac{\partial u}{\partial t}=0 \quad , \qquad \text{(С)}$$

при условиях

$$u\big|_{t=0}=\Phi(x), \qquad \frac{\partial u}{\partial t}\Big|_{t=0}=\varphi(x),$$

$$u\big|_{x=0}=0, \qquad u\big|_{x=l}=0 \quad , \qquad \frac{\partial u}{\partial t}\Big|_{t=+\infty}=0 \qquad (19)$$

имеет решение .

$$u(x,t)=\sum_{m=1}^{\infty}\left(C'_m e^{q'_m t}+C''_m e^{q''_m t}\right)\sin\frac{m\pi x}{l} \qquad (32)$$

где постоянные $C'_m, C''_m, q'_m, q''_m$ имеют значения, указанные раньше.

Дифференцируя (32) по $t$ и $x$ столько раз, сколько нужно для составления левой части (С), подставляя найденные значения производных в последнюю и собирая члены с одинаковыми значениями $\sin\frac{m\pi x}{l}$, получим выражение

$$\sum_{m=1}^{\infty}\left\{\left[\rho q'^2_m+\left(\kappa_m\eta+H\right)q'_m+\kappa_m M\right]C'_m e^{q'_m t}+\left[\rho q''^2_m+\left(\kappa_m\eta+H\right)q''_m+\kappa_m M\right]C''_m e^{q''_m t}\right\},$$

где

$$\kappa_m=\frac{m^2\pi^2}{l^2}.$$

Так как $q'_m, q''_m$ – корни уравнения

$$\rho q_m^2+\left(\kappa_m\eta+H\right)q_m+\kappa_m M=0,$$

то каждый член суммы (33) равен нулю. Следовательно, (32) действительно удовлетворяет уравнению (C).

Вследствие того, что

$$\sin\frac{m\pi x}{l}\Big|_{x=0}=\sin\frac{m\pi x}{l}\Big|_{x=l}=0$$

для всех целых $m$ как выражение (32), так и выражение

$$\frac{\partial u}{\partial t}=\sum_{m=1}^{\infty}(C'_m q'_m e^{q'_m t}+C''_m q''_m e^{q''_m t})\sin\frac{m\pi x}{l} \qquad (34)$$

обращаются в нуль при $x=0$ и при $x=l$, каково бы ни было $t$.

Далее, полагая в (32) и (34) $t=0$, находим:

$$u\big|_{t=0}=\sum_{m=1}^{\infty}(C'_m+C''_m)\sin\frac{m\pi x}{l},$$

$$\frac{\partial u}{\partial t}\Big|_{t=0}=\sum_{m=1}^{\infty}(C'_m q'_m+C''_m q''_m)\sin\frac{m\pi x}{l};$$

вследствие выбора коэффициентов $C'_m, C''_m$:

$$C'_m+C''_m=\frac{2}{l}\int_0^l \Phi(x)\sin\frac{m\pi x}{l}dx,\quad C'_m q'_m+C''_m q''_m=\frac{2}{l}\int_0^l \varphi(x)\sin\frac{m\pi x}{l}dx$$

Эти две величины оказываются соответственно равными $\Phi(x), \varphi(x)$.

Наконец, положив в (34) $t=+\infty$, благодаря свойствам $q'_m, q''_m$ найдем:

$$\frac{\partial u}{\partial t}\Big|_{t=+\infty}=0,\ 0\le x\le l.$$

Следовательно, найденное решение (32) задачи удовлетворяет как уравнению (C), так и всем сопровождающим его условиям. Соображения § 3 позволяют утверждать, что другого решения того же уравнения с теми же условиями не существует.

**§ 7.** Каков же физический смысл полученного решения? Положим

$$\sigma_m=\frac{\dfrac{m^2\pi^2\eta}{l^2}+H}{2\rho},$$

$$(m=1,2,\dots) \qquad (35)$$

$$\nu_m=+\left[\frac{\left(\dfrac{m^2\pi^2\eta}{l^2}+H\right)^2}{4\rho^2}-\frac{m^2\pi^2 M}{\rho l^2}\right]^{\frac{1}{2}},$$

так что

$$q_m = -\sigma_m + \nu_m\,, \qquad q''_m = -\sigma_m - \nu_m \tag{36}$$

Тогда решение (32) можно представить в таком виде:

$$u(x,t) = \sum_{m=1}^{\infty} e^{-\sigma_m t}(C'_m e^{\nu_m t} + C''_m e^{-\nu_m t})\sin\frac{m\pi x}{l}\quad , \tag{37}$$

так что смещение $u_m(x,t)$, получающееся от обертонов $m$ -ого порядка, есть

$$u_m(x,t) = e^{-\sigma_m t}(C'_m e^{\nu_m t} + C''_m e^{-\nu_m t})\sin\frac{m\pi x}{l}\,. \tag{38}$$

Предположим сначала, что для выбранного $m$

$$\frac{\left(\dfrac{m^2\pi^2\eta}{l^2} + H\right)^2}{4\rho^2} - \frac{m^2\pi^2 M}{\rho l^2} < 0\,,$$

как это может случиться, если $M$ очень велико по сравнению с $\eta$ и $H$.

Тогда, положив

$$\nu_m = +i\,\nu'_m\,, \qquad \left[\nu'_m > 0, i^2 = -1\right],$$

мы делаем подстановку

$$C'_m + C''_m = A_m \sin B_m\,,\; C'_m - C''_m = -iA_m \cos B_m\,; \tag{39}$$

после этого оказывается, что для данного $m$

$$u_m(x,t) = A_m e^{-\sigma_m t}\sin\left(\nu'_m t + B_m\right)\sin\frac{m\pi x}{l}\,.$$

Отсюда видно, что в рассматриваемом случае обертон $m$-го порядка является таковым в узком смысле слова, будучи периодическим движением с периодом

$$\tau_m = \frac{2\pi}{\nu'_m} = \frac{2\pi}{\left|\left[\dfrac{m^2\pi^2\eta}{\rho l^2} - \sigma_m^2\right]^{\frac{1}{2}}\right|}\,, \tag{40}$$

где

$$\sigma_m = \left[\frac{m^2\pi^2\eta}{l^2} + H\right] : 2\rho \tag{41}$$

Этот период $\tau_m$ тем меньше, чем выше порядок $m$ обертона.

С другой стороны, амплитуда колебания частицы нити, с абсциссой $x$, имеет величину

$$a_m = A_m e^{-\sigma_m t}\sin\frac{m\pi x}{l}\,,$$

которая, вследствие положительности $\sigma_m$, неограниченно убывает с возрастанием $t$ и притом тем сильнее, чем выше порядок обертона $m$. Логарифмический декремент затухания дается формулой (41).

В том случае, когда для данного $m$ корни $q'_m, q''_m$ действительны и различны, соответствующий обертон апериодичен и может быть так назван лишь в расширенном толковании этого слова. Уравнения (35) и (36) показывают, что в этом случае

$$q'_m < 0, q''_m < 0, |q''_m| > |q'_m|.$$

Поэтому смещение $u_m(x,t)$, вызываемое таким обертоном, для всех частиц нити с течением времени стремится к нулю и тем более энергично, чем больше $m$. Аналогичное явление происходит и в том случае, когда $q'_m = q''_m$.

Может показаться на первый взгляд, что полученные результаты совершенно те же, к каким приводит теория колебаний упругих струн с учетом сопротивления, пропорционального первой степени скорости. В обоих случаях как периодические обертоны, так и обертоны, не имеющие периода, затухают с течением времени. Но в то время, как все обертоны, чисто упругой струны затухают равномерно, имея общий декремент, зависящий только от коэффициента сопротивления $H$, декременты затухания обертонов упруго-вязкой нити, завися от $H$ и от коэффициента вязкости $\eta$, существенно зависят от порядка обертона. Чем выше порядок обертона, тем быстрее этот обертон затухает, так что для достаточно большого промежутка времени, протекшего после начала движения, нить можно считать колеблющейся только в основном тоне (который может оказаться и апериодическим). Таким образом, благодаря исключительной роли внутреннего трения, в процессе колебания упруго-вязких нитей происходит явление, которое можно было бы назвать *элевтерозом* и которое состоит в постепенном освобождении основного тона от сопутствующих обертонов. Это явление интересно в том отношении, что оно позволяет писать решение уравнения (C) не в виде (32) или (37), а в более простом:

$$u(x,t) = (C'_1 e^{q'_1 t} + C'_1 e^{q''_1 t}) \sin \frac{\pi x}{l} \qquad (42)$$

если только $t$ достаточно велико для того, чтобы можно было пренебречь влиянием всех обертонов. И во всяком случае можно найти такое конечное целое $N > 0$, что сумма

$$\sum_{m=1}^{N} (C'_m e^{q'_m t} + C''_m e^{q''_m t}) \sin \frac{m\pi x}{l}$$

достаточно точно представляет величину смещения частиц нити, каковы бы ни $x$ и $t$, $0 \le x \le l$, $t \ge 0$.

**§ 8.** Результаты предыдущего параграфа позволяют поставить очень простой опыт, посредством которого можно установить значение коэффициента вязкости $\eta$ и модуля сдвига $\mu$ для вещества нити. В самом деле, благодаря элевтерозу можно без значительной ошибки принять, что по прошествии достаточно большого промежутка времени $t_0$ после начала движения при $t \ge t_0$ нить движется в основном тоне

$$u(x,t)=u_1(x,t)=(C_1'e^{q_1't}+C_1'e^{q_1''t})\sin\frac{\pi x}{l}\ , \qquad (42)$$

где $\quad q_1'=-\sigma_1+\nu_1, \quad q_1''=-\sigma_1-\nu_1$

и $$\sigma_1=\frac{\dfrac{\pi^2\eta}{l^2}+H}{2\rho}, \nu_1=\left[\sigma_1^2-\frac{\pi^2 M}{\rho l^2}\right]^{\frac{1}{2}}.$$

С другой стороны, надлежащим выбором натяжения $T$ нити можно достигнуть того, что основной тон (42) окажется периодическим. Для этого нужно только сделать $T$ таким большим, чтобы $M=T+\mu$ удовлетворяло неравенству

$$\frac{\pi^2 M}{\rho l^2}>\frac{\left(\dfrac{\pi^2\eta}{l^2}+H\right)^2}{4\rho^2},$$

что, по крайней мере при небольших коэффициентах трения, всегда может быть достигнуто. Тогда уравнение (42) приобретает вид:

$$u(x,t)=A_1e^{-\sigma_1 t}\sin(\nu_1 t+B_1)\sin\frac{\pi x}{l}, \qquad (43)$$

где положено

$$A_1\sin B_1=C_1'+C_1'', \qquad -iA_1\cos B_1=C_1'-C'', \quad i^2=-1.$$

*Фотографируя колебание (43) какой-либо частицы нити (для выбранного значения* $x$*), мы можем измерить на получившемся фильме как величину логарифмического декремента затухания* $\sigma_1$*, так и циклическую частоту* $\nu_1$. Тогда из уравнений

$$2\rho\sigma_1=\frac{\pi^2\eta}{l^2}+H\ , \ \nu_1^2=\frac{\pi^2(T+\mu)}{\rho l^2}-\sigma_1^2\ , \qquad (44)$$

зная плотность $\rho$ вещества нити, ее длину $l$, ее натяжение $T$ и коэффициент внешнего трения $H$, можем вычислить[1)] коэффициент вязкости $\eta$ и модуль сдвига $\mu$.

**§ 9.** До сих пор, занимаясь однородным уравнением (С) § 4, мы предполагали нить свободной от влияния внешних сил и, в частности, невесомой. Теперь мы рассмотрим движение частиц нити, подверженной силе тяжести.

*[1)]Это выражение для логарифмического декремента* $\sigma_1$ *является обобщением приведенного в статье R.Becker'a, см.Zts.f.Physik,* $33_3$*,S.191,Formel(10a),1925.-У А.Лява (loc.cit.,стр.128) читаем: «Если бы существовали сопротивления только этого типа (т.е. пропорциональные быстроте изменения деформации.– А.Г.), то, как указал Кельвин, относительное уменьшение амплитуды в единицу времени было бы обратно пропорционально квадрату периода колебания; но ряд опытов с крутильными колебаниями проволок показал, что этот закон не подтверждается». [См.W.Thomson, Elasticity, Encycl.Brit.] Как видно из данных в тексте формул для* $\sigma_1$ *и* $\nu_1$ *соотношение между относительным уменьшением амплитуды в единицу времени и периодом колебания отличается от кельвинского.*

Последнюю будем считать направленной против положительного отсчета смещений, как в § 1 и 2. Продолжая считать, что движение происходит в среде, сопротивляющейся пропорционально первой степени скорости, и называя, как и раньше, коэффициент внешнего трения буквой $H$, напишем уравнение движения нити так:

$$\eta\frac{\partial^3 u}{\partial x^2 \partial t}+M\frac{\partial^2 u}{\partial x^2}-\rho\frac{\partial^2 u}{\partial t^2}-H\frac{\partial u}{\partial t}=\rho g \qquad \text{(D)}$$

при условиях

$$u\big|_{t=0}=\Phi(x), \quad \frac{\partial u}{\partial t}\Big|_{t=0}=\varphi(x),$$

$$u\big|_{x=0}=0, \quad u\big|_{x=l}=0, \quad \frac{\partial u}{\partial t}\Big|_{t=+\infty}=0. \qquad (45)$$

С точки зрения физической, смысл последнего из условий (45) ясен; ведь вязкость нити, так сказать, стремится уничтожить наличие градиентов скорости, т.е. благодаря отсутствию движения на концах свести к нулю все скорости и в середине.

Заметим, что уравнение (D) имеет частный интеграл:

$$-\frac{\rho g}{2M}x(l-x), \qquad (46)$$

удовлетворяющий первой группе условий (45) на концах нити. Поэтому общее решение, удовлетворяющее этой группе условий, есть

$$u(x,t)=-\frac{\rho g}{2M}x(l-x)+\sum_{m=1}^{\infty}(C'_m e^{q'_m t}+C''_m e^{q''_m t})\sin\frac{m\pi x}{l} \quad ; \qquad (47)$$

В самом деле, почленное дифференцирование (47) по $t$ дает:

$$\frac{\partial u}{\partial t}=\sum_{m=1}^{\infty}(q'_m C'_m e^{q'_m t}+q''_m C''_m e^{q''_m t})\sin\frac{m\pi x}{l}, \qquad (48)$$

Если ряд справа в (48) сходится равномерно для $t\geq 0$. Из (48) видно, что $\frac{\partial u}{\partial t}$ обращается в нуль при всяком $x$, когда $t=+\infty$.

Коэффициенты $C'_m$, $C''_m$ пока не определены. Выбираем их так, чтобы

$$C'_m+C''_m=\frac{2}{l}\int^{l}\left[\Phi(x)+\frac{\rho g}{2M}x(l-x)\right]\sin\frac{m\pi x}{l}dx,$$

$$(49)$$

$$q'_m C'_m+q''_m C''_m=\frac{2}{l}\int_0^l\varphi(x)\sin\frac{m\pi x}{l}dx,$$

Мы удовлетворяем первой группе условий (45), так как тогда, с одной стороны, имеем:

$$\Phi(x)+\frac{\rho g}{2M}x(l-x)=\sum_{m=1}^{\infty}(C'_m+C''_m)\sin\frac{m\pi x}{l},$$

$$\varphi(x)=\sum_{m=1}^{\infty}(q'_m C'_m+q''_m C''_m)\sin\frac{m\pi x}{l}, \tag{50}$$

а с другой – из (47) и (48) при $t=0$ находим:

$$u\big|_{t=0}=-\frac{\rho g}{2M}x(l-x)+\sum_{m=1}^{\infty}(C'_m+C''_m)\sin\frac{m\pi x}{l}$$

$$\frac{\partial u}{\partial t}\Big|_{t=0}=\sum_{m=1}^{\infty}(q'_m C'_m+q''_m C''_m)\sin\frac{m\pi x}{l}.$$

Далее, если $\Phi(x)$ и $\varphi(x)$ непрерывны вместе с их первыми производными и имеют вторые производные, удовлетворяющие условиям Дирихле для всех значений $x$ из интервала $0\le x\le l$ то, как мы уже видели в § 5. Можно путем почленного дифференцирования (47) нужное количество раз составить равномерно сходящиеся ряды для всех производных, входящих в уравнение (D).

Наконец, из (47) при $t=+\infty$ получается:

$$u\big|_{t=+\infty}=-\frac{\rho g}{2M}x(l-x).$$

Следовательно, упруго-вязкая тяжелая нить, приходя в состояние покоя спустя бесконечно-большой промежуток времени после начала движения, принимает форму параболы, а не цепной линии, как это можно было бы думать.

**§ 10.** Здесь мы рассмотрим еще более общий случай вынужденных колебаний упруго-вязкой нити. Относящееся к этому вопросу дифференциальное уравнение таково:

$$\eta\frac{\partial^3 u}{\partial x^2\partial t}+M\frac{\partial^2 u}{\partial x^2}-\rho\frac{\partial^2 u}{\partial t^2}-H\frac{\partial u}{\partial t}=-f(x,t), \tag{E}$$

где $f(x,t)$ – фактор вынуждающей силы. Этот фактор предположим заданным в виде ряда:

$$f(x,t)=\sum_{m=1}^{\infty}w_m(t)\sin\frac{m\pi x}{l}. \tag{51}$$

Мы будем искать решение уравнения (E) также в виде ряда:

$$u(x,t)=\sum_{m=1}^{\infty}\left(C'_m e^{q'_m t}+C''_m e^{q''_m t}+a_m(t)\right)\sin\frac{m\pi x}{l}, \tag{52}$$

так что, предполагая возможным почленное дифференцирование этого ряда нужное число раз, имеем:

$$\frac{\partial u}{\partial t}=\sum_{m=1}^{\infty}\left(q'_m C'_m e^{q'_m t}+q''_m C''_m e^{q''_m t}+\dot{a}_m(t)\right)\sin\frac{m\pi x}{l} \ ,$$

$$\frac{\partial^2 u}{\partial t^2}=\sum_{m=1}^{\infty}\left(q'^2_m C'_m e^{q'_m t}+q''^2_m C''_m e^{q''_m t}+\ddot{a}_m(t)\right)\sin\frac{m\pi x}{l} \ ,$$

$$\frac{\partial^2 u}{\partial x^2}=\sum_{m=1}^{\infty}\left(C'_m e^{q'_m t}+C''_m e^{q''_m t}+a_m(t)\right)\left(-\frac{m^2\pi^2}{l^2}\right)\sin\frac{m\pi x}{l} \tag{53}$$

$$\frac{\partial^3 u}{\partial x^2\partial t}=\sum_{m=1}^{\infty}\left(q'_m C'_m e^{q'_m t}+q''_m C''_m e^{q''_m t}+\dot{a}_m(t)\right)\left(-\frac{m^2\pi^2}{l^2}\right)\sin\frac{m\pi x}{l}$$

где $q'_m, q''_m$ – по-прежнему корни уравнения

$$\rho q_m^2+\left(\kappa_m\eta+H\right)q_m+\kappa_m M=0,$$

$$\kappa_m=\frac{m^2\pi^2}{l^2}, m=1,2,... \tag{54}$$

Подставляя (51) и (53) в (Е), приводим последнее к виду:

$$\sum_{m=1}^{\infty}\{\ \left[\rho q'^2_m+\left(\kappa_m\eta+H\right)q'_m+\kappa_m M\right]C'_m e^{q'_m t}+\left[\rho q''^2_m+\left(\kappa_m\eta+H\right)q''_m+\kappa_m M\right]C''_m e^{q''_m t}+$$

$$+[\rho\ddot{a}_m(t)+\left(\kappa_m\eta+H\right)\dot{a}_m(t)+\kappa_m M a_m(t)\ ]\ \}\sin\frac{m\pi x}{l}=$$

$$=\sum_{b=1}^{\infty}w_m(t)\sin\frac{m\pi x}{l}\ .$$

На основании (54) два члена под знаком суммы исчезают. Остальные члены отождествим, полагая для всех $m=1,2,...$

$$\rho\ddot{a}_m(t)+\left(\kappa_m\eta+H\right)\dot{a}_m(t)+\kappa_m M a_m(t)=w_m(t) \tag{55}$$

Корни характеристического уравнения для (55) определяются из (54). Они равны $q'_m$ и $q''_m$.

Заменив $\frac{d}{dt}$ на $D$, приведем (55) к виду[1]

$$\frac{1}{\rho}a'_m(t)+\frac{1}{\rho}a''_m(t)=a_m(t)=\left[\rho D^2+\left(\kappa_m\eta+H\right)D+\kappa_m M E\right]^{-1}w_m(t)\ ,$$

где $E$ – оператор тождественного преобразования и

$$a'_m(t)=\frac{1}{q''_m-q'_m}\left(D-q'_m E\right)^{-1}w_m(t),\quad a''_m(t)=\frac{1}{q'_m-q''_m}\left(D-q''_m E\right)^{-1}w_m(t), \tag{56}$$

т.е. $a'_m(t), a''_m(t)$ являются решениями неоднородных линейных уравнений первого порядка:

$$\frac{da'_m}{dt} - q'_m a'_m = \frac{w_m(t)}{q''_m - q'_m}, \quad \frac{da''_m}{dt} - q''_m a''_m = \frac{w_m(t)}{q'_m - q''_m}.$$

Легко видеть, что интегралы этих уравнений, обращающиеся в нуль при $t = 0$, соответственно равны:

$$a'_m(t) = \frac{-1}{q''_m - q'_m} e^{q'_m t} \int_0^t e^{-q'_m \theta} w_m(\theta) d\theta,$$

$$a''_m(t) = \frac{-1}{q'_m - q''_m} e^{q''_m t} \int_0^t e^{-q''_m \theta} w_m(\theta) d\theta,$$

так что

$$a_m(t) = \frac{1}{\rho}\left\{ \frac{-1}{q''_m - q'_m} \int_0^t e^{-q'_m (t-\theta)} w_m(\theta) d\theta + \frac{-1}{q'_m - q''_m} \int_0^t e^{-q''_m (t-\theta)} w_m(\theta) d\theta \right\}. \qquad (57)$$

Здесь существенным было предположение, что $q'_m \neq q''_m$.Однако и в случае $q'_{m_0} = q''_{m_0} = q_{m_0}$ применение символического способа интеграции (55) оказывается возможным; нужно только употребить подстановку:

$$\left(D - q_{m_0} E\right) a_{m_0}(t) = z_{m_0}, \quad \left(D - q_{m_0} E\right) z_{m_0} = \frac{w_{m_0}}{\rho},$$

чтобы получить для $a_{m_0}$ выражение:

$$a_{m_0}(t) = \frac{1}{\rho} e^{q_{m_0} t} \int_0^t e^{q_{m_0} \theta} d\theta \int_0^\theta e^{-q_{m_0}(\theta - \vartheta)} w_m(\vartheta) d\vartheta,$$

[1]Здесь штрихи не обозначают производных

обращающееся в 0 при $t = 0$.

На основании этих соображений напишем решение уравнения (E) в такой форме:

$$u(x,t) = \sum_{m=1}^{\infty} \left\{ \begin{aligned} &\left[ C'_m - \frac{1}{\rho\left(q''_m - q'_m\right)} \int_0^t e^{-q'_m \theta} w_m(\theta) d\theta \right] e^{q'_m t} + \\ &+ \left[ C''_m - \frac{1}{\rho\left(q'_m - q''_m\right)} \int_0^t e^{-q''_m \theta} w_m(\theta) d\theta \right] e^{q''_m t} \end{aligned} \right\} \sin\frac{m\pi x}{l} +$$

$$+ \frac{1}{\rho} \sin\frac{m_0 \pi x}{l} e^{q_{m_0} t} \int_0^t e^{q_{m_0} \theta} d\theta \int_0^\theta e^{-q_{m_0}(\theta - \vartheta)} w_{m_0}(\vartheta) d\vartheta. \qquad (58)$$

Здесь сумма в правой части распространяется лишь на те $m$, для которых $q'_m \neq q''_m$. Последний член, выделенный из суммы, относится к такому $m = m_0$, для которого $q''_{m0} = q'_{m_0} = q_{m_0}$. Если для всех $m$ корни $q'_m, q''_m$ характеристических уравнений (54) попарно различны, то этот последний член в выражении для $u(x,t)$ отсутствует.

**§ 11.** Найденное решение (58) уравнения (E) вынужденных колебаний обращается вместе со своей производной $\frac{\partial u}{\partial t}$ в нуль для всех $t$ при $x = 0$ и при $x = l$, так как

$$\sin\frac{m\pi 0}{l} = \sin\frac{m\pi l}{l} = 0$$

для всех целых $m$. Кроме того, легко убедиться в том, что

$$u\big|_{t=0} = \Phi(x), \quad \frac{\partial u}{\partial t}\Big|_{t=0} = \varphi(x),$$

если положить

$$C'_m + C''_m = \frac{2}{l}\int_0^l \Phi(x)\sin\frac{m\pi x}{l}dx, \; q'_m C'_m + q''_m C''_m = \frac{2}{l}\int_0^l \varphi(x)\sin\frac{m\pi x}{l}dx,$$

как мы это делали § 4. Действительно, при $t = 0$ интегралы исчезают как в $(58)$, так и выражении, получающемся из (58) дифференцированием по $t$.

Остальное дает разложение в ряды Фурье для $u$ и $\frac{\partial u}{\partial t}$ при $t = 0$.

Таким образом остается открытым лишь вопрос о сходимости ряда (58) и тех, которые получаются из него почленным дифференцированием столько раз, сколько нужно для составления левой части (E). Мы уже видели в § 5, что ряды

$$\sum_m^\infty \left(C'_m q'_m e^{q'_m t} + C''_m q''_m e^{q''_m t}\right)\sin\frac{m\pi x}{l},$$

$$\sum_m^\infty \left(C'_m q'^2_m e^{q'_m t} + C''_m q''^2_m e^{q''_m t}\right)\sin\frac{m\pi x}{l},$$

$$\sum_m^\infty \left(C'_m e^{q'_m t} + C''_m e^{q''_m t}\right)\left(-\frac{m^2\pi^2}{l^2}\right)\sin\frac{m\pi x}{l}$$

$$\sum_m^\infty \left(C'_m q'_m e^{q'_m t} + C''_m q''_m e^{q''_m t}\right)\left(-\frac{m^2\pi^2}{l^2}\right)\sin\frac{m\pi x}{l}$$

сходятся равномерно для всех $t \ge 0$ соответственно к функциям

$$\frac{\partial u}{\partial t}, \frac{\partial^2 u}{\partial t^2}, \frac{\partial^2 u}{\partial x^2}, \frac{\partial^3 u}{\partial x^2 t},$$

где $u$ обозначает смещение в случае отсутствия вынуждающей силы, если только $\Phi(x)$ и $\varphi(x)$ непрерывны с их первыми производными и если их вторые производные удовлетворяют условиям Дирихле на интервале $0 \le x \le l$. Следовательно, в данном случае подлежат рассмотрению ряды:

$$\sum_{m=1}^{\infty} \frac{\sin\frac{m\pi x}{l}}{\rho(q''_m - q'_m)} \int_0^t e^{q'_m(t-\theta)} w_m(\theta) d\theta,$$

$$\sum_{m=1}^{\infty} \frac{\sin\frac{m\pi x}{l}}{\rho(q'_m - q''_m)} \int_0^t e^{q''_m(t-\theta)} w_m(\theta) d\theta, \qquad (60)$$

А также те, которые получаются из них дифференцированием для составления всех производных в (E).

Пусть для всех $t \ge 0$ величина силы

$$f(x,t) = \sum_{m}^{\infty} w_m(t) \sin\frac{m\pi x}{l}, \qquad (61)$$

рассматриваемая как функция $x$, непрерывна в промежутке $0 \le x \le l$ и имеет там непрерывную первую производную. Ее вторая производная пусть удовлетворяет в этом промежутке условиям Дирихле. Тогда для всех $t \ge 0$ имеем

$$\frac{\partial f}{\partial x} = \sum_{m=1}^{\infty} w_m(t) \frac{m\pi}{l} \cos\frac{m\pi x}{l}, \qquad (62)$$

причем коэффициенты $w_m(t)\frac{m\pi}{l}$ этого ряда имеют порядок малости не ниже $\frac{1}{m^2}$.

Тогда, во-первых, ряд

$$\sum_{m=1}^{\infty} w_m(t)\left(-\frac{m^2\pi^2}{l^2}\right)\sin\frac{m\pi x}{l} \qquad (63)$$

сходится равномерно к функции $\frac{\partial^2 f}{\partial x^2}$ для всех $t \ge 0$, а, во-вторых, порядок малости $w_m(t)$ должен быть по крайней мере равен $m^{-3}$.

По этой же причине и ряды:

$$\sum_{m=1}^{\infty} \frac{\sin\frac{m\pi x}{l}}{\rho(q''_m - q'_m)} \int_0^t e^{q'_m(t-\theta)} w_m(\theta) d\theta,$$

$$\sum_{m=1}^{\infty} \frac{\sin\frac{m\pi x}{l}}{\rho(q''_m - q'_m)} \left[ w_m(t) + q'_m \int_0^t e^{q'_m(t-\theta)} w_m(\theta) d\theta \right],$$

$$\sum_{m=1}^{\infty} \frac{\sin\frac{m\pi x}{l}}{\rho(q''_m - q'_m)} \left(-\frac{m^2\pi^2}{l^2}\right) \int_0^t e^{q'_m(t-\theta)} w_m(\theta) d\theta, \qquad (64)$$

$$\sum_{m=1}^{\infty}\frac{\sin\frac{m\pi x}{l}}{\rho(q''_m-q'_m)}\left(-\frac{m^2\pi^2}{l^2}\right)\left[w_m(t)+q'_m\int_0^t e^{q'_m(t-\theta)}w_m(\theta)d\theta\right],$$

$$\sum_{m=1}^{\infty}\frac{\sin\frac{m\pi x}{l}}{\rho(q''_m-q'_m)}\left[\dot w_m(t)+q'_m w_m(t)+q'^2_m\int_0^t e^{q'_m(t-\theta)}w_m(\theta)d\theta\right],$$

а также ряды, получаемые из этих заменой $q''_m$ на $q'_m$ и обратно, равномерно сходятся при всех $t\ge 0$. Относительно первых четырех рядов это утверждение очевидно, так как порядок малости их коэффициентов при $\sin\frac{m\pi x}{l}$ не ниже порядка малости коэффициентов равномерно сходящегося ряда (63). Равномерная сходимость последнего ряда из (64) при всех $t\ge 0$ обеспечивается ограниченностью всех производных $\dot w_m(t)$, тем, что порядок малости $\dot w_m(t)$ не ниже $m^{-3}$, и, наконец, тем, что ряд (63) сходится равномерно.

Но (64) и аналогичные им ряды – это как раз те чести разложений функций

$$u(x,t),\ \frac{\partial u}{\partial t},\frac{\partial^2 u}{\partial t^2},\frac{\partial^2 u}{\partial x^2},\frac{\partial^3 u}{\partial x^2 t},$$

которые входят в уравнение (E) и которые получились бы путем почленного дифференцирования выражения (58) для $u(x,t)$ соответствующее количество раз. Принимая во внимание и равномерную сходимость рядов (59), приходим к заключению, что при предположениях, сделанных относительно функций $\Phi,\varphi$ и $f$, решение (58) уравнения (E) может быть почленно дифференцировано столько раз, сколько требуется для составления левой части (E). Что (58) есть действительно решение (E), можно было бы теперь непосредственно проверить.

Мы видели в начале этого параграфа, что это решение удовлетворяет условиям:

$$u\big|_{t=0}=\Phi(x),\qquad \frac{\partial u}{\partial t}\Big|_{t=0}=\varphi(x),$$

$$u\big|_{x=0}=0,\qquad u\big|_{x=l}=0.$$

Что касается четвертой степени стороны контура С, рассмотренного в § 3, то выражение (58) дифференцированием по $t$ при $t=+\infty$ дает:

$$\frac{\partial u}{\partial t}\Big|_{t=+\infty}=\sum_{m=1}^{\infty}\left(A'_m q'_m+A''_m q''_m\right)\sin\frac{m\pi x}{l}, \tag{65}$$

где

$$A'_m=\lim_{t\to+\infty}\left[\frac{1}{\rho(q''_m-q'_m)}\int_0^t e^{q'_m(t-\theta)}w_m(\theta)d\theta\right],$$

$$A''_m=\lim_{t\to+\infty}\left[\frac{1}{\rho(q'_m-q''_m)}\int_0^t e^{q''_m(t-\theta)}w_m(\theta)d\theta\right],$$

то $\quad \frac{\partial u}{\partial t}\Big|_{t=+\infty}=0.$

На основании теорем § 3 уравнение (58) §10 дает *единственное решение задачи.*

На этом заканчиваем первую часть своей работы, посвященную поперечным колебаниям нити, с тем, чтобы во второй части рассмотреть вопрос о *крутильных колебаниях*, когда ось нити во все время движения сохраняет свою прямолинейную форму.

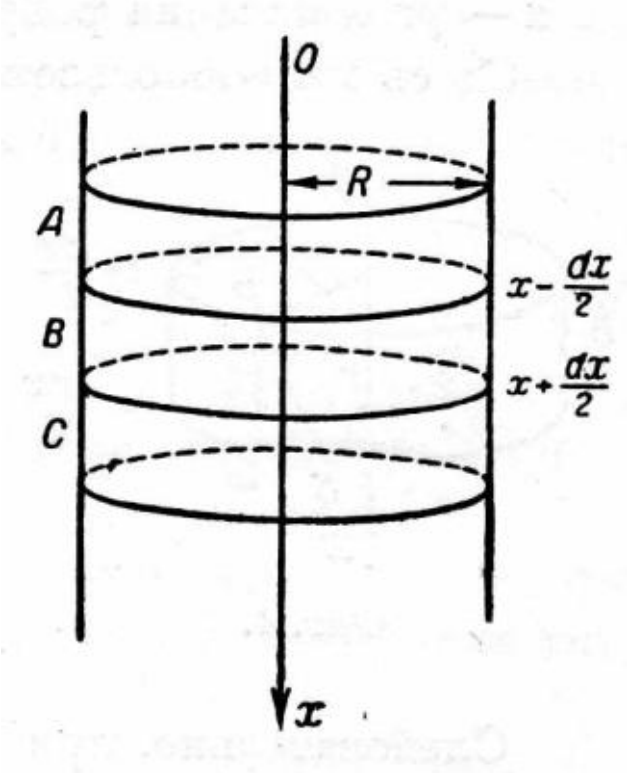

Фиг. 3.

**§12.** Пусть нить в форме прямого круглого цилиндра с радиусом сечения $R$ лежит своей осью вдоль оси $OX$, направленной вертикально вниз (фиг.3).

Пусть рядом бесконечно близких сечений нить, предполагаемая однородной, разбита на ряд бесконечно тонких слоев…,A,B,C,… Пусть слой B, ограниченный от соседних слоев A и C сечениями с абсциссами $x-\frac{\partial x}{2}, x+\frac{\partial x}{2}$, в момент времени $t$ оказывается повернутым вокруг OX на угол $\varphi$, отсчитываемый от некоторого выбранного для всех слоев начального положения. Угол $\varphi$ есть функция абсциссы $x$ и рассматриваемого момента времени $t$.

Выделим из слоя B элемент объема, ограниченный двумя цилиндрическими поверхностями с общей осью OX и радиусами сечений $r, r+dr$, горизонтальными плоскостями $x-\frac{\partial x}{2}, x+\frac{\partial x}{2}$ и двумя объема соприкасается по горизонталь ным плоскостям пересекающимися по OX плоскостями с азимутами $\psi, \psi+d\psi$. Этот элемент $x-\frac{\partial x}{2}, x+\frac{\partial x}{2}$ с аналогичными элементами из слоев A и C. Вследствие неодинаковости мгновенных линейных скоростей у этих элементов, на их поверхности раздела должны возникать из-за вязкости вещества нити силы внутреннего трения, которые направлены в сторону уменьшения относительных скоростей. По закону Ньютона абсолютная величина силы трения пропорциональна площади $rd\psi dr$.

Следовательно, сила действия элемента из слоя A на элемент из слоя B, соседней с первым, есть

$$-\eta r dr d\psi r \left(\frac{\partial \dot{\varphi}}{\partial x}\right)_{x-\frac{dx}{2}}.$$

Аналогично, действие элемента из слоя C на рассматриваемый элемент из B равно

$$+\eta r dr d\psi r \left(\frac{\partial \dot{\varphi}}{\partial x}\right)_{x+\frac{dx}{2}}.$$

Их совместное действие оказывается равным

$$\eta r^2 drd\psi \frac{\partial^2 \varphi}{\partial x^2} dx,$$

и его момент относительно $OX$ есть

$$\eta r^3 drd\psi \frac{\partial^2 \varphi}{\partial x^2} dx.$$

Поэтому полный момент сил внутреннего трения, действующих на весь слой В, равен

$$\eta \frac{\partial^2 \dot{\varphi}}{\partial x^2} dx \int_0^R r^3 dr \int_0^{2\pi} d\psi = +\frac{\pi R^4}{2} dx \eta \frac{\partial^2 \dot{\varphi}}{\partial x^2}. \qquad (66)$$

Далее, момент силы инерции для слоя В, очевидно, есть

$$-\frac{MR^2}{2}\ddot{\varphi} = -\frac{\pi R^4}{2} dx \rho \ddot{\varphi}; \qquad (67)$$

Здесь $M$ обозначает массу слоя В. Наконец, рассмотрим, как действуют на слой В силы упругости. Фиг.4 представляет коаксильную цилиндрическую часть слоя В той же высоты $dx$, но меньшего радиуса $r$.Благодаря кручению в верхнем сечении $x - \frac{dx}{2}$ этой части по закону Гука на элемент части, имеющий вид призмы с основанием $rdrd\psi$ и с высотой $dx$, действует сила упругости

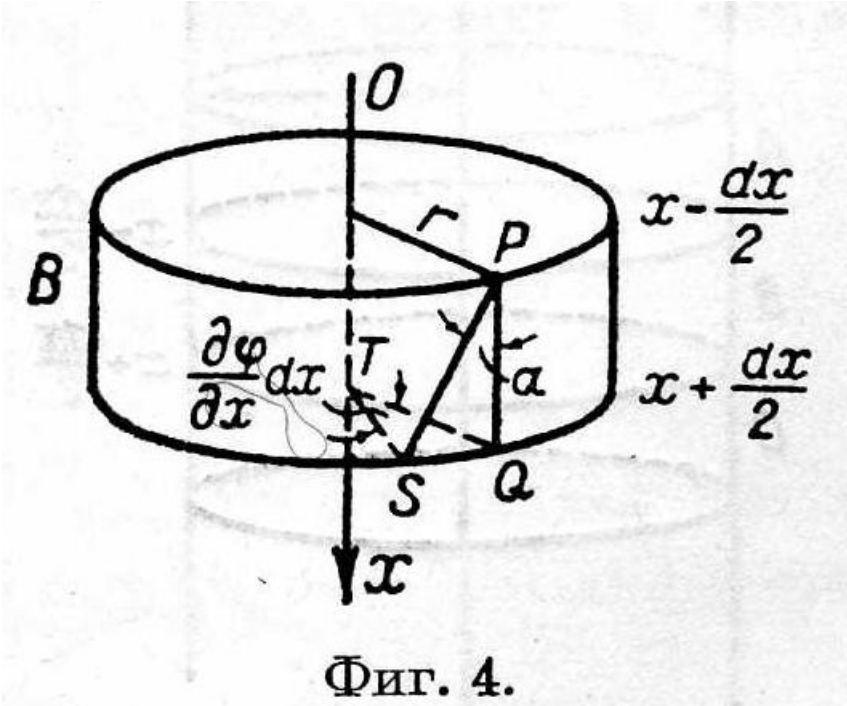

Фиг. 4.

$$\mu r dr d\psi tg\alpha = \mu r^2 drd\psi \left(\frac{\partial \varphi}{\partial x}\right)_{x-\frac{dx}{2}},$$

где $\alpha$ – угол наклона ребра призмы и $\mu$ – модуль сдвига.

Здесь мы воспользовались очевидным соотношением:

$$tg\alpha = \frac{SQ}{PQ} = \frac{r \cdot \angle QTS}{dx} = \frac{r}{dx}\left(\frac{\partial \varphi}{\partial x}\right)_{x-\frac{dx}{2}} \cdot dx.$$

Аналогично, в нижнем сечении рассматриваемой призмы действует в ту же сторону сила упругости

$$\mu r dr d\psi tg\alpha' = \mu r^2 drd\psi \left(\frac{\partial \varphi}{\partial x}\right)_{x+\frac{dx}{2}}.$$

Следовательно, призма оказывается подверженной деформирующему усилию, равному разности приведенных выше сил:

$$\mu r^2 drd\psi \left[\left(\frac{\partial \varphi}{\partial x}\right)_{x+\frac{dx}{2}} - \left(\frac{\partial \varphi}{\partial x}\right)_{x-\frac{dx}{2}}\right] = \mu r^2 drd\psi \frac{\partial^2 \varphi}{\partial x^2} dx.$$

Момент этого усилия относительно $OX$ есть

$$\mu r^3 dr d\psi \frac{\partial^2\varphi}{\partial x^2} dx;$$

полный момент упругих сил, действующих на весь слой В, ок5азывается равным

$$\mu\int_0^R r^3 dr \int_0^{2\pi} d\psi \frac{\partial^2\varphi}{\partial x^2} dx = +\frac{\pi R^4}{2} dx \mu \frac{\partial^2\varphi}{\partial x^2} \quad . \tag{68}$$

По принципу Даламбера сумма моментов (66),(67) и (68) должна равняться нулю, откуда получаем дифференциальное уравнение крутильных колебаний упруго-вязкой нити:

$$\eta \frac{\partial^3\varphi}{\partial x^2 \partial t} + \mu \frac{\partial^2\varphi}{\partial x^2} - \rho \frac{\partial^2\varphi}{\partial t^2} = 0. \tag{F}$$

Оно имеет тот же вид, что и уравнение (А) §1 для поперечных колебаний, если в последнем положить $T = 0, g = 0$. И физический смысл его тот же самый, только роль смещения $u$ здесь играет угол закручивания $\varphi$.

**§ 13.** Уравнение (F) должно быть сопровождено граничными условиями, определяемыми из постановки задачи. Здесь рассматривается случай, когда верхний конец нити закреплен, а нижний несет некоторый –дополнительный момент инерции, оставаясь свободным. Нить предполагается расположенной по-прежнему вертикально и движение ее происходит в вакууме, так что внешнее трение отсутствует. На верхнем конце при всех $t \geq 0$ должно быть

$$\varphi|_{x=0} = 0. \tag{69}$$

Посмотрим, какие условия должны удовлетворяться на нижнем, свободном, конце с абсциссой $x = l.$ Пусть дополнительный момент инерции обусловлен диском, прикрепленным к концу $l$ нити своим центром. Пусть $R_0$ – радиус диска, $d$ – его толщина, $\rho_0$ – плотность его вещества.

Тогда мгновенный момент сил инерции диска есть

$$-\frac{\pi d R_0^4 \rho_0}{2}\left(\frac{\partial^2\varphi}{\partial t^2}\right)_{x=l} . \tag{70}$$

Этот момент должен уравновешиваться с моментами сил упругости и вязкости нити в месте прикрепления диска.

Рассмотрим сначала действие упругих сил. В элементе сечения нити, который имеет площадь $r dr d\psi$, из-за действия диска развивается сила упругости, равная

$$\mu_0 r^2 dr d\psi \left(\frac{\partial \varphi}{\partial x}\right)_{x=l},$$

как это мы видели в § 12. Следовательно, наоборот, этот элемент конечного сечения нити действует на диск с силой

$$-\mu_0 r^2 dr d\psi \left(\frac{\partial \varphi}{\partial x}\right)_{x=l}.$$

Ее момент относительно $OX$ равен

$$-\mu_0 r^3 dr d\psi \left(\frac{\partial \varphi}{\partial x}\right)_{x=l},$$

и, следовательно, полный момент упругой реакции конца нити на диск есть:

$$-\frac{\mu_0 \pi R^4}{2}\left(\frac{\partial \varphi}{\partial x}\right)_{x=l},$$

где $R$, как и раньше, радиус сечения нити; $\mu_0$ есть модуль сдвига для конца нити.

Во-вторых, благодаря вязкости, тот же элемент $rdrd\psi$ нижнего сечения нити действует на диск по закону внутреннего трения с силой

$$-\eta_0 r^2 dr d\psi \left(\frac{\partial \varphi}{\partial x}\right)_{x=l},$$

где $\eta_0$, как и $\mu_0$, относится к веществу, прикрепляющему диск к нити. Момент этой силы относительно $OX$ есть

$$-\eta_0 r^3 dr d\psi \left(\frac{\partial \dot{\varphi}}{\partial x}\right)_{x=l},$$

и полный момент сил вязкости, действующих на диск, нужно считать равным

$$-\frac{\eta_0 \pi R^4}{2}\left(\frac{\partial \dot{\varphi}}{\partial x}\right)_{x=l}. \qquad (72)$$

Из (70) - (72) по принципу Даламбера получаем условие на свободном конце:

$$\rho_0 d\left(\frac{\partial^2 \varphi}{\partial t^2}\right)_{x=l} + \mu_0 \delta^4 \left(\frac{\partial \varphi}{\partial x}\right)_{x=l} + \eta_0 \delta^4 \left(\frac{\partial^2 \varphi}{\partial x \partial t}\right)_{x=l} = 0, \qquad (73)$$

где $\delta$ есть отношение радиуса нити к радиусу диска.

Наконец, должно быть дано начальное распределение углов закручивания $\varphi$ и угловых скоростей $\dot{\varphi}$:

$$\varphi|_{t=0} = \Phi(x), \quad \dot{\varphi}|_{t=0} = \varphi(x). \qquad (74)$$

**§ 14.** Уравнение

$$\eta\frac{\partial^3\varphi}{\partial x^2\partial t}+\mu\frac{\partial^2\varphi}{\partial x^2}-\rho\frac{\partial^2\varphi}{\partial t^2}=0 \qquad (F)$$

для крутильных колебаний интегрируем по способу Д. Бернулли, положив

$$\varphi = XT$$

где $X$ есть функция только $x$, а $T$ зависит только от $t$.Тогда уравнение (F) представляется в виде

$$X'':X=\rho T''':(\eta T'+\mu T),$$

и так как обе части суть функции переменных $x,t$, не зависящих друг от друга, то каждая из них есть постоянная величина. Обозначив последнюю $-\kappa$, получим два обыкновенных уравнения второго порядка:

$$X''+\kappa X=0\,. \qquad (75)$$

$$\rho T'''+\kappa\eta T'+\kappa\mu T=0, \qquad (76)$$

Причем интеграл (75), исчезающий при $x=0$, есть

$$\overline{X}=\sin\sqrt{x}x,$$

а (76) имеет общим интегралом выражение:

$$T=C'e^{q't}+C''e^{q''t},$$

где $C',C''$ – произвольные постоянные и $q',q''$ – корни уравнения

$$\rho q^2+\kappa\eta q+\kappa\mu=0\,.$$

Вследствие этого мы получили частное решение уравнения(F) в виде:

$$\varphi=\left(C'y^{q't}+C''y^{q''t}\right)\sin\sqrt{x}x; \qquad (77)$$

при $x=0$ это решение исчезает вместе со своей производной по $t$ для всех $t\ge 0$.

Параметр $\kappa$ должен принимать такое значение, чтобы удовлетворялось условие (73) на свободном конце нити. Для простоты допускаем, что $\mu_0=\mu$ и $\eta_0=\eta$, т.е. что способ прикрепления к диску не изменил упруго-вязких свойств в конце нити. На основании (77) условие (73) представляется в виде:

$$C'e^{q't}\left\{\rho_0 dq'^2\sin\sqrt{x}l+\sqrt{x}\delta^4\eta q'\cos\sqrt{x}l+\sqrt{x}\delta^4\mu\cos\sqrt{x}l\right\}+$$

$$+C''e^{q''t}\left\{\rho_0 dq''^2\sin\sqrt{x}l+\sqrt{x}\delta^4\eta q''\cos\sqrt{x}l+\sqrt{x}\delta^4\mu\cos\sqrt{x}l\right\}=0\,.$$

Чтобы это было справедливо при всех $t$, нужно чтобы $q',q''$ удовлетворяли уравнению:

$$(\rho_0 d \sin\sqrt{x}l)q^2 + (\sqrt{x}\delta^4\eta\cos\sqrt{x}l)q + (\sqrt{x}\delta^4\mu\cos\sqrt{x}l) = 0.$$

С другой стороны, $q'$ и $q''$ по их определению являются корнями для

$$\rho q^2 + \kappa\eta q + \kappa\mu = 0.$$

Следовательно, коэффициенты обоих последних уравнений должны быть пропорциональны, а это приводит к условию:

$$\sqrt{\kappa} = \frac{\delta^4\rho}{\rho_0 d} ctg\sqrt{x}l. \qquad (78)$$

Итак, для того, чтобы удовлетворялось условие (73) на свободном конце нити, $\kappa$ должны принимать лишь те значения $\sqrt{\kappa_1}, \sqrt{\kappa_2}, ..., \sqrt{\kappa_m}, ...$ которые служат корнями уравнения (78). Эти значения являются собственными значениями (Eigenwerte) рассматриваемой задачи. Не нарушая общности рассуждений, можно считать их положительными и занумерованными в порядке их возрастания(легко убедиться, что они различны, что их бесконечно много и что разность $\sqrt{\kappa_m} - \sqrt{\kappa_{m-1}}$ стремится к некоторому пределу, когда $m \to \infty$).

Вычисляя для каждого $\kappa_m, m = 1,2,...$соответствующие значения $q', q''$, видим, что

$$\varphi_m = \left(C'_m e^{q'_m t} + C''_m e^{q''_m t}\right)\sin\sqrt{\kappa_m}x$$

при любых постоянных $C'_m, C''_m$ для всех $m = 1,2,...$, являются решениями уравнения (F), удовлетворяющими условиям(69) и (73) на концах нити. Вследствие линейности уравнений (F), и (73), сумма

$$\varphi(x,t) = \sum_{m=1}^{\infty}\left(C'_m e^{q'_m t} + C''_m e^{q''_m t}\right)\sin\sqrt{\kappa_m}x \qquad (79)$$

при любых $C'_m, C''_m$ также является решением (F) при соблюдении обоих условий (69) и (73).

**§15.** Выберем теперь $C'_m, C''_m, m = 1,2,...$ так, чтобы соблюдались их начальные условия (74). Пусть данные функции $\Phi(x)$ и $\varphi(x)$ непрерывны с их первыми производными для всех $x$ из интервала $0 \le x \le l$ и пусть их вторые производные для тех же $x$ удовлетворяют условиям Дирихле.

Положим

$$\Phi(x) = \sum_{m=1}^{\infty} A_m \sin\sqrt{\kappa_m}x, \quad \varphi(x) = \sum_{m=1}^{\infty} B_m \sin\sqrt{\kappa_m}x \qquad (80)$$

Дальше будет показано, что ряд (79) может быть почленно дифференцирован по $x$ и по $t$ столько раз, сколько требуется для составления уравнения (F). Дифференцируя (79) по $t$, находим:

$$\frac{\partial\varphi}{\partial t} = \sum_{m=1}^{\infty}\left(q'_m C'_m e^{q'_m t} + q''_m C''_m e^{q''_m t}\right)\sin\sqrt{\kappa_m}\,x\,. \tag{81}$$

Подстановка $t = 0$ дает:

$$\varphi\big|_{t=0} = \sum_{m=1}^{\infty}\left(C'_m + C''_m\right)\sin\sqrt{\kappa_m}\,x\,,$$

$$\frac{\partial\varphi}{\partial t}\Big|_{t=0} = \sum_{m=1}^{\infty}\left(q'_m C'_m + q''_m C''_m\right)\sin\sqrt{\kappa_m}\,x\,.$$

Следовательно, чтобы соблюдались начальные условия (74), достаточно выбрать $C'_m$ и $C''_m$ согласно уравнениям:

$$C'_m + C''_m = A_m,\quad q'_m C'_m + q''_m C''_m = B_m,\qquad m = 1,2,\ldots \tag{82}$$

Остается лишь выяснить вопрос о сходимости рядов. В предыдущем параграфе было указано, что последовательность собственных значений

$$\sqrt{\kappa_1}, \sqrt{\kappa_2}, \ldots, \sqrt{\kappa_m}, \ldots$$

при возрастании $m$ стремится стать некоторой арифметической прогрессией. Это значит, что существует настолько большое целое $M > 0$, что ряды (80), в которых суммирование мы производим теперь от $m = M$ до $m = \infty$, состоят из членов, каждый из которых как угодно мало отличается от соответственных членов тригонометрических рядов. Так как сходимость ряда зависит именно от его остатка, то в данном случае могут быть применены основные теоремы рядов Фурье. Но тогда по данному поводу нам придется повторить те рассуждения, которые были приведены по аналогичному поводу в § 5. Следовательно, высказанное утверждение о равномерной сходимости всех рассматриваемых в этой задаче рядов можно считать доказанным.

Уравнение (79) дает решение для (F), которое при $t = 0$, при $x = 0$ и при $x = l$ определяется вместе со своей производной условиями (74),(69) и (73).

Последнее, будучи дифференциальным уравнением второго порядка, при помощи (74) выделяет из класса всевозможных интегралов (F) те, которые, удовлетворяя (73), должны одновременно удовлетворять условию (74) при

$$x = l:\qquad \varphi\Big|_{\substack{t=0\\x=l}} = \Phi(l),\quad \dot\varphi\Big|_{\substack{t=0\\x=l}} = \varphi(l)\,.$$

Кроме того, из (91) видно, что $\dot\varphi\big|_{t=\infty} = 0$.

На основании теорем § 3 можно утверждать, что решение (79) единственное.

**§ 16.** Подойдем к изучаемому вопросу с иной точки зрения. В § 10 уравнение вынужденных колебаний упруго-вязкой нити написано в следующей форме:

$$\eta \frac{\partial^3 u}{\partial x^2 \partial t} + M \frac{\partial^2 u}{\partial x^2} - \rho \frac{\partial^2 u}{\partial t^2} - H \frac{\partial u}{\partial t} = -f(x,t), \qquad \text{(E)}$$

где $f(x,t)$ - фактор вынуждающей силы.

Предварительно заметим следующее свойство функции

$$y(x) = f(x,0) - \rho \ddot{u}(x,0), \qquad (83)$$

представляющей величину, пропорциональную потерянной силе у частицы нити с абсциссой $x$ в начальный момент времени $t = 0$.Здесь, как всегда, точки над буквой обозначают порядок производной по времени. Пусть $x_0$ есть корень $y(x)$. Если положить в (E) $H = 0$,ограничиваясь для простоты случаем колебаний в вакууме, то при $x = x_0$ при помощи (83) из (E) найдем:

$$\left[ \eta \frac{\partial^3 u}{\partial x^2 \partial t} + M \frac{\partial^2 u}{\partial x^2} \right]_{\substack{x=x_0 \\ t=0}} = 0,$$

или, так как $\left.\frac{\partial u}{\partial t}\right|_{\substack{t=0 \\ x=x_0}} = \varphi(x_0)$, $\left.u\right|_{\substack{t=0 \\ x=x_0}} = \Phi(x_0)$,

то

$$\left[ \frac{\partial^2}{\partial x^2} \{\eta\varphi(x) + M\Phi(x)\} \right]_{x=x_0} = 0. \qquad (84)$$

Это значит, что те абсциссы $x$, для которых функция (83) обращается в нуль, являются точками перегиба кривой

$$Y = \eta\varphi(x) + M\Phi(x). \qquad (85)$$

Вместе с тем соотношение

$$\left[ \eta \frac{\partial^3 u}{\partial x^2 \partial t} + M \frac{\partial^2 u}{\partial x^2} \right]_{\substack{x=x_0 \\ t=0}} = 0$$

позволяет видеть, что в точках нити с такими абсциссами силы упругости уравновешиваются силами внутреннего трения в момент времени $t = 0$. В частности, если $\Phi(x_0) = 0$ , то при $t = 0$ для точки $x = x_0$ сила упругости отсутствует. Тогда должна отсутствовать т сила вязкости. Следовательно, соседние с $x_0$ слои в момент $t = 0$ тогда должны иметь скорости, равные между собой относительно скорости в $x_0$ по величине и противоположно направленные. Кроме того, из-за отсутствия силы упругости, по закону Гука, должно быть: $u(x_0,0) = 0$.

Существенно отметить, что кривая (85)произвольна в той же мере, что и функции $\Phi(x)$ и $\varphi(x)$.

В точке, подобной рассмотренной $x_0$, в начале движения силы упругости и вязкости взаимно уравновешиваются, и движение в этот момент таково, как если бы частица $x_0$,не будучи связана внутренними силами нити, была подвержена действию лишь вынуждающей силы $f(x_0,0)$. Начальные условия движения этой частицы соблюдались бы при любых предположениях относительно сил упругости и вязкости в этот момент и требовалось бы только наличие равновесия между ними. Такие точки $x_0$ являются в некотором смысле особыми, мы будем исключать их из рассмотрения при изучении явления в общем. Дальше, в § 18, мы увидим математический смысл этого исключения, заметив сейчас, что благодаря произволу в выборе функции $Y(x)$, можно считать, что выбранная кривая (85) не имеет точек перегиба, что и будет показано в § 19.

**§ 17.** Если бы на частицы нити, кроме внешних сил, действовали только силы упругости, то закон малых(поперечных или крутильных) колебаний нити выражался бы уравнением

$$f(x,t)-\rho\ddot{u}(x,t)=Mu(x,t),$$

где $u(x,t)$, смотря по обстоятельствам, есть или линейное перемещение или угол закручивания. Следовательно, в этом случае было бы:

$$u(x,t)=\frac{1}{M}\left[f(x,t)-\rho\ddot{u}(x,t)\right].$$

Это уравнение чисто упругих перемещений является выражением закона Гука на основании принципа Даламбера.

Если же в рассмотрение вводятся такие силы, как, например, внутреннее трение, то уравнение (86),благодаря аддитивным свойствам малых деформаций, должно быть дополнено. Дополнительный к (86) член можно было бы составить следующим образом. Пусть частица нити с абсциссой $x$ в момент времени $\theta$ получает мгновенный импульс, равный единице. Пусть вследствие действия этого импульса та же частица получает в более поздний момент $t\geq 0$ дополнительное к (86) смещение, которое мы обозначим $K(x;t,\theta)$. Благодаря аддидивности перемещений можем тогда утверждать, что если на частицу $x$ в момент $\theta$ подействовал импульс $F(x,\theta)d\theta$, то в момент $t>0$ она получит дополнительное к (86) смещение, равное $K(x;t,\theta)F(x,\theta)d\theta$.

Вследствие той же аддитивности дополнительное к (86) смещение нужно считать равным:

$$\int_0^t K(x;t,\theta)F(x,\theta)d\theta,$$

если на частицу нити $x$ действовала в течение всего предшествующего времени переменная сила $F(x,\theta)$, $0 \le \theta \le t$.

Так как влияние дополнительных к упругим внутренних сил(например, сил вязкости) считается учтенным посредством функции $K(x;t,\theta)$, то по принципу Даламбера $F(x,\theta)$ следует полагать равным

$$f(x,\theta) - \rho\ddot{u}(x,\theta).$$

Подставляя все это в (86), мы приходим к уравнению:

$$u(x,t) = \frac{1}{M}\left[f(x,t) - \rho\ddot{u}(x,t)\right] + \int_0^t K(x;t,\theta)\left[f(x,\theta) - \rho\ddot{u}(x,\theta)\right]d\theta \qquad \text{(G)},$$

данному В.Вольтерра, а также Больцманом.[1)]

Здесь ядро $K(x;t,\theta)$ предполагается обычно известным из опыта. Тогда (G), дает $u(x,t)$.

§ **18**. Поставим иную задачу. Именно, пользуясь тем, что закон поперечных или крутильных колебаний известен на основании изложенного в предыдущих параграфах, попытаемся определять вид ядра $K(x;t,\theta)$ для упруго-вязкой нити.

Умножив (G) на $M$, введем обозначения:

$$Mu(x,t) = U(x,t), \qquad f(x,t) - \rho\ddot{u}(x,t) = V(x,t).$$

Для краткости положим $U(x,t),-V(x,t) = W(x,t)$. Тогда

$$w(x,t) = \int_0^t K(x,t-\theta)V(x,\theta)d\theta \qquad (87)$$

Под интегралом правой части выполним преобразование

$$t - \theta = \vartheta, d\theta = -d\vartheta.$$

Тогда $\int_0^t K(x,t-\theta)V(x,\theta)d\theta = -\int_t^0 K(x,\vartheta)V(x,t-\vartheta)d\vartheta$.

Уравнение (87) принимает вид

$$\int_0^t V(x,t-\vartheta)K(x,\vartheta)d\vartheta = W(x,t). \qquad (88)$$

Это - уравнение Вольтерра первого рода с ядром $V(x,t-\vartheta)$ циклического типа, изучением которого занимался в числе других и Пуассон, благодаря чему его называют именем последнего. Из теории уравнений Вольтерра следует, что решение $K(x,\vartheta)$ уравнения (88), если оно существует, непременно единственно.

Пусть $x$ не есть точка перегиба кривой

$$Y = \eta\varphi(x) + M\Phi(x).$$

[1)] *Ann.d.Phys.,7, 624,1876; V.Volterra,Drei Vorles.u.s.w. Leipz.,*

Рассмотрим множество моментов времени $t = 0, \tau, 2\tau, ...,$ отделенных друг от друга достаточно малым интервалом $\tau$. Обозначим:

$$W(x, i\tau) = W_i, \quad i = 0,1,2,...$$

$$\tau V(x, (i-j)\tau) = V_{i-j}, j \le i,$$

$$K(x, j\tau) = K_j, j = 0,1,2,...$$

Тогда уравнение (88) приближенно может быть заменено следующей системой:

$$\sum_{j=0}^{i} V_{i-j} K_j = W_i \;, \; i = 0,1,2,..., \qquad (89)$$

или в развернутом виде:

$$V_0 K_0 = W_0$$

$$V_1 K_0 + V_0 K_1 = W_1$$

$$V_2 K_0 + V_1 K_1 + V_0 K_2 = W_2$$

……………………………..

Определитель этой системы есть:

$$V_0^{n+1} = \left[ f(x,0) - \rho \ddot{u}(x,0) \right]^{n+1},$$

где $n$–любое целое неотрицательное число. Так как этот определитель не равен нулю при сделанном предположении относительно точки $x$, то из (89) можно найти сколь угодно большое число $n$ неизвестных $K_0, K_1, ....,$ и притом единственным образом определяемых при данном $\tau$. Находя достаточно большое количество $K_i$ при выбранном $x$ и $\tau$, мы *по способу Прони* [1)] можем найти приближенное выражение для $K(x, t - \vartheta)$ в виде суммы показательных функций:

$$K(x, t-\theta) = \sum_{p=1}^{N} A_p(x) e^{\alpha_p(x)[t-\theta]}, \qquad (90)$$

где $A_p(x), \alpha_p(x)$ – надлежащим способом подобранные для данного $x$ числа.

**§ 19.** Пусть теперь $V(x,0) = 0$ для данного $x$, являющегося особой точкой в смысле, установленном в § 16, и пусть производная от $V(x,t)$ по $t$ при $t = 0$ не исчезает. Хотя система (89) предыдущего параграфа и не разрешима в этом случае, тем не менее и для таких значений $x$ ядро $K$ может быть найдено при помощи следующего общего приема. Дифференцируя (88) по $t$, находим:

$$K(x,t)V(x,0) + \int_0^t \frac{\partial V(x, t-\vartheta)}{\partial t} K(x,\vartheta) d\vartheta = \frac{\partial W(x,t)}{\partial t},$$

или, так как $V(x,0) = 0$:

[1] *Витекер Э., Робинсон Г., Матем. обработка результатов наблюдений, ОНТИ,стр.343,1933.*

$$\int_0^t \frac{\partial V(x,t-\vartheta)}{\partial t} K(x,\vartheta)d\vartheta = \frac{\partial W(x,t)}{\partial t} \quad , \qquad (91)$$

Т.е. опять уравнение типа (88), но с ядром $\frac{\partial V}{\partial t}$, не исчезающим при $t = 0$.Тогда остается лишь повторить все сказанное по этому поводу в § 18.

Если бы для данного $x$ производная $\frac{\partial V}{\partial t}$ обращалась в нуль при $t = 0$, но $\frac{\partial^2 V}{\partial t^2}\Big|_{t=0} \neq 0$, то, дифференцируя (91) по $t$, опять получили бы уравнение типа (88):

$$\int_0^t \frac{\partial^2 V(x,t-\vartheta)}{\partial t^2} K(x,\vartheta)d\vartheta = \frac{\partial^2 W(x,t)}{\partial t^2}$$

с ядром, не исчезающим при $t = 0$. Общность приведенных рассуждений очевидна.

Решающее слово, как всегда, принадлежит опыту. Здесь может быть одно из двух: или наблюдения над малыми колебаниями однородных нитей в самом деле приводят к тем законам движения, о которых говорилось выше, или же действительные колебания происходят не так, как это следовало бы при сделанных предположениях по теории. В первом случае было бы возможно заключить, что упругое последействие, дающее себя знать в рассматриваемом процессе движения нити, является следствием совместного влияния упругих и вязких сил. Во втором случае можно было бы по меньшей мере утверждать, что наличие сил упругости и вязкости не исчерпывает проблемы упругого последействия. Однако и во втором случае проведенное выше исследование вопроса все-таки является некоторым приближением к разгадке истинной природы упругого последействия.

*Московский текстильный институт.*

## LE PROBLÈME DE L'ACTION POSTÉRIEURE ÉLASTIQUE ET LA FRICTION INTÉRIEURE

A.N.GUÉRASSIMOV

L'ouvrage proposé est un essai théorique à appliquer les équations de la théorie de l'élasticité et de la hydrodynamique classique à la question de petits déplacements de fils homogènes, appelés ici la délastovisqueux.

Dans cet ouvrage est donnée une equation aux dérivées partielles du troisième ordre et il a été démontré que cette equation ne peut avoir qu'une seule solution qui prend avec sa dérivée première par rapport à t les valeurs données sur la borne du domaine fundamental. Cette equation peut être résolue en appliquant par exemple la méthode connue de D.Bernoulli.

Les mouvements des fils visqueux et en même temps élastiques peuvent être décrits moyennant l'équation intégrale de L.Boltzmann, la fonction influente de cette équation étant donnée sous la forme d’une combinaison linéaire de fonctions exponentielles don’t les coefficients ont été trouvés par la méthode de Prony.

Les résultats de la théorie exposée s’appliquent actuellement à la méthode acoustique pour la déterminaison pratique et assez simple du coefficient de la friction intérieure et du module de l’élasticité.

Institut textile de Moscou.

*ПММ, 1938,т.II, в. 1, с.137*

4.ПИСЬМО В РЕДАКЦИЮ (Поправка к статье)

В §3 моей статьи « Проблема упругого последействия и внутреннее трение»(сборн. «Прикл. мат. и мех.», т. I, вып. 4, 1937) намечено доказательство единственности решения данного там уравнения колебаний упругой вязкой нити. Приведенная схема доказательства вряд ли оказалась убедительной, и к тому же отсутствуют указания на способ перехода к пределу. Ниже намечен путь к ее исправлению.

Речь идет об интеграле $u(x,t)$ уравнения

$$\eta\frac{\partial^3 u}{\partial x^2\partial t}+M\frac{\partial^2 u}{\partial x^2}-\rho\frac{\partial^2 u}{\partial t^2}-H\frac{\partial u}{\partial t}=0\ , \qquad (1)$$

удовлетворяющего при $T>0$ условиям:

$$u\big|_{x=0}=0,\ \ u\big|_{x=l}=0,\ \ u\big|_{t=0}=0,\ \ \frac{\partial u}{\partial t}\Big|_{t=0}=0\ ,\ u\big|_{t=T}=0,\ \ \frac{\partial u}{\partial t}\Big|_{t=T}=0\ . \qquad (2)$$

Расширяем область S $(0\le x\le l, 0\le t\le T)$, восстанавливая перпендикуляры равной длины $\varepsilon$ в плоскости $XT$ к границе $C$ этой области (направленные наружу). Концы этих перпендикуляров определяют замкнутый контур $\Sigma$, содержащий внутри себя $C$ и нигде не касающийся последнего. Можно потребовать, чтобы $u$, которое по условию исчезает вместе с $\frac{\partial u}{\partial t}$ на $C$, исчезало бы тождественно внутри области $S_+$ между $C$ и $\Sigma$ и на $\Sigma$.

Применение формулы Гаусса к $S$ и $S_+$ приводит к выводу, что

$$\iint_S\left(\frac{\partial v}{\partial x}\right)^2 dxdt=0 \ \text{(где } v=\eta\frac{\partial u}{\partial t}+Mu\text{ ).}$$

Отсюда следует, что $u\equiv 0$ в $S$, как это и было показано. Результат остается в силе, каково бы ни было $\varepsilon>0$.

Заставляя $\varepsilon$ стремиться к нулю, получим, что решение уравнения (1) при условиях (2) равно тождественно нулю в $S$.

Так как это верно при всяком $T>0$, то при $T\to+\infty$ получим доказываемую теорему. При этом нет никакой надобности рассматривать знак $\eta H-\sigma M$. Но зато существенным является требование, чтобы не только $\frac{\partial u}{\partial t}$, но и $u$ обращалось в 0 при $t=+\infty$. Это значит, что потенциальная и кинетическая энергия колебаний с течением времени переходит в теплоту.

Далее в цитированной работе (§ 11) идет речь о сходимости ряда

$$f(x,t)=\sum_m w_m(t)\sin\frac{m\pi x}{l}$$

для $u(x,t)$ уравнения

$$\eta\frac{\partial^3 u}{\partial x^2\partial t}+M\frac{\partial^2 u}{\partial x^2}-\rho\frac{\partial^2 u}{\partial t^2}-H\frac{\partial u}{\partial t}=-f(x,t).$$

Чтобы ряды, составленные формально для

$$u,\frac{\partial u}{\partial t},\frac{\partial^2 u}{\partial t^2},\frac{\partial^2 u}{\partial x^2},\frac{\partial^3 u}{\partial x^2\partial t},$$

равномерно и абсолютно сходились, достаточно предположить, что порядок малости $w_m(t)$ равномерно для всех $t\ge 0$ с возрастанием $m$ был не ниже, чем $m^{-2}$, т.е. что производная первого порядка по $x$ от $f$ удовлетворяет в $[0,l]$ условиям Дирихле. Значит, нет надобности считать $f'_x$ непрерывной, как это указывается в тексте.

## 5. ОСНОВАНИЯ ТЕОРИИ ДЕФОРМАЦИЙ УПРУГО-ВЯЗКИХ ТЕЛ

### § 1. Элементарные деформации для однородных изотропных упругих тел

Пусть упругий однородный и изотропный стержень длиной $dx$ под действием растягивающего усилия приобретает удлинение $\partial u_x$.

По закону Гука, верному в пределах упругости, при такой деформации стержня в последнем возникает напряжение $v'_x$, равное:

$$v'_x = E\frac{\partial^2 u_x}{\partial x}, \qquad (1)$$

где $E$ есть модуль Юнга. Если обозначить абсолютные величины относительных сокращений линейных размеров в направлениях осей координат $OY$ и $OZ$ соответственно через $\frac{\partial u_y}{\partial y}$ и $\frac{\partial u_z}{\partial z}$, то

$$v'_x = -m'E\frac{\partial u_y}{\partial y} = -m'E\frac{\partial u_z}{\partial z}. \qquad (1')$$

Минус в правой части указывает на сокращение поперечных размеров, а $m'$ есть положительная постоянная, обратная коэффициенту Пуассона.[1)]

Пусть ABCD (фиг.1) есть элементарный параллелепипед, выделенный внутри однородного и изотропного упругого тела; ребра параллелепипеда, параллельные осям прямоугольной системы координат, имеют соответственно длины $\partial x, \partial y, \partial z$. Считая грань $AD$ неподвижной, представим себе, что

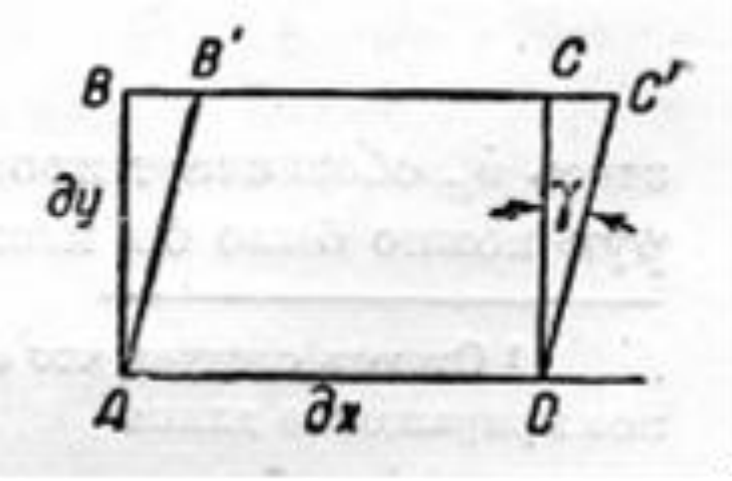

Фиг.1

противоположная грань $BC$ вынуждена сдвинуться вдоль $OX$ до положения $B'C'$. Обозначим $BB'$ через $\partial u_x$. По закону Гука развивающееся при такой деформации напряжение $\tau'_{xy}$ должно быть пропорционально углу сдвига $\gamma$, так что

$$\tau'_{xy} = H\gamma = H\frac{\partial u_x}{\partial x}; \qquad (2)$$

здесь $H$ есть модуль сдвига.

Формулы, аналогичные (1) и (2), могут быть записаны и для двух других направлений $OY$ и $OZ$.

### §2. Деформации скоростей в упруго-вязком теле

Теперь рассмотрим это явление с гидромеханической точки зрения.

Когда стержень растягивается, то частицы вещества перемещаются вдоль

[1)]*А. и Л. Феппль. Сила и деформация, I, ГТТИ, стр. 33, 1933.*

оси стержня. Вещество как бы течет по стержню и последний представляет нечто подобное тому, что в гидродинамике называют трубкой тока. Однако площадь поперечного сечения стержня при растяжении сокращается, следовательно, для соблюдения условия неразрывности скорость движения частиц становится больше той, которая была бы при растяжении с постоянным сечением.[1)] Появление ускорений, связанных с относительным сокращением поперечных размеров, может быть истолковано как результат действия сил; роль последних играют силы внутреннего трения (или вязкости) и можно было бы говорить о «напряжении» $v''_x$ сил внутреннего трения, перпендикулярном к поперечному сечению стержня. Естественно предположить,[1)] что это напряжение пропорционально изменению скорости при переходе от одного сечения к другому, бесконечно близкому, так что

$$v''_x = \varepsilon \frac{\partial \dot{u}_x}{\partial x}; \qquad (3)$$

здесь $\dot{u}_x$ обозначает проекцию скорости на ось $OX$, а $\varepsilon$ есть постоянная, которую можно было бы назвать продольным модулем вязкости.[2)]

Формула (3) является аналитическим выражением закона Ньютона, по которому сила внутреннего трения пропорциональна градиенту скорости. Другой случай применения этого закона можно проследить на деформации сдвига.

В самом деле, возвращаясь к фиг.1, видим, что слои тела, параллельные оси $OX$, смещаясь, вынуждены скользить друг по другу. При этом они обмениваются количествами движения, т.е. приобретают ускорения единственно вследствие взаимодействий. Идя дальше, можно ввести в рассмотрение «тангенциальные напряжения $\tau''_{yx}$ сил внутреннего трения», действующие на единичную площадку, перпендикулярную в $OY$, и направленные параллельно $OX$. По закону Ньютона напряжение $\tau''_{yx}$ сил внутреннего трения пропорционально градиенту скорости $\dot{u}_x$ при переходе от одного слоя к другому, с ним соприкасающемуся, т.е.

[1)] *Отсюда следует, что если при чисто упругой деформации мы имели бы относительное приращение длины* $\frac{\partial u_x}{\partial x} = \frac{v_x}{E}$, *то благодаря текучести вещества действительное приращение длины будет больше, а именно:* $\frac{\partial u_x}{\partial x} = \frac{v_x}{E - \varepsilon^2}$, *отсюда* $v_x = E\frac{\partial u_x}{\partial x} - \varepsilon^2 \frac{\partial u_x}{\partial x} = v'_x + v''_x$, *где* $v''_x = -\varepsilon^2 \frac{\partial u_x}{\partial x}$. *Все дальнейшее рассуждение основывается на соотношении* $-\varepsilon \frac{\partial u_x}{\partial x} = \frac{\partial}{\partial t}\left(\frac{\partial u_x}{\partial x}\right)$, *лежащем в основе релаксационной теории Максвелла, откуда получаем* $v''_x = +\varepsilon \frac{\partial u_x}{\partial x}$.

*Таким образом, полное напряжение* $v_x$ *складывается из чисто упругого* $v'_x$ *и вязкого* $v''_x$.

[2)] *Так как* $\frac{\partial \dot{u}_x}{\partial x} = \frac{\partial}{\partial t}\left(\frac{\partial u_x}{\partial x}\right)$, *то (3) обозначает, что* $v''_x$ *пропорционально скорости деформации.*

*См. статью «Elasticity», W.Thomson в IX издании «Британской энциклопедии».*

$$\tau''_{yx} = \eta \frac{\partial \dot{u}_x}{\partial y}, \qquad (4)$$

где постоянная $\eta$ может быть названа поперечным модулем вязкости.[1)]

Аналогично (3) и (4) можно было бы написать подобные формулы и для других осей координат.

Сопоставляя (1) с (3) и (2) с (4), мы видим формальное сходство элементарных законов Ньютона и Гука. Это сходство может быть сформулировано таким образом: скорость деформации по отношению к закону Ньютона играет такую же роль, как сама деформация по отношению к закону Гука. Однако коэффициенты пропорциональности в обоих законах существенно различны.

§ 3.Общая деформация упруго-вязких аморфных тел

Обычно считают, что до предела упругости тело ведет себя строго по закону Гука и разумеют под пределом упругости верхнюю границу точности этого закона. В этих пределах тело должно считаться вполне упругим. С достижением предела упругости закон Гука перестает быть верным и отступление от него выражается у аморфных тел прежде всего в появлении признаков текучести. При достаточно больших нагрузках текучесть вещества настолько дает себя чувствовать, что аморфное тело ведет себя как вязкая жидкость.

Мы полагаем, что для некоторых веществ такое разделение свойств вряд ли может быть оправдано.

Во-первых, едва ли можно предположить возможность перерождения упругих сил в силы внутреннего трения. В самом деле, природа упругих сил Существенно отличается от природы сил вязкости. Упругие силы зависят от смещений, т.е. от взаимного расположения частиц тела. И потому консервативны, в то время как силы внутреннего трения, будучи обусловлены различием в скоростях отдельных частиц, не обладают силовой функцией. Наличию упругих сил соответствует наличие некоторых обратимых процессов, сопровождающихся переходом потенциальной энергии в кинетическую и обратно. Наоборот, дифференциальные уравнения движения вязких тел, содержа члены с первой степенью дифференциала времени, описывают принципиально необратимые процессы.

[1)] *Легко видеть, что* $\frac{\partial \dot{u}_x}{\partial y} = \frac{\partial \gamma}{\partial t}$, *где* $\gamma$ *– угол сдвига (фиг.1). Следовательно, и в этом случае напряжение* $\tau''_{yx}$ *пропорционально скорости деформации.*

Во-вторых, можно назвать вещество, которое, несомненно, одновременно обладает как упругими, так и вязкими свойствами. Достаточно назвать, например, хлебное тесто. При разрыве оно дает края, подобные тем, которые получаются при разрыве куска металла. Но с течением времени резкость краев сглаживается и кусок теста постепенно приобретает округлость формы, свойственную жидким телам.

Все это заставляет нас полагать, что при изучении деформации аморфных тел нужно принимать во внимание, вообще говоря, как наличие сил упругости, так и сил внутреннего трения.

На основании результатов §§ 1 и 2 и сказанного здесь заключаем, что при продольном растяжении упруго-вязкого тела возникает напряжение, которое благодаря независимости и аддитивности величин (1) и (2) нужно считать равным:

$$v_x = v'_x + v''_x = \frac{\partial}{\partial x}\left(Eu_x + \varepsilon\dot{u}_x\right). \qquad (5)$$

Аналогично при чистом сдвиге в упруго-вязком теле появляется тангенциальное напряжение

$$\tau_{yx} = \tau'_{yx} + \tau''_{yx} = \frac{\partial}{\partial y}\left(Hu_x + \eta\dot{u}_x\right). \qquad (6)$$

Подобно (5) и (6) можно написать формулы для нормальных и тангенциальных напряжений по прочим направлениям.

## § 4. Тензор упругой деформации

Рассмотрим[1)] две бесконечно близкие частицы однородного и изотропного упругого тела, положения которых определяются векторами $r$ и $r+dr$. После деформации положения этих частиц будут определены векторами

$$r + u(r), r + dr + u(r) + du(r),$$

где $u(r)$ есть вектор смещения частицы с начальным положением $r$. Поэтому смещение частицы с начальным положением $r+dr$ есть

$$u(r) + du(r) = u + \left(\frac{du}{dr}, dr\right).$$

Здесь $\frac{du}{dr}$ есть тензор, представляющий собой производную от вектора $u$ и по вектору $r$. Обычным приемом его можно разложить на симметричную часть $Defu$, равную:

[1)] *Н.Е.Кочин. Векторное исчисление и начала векторного исчисления, ГТТИ, стр.364, 1934.*

$$Defu = \begin{Bmatrix} \frac{\partial u_x}{\partial x} & \frac{1}{2}\left(\frac{\partial u_x}{\partial y} + \frac{\partial u_y}{\partial x}\right) & \frac{1}{2}\left(\frac{\partial u_x}{\partial z} + \frac{\partial u_z}{\partial x}\right) \\ \frac{1}{2}\left(\frac{\partial u_y}{\partial x} + \frac{\partial u_x}{\partial y}\right) & \frac{\partial u_y}{\partial y} & \frac{1}{2}\left(\frac{\partial u_y}{\partial z} + \frac{\partial u_z}{\partial y}\right) \\ \frac{1}{2}\left(\frac{\partial u_z}{\partial x} + \frac{\partial u_x}{\partial z}\right) & \frac{1}{2}\left(\frac{\partial u_z}{\partial y} + \frac{\partial u_y}{\partial z}\right) & \frac{\partial u_z}{\partial z} \end{Bmatrix} \qquad (7)$$

и антисимметричную часть

$$A' = \begin{Bmatrix} 0 & -\omega_z & +\omega_y \\ +\omega_z & 0 & -\omega_x \\ -\omega_y & +\omega_x & 0 \end{Bmatrix}, \qquad (8)$$

где вектор $\omega$ определяется из уравнения

$$\omega = \frac{1}{2} rotu; \qquad (9)$$

поэтому

$$du(r) = \left(\frac{du}{dr}, dr\right) = (Defu, dr) + (A', dr)$$

и, следовательно, смещение частицы $r + dr$ после деформации оказывается равным

$$u + (Defu, dr) + \frac{1}{2}[rotu, dr], \qquad (10)$$

так как легко проверить, что $(A', dr) = [\omega, dr]$

Из (10) видно, что перемещение частицы состоит из поступательного перемещения, из деформации этой частицы состоит из поступательного перемещения, из деформации этой частицы и из вращения ее как целого вокруг оси с направлением вектора $rotu$ на угол $\frac{1}{2}|rotu|$.

Симметричный тензор $Defu$, определяемый уравнением (7), носит название тензора деформаций смещений. Сумма его компонентов, стоящих на главной диагонали (первый инвариант), равна объемному расширению:

$$\frac{\partial u_x}{\partial x} + \frac{\partial u_y}{\partial y} + \frac{\partial u_z}{\partial z} = divu = \theta'.$$

## § 5. Обобщенный закон Гука

На основании формул (1) и (1′) можем утверждать, что напряжению $v_x$ соответствует относительное удлинение вдоль $OX$, равное $\frac{v'_x}{E}$. С другой стороны, вследствие растягивающих напряжений вдоль $OY$ и $OZ$ получаются сокращения в направлении $OX$, соответственно равные $\frac{v'_y}{m'E}$ и $\frac{v'_z}{m'E}$. Следовательно, полное удлинение вдоль $OX$ благодаря аддитивности смещений оказывается равным:

$$\lambda_x{}' = \frac{v'_x}{E} - \frac{v'_y}{m'E} - \frac{v'_z}{m'E}. \tag{12}$$

Аналогично для двух других осей получим удлинения:

$$\lambda_y{}' = \frac{v'_y}{E} - \frac{v'_x}{m'E} - \frac{v'_z}{m'E}, \quad \lambda_z{}' = \frac{v'_z}{E} - \frac{v'_x}{m'E} - \frac{v'_y}{m'E}. \tag{12'}$$

Уравнение (12′) можно представить в виде:

$$\lambda'_\alpha = \frac{1}{m'E}\{(m'+1)v'_\alpha - \pi'\}, \quad (\alpha = x, y, z)$$

где $\pi' = v'_x + v'_y + v'_z$.

Наряду с симметричным тензором $Defu$ и упругих деформаций рассмотрим теперь тензор $\Pi'$ упругих напряжений, который, будучи отнесен к главным осям тензора $Defu$, имеет вид:

$$\Pi' = \begin{Bmatrix} v'_x & 0 & 0 \\ 0 & v'_y & 0 \\ 0 & 0 & v'_z \end{Bmatrix}. \tag{13}$$

Если внести еще отнесенный к тем же осям единичный тензор

$$I = \begin{Bmatrix} 1 & 0 & 0 \\ 0 & 1 & 0 \\ 0 & 0 & 1 \end{Bmatrix}. \tag{14}$$

то (12) может быть представлено в следующей тензорной форме, служащей выражением обобщенного закона Гука:

$$Defu = \frac{1}{m'E}\{(m'+1)\Pi' - \pi' I\}. \tag{15}$$

Решая (15) относительно $\Pi'$ и вводя упругие постоянные Ламé по формулам:

$$\frac{2\mu'}{m'E} = \frac{1}{m'+1}, \qquad \frac{2\mu}{m'-2} = \lambda, \tag{16}$$

Получим обобщенный закон Гука в виде:

$$\Pi' = 2\mu' Defu + \lambda'\theta' I, \tag{17}$$

откуда непосредственно могут быть получены компоненты тензора упругих напряжений в любой системе координат.

§ 6. Обобщение закона Ньютона

Рассуждения двух последних параграфов почти буквально применяются к рассмотрению деформации скорости и к обобщению закона Ньютона.

Если рассматриваемые частицы приобретут некоторое отрицательное приращение скорости в направлении $OX$, то вследствие увеличения поперечных размеров возникнут напряжения вязкости $v''_y$ и $v''_z$, действующие в направлении $OY$ и $OZ$.

Поэтому деформации скорости в направлениях осей координат выразятся на основании (3) при помощи уравнений, аналогичных $(12')$. Именно, градиент скорости в направлении $OX$ оказывается равным:

$$\lambda''_x = \frac{v''_x}{\varepsilon} - \frac{v''_y}{m''\varepsilon} - \frac{v''_z}{m''\varepsilon}, \tag{18}$$

$v''_y, v''_z$ – «напряжения» сил внутреннего трения, действующие в направлении осей $OX, OY, OZ; \varepsilon$ –продольный модуль вязкости и $\frac{1}{m''}$ –число, аналогичное коэффициенту Пуассона. Для двух других осей получим:

$$\lambda''_y = \frac{v''_y}{\varepsilon} - \frac{v''_x}{m''\varepsilon} - \frac{v''_z}{m''\varepsilon} \qquad \lambda''_z = \frac{v''_z}{\varepsilon} - \frac{v''_x}{m''\varepsilon} - \frac{v''_y}{m''\varepsilon} \tag{18'}$$

Эти уравнения могут быть представлены в более сокращенном виде:

$$\lambda''_\alpha = \frac{1}{m''\varepsilon}\{(m''+1)v''_\alpha - \pi''\} \qquad (\alpha = x, y, z) \tag{18}$$

если $\pi''$ определить аналогично $\pi'$:

$$\pi'' = v''_x + v''_y + v''_z .$$

Введение тензоров $\Pi''$ и $\mathrm{I}$, аналогичных (13) и (14), позволит написать обобщенный закон Ньютона в виде:

$$Def\ddot{u} = \frac{1}{m''\varepsilon}\{(m''+1)\Pi'' - \pi''\mathrm{I}\}, \tag{19}$$

где $Def\ddot{u}$ есть симметричный тензор деформации скорости, равный, подобно (7):

$$Def\dot{u} = \begin{Bmatrix} \frac{\partial \dot{u}_x}{\partial x} & \frac{1}{2}\left(\frac{\partial \dot{u}_x}{\partial y} + \frac{\partial \dot{u}_y}{\partial x}\right) & \frac{1}{2}\left(\frac{\partial \dot{u}_x}{\partial z} + \frac{\partial \dot{u}_z}{\partial x}\right) \\ \frac{1}{2}\left(\frac{\partial \dot{u}_y}{\partial x} + \frac{\partial \dot{u}_x}{\partial y}\right) & \frac{\partial \dot{u}_y}{\partial y} & \frac{1}{2}\left(\frac{\partial \dot{u}_y}{\partial z} + \frac{\partial \dot{u}_z}{\partial y}\right) \\ \frac{1}{2}\left(\frac{\partial \dot{u}_z}{\partial x} + \frac{\partial \dot{u}_x}{\partial z}\right) & \frac{1}{2}\left(\frac{\partial \dot{u}_z}{\partial y} + \frac{\partial \dot{u}_y}{\partial z}\right) & \frac{\partial \dot{u}_z}{\partial z} \end{Bmatrix}. \qquad (20)$$

Вводя вместо $\varepsilon$ и $m''$ коэффициенты вязкости $\mu''$ и $\lambda''$, аналогичные коэффициентам Ламé $\mu'$ и $\lambda'$, при помощи формул:

$$\frac{2\mu''}{m''\varepsilon} = \frac{1}{m''+1} \quad , \quad \frac{2\mu''}{m''-2} = \lambda'', \qquad (21)$$

и решая (19) относительно $\Pi''$, получим обобщенный закон Ньютона в виде:

$$\Pi'' = 2\mu'' Def\dot{u} + \lambda''\theta'' \mathrm{I}, \qquad (22)$$

где $\theta'' = div\dot{u}$ представляет собой скоростной коэффициент расширения.

Наконец, если принять во внимание наличие гидромеханического давления $p$, то благодаря аддитивности всех компонентов сил для тензора напряжений вязкости получим:[1)]

$$\Pi'' = 2\mu'' Def\dot{u} + (\lambda''\theta'' - p)\mathrm{I}$$

§ 7. Обобщенный закон Гука-Ньютона для упруго-вязких деформаций.

Принимая во внимание независимость упругих и вязких сил друг от друга и их аддитивность для малых смещений, мы можем теперь построить тензор упруго-вязких «напряжений» $\Pi$, выразив его в функции тензоров деформаций смещений и скоростей и их инвариантов $\theta'$ и $\theta''$.

Именно, мы должны иметь:

$$\Pi = \Pi' + \Pi''. \qquad (24)$$

По правилу сложения тензоров, найдем:

$$\Pi = 2\mu' Defu + \lambda'\theta'\mathrm{I} + 2\mu'' Def\dot{u} + (\lambda''\theta'' - p)\mathrm{I} =$$

*1) Страхович, Прикладная математика и механика, 11, 1934.*

$$
= \begin{Bmatrix} 2\frac{\partial T_x}{\partial x} & \frac{\partial T_x}{\partial y}+\frac{\partial T_y}{\partial x} & \frac{\partial T_x}{\partial z}+\frac{\partial T_z}{\partial x} \\ \frac{\partial T_y}{\partial x}+\frac{\partial T_x}{\partial y} & 2\frac{\partial T_y}{\partial y} & \frac{\partial T_y}{\partial z}+\frac{\partial T_z}{\partial y} \\ \frac{\partial T_z}{\partial x}+\frac{\partial T_x}{\partial z} & \frac{\partial T_z}{\partial y}+\frac{\partial T_y}{\partial z} & 2\frac{\partial T_z}{\partial z} \end{Bmatrix} + \qquad (25)
$$

$$
+ \begin{Bmatrix} \frac{\partial(\lambda' u_x + \lambda'' \dot{u}_x)}{\partial x} - p & 0 & 0 \\ 0 & \frac{\partial(\lambda' u_y + \lambda'' \dot{u}_y)}{\partial y} - p & 0 \\ 0 & 0 & \frac{\partial(\lambda' u_z + \lambda'' \dot{u}_z)}{\partial z} - p \end{Bmatrix},
$$

где $T_x = \mu' u_x + \mu'' \dot{u}_x, \quad T_y = \mu' u_y + \mu'' \dot{u}_y, \quad T_z = \mu' u_z + \mu'' \dot{u}_z.$

Отсюда видно, что $\Pi$ есть симметричный тензор.

§ 8. Вычисление $div\Pi$

Обозначим тензоры, стоящие в (25) справа, через $\Sigma$ и $T$. Тогда

$$
div\Pi = div\Sigma + divT \quad . \qquad (26)
$$

Легко убедиться, что проекция вектора $div\Sigma$ на ось $OX$ есть

$$
(div\Sigma)_x = \Delta(\mu' u_x + \mu'' \dot{u}_x) + \frac{\partial}{\partial x} div(\mu' u + \mu'' \dot{u}),
$$

где $\Delta$ – оператор Лапласа. Аналогично найдем:

$$
(div\Sigma)_y = \Delta(\mu' u_y + \mu'' \dot{u}_y) + \frac{\partial}{\partial y} div(\mu' u + \mu'' \dot{u}),
$$

$$
(div\Sigma)_z = \Delta(\mu' u_z + \mu'' \dot{u}_z) + \frac{\partial}{\partial z} div(\mu' u + \mu'' \dot{u}).
$$

Далее, для компонентов $divT$ получим:

$$
(divT)_x = \frac{\partial^2}{\partial x^2}(\mu' u_x + \mu'' \dot{u}_x) - \frac{\partial p}{\partial x},
$$

$$
(divT)_y = \frac{\partial^2}{\partial y^2}(\mu' u_y + \mu'' \dot{u}_y) - \frac{\partial p}{\partial y},
$$

$$
(divT)_z = \frac{\partial^2}{\partial z^2}(\mu' u_z + \mu'' \dot{u}_z) - \frac{\partial p}{\partial z}.
$$

Следовательно,

$$(div\Pi)_x = \Delta(\mu' u_x + \mu'' \dot{u}_x) + \frac{\partial}{\partial x} div(\mu' u + \mu'' \dot{u}) + \frac{\partial^2}{\partial x^2}(\mu' u_x + \mu'' \dot{u}_x) - \frac{\partial p}{\partial x},$$

$$(div\Pi)_y = \Delta(\mu' u_y + \mu'' \dot{u}_y) + \frac{\partial}{\partial y} div(\mu' u + \mu'' \dot{u}) + \frac{\partial^2}{\partial y^2}(\mu' u_y + \mu'' \dot{u}_y) - \frac{\partial p}{\partial y},$$

$$(div\Pi)_z = \Delta(\mu' u_z + \mu'' \dot{u}_z) + \frac{\partial}{\partial z} div(\mu' u + \mu'' \dot{u}) + \frac{\partial^2}{\partial z^2}(\mu' u_z + \mu'' \dot{u}_z) - \frac{\partial p}{\partial z}.$$

Отсюда непосредственно видно, что

$$div\Pi = \Delta(\mu' u + \mu'' \dot{u}) + [(\mu' + \lambda') grad\theta' + (\mu'' + \lambda'') grad\theta''] - gradp. \qquad (27)$$

§ 9. Уравнения Ламе-Навье-Стокса для упруго-вязких деформаций

Пусть $\rho$ есть масса в единице объема данного упруго-вязкого тела. Если на единицу массы действует внешняя сила $F$, то действующая на элемент $dv$ объема тела внешняя сила есть $\rho F dv.$ Далее, очевидно, что $-\rho\frac{d\dot{u}}{dt}dv$ представляет собой силу инерции рассматриваемого элемента $dv$. Если, наконец, на элемент $ds$ поверхности рассматриваемого объема $V$ действует напряжение по нормали $n$, равное $v_n$, то по принципу Д'Аламбера условие динамического равновесия объема представится в виде:

$$\int_V \rho\left(F - \frac{d\dot{u}}{dt}\right)dv + \int_S v_n ds = 0.$$

Заменяя, по формуле Гаусса, интеграл по поверхности объемным интегралом, получаем:

$$\int_V \left\{\rho\left[F - \frac{d\dot{u}}{dt}\right] + div\Pi\right\}dv = 0,$$

откуда вследствие произвола в выборе $V$ должно быть:

$$\rho\left(\frac{d\dot{u}}{dt} - F\right) = div\Pi\,. \qquad (28)$$

Подставляя найденное для $div\Pi$ значение (27), мы получим закон изменения деформаций в упруго-вязком теле с течением времени:

$$\rho\left(\frac{d\dot{u}}{dt} - F\right) = \Delta(\mu' u + \mu'' \dot{u}) + [(\mu' + \lambda') grad\theta' + (\mu'' + \lambda'') grad\theta''] - gradp. \qquad (29)$$

При $\mu'' = \lambda'' = 0$ и при $p = 0$ получаем уравнения Ламé, лежащие в основе теории упругости:

$$\rho\left(\frac{d\dot{u}}{dt}-F\right)=\mu'\Delta u+.(\mu'+\lambda')grad(divu).$$

Наоборот, при $\mu'=\lambda'=0$ из (29) получается система уравнений Навье-Стокса для движения вязких жидкостей:

$$\rho\left(\frac{d\dot{u}}{dt}-F\right)=\mu''\Delta\dot{u}+.(\mu''+\lambda'')grad(divu)-gradp,$$

где обычно полагают $\mu''+\lambda''=\frac{1}{3}\mu''$.

Наконец, применение (29) к одномерной задаче о колебаниях упруго-вязких нитей приводит к уравнению:

$$\mu'\frac{\partial^3 u}{\partial x^2\partial t}+\mu'\frac{\partial^2 u}{\partial x^2}-\rho\frac{\partial^2 u}{\partial t^2}=-\rho F,$$

которое подробно изучалось в нашей работе «Проблема упругого гистерезиса и внутреннее трение».[1]

LES PRINCIPES DE LA THÉORIE DE DÉFORMATIONS DE CORPS ÉLASTO-VISQUEUX

A.N.GUÉRASSIMOV

(Résumé)

Dans les considerérations proposes ici nous avons donné un système d'équations linéaires, contenant comme un cas particulier les équations de Lamé de la théorie d'élasticité classique et les équations de Navier-Stokes de l'hydromécanique d'un liquid visqueuix à la question des deformations de tous les corps amorphes.

[1] *См. Прикл.мат. и мех., нов. Серия, вып. 4, т. I, 1938.*

*ПММ,1939,т.II,в.4, 467-486.*

## **6.**К ВОПРОСУ О МАЛЫХ КОЛЕБАНИЯХ УПРУГО-ВЯЗКИХ МЕМБРАН

Вообразим тонкую пластинку, вещество которой обладает одновременно и упругими и вязкими свойствами. Возьмем прямоугольную систему $OXYZ$ декартовых координат так, чтобы $OX, OY$ лежали в плоскости пластинки, и разобьем последнюю на призматические столбики системой равноотстоящих друг от друга плоскостей, параллельных $YOZ$ и $ZOX$ и достаточно близких друг к другу. Будем считать, что во все время движения каждый из таких столбиков остается параллелепипедом, параллельным OZ, и не деформирует своего сечения плоскостью, параллельной $XOY$.

Такое предположение избавит от необходимости рассматривать моменты, повертывающие боковые грани столбиков. Следовательно, речь будет идти лишь о достаточно тонкой пластинке, которую в дальнейшем будем называть мембраной.[1)]

Рассмотрим три столбика $A, B, C$ находящиеся в одном ряду, параллельном $OX$, и последовательно соприкасающиеся друг с другом гранями, перпендикулярными к $OX$ (фиг.1). Если стороны основания столбика $B$ имеют длины $\Delta x, \Delta y,$ высота этого столбика есть $l$, а плотность вещества в $B$ есть $\rho$, то

$$-\rho l \Delta x \Delta y \frac{\partial^2 z}{\partial t^2} \qquad (1)$$

есть сила инерции для $B$ в тот момент, когда смещение центра его массы равно $z$.

Во-вторых, если натяжение мембраны на гранях с абсциссами $x-\frac{\Delta x}{2}, x+\frac{\Delta x}{2}$, рассчитанное на единицу площади, соответственно равно $T$ и $T'$ и образует с $OX$ углы $\alpha$ и $\alpha'$, то на $B$ действует вертикальная (параллельная $OZ$) слагающая силы натяжения, равная:

$$T'l'\Delta y \sin\alpha' - Tl\Delta y \sin\alpha \approx \Delta y\left[T'l' tg\alpha' - Tl tg\alpha\right] =$$

$$= \Delta y\left[\left(Tl\frac{\partial z}{\partial x}\right)_{x+\frac{\Delta x}{2}} - \left(Tl\frac{\partial z}{\partial x}\right)_{x-\frac{\Delta x}{2}}\right] = \Delta x \Delta y \frac{\partial}{\partial x}\left(Tl\frac{\partial z}{\partial x}\right). \qquad (2)$$

Аналогично из-за натяжения вдоль $OY$ на $B$ действует вертикальная слагающая натяжения, имеющая величину

$$\Delta x \Delta y \frac{\partial}{\partial y}\left(Tl\frac{\partial z}{\partial y}\right). \qquad (2')$$

[1]*Что касается контура, ограничивающего мембрану, то для него предполагается существование функции Грина в обычном, узком, смысле.*

Пусть, далее, на гранях $x-\dfrac{\Delta x}{2}, x+\dfrac{\Delta x}{2}$ столбика $B$ значения модуля сдвига равны соответственно $\mu$ и $\mu'$. В таком случае на $B$ действует вертикальная слагающая силы упругости, равная

$$\mu' l' \Delta y \sin\alpha' - \mu l \Delta y \sin\alpha = \Delta x \Delta y \frac{\partial}{\partial x}\left(\mu l \frac{\partial z}{\partial z}\right). \tag{3}$$

Подобным образом благодаря сдвигам в направлении $OY$ на $B$ действует сила упругости с вертикальной слагающей

$$\Delta x \Delta y \frac{\partial}{\partial y}\left(\mu l \frac{\partial z}{\partial y}\right). \tag{3'}$$

Исследуем, наконец, действие сил внутреннего трения; при этом будем исходить из закона Ньютона, по которому сила трения между двумя слоями пропорциональна градиенту скорости и площади соприкосновения слоев.

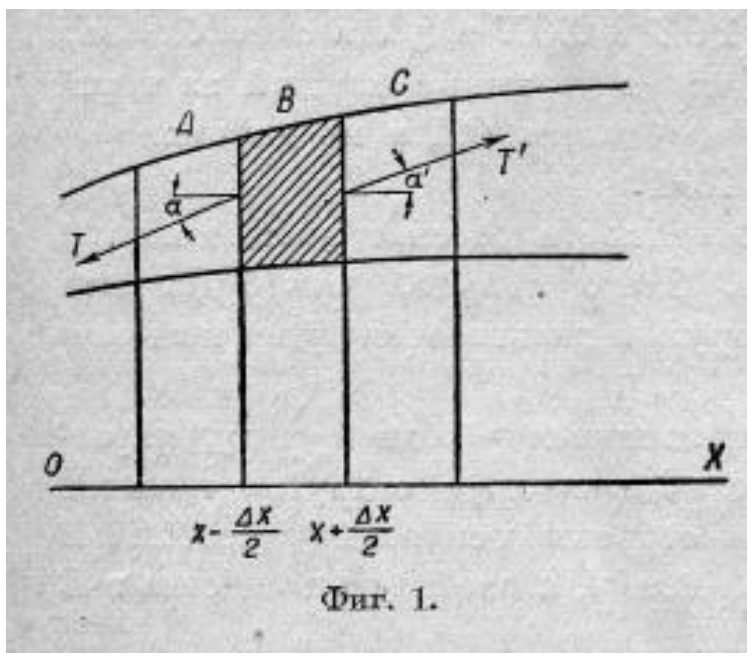

Фиг. 1.

Рассмотрим ряд столбиков, параллельный $OX$ (фиг. 1). На поверхности соприкосновения столбиков $A$ и $B$ возникает сила внутреннего трения, которая, будучи отнесена к $B$, имеют величину: $-\Delta y\left(\eta l \dfrac{\partial \dot z}{\partial x}\right)_{x-\frac{\Delta x}{2}}$. Здесь точка обозначает производную по времени, а $\eta$ есть коэффициент вязкости. Аналогично, $B$ действует на $C$ с силой: $-\Delta y\left(\eta l \dfrac{\partial \dot z}{\partial x}\right)_{x+\frac{\Delta x}{2}}$. Следовательно, столбик $C$ действует на $B$ с силой, равной

$$+\Delta y\left(\eta l \frac{\partial \dot z}{\partial x}\right)_{x+\frac{\Delta x}{2}}.$$

Поэтому в ряде столбиков, параллельном $OX$, на $B$ действует вертикальная сила внутреннего трения:

$$\Delta y\left[\left(\eta l \frac{\partial \dot z}{\partial x}\right)_{x+\frac{\Delta x}{2}} \cdot - \left(\eta l \frac{\partial \dot z}{\partial x}\right)_{x-\frac{\Delta x}{2}} \cdot\right] = \Delta x \Delta y \frac{\partial}{\partial x}\left(\eta l \frac{\partial \dot z}{\partial x}\right). \tag{4}$$

Совершенно так же найдем, что в ряде столбиков, параллельном $OY$, на $B$ действует сила внутреннего трения, имеющая величину:

$$\Delta x \Delta y \frac{\partial}{\partial y}\left(\eta l \frac{\partial \dot z}{\partial y}\right). \tag{4'}$$

На основании принципа Даламбера уравнение движения получится, если приравнять нулю сумму сил (1)-(4′). Итак, имеем:

$$\frac{\partial}{\partial x}\left[l(T+\mu)\frac{\partial z}{\partial x} + l\eta \frac{\partial \dot z}{\partial x}\right] + \frac{\partial}{\partial y}\left[l(T+\mu)\frac{\partial z}{\partial y} + l\eta \frac{\partial \dot z}{\partial y}\right] = \rho l \frac{\partial^2 z}{\partial t^2},$$

или, положив

$$T + \mu = M, \tag{5}$$

напишем уравнение колебаний пластинки в виде:

$$\frac{\partial}{\partial x}\left[\left(\eta\frac{\partial \dot z}{\partial x}+M\frac{\partial z}{\partial x}\right)l\right]+\frac{\partial}{\partial y}\left[\left(\eta\frac{\partial \dot z}{\partial y}+M\frac{\partial z}{\partial y}\right)l\right]=\rho l\frac{\partial^2 z}{\partial t^2}. \qquad (A')$$

В частности, если толщина $l$ пластинки остается постоянной, то $(A')$ принимает более простой вид:

$$\frac{\partial^2}{\partial x^2}\left[\left(\eta\frac{\partial \dot z}{\partial t}+Mz\right)\right]+\frac{\partial^2}{\partial y^2}\left[\left(\eta\frac{\partial \dot z}{\partial t}+Mz\right)\right]=\rho\frac{\partial^2 z}{\partial t^2}. \qquad \text{(A)}$$

В дальнейшем мы ограничимся лишь случаем однородной пластинки и уравнение колебаний будем писать в виде (A). В одной из –предыдущих работ [1)] мной было показано, что малые колебания упруго-вязких тел описываются при помощи обобщенной системы уравнений Навье-Стокса-Ламé, которая, будучи представлена в векторной форме, имеет вид:

$$\rho\left(F-\frac{d\dot u}{dt}\right)+\Delta(\mu' u+\mu''\dot u)+(\mu'+\lambda')grad\theta'+(\mu''+\lambda'')grad\theta''=gradp; \qquad (6)$$

Здесь $F$ – внешняя объемная сила на единицу массы, $\mu'$ и $\lambda'$ –упругие постоянные по Ламе, $\mu''$ и $\lambda''$ – аналогичные им постоянные для сил внутреннего трения, $p$ – давление,

$$\theta'=divu=\frac{\partial u_x}{\partial x}+\frac{\partial u_y}{\partial y}+\frac{\partial u_z}{\partial z}, \qquad \theta''=div\dot u=\frac{\partial \dot u_x}{\partial x}+\frac{\partial \dot u_y}{\partial y}+\frac{\partial \dot u_z}{\partial z},$$

$$\Delta=\frac{\partial^2}{\partial x^2}+\frac{\partial^2}{\partial y^2}+\frac{\partial^2}{\partial z^2}.$$

Применим (6) в случае колебания ненатянутой $(T=0)$ однородной упруго-вязкой мембраны, положив

$$F_x=F_y=F_z=0,$$

$$u=z(x,y,t),\ u_x=u_y=0,\ u_z=z,\ \frac{\partial u_x}{\partial x}=\frac{\partial u_x}{\partial y}=\frac{\partial u_x}{\partial z}=0,$$

$$\frac{\partial u_y}{\partial x}=\frac{\partial u_y}{\partial y}=\frac{\partial u_y}{\partial z}=0,\ \frac{\partial u_z}{\partial x}=\frac{\partial z}{\partial x},\ \frac{\partial u_z}{\partial y}=\frac{\partial z}{\partial y},\ \frac{\partial u_z}{\partial z}=1,\ p=0.$$

Легко найдем, что $\theta'=1, \theta''=0,$ следовательно,

$$grad\theta'=grad\theta''=0.$$

Далее вычисляем:

[1] *«Основания теории деформаций упруго-вязких тел», Прикл. мат. и мех., 1938,т.II,вып.3,379-388*

$$\Delta u_x = 0,\ \Delta u_y = 0,\ \Delta u_z = \frac{\partial^2 z}{\partial x^2} + \frac{\partial^2 z}{\partial y^2},$$

$$\Delta \dot{u}_x = 0,\ \Delta \dot{u}_x = 0,\ \Delta \dot{u}_z = \frac{\partial^2 \dot{z}}{\partial x^2} + \frac{\partial^2 \dot{z}}{\partial y^2},$$

$$\frac{du_x}{dt} = 0,\quad \frac{d\dot{u}_y}{dt} = 0,$$

$$\frac{d\dot{u}_z}{dt} = \frac{\partial \dot{z}}{\partial t} + \frac{\partial \dot{z}}{\partial z}\dot{z} = \frac{\partial \dot{z}}{\partial t} + \dot{z}\frac{\partial}{\partial t}\left(\frac{\partial z}{\partial z}\right) = \frac{\partial \dot{z}}{\partial t}.$$

Из трех уравнений системы, соответствующей векторному уравнению (6), остается лишь одно последнее, принимающее вид:

$$\frac{\partial^2}{\partial x^2}\left(\mu''\dot{z} + \mu' z\right) + \frac{\partial^2}{\partial y^2}\left(\mu''\dot{z} + \mu' z\right) = \rho\frac{\partial^2 z}{\partial t^2}.$$

Это уравнение совпадает с (A), так как $\mu' = \mu,$ а $\mu'' = \eta$ есть поперечный коэффициент трения.[1)] Если бы мембрана была натянута $(T \neq 0)$, то в таком случае создавалось бы давление $p$ с градиентом $gradp = -T\left(\frac{\partial^2 z}{\partial x^2} + \frac{\partial^2 z}{\partial y^2}\right)$, и в этом случае мы опять пришли бы к уравнению (A).

Положив

$$\eta\frac{\partial z}{\partial t} + Mz = u, \tag{7}$$

можем представить (A) в виде

$$\frac{\partial^2 u}{\partial x^2} + \frac{\partial^2 u}{\partial y^2} = \rho\frac{\partial^2 z}{\partial t^2}. \tag{8}$$

Пусть $z$ есть решение уравнения (A), удовлетворяющее условиям:

$$z\big|_{t=0} = 0;\ \frac{\partial z}{\partial t}\Big|_{t=0} = 0;$$

$$z\big|_C = 0 \text{ для всех } t; \tag{9}$$

$$\frac{\partial z}{\partial t}\Big|_{t=T} = 0;$$

[1)] *«Основания теории деформаций упруго-вязких тел», loc.cit.*

Здесь C обозначает замкнутую кривую, лежащую в плоскости $XOY$, и $T$ есть положительная величина, как угодно большая, но конечная и постоянная. Из условия на контуре C: $z|_C = 0$ следует, что для любого $t$

$$\frac{\partial z}{\partial t}|_C = 0 .$$

При этих предположениях относительно $z$, $u$ есть решение уравнения (8), удовлетворяющее условиям:

$$u|_{t=0} = 0;\ u|_C = 0 \text{ для всех } t \geq 0. \qquad (10)$$

В таком случае $u, \frac{\partial z}{\partial t}$ для всех значений $t \geq 0$. тождественно исчезают, а следовательно, $z(x, y, t) \equiv 0$.

Через точки контура $C$ проведем образующие цилиндрической поверхности параллельно оси $OT$ и рассмотрим поверхность $S$, состоящую на боковой поверхности цилиндра с направляющей $C$, из части плоскости $XOY$, лежащей внутри $C$, и на части плоскости $t = T$, ограниченной $C$. Обозначим $D$ область, ограниченную поверхностью $S$ и лежащую внутри $S$.

Пусть $P(x, y, t)$ непрерывна вместе со своими производными до второго порядка по $x$ и $y$, а $Q(x, y, t)$– до первого порядка по $t$. Тогда должно быть:[1)]

$$\iiint_D \left\{ \frac{\partial^2 P}{\partial x^2} + \frac{\partial^2 P}{\partial y^2} - \frac{\partial Q}{\partial t} \right\} dxdydt = \iint_S \frac{\partial P}{\partial x} dydt + \frac{\partial P}{\partial y} dtdx - Qdxdy. \qquad (11)$$

Положим здесь

$$P = \frac{1}{2} u^2, \qquad Q = \int_0^t \rho u \frac{\partial^2 z}{\partial t^2} dt . \qquad (12)$$

Тогда (11) принимает вид:

$$\iiint_D \left\{ \left( \frac{\partial u}{\partial x} \right)^2 + \left( \frac{\partial u}{\partial y} \right)^2 \right\} dxdydt =$$

$$\iint_S \left\{ u \frac{\partial u}{\partial x} dydt + u \frac{\partial u}{\partial y} dtdx - dxdy \int_0^t \rho u \frac{\partial^2 z}{\partial t^2} \right\} . \qquad (13)$$

Но

$$\iint_S u \frac{\partial u}{\partial x} dydt + u \frac{\partial u}{\partial y} dtdx = \int_0^T dt \int_C u \left( \frac{\partial u}{\partial x} dy + \frac{\partial u}{\partial y} dx \right).$$

[1] *Излагаемый здесь способ доказательства единственности решения является обобщением применяемого автором рассуждения для линейной задачи. См. «Проблем. упр. последейств.», Прикл. мат. и мех., т.I,в.4, 1938.*

Так как, по предположению, $u|_C = 0$, то этот интеграл исчезает. Далее, имеем

$$+\iint_S dxdy\int_0^t \rho u\frac{\partial^2 z}{\partial t^2}dt = -\iint_S dxdy\int_0^0 \rho u\frac{\partial^2 z}{\partial t^2}dt + \iint_S dxdy\int_0^T \rho u\frac{\partial^2 z}{\partial t^2}dt.$$

Здесь первый интеграл справа исчезает, второй же преобразуется на основании (7) так:

$$\int_0^T \rho u\frac{\partial^2 z}{\partial t^2}dt = \int_0^T \frac{\partial}{\partial t}\left\{\frac{\eta\rho}{2}\left(\frac{\partial z}{\partial t}\right)^2 + M\rho z\frac{\partial z}{\partial t}\right\}dt - \int_0^T M\rho\left(\frac{\partial z}{\partial t}\right)^2 dt =$$

$$= \left[\frac{\eta\rho}{2}\left(\frac{\partial z}{\partial t}\right)^2 + M\rho z\frac{\partial z}{\partial t}\right]_{t=0}^{t=T} - \int_0^T M\rho\left(\frac{\partial z}{\partial t}\right)^2 dt.$$

По условию $\frac{\partial z}{\partial t}\Big|_{t=0} = \frac{\partial z}{\partial t}\Big|_{t=T} = 0$, следовательно, вставка в правой части исчезает. Окончательно имеем

$$\iint_S dxdy\int_0^t \rho u\frac{\partial^2 z}{\partial t^2}\,dt = -\iint_S dxdy\int_0^T M\rho\left(\frac{\partial z}{\partial t}\right)^2 dt,$$

Так что (13) принимает вид:

$$\iiint_D \left\{\left(\frac{\partial u}{\partial x}\right)^2 + \left(\frac{\partial u}{\partial y}\right)^2\right\}dxdydt = +\iiint_D M\rho\left(\frac{\partial z}{\partial t}\right)^2 dxdydt. \quad (14)$$

Восстановим к $S$ направленные наружу нормали равной длины $\varepsilon$. Место концов этих нормалей есть замкнутая поверхность $\Sigma$, содержащая $S$ внутри себя. Назовем $D_+$ двухсвязную область между $S$ и $\Sigma$.

Не нарушая общности рассуждений, можно считать, что $z$ тождественно исчезает внутри $D_+$ и на $\Sigma$. Тогда применение формулы Гаусса к $D_+$ возможно, и мы имеем:

$$0 = \iiint_{D_+} \left\{\left(\frac{\partial u}{\partial x}\right)^2 + \left(\frac{\partial u}{\partial y}\right)^2\right\}dxdydt = -\iiint_D M\rho\left(\frac{\partial z}{\partial t}\right)^2 dxdydt. \quad (15)$$

Складывая почленно (14) и (15), получим:

$$\iiint_D \left\{\left(\frac{\partial u}{\partial x}\right)^2 + \left(\frac{\partial u}{\partial y}\right)^2\right\}dxdydt = 0, \quad (16)$$

Откуда следует, что $\frac{\partial u}{\partial x} = \frac{\partial u}{\partial y} = 0$. Значит, $u$ может зависеть лишь от $t$. Однако по условию, когда точка $(x, y)$ лежит на боковой поверхности цилиндра $S$, имеем $u = 0$. Таким образом, должно быть $u \equiv 0$. Обращаясь к уравнению (8) и замечая, что $\left(\frac{\partial u}{\partial x}\right)^2 = \left(\frac{\partial u}{\partial y}\right)^2 \equiv 0$, находим, что $\frac{\partial u}{\partial t}$ может быть лишь функцией координат $x, y$.

Но при $t=0, \frac{\partial z}{\partial t}=0$. Отсюда $\frac{\partial z}{\partial t}\equiv 0$. Из (7) видно, что при $u\equiv 0, \frac{\partial z}{\partial t}\equiv 0$ необходимо следует $z\equiv 0$ в $D$. Наоборот из $z\equiv 0$ в $D$ вытекает, что $z\equiv 0$ в $D_+$, если $z$ удовлетворяет условиям (9) на $S$ и исчезает вместе с $\dot z$ на $\Sigma$. Полученный результат верен при любых конечных положительных $T,\varepsilon$. В частности, он остается в силе и при $\varepsilon\to 0$. Заставляя затем $T$ стремиться к $+\infty$, мы придем к выводу, что существует лишь тождественно равное нулю решение $z$ уравнения

$$\Delta u=\rho\frac{\partial^2 z}{\partial t^2},\quad u=\eta\dot z+Mz,$$

удовлетворяющее условиям

$$z\big|_{t=0}=0,\quad \frac{\partial z}{\partial t}\Big|_{t=0}=0,\ z\big|_C=0,\qquad \frac{\partial z}{\partial t}\Big|_{t=+\infty}=0.$$

Пусть теперь имеем неоднородное уравнение

$$\Delta(\eta\dot z+Mz)=\rho\ddot z-f(x,y,t), \tag{B}$$

причем $z$ должно удовлетворять условиям:

$$z\big|_{t=0}=\Phi(x,y);\quad \dot z\big|_{t=0}=\varphi(x,y);\qquad z\big|_C=0;\ \dot z\big|_{t=+\infty}=0. \tag{$B'$}$$

Если бы уравнение (B) при условиях $(B')$ имело два различных решения $z_1$ и $z_2$, то их разность $z=z_2-z_1$ была бы решением однородного уравнения (A) с условиями

$$z\big|_{t=0}=0;\quad \dot z\big|_{t=0}=0;\quad z\big|_C=0;\quad \dot z\big|_{t=+\infty}=0,$$

Причем, вопреки доказанному, $z$ не исчезало бы тождественно. Поэтому (B) может иметь не больше одного решения $z(x,y,t)$, обладающего непрерывными производными до второго порядка включительно и удовлетворяющего условиям $(B')$. Смысл условия $\dot z\big|_{t=+\infty}=0$ тот, что как бы велик ни был конечный промежуток времени, в течение которого мембрана подвергается воздействию внешней силы $f$, кинетическая энергия мембраны при $t=+\infty$ должна неизбежно обратиться в нуль.

Вернемся к однородному уравнению:

$$\Delta(\eta\dot z+Mz)=\rho\frac{\partial^2 z}{\partial t^2}. \tag{A}$$

Было показано, что это уравнение может иметь не больше одного решения $z(x,y,t)$, удовлетворяющего условиям:

$$z\big|_{t=0}=\Phi(x,y);\quad \frac{\partial z}{\partial t}\Big|_{t=0}=\varphi(x,y);\qquad z\big|_C=0;\ \frac{\partial z}{\partial t}\Big|_{t=+\infty}=0 \tag{$A'$}$$

Такое решение мы найдем, применяя способ Д.Бернулли. Для этого будем искать интеграл уравнения (A) в виде:

$$z=ZT,$$

где $Z$ есть функция координат $x,y,$ а $T$ зависит только от времени. Тогда

$$\eta\frac{\partial z}{\partial t}+Mz=(\eta T'+MT)Z,$$

$$\Delta\left(\eta\frac{\partial z}{\partial t}+Mz\right)=(\eta T'+MT)\Delta Z\,,$$

$$\rho\frac{\partial^2 z}{\partial t^2}=\rho ZT''\,.$$

Подставляя это в (А), приводим его к виду:

$$\frac{\Delta Z}{Z}=\frac{\rho T''}{\eta T'+MT}\,.$$

Так как между положением точки $(x, y)$ на мембране и временем $t$ не может быть никакой зависимости, то каждая часть предыдущего равенства должна быть постоянной. Обозначая последнюю $(-\lambda)$, получаем два независимых уравнения:

$$\Delta Z+\lambda Z=0, \tag{17}$$

$$\rho T''+\lambda\eta T'+\lambda MT=0\,. \tag{18}$$

Рассмотрим подробнее уравнение (17). Если $\lambda=0,$ то это уравнение не удовлетворяется никаким $Z(x, y)$, исчезающим на контуре C, кроме тождественного нуля. В самом деле, гармоническая функция $Z$, исчезающая на границе области, тождественно равна нулю в этой области. Более того, ни при каком отрицательном $\lambda$ уравнение (17) не может иметь отличного от нуля решения, обращающегося в нуль на C. Действительно, если бы такое решение $Z(x, y)$ существовало, то оно необходимо должно было бы достигать экстремума в некоторой точке $P(x_0, y_0)$ внутри C. Этот экстремум был бы положительным в случае максимума или отрицательным в случае минимума. В случае максимума должна было бы иметь место

$$\left(\frac{\partial^2 Z}{\partial x^2}\right)_{x_0,y_0}<0\,,\ \left(\frac{\partial^2 Z}{\partial y^2}\right)_{x_0,y_0}<0\,,\ \lambda Z<0\,.$$

Наоборот, для минимума

$$\left(\frac{\partial^2 Z}{\partial x^2}\right)_{x_0,y_0}>0\ \left(\frac{\partial^2 Z}{\partial y^2}\right)_{x_0,y_0}>0\,,\ \lambda Z>0\,.$$

Очевидно, ни тот, ни другой случай несовместим с уравнением (17). Поэтому краевая задача $\Delta Z+\lambda Z=0,\ Z\big|_C=0$ может допускать не исчезающее тождественно решение лишь для существенно положительных значений $\lambda$.

Мы рассмотрим сначала уравнение

$$\Delta Z=f(x, y)\,, \tag{19}$$

где $f(x, y)$ непрерывна внутри C. Найдем решение (19), обращающееся в нуль на C. Пусть $u, v-$ две функции от $x, y$, непрерывные в области $D$ с их частными производными до второго порядка. Тогда должно быть:

$$\iint\limits_D\{v\Delta u-u\Delta v\}d\xi d\eta=\int\limits_{\overset{C}{\mapsto}}\left\{u\frac{\partial v}{\partial n}-v\frac{\partial u}{\partial n}\right\}ds\,, \tag{20}$$

где $D$–область, ограниченная контуром С, пробегаемым в положительном направлении; на последнее указывает стрелка под знаком С. Здесь $n$–направление внутренней к С нормали. Положим в (20) $u = Z, v = G(x,y,\xi,\eta)$, причем $(x,y) \subset D$.

$$G(x,y,\xi,\eta) = \ln\frac{1}{r} + g(x,y,\xi,\eta) \quad ,$$

где $r^2 = (x-\xi)^2 + (y-\eta)^2$, а $g$– гармоническая в $D$ функция, принимающая на С значения, равные $-\ln\frac{1}{r}$. Тогда $G(x,y,\xi,\eta)$– функция гармоническая в $D$ всюду, кроме точки $\xi = x, \eta = y$ и обращающаяся в нуль на С. Описав вокруг точки $(x,y) \subset D$ окружность Г радиусом $\varepsilon$ так, чтобы Г лежала внутри С, применим (20) к области $D'$, ограниченной С и Г. Так как на С имеем $u = Z = 0, v = G = 0$ и так как в $D' : \Delta v = 0, \Delta Z = f(\xi,\eta)$, то (20) принимает вид:

$$\iint_{D'} f(\xi,\eta) G(x,y;\xi,\eta) d\xi d\eta =$$

$$=- \int_{\overset{\Gamma}{\mapsto}} \left\{ Z\frac{\partial \ln\frac{1}{r}}{\partial n} + Z\frac{\partial g}{\partial n} - \frac{\partial Z}{\partial n}\ln\frac{1}{r} - \frac{\partial Z}{\partial n} g \right\} ds. \qquad (21)$$

Вследствие того, что $Z$ и $g$ с их производными по нормали ограничены на Г, то второй и четвертый члены правой части при достаточно малом $\varepsilon$ могут быть сделаны как угодно малы по абсолютной величине. Это же относится и к третьему члену, потому что

$$\left| \int_\Gamma \frac{\partial Z}{\partial n} \ln\frac{1}{r} ds \right| \le \max\left|\frac{\partial Z}{\partial n}\right|_\Gamma \cdot \left|\ln\frac{1}{r}\right| \cdot \varepsilon \cdot 2\pi = \max\left|\frac{\partial Z}{\partial n}\right|_\Gamma \cdot 2\pi\frac{\ln N}{N},$$

где $N > 1$ тем больше, чем меньше $\varepsilon$. А так как $\frac{\ln N}{N} \to 0$, когда $N$ возрастает, то высказанное утверждение становится очевидным.

Рассмотрим теперь первый член справа. Имеем:

$$\int_\Gamma Z\frac{\partial \ln\frac{1}{r}}{\partial r} ds \le \max(Z)_\Gamma \cdot \int_\Gamma -\frac{\partial \ln\frac{1}{r}}{\partial r} r d\varphi = \max(Z)_\Gamma \cdot 2\pi .$$

Если $\varepsilon \to 0$, то $\max(Z)_\Gamma \to Z(x,y)$. Поэтому (21) дает:

$$Z(x,y) = -\frac{1}{2\pi}\iint_D f(x,y;\xi,\eta) d\xi d\eta . \qquad (22)$$

В частности, если положить $f(x,y) = -\lambda Z(x,y)$, то (22) принимает вид:

$$Z(x,y) = \lambda \iint_D K(x,y;\xi,\eta) Z(\xi,\eta) d\xi d\eta, \qquad (23)$$

где

$$K = \frac{1}{2\pi} G, \quad \lambda > 0 .$$

Ядро $K$ интегрального уравнения (23) при $\xi = x, \eta = y$ логарифмически обращается в бесконечность. Однако это обстоятельство не является существенным, так как $K(x, y; \xi, \eta)$ обладает интегрируемым в $D$ квадратом.[1)] Также несущественно для дальнейшего и то, что оно не принадлежит к числу регулярных, так как после первой итерации получается ограниченное ядро.

Далее, можно доказать симметричность $K(x, y, \xi, \eta)$. Возьмем внутри $D$ две точки $A'(x', y')$ и $A''(x'', y'')$ и опишем около них как центров окружности $\Gamma', \Gamma''$ радиусами, соответственно равными $\varepsilon'$ и $\varepsilon''$ и настолько малыми, что $\Gamma', \Gamma''$ целиком лежат внутри $D$. Применим (20) к области $D'$, ограниченной $C, \Gamma'$ и $\Gamma''$, положив:

$$v = G'(x', y'; \xi, \eta) = \ln\frac{1}{r'} + g'(x', y'; \xi, \eta) = G(x', y'; \xi, \eta),$$

$$v = G''(x'', y''; \xi, \eta) = \ln\frac{1}{r''} + g''(x'', y''; \xi, \eta) = G(x'', y''; \xi, \eta),$$

где $\quad r'^2 = (x' - \xi)^2 + (y' - \eta)^2, \quad r''^2 = (x'' - \xi)^2 + (y'' - \eta)^2.$

Так как в $D'$ имеем $\Delta G' = \Delta G'' = 0$, причем $G'$ и $G''$ на С обращаются в нуль, то (20) дает:

$$-\int\limits_{\substack{\Gamma' \\ \downarrow\mapsto\uparrow}} \left\{G''\frac{dG'}{dn} - G'\frac{dG'}{dn}\right\}ds - \int\limits_{\substack{\Gamma'' \\ \downarrow\mapsto\uparrow}} \left\{G''\frac{dG'}{dn} - G'\frac{dG''}{dn}\right\}ds = 0.$$

Мы видели уже, что первый интеграл стремился к $-2\pi\Pi(x'', y'', x', y')$ при $\varepsilon \to 0$. Аналогично найдем, что второй стремится к $+2\pi G(x', y'; x'', y'')$ при $\varepsilon'' \to 0$. Предыдущее равенство позволяет видеть, что $G(x, y; \xi, \eta) = G(\xi, \eta; x, y)$, если $(x, y)$ и $(\xi, \eta)$ лежат внутри $D$.

Так что ядро $K(x, y; \xi, \eta)$ обладает интегрируемым в $D$ квадратом и симметрично, то к нему применимы результаты теории, развитой Гильбертом

[1] *Интегрируемость $K^2$ прежде всего установить, если убедиться, что* $\iint (\ln r)^2 d\xi d\eta$ *остается конечным в круге радиуса $R$, описанном из точки $(x, y)$ как центра. Не нарушая общности, можно принять $x = y = 0$. Тогда*

$$\iint (\ln r)^2 d\xi d\eta = 2\pi\int\limits_0^R r(\ln r)^2 dr = 2\pi\left[e^u\left(u^2 - 2u + 2\right)\right]_{-\infty}^{\ln R},$$

*что представляет собой конечную величину. Аналогично убедимся в конечности интеграла* $\iint \ln r d\xi d\eta = 2\pi\int\limits_0^R r\ln r dr$ *,распространенного на ту же площадь, а следовательно, и интеграла*

$$\iint\limits_D [G(x, y; \xi, \eta)]^2 d\xi d\eta = \iint\limits_D \left[\left(\ln\frac{1}{r}\right)^2 + 2g(x, y; \xi, \eta)\ln\frac{1}{r} + \{g(x, y; \xi, \eta)\}^2\right] d\xi d\eta.$$

*Так как функция* $\iint dxdy \iint\limits_D [Kx, y; \xi, \eta)]^2 d\xi d\eta \le M$ *,где $M$ – конечное число, что и требовалось доказать. (см. Курант-Гильберт, Методы математической физики, ГТИ, стр. 347, 1933.)*

и Шмидтом.[1] В частности, если $Z(x,y)$ есть непрерывное в $D$ вместе со своими производными первого порядка решение уравнения (23), то оно имеет непрерывные в $D$ производные второго порядка, удовлетворяя уравнению $\Delta Z + \lambda Z |= 0$, причем $Z|_C = 0$.[2] Далее, можно утверждать, что существует по крайней мере одно положительное число $\lambda$, для которого (23) имеет не равное тождественное нулю решение $Z(x,y)$. Более того, таких чисел $\lambda_k$ должно быть бесконечное, именно счетное, множество, так как в противном случае разложение ядра $K(x,y;\xi,\eta)$ в билинейный ряд по фундаментальным функциям не могло бы дать требуемой особенности $K(x,y;\xi,\eta)$ в точке $\xi = x$, $\eta = y$. Итак, существует счетное множество положительных фундаментальных[3] чисел $\lambda_k$ и соответствующих фундаментальных функций $Z_k(x,y)$, причем все числа $\lambda_k$, $k = 1,2,...$, можно считать занумерованными в порядке их возрастания. При этом ряд

$$\sum_{k=1}^{\infty} \frac{1}{\lambda_k^2} \tag{24}$$

должен непременно сходиться, следовательно, множество чисел $\lambda_k$ не может иметь предельной точки на конечном расстоянии. Действительно, в противном случае в написанном выше ряде было бы бесконечно большое количество конечных членов, как угодно мало отличающихся друг от друга. Но тогда этот ряд не мог бы сходиться. Таким образом можно указать такое конечное число $\Lambda, \Lambda > 0$, что для любого $k$ имеет место неравенство:

$$\lambda_k - \lambda_{k-1} \geq \Lambda. \tag{25}$$

Более того, можно показать, что с возрастанием номера $k$ фундаментальное число $\lambda_k$ становится бесконечно большой величиной того же порядка, что и $k$, какова бы ни была конечная замкнутая область $D$. Для доказательства рассмотрим два квадрата $Q', Q''$ со сторонами $\alpha', \alpha'', \alpha' > \alpha''$, из которых первый заключает в себе область $D$, а другой содержится в этой области. Пусть $\lambda_{k'}, \lambda_k$ и $\lambda_{k''}$ - фундаментальные числа номера k рассматриваемой краевой задачи, но вычисленные соответственно для $Q', D, Q''$. Так как частоты обертона порядка $k$ для мембраны должны быть тем выше, чем меньше площадь мембраны, то имеем очевидное соотношение: $\lambda'_{k_k} \leq \lambda_k \leq \lambda''_k$.

[1] *См. И.И. Привалов. Интегральные уравнения, ОНТИ, стр.166-169, 1935.*

[2] *Курант-Гильберт. Методы матем. физики,I, ГТТИ, стр. 343-344, 1933.*

[3] *Согласно принятой терминологии, фундаментальными числами рассматриваемой задачи являются не $\lambda_k$, а противоположные им числа -$\lambda_k$. См. Е.Гурса, Курс матем. анализа, т. III, 2, ОНТИ, стр. 189, 1934.*

Но исследование краевой задачи $\Delta Z + \lambda' Z = 0, Z\big|_{q'} = 0$ для квадрата 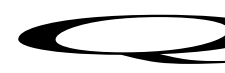
со стороной $\alpha'$ может быть выполнено путем разделения переменных. Положив $Z = XY$, где $X, Y$ соответственно функции только $x$ и только $y$, вместо $\Delta Z + \lambda' Z = 0$ получим:

$$X''Y + XY'' + \lambda' XY = 0\,,$$

или

$$\frac{X''}{X} = -\frac{Y'' + \lambda' Y}{Y} = -\kappa'\,,$$

где $\kappa'$ должна быть постоянной. Тогда необходимо имеем:

$$X = \sin\sqrt{\kappa'}x\,, \quad Y = \sin\sqrt{\lambda' - \kappa'}y$$

для того, чтобы $Z = XY$ исчезала на сторонах квадрата $Q'$; последний мы предполагаем расположенным так, что его вершины имеют координаты $(0,0), (\alpha', 0), (\alpha', \alpha'), (0, \alpha')$.

Отсюда следует, что

$$\kappa'_m = \frac{m^2\pi^2}{\alpha'^2}, \quad \lambda'_{mn} = \frac{n^2 + m^2}{\alpha'^2}\pi^2,$$

где $m, n$ – любые целые положительные числа. Определим число $\Lambda'$ фундаментальных чисел $\lambda''_{mn}$ этой задачи, не превышающих $\lambda'$. Очевидно, оно равно числу точек $(m, n)$ с целыми и положительными $m, n$, лежащих внутри квадранта:

$$m^2 + n^2 \le \frac{\lambda'\alpha'^2}{\pi^2}, m > 0, n > 0.$$

По мере того как $\lambda'$ растет, число $\Lambda'$ все более и более приближается к площади этого квадранта, равной $\dfrac{\lambda'\alpha'^2}{4\pi}$. Поэтому

$$\Lambda' \approx \frac{\lambda'\alpha'^2}{4\pi},$$

и, следовательно,

$$\frac{\lambda'}{\Lambda'} \approx \frac{4\pi}{\alpha'^2}.$$

В частности, если расположить фундаментальные числа $\lambda'$ в возрастающей последовательности, то

$$\frac{\lambda'_k}{k} \approx \frac{4\pi}{\alpha'^2}.$$

Аналогично, для квадрата $Q''$ со стороной $\alpha''$ найдем: $\frac{\lambda_k''}{k} \approx \frac{4\pi}{\alpha''^2}.$

Поэтому благодаря неравенству $\lambda_k' \le \lambda_k \le \lambda_k''$ видим, что $\lambda_k$ при достаточно больших $k$ оказывается заключенным между $k \cdot \frac{4\pi}{\alpha'^2}$ и $k \cdot \frac{4\pi}{\alpha''^2}$, так что имеет место равенство

$$\lambda_k = k \cdot \left[ \frac{4\pi}{\alpha'^2} + \omega \cdot 4\pi \left( \frac{1}{\alpha''^2} - \frac{1}{\alpha'^2} \right) \right], \tag{26}$$

где

$$0 \le \omega \le 1.$$

Покажем теперь, что ранг каждого из чисел $\lambda_k, k = 1,2,...,$не больше единицы, т.е. что каждому $\lambda_k$ соответствует одна и только одна фундаментальная функция $Z_k(x,y)$ данной краевой задачи $\Delta Z + \lambda_k Z = 0,\ Z\big|_C = 0.$ Для этого достаточно показать, что два исчезающих на контуре $C$ решения $Z_k', Z_k''$ уравнения $\Delta Z + \lambda_k Z = 0$ связаны между собой соотношением $Z_k'' = \alpha Z_k'$ ,где $\alpha = const.$ Пусть $\Gamma -$ произвольный замкнутый контур, целиком лежащий в $D.$ Применим формулу Грина:

$$\iint\limits_D \{Z_k'' L Z_k' - Z_k' L Z_k''\} dxdy =$$

$$= \int\limits_{\substack{C \\ \downarrow\mapsto\uparrow}} \left\{ Z_k' \frac{\partial Z_k''}{\partial n} - Z_k'' \frac{\partial Z_k'}{\partial n} \right\} ds - \int\limits_{\substack{\Gamma \\ \downarrow\mapsto\uparrow}} \left\{ Z_k' \frac{\partial Z_k''}{\partial n} - Z_k'' \frac{\partial Z_k'}{\partial n} \right\} ds$$

к области $D'$ ,лежащей между $C$ и $\Gamma$.

Здесь положено $LZ \equiv \Delta Z + \lambda_k Z.$

Так как $LZ_k' = LZ_k'' = 0$ в $D'$ и $Z_k' = Z_k'' = 0$ на $C$,

то остается:

$$\int\limits_{\substack{\Gamma \\ \downarrow\mapsto\uparrow}} \left\{ Z_k' \frac{\partial Z_k''}{\partial n} - Z_k'' \frac{\partial Z_k'}{\partial n} \right\} ds = 0,$$

где $n -$ любое направление; следовательно,

$$\frac{\partial}{\partial n} \ln Z_k'' = \frac{\partial}{\partial n} \ln Z_k',$$

откуда

$$Z_k'' = \alpha Z_k',\ \ \alpha = const,$$

*что и требовалось доказать.*

Фундаментальные функции

$$Z_k(x,y) \quad k=1,2,... \tag{26'}$$

попарно ортогональны друг другу, т.е.

$$\iint_D Z_i(x,y)Z_j(x,y)dxdy=0\,,\ i\neq j, i, j=1,2,... \tag{27}$$

Кроме того, благодаря произволу в выборе постоянных, подобных $\alpha$, эти функции можно считать нормированными, так что

$$\iint_D [Z_k(x,y)]^2 dxdy=1,\ k=1,2,... \tag{28}$$

Наконец, заметим, что система фундаментальных функций $Z_k(x,y)$ непременно должна быть замкнутой в классе непрерывных с частными производными 1-го порядка функций; это значит, что не существует никакой функции $\zeta(x,y)$, обладающей в $D$ непрерывными производными 1-го порядка и ортогональной ко всем $Z_k$, то она должна быть ортогональной и к ядру $K(x,y;\xi,\eta)$, так что

$$\Xi(x,y)=\iint_D K(x,y;\xi,\eta)\zeta(\xi,\eta)d\xi d\eta\equiv 0\,.$$

С другой стороны, если $\xi(x,y)$ непрерывна в $D$ вместе со своими производными 1-го порядка, то $\Xi(x,y)$ должна иметь непрерывные вторые производные и удовлетворять условиям:

$$\Delta\Xi(x,y)=\zeta(x,y),\quad \Xi(x,y)\big|_C=0\,.$$

Так как $\Xi(x,y)\equiv 0$, то первое их этих условий дает $\zeta(x,y)\equiv 0$, *что и требовалось доказать.*

Для краевой задачи штурм-лиувиллевского типа с одним переменным можно доказать ограниченность всех фундаментальных функций в их совокупности. Но для рассматриваемого случая краевой задачи с двумя переменными это положение перестает быть верным.[1)] Однако можно показать, что каждая из рассматриваемых фундаментальных функций из рассматриваемых фундаментальных функций $Z_k(x,y)$ не превосходит по модулю некоторого положительного не зависящего от $k$ и конечного числа, умноженного на соответствующее фундаментальное число $\lambda_k$. Это обстоятельство является простым следствием неравенства Шварца. Именно, так как

$$Z_k(x,y)=\lambda_k\iint_D K(x,y;\xi,\eta)Z_k(\xi,\eta)d\xi d\eta\,, \tag{29}$$

то, применяя названное неравенство, имеем:

[1] W.Sternberg, Üb. die asyptot.Integration part. Diffgl., II, Math. Ann.86.

$$|Z_k(x,y)|^2 = \left|\lambda_k \iint_D K(x,y;\xi,\eta) Z_k(\xi,\eta)\, d\xi d\eta\right|^2 \le$$

$$\le |\lambda_k|^2 \iint_D [K(x,y;\xi,\eta)]^2 d\xi d\eta \cdot \iint_D [Z_k(\xi,\eta)]^2 d\xi d\eta =$$

$$= |\lambda_k|^2 \iint_D [K(x,y;\xi,\eta)]^2 d\xi d\eta\ .$$

Но функция от $x$ и $y$,

$$= \iint_D [K(x,y;\xi,\eta)]^2 d\xi d\eta,$$

как мы видели, непрерывна в замкнутой области $D+C$, следовательно, ограничена там и снизу и сверху. Поэтому существует такое не зависящее от $k$, $x$ и $y$ положительное конечное число $M^2$, настолько большое, что для данной ограниченной области $D$

$$\iint_D [K(x,y;\xi,\eta)]^2 d\xi d\eta \le M^2.$$

поэтому

$$[Z_k(\xi,\eta)]^2 \le |\lambda_k|^2 \cdot M^2,$$

*что и требовалось доказать.*

Отсюда вытекает важное *следствие.* Рассмотрим ряд Фурье:

$$\sum_{k=1}^{\infty} \zeta_k Z_k(x,y), \tag{31}$$

расположенный по функциям $Z_k(x,y)$ рассматриваемой нормированной и ортогональной системы. Если порядок малости его коэффициентов $\zeta_k$ при возрастании $k$ делается не ниже, чем $k^{-3}$, то ряд (31) сходится абсолютно и равномерно.

Пусть $\zeta_k$ при $k \to \infty$ имеет порядок малости не ниже, чем $k^{-4}$. Тогда по отношению к ряду (31) возможно почленное применение операции $\Delta$. В самом деле, так как

$$\Delta Z_k(x,y) = -\lambda_k Z_k(x,y),$$

то ряд

$$\sum_{k=1}^{\infty} \zeta_k \Delta Z_k(x,y) = -\sum_{k=1}^{\infty} \zeta_k \lambda_k Z_k(x,y) \tag{32}$$

сходится равномерно и абсолютно, если только его коэффициенты $\zeta_k, \lambda_k$ имеют порядок малости не ниже чем $k^{-3}$. А так как $\lambda_k$ при $k \to \infty$ растет так же, как $k$, то для указанной сходимости ряда (32) достаточно, чтобы порядок $\zeta_k$ был по крайней мере равен $k^{-4}$.

Обратимся теперь к уравнению (18), положив в нем $\lambda = \lambda_k$:

$$\rho T'' + \lambda_k \eta T' + \lambda_k M T = 0. \tag{33}$$

Подстановка $T = e^{q_k t}$ приводит к характеристическому уравнению:

$$\rho q_k^2 + \eta \lambda_k q_k + M\lambda_k = 0, \tag{34}$$

откуда

$$q_{k'_k}, q_{k''} = \frac{-\lambda_{k''} \pm \sqrt{\lambda_k^2 \eta^2 - 4M\rho\lambda_k}}{2\rho}, \ k = 1,2,... \tag{35}$$

Следовательно, общий интеграл уравнения (18) есть:

$$T_k = C_{k'} e^{q'_k t} + C_{k''} e^{q''_k t}, \ k = 1,2,..., \tag{36}$$

где $C_{k'}, C_{k''}$ – постоянные, значения которых пока не определены.

Из (26) и (36) можно построить решение $z = \sum Z_k T_k$ уравнения (А), обращающееся в нуль на контуре $C$. Для этого достаточно положить:

$$z(x,y,t) = \sum_{k=1}^{\infty} \left\{ C_{k'} e^{q'_{k'} t} + C_{k''} e^{q''_{k'} t} \right\} Z_k(x,y), \tag{37}$$

разумеется, в предположении, что ряд (37) сходится.

Легко видеть, что (37) удовлетворяет условию:

$$\frac{\partial z}{\partial t}\Big|_{t=+\infty} = 0,$$

так как при $t \to +\infty$ производная по $t$ от каждого члена суммы (37) вследствие свойства $q_{k'}, q_{k''}$ стремится к нулю независимо от $(x, y)$.

Чтобы соблюсти первые два условия $(A')$, нужно только выбрать надлежащим образом коэффициенты $C_{k'}, C_{k''}$.

Пусть

$$\Phi(x,y) = \sum_{k=1}^{\infty} A_k Z_k(x,y), \quad \varphi(x,y) = \sum_{k=1}^{\infty} B_k Z_k(x,y) \tag{38}$$

определяет разложение данных функций $\Phi, \varphi$ по фундаментальным функциям $Z_k(x,y)$ ядра $K(x,y;\xi,\eta)$. Чтобы удовлетворить, по крайней мере формально, условиям

$$z\big|_{t=0} = \Phi(x,y), \quad \frac{\partial z}{\partial t}\Big|_{t=0} = \varphi(x,y),$$

достаточно положить

$$C_{k'} + C_{k''} = A_k,$$

$$q_{k'}C_{k'} + q_{k''}C_{k''} = B_k. \qquad k = 1,2,... \tag{39}$$

При выбранных таким образом $C_{k'}, C_{k''}$ выражение (37) с формальной стороны удовлетворит всем условиям задачи, и нужно только показать, что ряд (37), а также те ряды, которые получаются из него дифференцированием столько раз, сколько нужно для составления уравнения (A), равномерно и абсолютно сходятся для рассматриваемых значений $(x,y)$ и для всех $t \geq 0$.

Пусть коэффициенты $A_k, B_k$ разложений (38) функций $\Phi, \varphi$ в ряды Фурье имеют при $k \to \infty$ порядок малости не ниже, чем $k^{-4}$. Так как $q_{k'}, q_{k''}$ при

$k \to \infty$ растут по модулю так же, как и $k$, то для этого достаточно предположить, что в таком случае ряды

$$\sum_{k=1}^{\infty} A_k Z_k(x,y), \quad \sum_{k=1}^{\infty} B_k Z_k(x,y)$$

равномерно и абсолютно сходятся к функциям $\Phi(x,y), \varphi(x,y)$. Более того, по отношению к этим рядам допустимо почленное применение операции $\Delta$, причем получившиеся от этого ряды сходятся равномерно и абсолютно к функциям:

$$\Delta\Phi(x,y), \Delta\varphi(x,y). \, (t \geq 0)$$

Так как умножение чисел $C_{k'}, C_{k''}, C_{k'}q_{k'}, C_{k''}q_{k''}$ на величины $e^{q_{k'}t}, e^{q_{k''}t}$ $(t \geq 0)$

не увеличивает модуля этих чисел, то ряды

$$\sum_{л=1}^{\infty} \left\{ C_{k'}e^{q_{k'}t} + C_{k''}e^{q_{k''}t} \right\} Z_k(x,y),$$

$$\sum_{л=1}^{\infty} \left\{ C_{k'}q_{k'}e^{q_{k'}t} + C_{k''}q_{k''}e^{q_{k''}t} \right\} Z_k(x,y),$$

$$-\sum_{л=1}^{\infty} \left\{ C_{k'}e^{q_{k'}t} + C_{k''}e^{q_{k''}t} \right\} \lambda_k Z_k(x,y),$$

$$-\sum_{л=1}^{\infty} \left\{ C_{k'}q_{k'}e^{q_{k'}t} + C_{k''}q_{k''}e^{q_{k''}t} \right\} \lambda_k Z_k(x,y)$$

равномерно и абсолютно сходятся к функциям, соответственно равным

$$z(x,y);\ \frac{\partial z}{\partial t};\ \Delta z(x,y);\ \Delta\frac{\partial z}{\partial t}.$$

Кроме того, из сходимости последнего из четырех написанных выше рядов благодаря совпадению порядков $q_{k'}$ и $q_{k''}$ с порядком $\lambda_k$ при $k\to\infty$ получается, что и ряд

$$\sum_{л=1}^{\infty}\left\{C_{k'}q_{k''}^2e^{q_{k'}t}+C_{k''}q_{k''}^2e^{q_{k''}t}\right\}Z_k(x,y)$$

сходится равномерно и абсолютно. Но последний ряд получается путем почленного двукратного дифференцирования по $t$ ряда для $z(x,y,t)$; следовательно, его сумма равна $\frac{\partial^2 z}{\partial t^2}$.

Из этих соображений следует, что если порядок малости коэффициентов $A_k, B_k$ в разложениях для $\Phi,\varphi$ не ниже, чем $k^{-4}$, то все необходимые действия, формально выполняемые над найденным решением $z(x,y,t)$; для составления уравнения (A), вполне законны, и (37) действительно есть решение задачи о свободных колебаниях упруго-вязкой мембраны, закрепленной в точках контура $C$, на котором она натянута. Мы видели уже, что другого решения этой задачи быть не может.

Рассмотрим теперь вопрос *о вынужденных колебаниях упруго-вязкой мембраны, закрепленной в точках замкнутого контура $C$, на который она натянута.* С математической точки зрения этот вопрос сводится к изысканию решения $z$ неоднородного уравнения

$$\frac{\partial^2}{\partial x^2}\left(\eta\dot z+Mz\right)+\frac{\partial^2}{\partial y^2}\left(\eta\dot z+Mz\right)=\rho\frac{\partial^2 z}{\partial t^2}-f(x,y,t), \qquad (B)$$

которое удовлетворяло бы условиям:

$$z\big|_{t=0}=\Phi(x,y);\ \ \dot z\big|_{t=0}=\varphi(x,y);\ z\big|_C=0;\ \ \dot z\big|_{t=+\infty}=0. \qquad (B')$$

Мы уже видели, что если такое решение существует, то оно единственное. А так как решение однородного уравнения (A) уже найдено, то остается показать способ определения частного решения уравнения (B). Положим

$$f(x,y,t)=\sum_{k=1}^{\infty}w_k(t)Z_k(x,y), \qquad (40)$$

где $Z_k(x,y),\ k=1,2,...,-$ функции, принадлежащие нормальной ортогональной системе ядра $K(x,y;\xi,\eta)$., рассмотренные выше.

Решение задачи (B), $(B')$ будем искать в виде:

$$z(x,y,t)=\sum_{k=1}^{\infty}\{C_{k'}e^{q_{k'}t}+C_{k'}e^{q_{k''}t}+a_k(t)\}Z_k(x,y) \qquad (41)$$

с определенными раньше значениями $C_{k'}, C_{k''}, q_{k'}, q_{k''}$. Здесь $a_k(t)$ являются пока не определенными функциями.

Подставляя (41) в (В), найдем:

$$\sum_{k=1}^{\infty}\{-\left[\rho q_{k'}^2+\eta\lambda_k q_{k'}+M\lambda_k\right]C_{k'}e^{q_{k'}t}-\left[\rho q_{k''}^2+\eta\lambda_k q_{k''}+M\lambda_k\right]C_{k''}e^{q_{k''}t}-$$

$$-\left[\rho\ddot{a}_k+\eta\lambda_k\dot{a}_k+M\lambda_k a_k\right]\}Z_k(x,y)=-\sum_{k=1}^{\infty}w_k(t)Z_k(x,y).$$

По определению чисел $q_{k'}q_{k''}$ выражения в квадратных скобках под знаком суммы в левой части исчезают, следовательно, написанное тождество будет удовлетворено, если выбрать $a_k$ согласно уравнениям:

$$\rho a_{k''}(t)+\eta\lambda_k a_{k'}(t)+M\lambda_k a_k(t)=w_k(t) \qquad (k=1,2,...). \qquad (42)$$

Уравнение (42) может быть представлено в виде:

$$\left(\rho D^2+\eta\lambda_k D+M\lambda_k E\right)a_k(t)=w_k(t),$$

где $D=\dfrac{d}{dt}$ и $E$ – оператор тождественного преобразования. Отсюда, наоборот,

$$a_k(t)=\left(\rho D^2+\eta\lambda_k D+M\lambda_k E\right)^{-1}w_k(t).$$

Рассмотрим случай, когда $q_{k'}\neq q_{k''}$. Тогда можно написать:

$$a_k(t)=\frac{1}{\rho\left(q_{k'}-q_{k''}\right)}\left(D-q_{k'}E\right)^{-1}w_k(t)+\frac{1}{\rho\left(q_{k''}-q_{k'}\right)}\left(D-q_{k''}E\right)^{-1}w_k(t)=$$

$$=a_{1k}(t)+a_{21k}(t),$$

где положено

$$a_{1k}(t)=\frac{1}{\rho\left(q_{k'}-q_{k''}\right)}\left(D-q_{k'}E\right)^{-1}w_k(t),\quad a_{2k}(t)=\frac{1}{\rho\left(q_{k''}-q_{k'}\right)}\left(D-q_{k''}E\right)^{-1}w_k(t).$$

Отсюда легко найти, что исчезающее при $t=0$ значение функции $a_k(t)$ должно быть таково:

$$a_k(t)=\frac{1}{\rho\left(q_{k'}-q_{k''}\right)}\int_0^t w_k(\theta)e^{q_{k'}(t-\theta)}d\theta+\frac{1}{\rho\left(q_{k''}-q_{k'}\right)}\int_0^t w_k(\theta)e^{q_{k''}(t-\theta)}d\theta, \qquad (43)$$

$$k=1,2,..., \quad q_{k'}\neq q_{k''}.$$

Если бы оказалось, что для некоторого $k = k_0$

$$\lambda_{k_0} = \frac{4M\rho}{\eta}, \tag{44}$$

то $q_{k_0'} = q_{k_0} = q_{k_0}$. Изложенный выше способ неприменим в этом случае. Но тогда

$$\rho D^2 + \eta\lambda_{k_0} D + M\lambda_{k_0} = \rho\left(D - q_{k_0}E\right)^2,$$

следовательно,

$$a_{k_0}(t) = \frac{1}{\rho}\left(D - q_{k_0}E\right)^{-1}\left(D - q_{k_0}E\right)^{-1} w_{k_0}(t) =$$

$$= \frac{1}{\rho}\left(D - q_{k_0}E\right)^{-1}\int_0^{\theta} w_{k_0}(\vartheta)e^{q_{k_0}(\theta-\vartheta)}d\vartheta = \frac{1}{\rho}\int_0^{t} e^{q_{k_0}(t-\theta)}d\theta\int_0^{\theta} w_{k_0}(\vartheta)e^{q_{k_0}(\theta-\vartheta)}d\vartheta =$$

$$= \frac{1}{\rho}\int_0^{t} d\theta\int_0^{\theta} w_{k_0}(\vartheta)e^{q_{k_0}(\theta-\vartheta)}d\vartheta. \tag{45}$$

Таким образом для этого случая соответствующий коэффициент $a_{k_0}(t)$ представится в виде одного двойного интеграла. Как видно, таких коэффициентов не может быть больше одного.

Итак, мы получили решение уравнения (B):

$$\Delta u = \rho\frac{\partial^2 z}{\partial t^2} - f(x,y,t), \qquad u = \eta\frac{\partial z}{\partial t} + Mz$$

при условиях $(B')$:

$$z\big|_{t=0} = \Phi(x,y); \quad \frac{\partial z}{\partial t}\Big|_{t=0} = \varphi(x,y); \quad z\big|_C = 0; \quad \frac{\partial z}{\partial t}\Big|_{t=+\infty} = 0$$

в виде: $z = z_0 + \zeta,$ (46)

где

$$z_0 = \sum_{k=1}^{\infty}\left\{C_{k'}e^{q_{k'}t} + C_{k''}e^{q_{k'}t}\right\}Z_k(x,y) \tag{47}$$

есть решение соответствующего однородного уравнения (A) с надлежащим способом выбранными коэффициентами $C_{k'}, C_{k''}$, а

$$\zeta = \frac{1}{\rho}\sum_{k=1}^{\infty}\left\{\frac{1}{q_{k'} - q_{k''}}\int_0^{t} w_k(\theta)e^{q_{k'}(t-\theta)}d\theta + \frac{1}{q_{k''} - q_{k'}}\int_0^{t} w_k(\theta)e^{q_{k'}''(t-\theta)}d\theta\right\}Z_k(x,y). \tag{48}$$

Здесь $Z_k(x, y)$, $k = 1,2,...,$ есть замкнутая система нормированных ортогональных друг к другу линейно независимых функций, каждая из которых, являясь решением краевой задачи

$$\Delta Z + \lambda Z = 0, \qquad Z\big|_C = 0,$$

соответствует одному и только одному определенному положительному числу $\lambda_k$. Каждая из функций $Z_k(x, y)$ непрерывна в $D$ и обладает непрерывными производными 1-го и 2-го порядка.

Докажем, что ряд (48) может быть почленно дифференцирован столько раз, сколько требуется для составления уравнения (B). При этом является существенным следующее предположение. Пусть $w_k$ есть максимум функции $|w_k(\theta)|, k = 1,2,...,$ когда $\theta$ меняется от 0 до $+\infty$. В таком случае легко получить, что

$$\left|\int_0^t w_k(\theta)e^{q_{k'}(t-\theta)}d\theta\right| \le \frac{w_k}{|q_{k'}|},$$

$$\left|\int_0^t w_k(\theta)e^{q_{k''}(t-\theta)}d\theta\right| \le \frac{w_k}{|q_{k''}|}.$$

Высказанное выше утверждение имеет место, если $w_k$ при возрастании $k$ представляет малую величину порядка не ниже, чем $k^{-3}$. Действительно, так как при возрастании $k$ числа $|q_{k'}|, |q_{k'}|, |q_{k'} - q_{k''}|, \lambda_k$ являются большими того же порядка, что и $k$, и так как отношения

$$\frac{Z_k(x, y)}{q_{k'} - q_{k''}}$$

при $k \to \infty$ равномерно ограничены, то суммы

$$\left|\sum_{k=n+1}^{n+p}\left\{\frac{Z_k(x, y)}{\rho(q_{k'} - q_{k''})}\int_0^t w_k(\theta)e^{q_{k'}(t-\theta)}d\theta + \frac{Z_k(x, y)}{\rho(q_{k''} - q_{k'})}\int_0^t w_k(\theta)e^{q_{k'}(t-\theta)}d\theta\right\}\right|,$$

$$\left|M\sum_{k=n+1}^{n+p}\left\{\frac{-\lambda_k Z_k(x, y)}{\rho(q_{k'} - q_{k''})}\int_0^t w_k(\theta)e^{q_{k'}(t-\theta)}d\theta + \frac{-\lambda_k Z_k(x, y)}{\rho(q_{k''} - q_{k'})}\int_0^t w_k(\theta)e^{q_{k'}(t-\theta)}d\theta\right\}\right|,$$

$$\left|\eta\sum_{k=n+1}^{n+p}\left\{\frac{-\lambda_k q_{k'} Z_k(x, y)}{\rho(q_{k'} - q_{k''})}\int_0^t w_k(\theta)e^{q_{k'}(t-\theta)}d\theta + \frac{-\lambda_k q_{k''} Z_k(x, y)}{\rho(q_{k''} - q_{k'})}\int_0^t w_k(\theta)e^{q_{k'}(t-\theta)}d\theta\right\}\right|,$$

$$\left|\sum_{k=n+1}^{n+p}\left\{w_k(t)Z_k(x, y) + \frac{q_{k'}^2 Z_k(x, y)}{q_{k'} - q_{k''}}\int_0^t w_k(\theta)e^{q_{k'}(t-\theta)}d\theta + \frac{q_{k''}^2 Z_k(x, y)}{q_{k''} - q_{k'}}\int_0^t w_k(\theta)e^{q_{k'}(t-\theta)}d\theta\right\}\right|$$

при достаточно большом $n$ и каком угодно конечном $p \ge 1$ представляют собой малые порядков, соответственно равных

$n^{-4}, n^{-3}, n^{-2}, n^{-1}$.

Это значит, что ряды

$$\sum_{k=n+1}^{n+p} \left\{ \frac{Z_k(x,y)}{\rho(q_{k'} - q_{k''})} \int_0^t w_k(\theta) e^{q_{k'}(t-\theta)} d\theta + \frac{Z_k(x,y)}{\rho(q_{k''} - q_{k'})} \int_0^t w_k(\theta) e^{q_{k''}(t-\theta)} d\theta \right\},$$

$$M \sum_{k=n+1}^{n+p} \left\{ \frac{\Delta Z_k(x,y)}{\rho(q_{k'} - q_{k''})} \int_0^t w_k(\theta) e^{q_{k'}(t-\theta)} d\theta + \frac{\Delta Z_k(x,y)}{\rho(q_{k''} - q_{k'})} \int_0^t w_k(\theta) e^{q_{k''}(t-\theta)} d\theta \right\},$$

$$\eta \sum_{k=n+1}^{n+p} \left\{ \frac{q_{k'} \Delta Z_k(x,y)}{\rho(q_{k'} - q_{k''})} \int_0^t w_k(\theta) e^{q_{k'}(t-\theta)} d\theta + \frac{q_{k''} \Delta Z_k(x,y)}{\rho(q_{k''} - q_{k'})} \int_0^t w_k(\theta) e^{q_{k''}(t-\theta)} d\theta \right\},$$

$$\rho \sum_{k=n+1}^{n+p} \left\{ \frac{w_k(t) Z_k(x,y)}{\rho} + \frac{q_{k'}^2 Z_k(x,y)}{\rho(q_{k'} - q_{k''})} \int_0^t w_k(\theta) e^{q_{k'}(t-\theta)} d\theta + \frac{q_{k''}^2 Z_k(x,y)}{\rho(q_{k''} - q_{k'})} \int_0^t w_k(\theta) e^{q_{k''}(t-\theta)} d\theta \right\}$$

равномерно и абсолютно сходятся к функциям, соответственно равным

$$\zeta(x,y,t), \ \Delta(M\zeta), \ \ \Delta(\eta \frac{\partial \zeta}{\partial t}), \ \ \rho \frac{\partial^2 \zeta}{\partial t^2},$$

*что и утверждалось.*

Итак, определяемая рядом (48) функция $\zeta(x,y,t)$ допускает почленное дифференцирование , необходимое для составления уравнения (B).

Что касается решения однородного уравнения

$$z_0(x,y,t) = \sum_{k=1}^{\infty} \left\{ C_{k'} e^{q_{k'} t} + C_{k''} e^{q_{k''} t} \right\} Z_k(x,y), \tag{49}$$

то мы видели еще раньше, что только что полученный вывод для $\zeta(x,y,t)$ применим и к $z_0(x,y,t)$.

Таким образом, путем почленного дифференцирования суммы рядов (47) и (48) могут быть составлены все производные, входящие в уравнение (A).

Наконец, заметим, что как в случае колебания упруго-вязких нитей, так и здесь имеет место явление, названное нами *элевтерозом* и состоящее в последовательном затухании обертонов, начиная с наиболее высоких. Это обстоятельство позволяет в найденном решении ограничиваться лишь конечным числом членов. В частности, можно было бы применить акустический метод к практическому определению модуля сдвига $\mu$ и коэффициента внутреннего трения $\eta$. А знание последнего могло бы послужить для определения продолжительности релаксации на основании максвелловских уравнений.

SUR LA QUESTION DES PETITES VIBRATIONS DE MEMBRANES ÉLASTO-VISQUESEWS

A.N.GUÉRASSIMOV

Le mémoire contient une application de la théorie générale donnée par l'auteur dans son travail:"Principes de la théorie de la deformation des corps élasto-visqueux" (Прикл. Мат. и мех.,т.II, dsg. 3, 1938)

*Диссертация на степень кандидата физико-математических наук, Москва, МГУ,1942*

# **7.**НЕКОТОРЫЕ ЗАДАЧИ ТЕОРИИ УПРУГОСТИ С УЧЕТОМ ПОСЛЕДЕЙСТВИЯ И РЕЛАКСАЦИИ ПО ЛИНЕЙНОМУ ЗАКОНУ

*ОГЛАВЛЕНИЕ*

## ПРЕДИСЛОВИЕ

В 1939-1940 году нами был опубликован ряд статей [30-32], посвященных колебаниям упруго-вязких тел. Эти работы предлагаются теперь нами в качестве диссертации на соискание степени кандидата физико-математических наук.

Нам казалось необходимым снабдить их общим предисловием, обзором литературы по упруго-вязким и пластическим деформациям, критикой сообщаемым в них результатов и дополнением, содержащим некоторые новые результат ы, полученные автором в 40-42 годах. Все это не могло быть сделано своевременно по ряду не зависящих от нас причин.

Сообразно с этим мы прилагаем здесь оттиски упомянутых выше работ и в настоящей статье располагаем материал по следующему плану.

Сначала мы даем обзор литературы по данному вопросу и по смежным с ним задачам математической теории упругости. Этот обзор состоит из двух частей, не разграниченных строго. Одна часть касается теории пластических деформаций, другая посвящена специально линейным законам упругого последействия и релаксаций.

Затем мы вкратце останавливаемся на содержание наших статей и на критике сообщаемых там выводов.

Наконец, последняя часть содержит изложение задачи о колебаниях упруго-вязкой плоской мембраны, релаксирующей и последействующей по линейному закону. В этой части сообщается то, что было получено автором этой диссертации уже после опубликования приложенных работ.

Москва, октябрь 1942 г.

## I. ОБЗОР ЛИТЕРАТУРЫ

Вопросом о внутреннем трении в жидкостях интересовался еще Ис.Ньютон [41].В то время как Гук дал свой знаменитый закон, легший в основу классической теории упругости, Ньютон в отделе IX книги II своих «Philosophia naturalis principia mathematica» loci cito, посвященном вопросу «О круговом движении жидкостей» писал/ «Предположение» и «Предложение LI. X Теорема XXXIX.»/: « Сопротивление, происходящее от недостатка скользкости жидкости при прочих равных условиях предполагается пропорциональным скорости, с которой частицы жидкости разъединяются друг от друга». И далее: «Если в однородной и беспредельной жидкости вращается равномерно около постоянной своей оси твердый бесконечно длинный цилиндр и жидкость приводится в движение единственно только этим натиском, причем всякая ее частица продолжает сохранять свое равномерное движение, то я утверждаю, что времена обращений частиц жидкости пропорциональны их расстояниям до оси цилиндра»./стр. 436-448 loci cito, по русск. Перев. Акад. А.Н.Крылова/. С геометрической ясностью Ньютон доказывает эту « Теорему XXXIX», опираясь на приведенное выше « Предположение». Последнее известно в настоящее время, как « Закон Ньютона о внутреннем трении», по которому сила трения между двумя жидкими слоями пропорциональна поверхности раздела трущихся слоев и градиенту скорости.

В своих рассуждениях Ньютон, однако, допустил ошибку, отразившуюся на всех выводах отдела IX второй книги его « Начал». Именно, он ошибочно рассматривал действие сил трения вместо того, чтобы рассматривать действие *моментов* этих сил относительно оси вращения цилиндра. Эта ошибка была замечена Стоксом[42]. Если, следуя Стоксу, внести исправление в рассуждения Ньютона, то последняя часть приведенной выше теоремы LI, должна быть заменена словами: «… то я утверждаю, что времена обращений частиц жидкости пропорциональны квадратам их расстояний от оси вращения». Аналогичное исправление нужно внести и в « Теорему LII» относительно шара, вращающегося в беспредельной однородной жидкости, и во все последующие предложения этого отдела.

На основе законов Гука и Ньютона начинают развиваться, с одной стороны, математическая теория упругости с теорией сопротивления материалов, а с другой – динамика подвижных средин, в частности, гидродинамика вязких и пластических тел. Эти отрасли естествознания завершаются, с одной стороны, уравнениями Lamé классической теории упругости, с другой стороны – системой уравнений Navier-Stoqes'а гидродинамики вязких жидкостей.

В связи с возрастанием интереса к внутренним свойствам вещества в промышленности и в технике появляется все большее и большее количество экспериментальных и теоретических работ, касающихся пластически и упруго-вязких свойств материалов. Эти вопросы стали особенно насущными в СССР и, пожалуй, особенно теперь, в период военного времени. Вот почему именно теперь появляется на свет такое большое количество работ по пластическим деформациям и по упругому последействию с релаксацией, выполненных нашими советскими экспериментаторами и теоретиками. Достаточно назвать такие отрасли тяжелой индустрии, как штамповка, прокатка, вальцовка металлов, целиком основанные на теории пластических деформаций. Но и в легкой промышленности нашего Союза затрагиваемая здесь тема приобретает значение исключительной важности, вследствие чего в течение последних 10 лет ряд советских институтов легкой промышленности методически ставил на разрешение вопросы механических свойств текстильных, кожевенных, резиновых и пр. материалов и изделий.

Учение о пластических деформациях и об упругом взаимодействии и релаксации до последних лет развивалось по двум направлениям.

Одно из этих направлений , обильно представленное рядом блестящих исследований, было намечено B.de-Saint-Venant’ ом [1] на основе явления, установленного Trésca[1] и заключающегося в постоянстве наибольшего напряжения при сдвиге. Предполагая, что направления наибольших скоростей сдвига совпадают с направлениями наибольших скалывающих напряжений и что пластически деформируемое тело несжимаемо, S.-Venant’ из гипотезы Trésca получает систему уравнений движения пластического тела для случая плоской задачи в виде

$$\frac{\partial\sigma_x}{\partial x}+\frac{\partial\tau_{xy}}{\partial y}+X=\rho\left\{\frac{\partial v_x}{\partial t}+v_x\frac{\partial v_x}{\partial x}+v_y\frac{\partial v_x}{\partial y}\right\} ,$$

$$\frac{\partial\tau_{xy}}{\partial x}+\frac{\partial\sigma_y}{\partial y}+Y=\rho\left\{\frac{\partial v_y}{\partial t}+v_x\frac{\partial v_y}{\partial x}+v_y\frac{\partial v_y}{\partial y}\right\} ,$$

$$(\sigma_x-\sigma_y)^2+4\tau_{xy}^2=const ,$$

$$\frac{\partial v_x}{\partial x}+\frac{\partial v_y}{\partial y}=0 ,$$

$$\frac{\partial v_y}{\partial y}-\frac{\partial v_x}{\partial x}=\frac{\sigma_y-\sigma_x}{2\tau_{xy}}\left\{\frac{\partial v_x}{\partial y}-\frac{\partial v_y}{\partial x}\right\} .$$

Здесь X и Y – компоненты внешней силы на единицу объема, ρ – плотность, σ и τ с индексами–компоненты тензора напряжения, $v_x$ и $v_y$ - компоненты скорости. Третье уравнение выражает результаты опытов Trésca, четвертое – уравнение несжимаемости, пятое– коллинеарности максимальных скоростей скольжения и скалывающих напряжени й. M.Lévy [2-3] распространил уравнения пластических деформаций S.-Venant’a на случай трехмерной задачи; вместо условия коллинеарности максимальных скоростей и напряжений он ввел естественное обобщение этого предположения: тензор скоростей деформаций должен быть линейно связан с тензором напряжений. Условие же постоянства скалывающих напряжений Lévy ошибочно заменил уравнением:

$$(8q+16K^2)^3-3q^2(4K^2+q)+27R^2=0 , K=const,$$

где положено:

$$\sigma_x+\sigma_y+\sigma_z=N,$$

$$\Delta_x = \sigma_x - \tfrac{1}{3}N \;,\; \Delta_y = \sigma_y - \tfrac{1}{3}N \;\;,\; \Delta_z = \sigma_z - \tfrac{1}{3}N$$

$$q = \Delta_x\Delta_y + \Delta_y\Delta_z + \Delta_z\Delta_x - \left(\tau_{xy}^2 + \tau_{yz}^2 + \tau_{zx}^2\right) ,$$

$$R = \Delta_x\tau_{yz}^2 + \Delta_y\tau_{zx}^2 + \Delta_z\tau_{zx}^2 - 2\tau_{yz}\tau_{zx}\tau_{xy} .$$

Однако S.-Venant [4] заметил, что из уравнения Lévy, приведенного сейчас, не получается как частный случай его аналог для двухмерной задачи.

Далее, к тому же классу работ следует отнести сообщения H.Hencky [5-6], который представляет тензор пластических напряжений в виде суммы тензора упругих напряжений и «тензора трения». Компоненты первого из них, согласно теории R.von-Mises'a[7-8] Hencky пишет в виде

$$\sigma_x = K\frac{\partial v_x}{\partial x} \;,\; \sigma_y = K\frac{\partial v_y}{\partial x} ,\; \sigma_z = K\frac{\partial v_z}{\partial x} \;,$$

$$\tau_{yz} = \frac{K}{2}\left(\frac{\partial v_y}{\partial z} - \frac{\partial v_z}{\partial y}\right) ,\quad \tau_{zx} = \frac{K}{2}\left(\frac{\partial v_z}{\partial x} - \frac{\partial v_x}{\partial z}\right) ,\; \tau_{xy} = \frac{K}{2}\left(\frac{\partial v_x}{\partial y} - \frac{\partial v_y}{\partial x}\right) .$$

Если через λ с соответствующими индексами обозначить слагающие тензора скоростей деформаций, то компоненты тензора трения по Hencky напишутся так:

$$\sigma'_x = 2\kappa\lambda_x \;,\; \sigma'_y = 2\kappa\lambda_y ,\; \sigma'_z = 2\kappa\lambda_{zx} \;;$$

$$\tau_{yz} = \kappa\lambda_{yz} \;,\; \tau_{zx} = \kappa\lambda_{zx} ,\; \tau_{xy} = \kappa\lambda_{xy} \;,$$

где $\kappa$ – постоянная, зависящая от рода и свойств вещества. В том же роде строятся уравнения пластических деформаций в теориях A.Haar'a, von.Karman'a [9] и L.Prandtl'a [10-11]. В нашу задачу не входит обзор всех работ-это очень хорошо сделано в маленькой книжке С.Г.Михлина [12]. Приведенного выше достаточно для характеристики общего направления исследований по пластическим деформациям; сюда же следует отнести еще интересную теоретическую работу C.Caratheodory , E.Schmidt'a[13], посвященную семействам кривых, введенных Hencky и Prandtl'ем, затем ряд исследований, выполненных берлинским профессором-женщиной H.Pollaczek-Geiringer [14], наконец, фундаментальную работу акад. С.Соболева [15], краткое переложение которой можно найти в указанной крижке Михлина, а также диссертации А.А.Ильюшина [16], и П.М.Огибалова [17]. Не претендуя на полноту списка, упомянем экспериментальные и теоретические исследования Н.С. Курнакова Жемчужного[18], Курнакова и Рапке [19],Г.Дж. Тэкелла[20], прекрасную книжку A.Nadai [21] и работы лаборатории С.И.Губкина [22].

Другой класс исследований, представленный не меньшим количеством работ, посвящен упругому последействию и тому явлению , которое Cl.Maxwell назвал *релаксацией*. Это явление, как известно, заключается в изменении напряжения со временем при постоянной деформации. Основания релаксационной теории даны Cl.Maxwell'ом в 1868 году[23] . Я позволю себе привести здесь выдержку из подлинного текста знаменитого сочинения, в которой и заключена сущность всей теории:

«…A distortion or strain of some kind, which we may call S, is produced in the body by displacement. A state of stress or elastic force which we may call F is thus excited. The relation between the stress and the strain be written F=ES, where E is the coefficient of elasticity, for that particular kind of strain. In a solid body free from viscosity, F will remain =ES and

$$\frac{dF}{dt} = E\frac{dS}{dt}.$$

If, however, the body is viscous, E will not remain constant, but will tend to disappear at a rate depending of the value of F and on the nature of the body.
If we suppose this rate proportional to F, the equation may be written

$$\frac{dF}{dt} = E\frac{dS}{dt} - \frac{F}{T} \quad ,$$

which will indicate the actual phenomena in an empirical manner. For if S be constant,

$$F - ESe^{-\frac{t}{T}} \quad ,$$

showing that F gradually disappears, go that if the body is left to itself it gradually loses any internal stress, and the pressures are finally distributed as in a fluid at rest. If $\frac{dS}{dt}$ is constant, that is, if free is a steady motion of the body which continually increases the displacement,

$$F = ET\frac{dS}{dt} + Ce^{-\frac{t}{T}} \quad ,$$

showing that F tend to a constant value depending on the rate of displacement…»

Итак, по Maxwell‘у, тело с течением времени теряет способность сопротивляться приложенным извне механическим воздействиям. Иными словами, напряжение есть функция времени, убывающая по показательному закону одного из двух приведенных в тексте видов. Это положение и лежит в основе релаксационной теориии. Постоянный промежуток времени T носит название периода релаксации, обратная ему величина

$$r = \frac{1}{T}$$

Называется обычно скоростью релаксации. Заметим, что Maxwell дает без вывода основное уравнение

$$\dot{F} = E\dot{S} - rF \quad ,$$

просто ссылаясь на экспериментальные данные. Мы увидим в дальнейшем, при изложении работы R.Becker’a [29], каким образом это уравнение может быть обосновано выражение теоретически. Здесь же обратим внимание на то, что выражение «if, however, the body is viscous» допускает значительную широту трактовки понятия «вязкого» тела.

В 1876 году появилась замечательная работа L.Boltzmann’a [24], в которой явление упругого последействия феноменологически истолковывается, как результат своего рода «памяти» /Erinnerungsvermögen/ тела по отношению к его предшествующим «переживаниям» /Erlebnisse/. Вкратце, теория Boltzmann'а сводится к следующему.

Представим себе упругое тело, на которое действует сила P(t) . Назовем y(t) изменение той координаты, дифференциал dy которой, будучи умножен на P(t), дает величину производимой на пути dy работы силы P(t). При этом координата и сила понимаются в обобщенном смысле слова. Согласно закону Гука, в статическом случае упругой деформации мы должны иметь:

$$y(t) = \frac{F}{E} P(t),$$

где F=const зависит от формы тела , E=const есть модуль упругости. В последействующем же теле y(t) зависит не только от мгновенного значения силы P(t), но и от совокупности всех предшествующих значений в их хронологическом порядке. Пусть $P(\vartheta)$ есть значение силы в момент $\vartheta \le t$. Ее влияние на значение y в момент t должно быть, по Boltzmann'у, прямо пропорционально импульсу $P(\vartheta)d\vartheta$ и быть тем меньше, чем больше промежуток времени между обоими рассматриваемыми моментами $\vartheta$ и t. Если обозначить через $\beta$ некоторую « постоянную последействия» и $\varphi(t-\vartheta)$– убывающую «функцию последействия» (Nachwirkungs funktion), то изменение y(t)вносимое к моменту t силой $P(\vartheta)$, может быть представлено в виде

$$\frac{F}{E}\beta P(\vartheta)\varphi(t-\vartheta)d\vartheta;$$

при этом теория размерности требует, чтобы

$$[\beta\varphi(t-\vartheta)] = \text{сек}^{-1},$$

чего можно достигнуть, положив β отвлеченным числом и $[\varphi]$=сек$^{-1}$.

Результирующее изменение рассматриваемой координаты получается суммированием всех предшествующих моменту t значений для предыдущего выражения, то – есть

$$y(t) = \frac{F}{E}\left\{P(t) + \beta\int_{-\infty}^{t} P(\vartheta)\varphi(t-\vartheta)d\vartheta\right\},$$

или, положив

$$t - \vartheta = \varpi,$$

$$y(t) = \frac{F}{E}\left\{P(t) + \beta\int_{0}^{\infty} P(t-\varpi)\varphi(\varpi)d\varpi\right\}.$$

L.Boltzmann дает это уравнение в виде, решенном относительно P(t), что можно легко сделать, ввиду малости постоянной последействия β . тогда первое из ряда последовательных приближений дает:

$$P(t) = \frac{E}{F}\left\{ y(t) - \beta\int_0^{\infty} y(t-\varpi)\varphi(\varpi)d\varpi \right\} \qquad \text{(B)}$$

Wiechert'ом[25], было показано, что это последнее уравнение есть необходимое следствие релаксационной теории Maxwell'a. Vito Volterra [26] вместо приведенного только что уравнения Boltzmann'a получает в своей « Теории наследственности» более общее выражение правой части в виде бесконечного ряда интегралов возрастающей кратности.

Таким образом, Boltzmann и Volterra описывают явление последействия посредством линейного интегрального уравнения приведенного выше вида.Tompson[36] показал, что упругое последействие, по крайней мере для простейших деформаций вроде простого растяжения, может быть описано при помощи линейной дифференциальной зависимости вида

$$\sigma = E\varepsilon + \mu\dot{\varepsilon},$$

где σ – напряжение, ε – деформация, Е и μ – физические постоянные.

Ядро $\varphi(\varpi)$ интегрального уравнения Boltzmann'a обычно определяется из опытных данных. В довольно широких пределах можно считать (Becker,[29]). $\varphi(\varpi) = \frac{1}{\varpi}$

Сам Becker в своей работе, как увидим, дает для $\varphi(\varpi)$ иное выражение, хорошо согласующееся с опытом. Работы Jordan'a [33] в этом направлении показали, что при $\varphi(\varpi)$, данном Becker' ом, уравнение (В) достаточно хорошо согласуется с действительностью, а в 1912 году Jordan'y удалось даже подметить [о 34], что (В) с успехом может быть применено к вопросам магнитного возбуждения железа и поляризации диэлектриков. Опыты с крутильными (Boltzmann) и поперечными (Bennewitz) колебаниями также подтверждают верность уравнения (И) в указанном Becker' ом предположении относительно вида функции $\varphi(\varpi)$. Частные случаи применения теории Boltzmann'a разобраны в цитируемой статье R.Becker'a [29] ( Belastungs,-Entlastungs,-Überlagerungsversuch, periodische Beanspruchung, u.a) с указанием на выводы из их экспериментальной проверки. Результаты такой проверки, как говорит Becker, как будто подтверждают верность (В) при Becker'овском виде ядра $\varphi(\varpi)$. С другой стороны, лаборатория испытания материалов при Московском Государственном Университете в институте Механики дает более сложное, но зато лучше согласующееся с опытом выражение для ядра. Так, Бронский [40] приводит следующий вид, установленный из многочисленных опытов, осуществленных в лабораториях МГУ:

$$\varphi(\varpi) = \bar{\sigma}\alpha e^{-\varpi^{\alpha}}\varpi^{\alpha-1} \qquad ,$$

где $\bar{\sigma}$ и α –положительные постоянные, α<1. Теории Wiechert'a[25] и Bennewitz'a [27], поскольку они феноменолочески верно отображает процесс последействия, приводят к тем же результатам, что и теория Boltzmann'a. Я позволю себе не останавливаться на них.

В 1913 году von-Karman'м была сделана первая попытка построить теорию упругого последействия, основанную на наличии определенных неоднородностей в теле.

Эта мысль о пластической неоднородности, как о причине упругого последействия, еще более ярко выражена в работ v.-Wartenberg'a[28];благодаря последней теории упругого последействия делает с 1913 года значительный шаг вперед в направлении от чистого описания количественной стороны к физическому существу явления.

Теория von-Wartenberg'a имеет ввиду поликристаллические системы( металлы) и исходит из свойства монокристаллов не показывать никакой склонности к последействию. Но в поликристаллических телах явно обнаруживается способность последействовать и Wartenberg причину этого видит в беспорядочности расположения кристаллитов с различной сопротивляемостью деформациям. Он полагает тело состоящим из множества легко деформируемых, пластических, частиц, каждая из которых в течение опыта остается в окружении ведущих себя как вполне упругие. Во время нагружения легко деформируемые частицы должны разнапрягаться (sich entsprannen) . При разгружении, вследствие упругого окружения, потерявших частично или целиком свое напряжение частица деформируется в обратном направлении. Однако, если первоначальное нагружение зашло так далеко, что разнапряглись значительные группы лежащих друг около друга частиц, то после снятия нагрузки действием окружающих упругих частиц эти группы вообще не вернутся к прежнему состоянию и будет иметь место пластическая деформация. В кристаллических телах пластичность может быть объяснена сдвигами вдоль плоскостей спайности или скольжением кристаллитов друг по другу своими гранями.

(?) склоняется к предположению, что как и в вязких жидкостях, при этом развиваются силы внутреннего трения, пропорциональные первой степени градиента скорости. Очевидно, это предположение, законное для аморфных тел, для поликристаллических тел вряд ли соответствует действительности. Но оно, как показано Wiechert'ом и Becker'ом /loc. cit./, укладывается в рамки уравнения Boltzmann'a при некоторой модификации ядра. Эта модификация осуществлена в работе Becker'a [29], так много раз уже названной. Вот рассуждения Becker'a по этому поводу.

Пусть тело состоит из N частей, каждая из которых пластически однородна и которые располагаются внутри тела без какого-либо порядка. Не нарушая общности рассуждений, эти части можно считать равными по объему. При упругом характере деформации, под действием нагрузки P(t) во всех частях тела создалось бы нормальное напряжение

$$\sigma_0 = f \cdot P \;, \quad f = const \;.$$

Рассмотрим одну из частиц, которая, благодаря своей способности разнапрягаться, находится под напряжением $\sigma < \sigma_0$. Разность $\sigma_0 - \sigma$ есть мера уменьшения участия данной частицы в общем несении нагрузки P, но обусловленное этим увеличение напряжения в других местах тела может не приниматься в расчет, если только N достаточно велико.

Кроме чисто упругой деформации $y = \frac{F}{E}P$, тело получит, вследствие сказанного, дополнительную деформацию

$$\Delta y = \frac{F}{E}\lambda(\sigma_0 - \sigma) \quad ,$$

где λ – некоторая физическая постоянная. Тогда полное мгновенное значение деформации в момент времени t должно быть принято равным

$$y(t) = \frac{F}{E}\left\{P(t) + \lambda\sum_1^N(\sigma_0 - \sigma)\right\} \ ,$$

где суммирование распространяется на все N частей тела. По Maxwell'ю (см. выше формулу $F = ET\dot{S} + Ce^{-\frac{t}{T}}$), при постоянной скорости деформации $\dot{S}$ напряжение F с возрастанием t стремится к постоянному значению, равному

$$\tilde{F} = ET\dot{S} = ETv\,,$$

где v-скорость деформации. Ограничиваясь одномерным случаем, что однако не суживает общности вывода, если смещение u есть функция координаты x и времени t, можно написать $S = \frac{\partial u}{\partial x}$ и следовательно

$$v = \dot{S} = \frac{\partial}{\partial t}\left(\frac{\partial u}{\partial x}\right) = \frac{\partial}{\partial x}\left(\frac{\partial u}{\partial t}\right) = \frac{\partial \dot{u}}{\partial x};$$

поэтому предыдущее представляется в виде

$$\tilde{F} = \eta\frac{\partial \dot{u}}{\partial x}.$$

Здесь $\eta = ET$ есть коэффициент внутреннего трения (Maxwell). Отсюда видно, что «вязкая» часть $\tilde{F}$ напряжения пропорциональна градиенту скорости $\frac{\partial \dot{u}}{\partial x}$, как это и должно быть по закону Ньютона. Предположение Wartenberg'а относительно сил внутреннего трения адекватно предположению верности релаксационной теории для поликристаллов. С другой стороны, приведенная выше формула этой теории по Becker'у может быть оправдана допущением, что пластическая разрядка напряженного состояния происходит так, что за время dt напряжение уменьшается на $d\sigma_1$,пропорциональное σ, то-есть

$$d\sigma_1 = -r\sigma dt,$$

где r– физическая постоянная, обратная по отношению к периоду релаксации. Если P(t)за время dt меняется на dP, то кроме пластического изменения напряжения войдет еще чисто упругое:

$$d\sigma_2 = fdP,$$

так что полное изменение $d\sigma = d\sigma_1 + d\sigma_2$ должно удовлетворять соотношению

$$\dot{\sigma} = -r\sigma + f\dot{P},$$

а это есть, с точностью до обозначений, не что иное, как цитированное выше уравнение Maxwell'а ( см. английский текст, P пропорционально S, $\sigma = F,\ r = \dfrac{1}{T}$):

$$\dot{F} = E\dot{S} - \frac{F}{T}.$$

Используя тождество

$$e^{-rt}\frac{d}{dt}(e^{rt}\sigma) = \dot{\sigma} + r\sigma,$$

Becker из предыдущего получает

$$\sigma(t) = f\left\{P(t) - \int_{-\infty}^{t} rP(\vartheta)e^{-r(t-\vartheta)}d\vartheta\right\} \quad .$$

Отсюда, при помощи подстановки $t - \vartheta = \omega$ и соотношения $fP(t) = \sigma_0(t)$, для величины одного из слагаемых в предыдущем общем выражении для деформации y(t) находим

$$\sigma_0(t) - \sigma(t) = f\int_{-\infty}^{t} rP(\vartheta)e^{-r(t-\vartheta)}d\vartheta = f\int_{0}^{\infty} rP(t-\omega)e^{-r\omega}d\omega \ .$$

Если через $F(r)dr$ обозначить число частей из N, обладающих пластичностью, лежащей между $r$ и $r + dr$, то для $\sum_{1}^{N}(\sigma_0 - \sigma)$ получим выражение

$$\sum_{1}^{N}(\sigma_0 - \sigma) = f\int_{0}^{\infty} P(t-\omega)d\omega)\int_{0}^{\infty} F(r)re^{-r\omega}dr \ ,$$

причем $\int_{0}^{\infty} F(r)dr = N$.

Положив еще

$$f\int_{0}^{\infty} F(r)re^{-r\omega}dr = \frac{\beta}{\lambda}\varphi(\omega),$$

придем к уравнению Boltzmann'а (B) .

При этом существенным обстоятельством является предположение, что тело пластически неоднородно, то-есть, что $r$ действительно принимает неотрицательные значения, как угодно отличающиеся друг от друга.

Замечая, что части тела с очень малыми значениями $r$ не оказывают влияния на ход деформации и принимая, что существует такое малое $r_* > 0$ и такое большое $R_* > 0$, для которых

$$F(r) = 0 \text{ при } r < r_* \text{ и } r > R_*,$$

$$F(r) = \frac{1}{r} \text{ при } r_* \le r \le R_* ,$$

Becker получает:

$$\varphi(\omega) = \frac{1 - e^{-R_*\omega}}{\omega} ,$$

$$\beta = \lambda f N \frac{\bar{\varepsilon}}{R_*},$$

где $\bar{\varepsilon}$ - среднее значение пластичности. При этом поставлено условие, что во время наблюдения наиболее твердые частицы пластически не деформируются, так что для всех $\omega$ величина $r_*\omega$ настолько мала, что можно положить $e^{-r_*\omega} = 1$.

В качестве следствия из этой теории, для наблюдаемого модуля упругости $E_*$ получается значение:

$$E_* = \frac{E}{1 + \beta[0{,}577 + \ln R_*] + \beta \ln t},$$

откуда видно, что $E_*$ убывает с течением времени, оставаясь меньше должного модуля $E$.

Было сказано уже, что опыты Boltzmann‘a, Jordan’a и Bennewitz’a подтверждают этот результат. Однако, работами Sieg’a (Phys. Review, 35, 1912 и D.Schenk’a (ZS.f.Phys., 1931) установлено, что затухание колебаний и период для металлических стержней несомненно зависят от амплитуды-факт, стоящий в явном противоречии с найденным выражением для кажущегося модуля $E_*$ упругости. С другой стороны, опыты с поликристаллическими телами показывают, что скорость течения растет не пропорционально силе, а много быстрее. По этому поводу были предложены зависимости, отличающиеся от закона Ньютона для внутреннего трения (Schönborn, ZS.f.Physik, 3, 1922; Geiss, ZS.f.Phys.,29,1924) , но эти зависимости трудно уложить в рамки какой-либо теории. Поэтому Becker делает попытку модифицировать вид ядра $\varphi(\omega)$, исходя из теории спонтанных разрядок местных напряжений.

Именно, пусть свободное от напряжений тело в момент $t = 0$ подвергается действию силы $P$, из-за чего создается во всем объеме напряжение $\sigma_0 = fP$. Рассмотрим малую однородную область, пластичность которой характеризуется постоянной к таким образом, что когда напряжение в этой области равно $\sigma$, вероятность $dw$ спонтанной разрядки напряжения в течение промежутка времени $dt$ в этой области равна

$$dw = S(\sigma)\kappa dt.$$

Относительно функции $S(\sigma)$ Becker предполагает, что она не отрицательна, растет с ростом $|\sigma|$, от которого она только и зависит, и при $\sigma = 0$ обращается в нуль. Далее предполагается, что опыт ограничивается только чистым последействием, когда $\sigma$ меняется скачком только в изолированных друг от друга малых частях тела, окруженных такими, где поддерживается «должное» напряжение (Sollspannung) $fP$.

Тогда, вводя функцию распределения $F(\kappa)$, определяемую аналогично выше рассмотренной $F(r)$, путем рассуждений, совершенно подобных изложенным по поводу $F(r)$, Becker находит для ядра $\varphi$ в уравнении Boltzmann'а

$$\varphi(\omega) = \frac{1 - e^{-S(\sigma_0)\omega}}{\omega},$$

где вместо постоянной $R_*$ в показателе правой части входит функция $S(\sigma_0)$, обеспечивающая должную зависимость периода колебаний, декремента затухания и кажущегося модуля упругости от амплитуды. Ограничиваясь, по недостатку места, этим, далеко не полным, обзором литературы, мы переходим теперь к разбору предлагаемых в качестве диссертации статей

## II. ОБЗОР РАБОТ, ПРЕДЛАГАЕМЫХ В КАЧЕСТВЕ ДИССЕРТАЦИИ

К настоящей статье приложены оттиски трех работ А.Н.Герасимова, посвященных упруго-вязким деформациям в твердых телах.

Первая из этих работ[30] носит название «Проблема упругого последействия и внутреннее трение» (Le problème du frottement intér. et de l'action postérieeux) ,была выполнена в 1937 году по распоряжению Научно-исследовательского Сектора Московского Текстильного института, где работал тогда автор в качестве исполняющего обязанности доцента при кафедре физики. Результаты этой работы и были опубликованы потом в сборнике «Прикладной математики и механики».

В цитируемой статье рассматривается задача о поперечных колебаниях закрепленной на концах нити под действием лежащих в одной плоскости сил, перпендикулярных к длине, так что частицы нити не испытывают смещений вдоль нити. Вещество последней предполагается одновременно обладающим и упругими/закон Гука/, и вязкими /закон Ньютона/ свойствами. Выводится уравнение *) движения нити в форме /стр.500/:

$$\eta \frac{\partial^3 u}{\partial x^2 \partial t} + M \frac{\partial^2 u}{\partial x^2} - \rho \frac{\partial^2 u}{\partial t^2} - H \frac{\partial u}{\partial t} + f(x,t) = 0 , \qquad \text{(B)}$$

сопровождаемое начальными условиями:

*) *Это уравнение (B) содержит в себе, как частный случай, то, которое получил П.М.Огибалов [38] для сдвига*

$$u\big|_{t=0}\Phi(x), \quad \frac{\partial u}{\partial t}\Big|_{t=0} = \varphi(x) \ ,$$

и условиями на концах:

$$u\big|_{x=0} = 0 \ , \qquad u\big|_{x=l} = 0 \ .$$

Дается доказательство единственности решения этой при двух существенных предположениях. Во-первых, допускается, что решение может быть продолжено за конец нити $x = 0$ и $x = 1$, и, во-вторых, что скорость частиц нити стремится к нулю при $t \to \infty$. Последнее допущение физически бесспорно, но, как показал А.Ю.Ишлинский, совершенно не необходимо для доказательства единственности. После этого разбирается решение задачи о колебаниях упруго-вязкой нити по способу разделения переменных. Решение получается, конечно, в виде бесконечного ряда:

$$u(x,t) = \sum_{m=1}^{\infty}\left\{\left[C'_m - \frac{1}{\rho(q''_m - q'_m)}\int_0^t e^{-q'_m\theta} w_m(\theta)d\theta\right]e^{q'_m t} + \left[C''_m - \frac{1}{\rho(q'_m - q''_m)}\int_0^t e^{-q''_m\theta} w_m(\theta)d\theta\right]e^{q''_m t}\right\}\sin\frac{m\pi x}{l} +$$

$$+\frac{1}{\rho}\sin\frac{m_0\pi x}{l}e^{q_{m_0}t}\int_0^t e^{q_{m_0}\theta}\int_0^\theta e^{-q_{m_0}(\theta-\vartheta)}w_0(\vartheta)d\vartheta \qquad /53/$$

/см.стр.521/. Здесь $q'_m$ и $q''_m$ – корни уравнения

$\rho q^2 + (\kappa\eta + H)q + \kappa M = 0$ , $\sqrt{\kappa} = \frac{m\pi}{l}$, m – целое, $m_0$ – таково, что $q'_{m_0} = q''_{m_0} = q_{m_0}$, если такое $q_{m_0}$ существует. Функции $w_m(t)$– коэффициенты при $\sin\frac{m\pi x}{l}$ в разложении вынуждающей силы: $f(x,t) = \sum_{m=1}^{\infty} w_m(t)\sin\frac{m\pi x}{l}$ в ряд Fourier.

Уравнение (В) отличается от хорошо известного уравнения для хорошо известного уравнения для затухающих вынужденных колебаний упругой струны дополнительным членом

$$\eta\frac{\partial^3 u}{\partial x^2 \partial t}$$

С производной третьего порядка. Его введение обусловлено наличием упругого последействия по линейному Ньютонову закону о пропорциональности между силой внутреннего трения и градиентом скорости. Как физический результат, получается явление, названное *элевтерозом* и состоящее в «освобождении» основного тона от обертонов, начиная с наиболее высоких. Это явление было подтверждено в 1940 году А.Ю. Ишлинским в задаче о продольных колебаниях стержней, релаксирующих и последействующих по линейному закону. Далее, следует оба начальных условия и что доказательство того, что при надлежащем выборе постоянных $C'_m$ и $C''_m$ могут быть удовлетворены оба начальных условия и что полученное решение задачи не только формально, но и фактически удовлетворяет основному дифференциальному уравнению.

Во второй части рассматриваемой работы решается задача о крутильных колебаниях упруго-вязкой нити, закрепленной на одном конце ее, несущей некоторый дополнительный момент инерции на другом.

Наконец, в третьей части автор показывает, что, пользуясь найденным решением о колебаниях нити, можно из уравнения Boltzmann'а /стр.533/:

$$u(x;t) = \frac{1}{M}\left[f(x,t) - \rho\ddot{u}(x;t)\right] + \int_0^t K(x;t-\vartheta)\left[f(x,\vartheta) - \rho\ddot{u}(x;\vartheta)\right]d\vartheta$$

Определить, для ряда дискретных моментов времени, значения ядра $K(x;n\tau)$, $n = 0,1,2,...$Тогда, пользуясь приемом, указанным Prony и описанным у Витекера и Робинсона в «Матем. обработке результатов наблюдений», можно построить аналитическое выражение для ядра в виде линейного агрегата из показательных функций.

Вторая из приложенных к диссертации работ [31]называется « К вопросу о малых колебаниях упруго-вязких мембран». Она была выполнена непосредственно после первой и напечатана в том же 39 году. Она посвящена распространению результатов первой на случай двухмерной задачи о мембране.

Дается общее уравнение в виде

$$\frac{\partial^2}{\partial x^2}\left[\eta\frac{\partial u}{\partial t} + Mu\right] + \frac{\partial^2}{\partial y^2}\left[\eta\frac{\partial u}{\partial t} + Mu\right] = \rho\frac{\partial^2 u}{\partial t^2}\ , \qquad /A/$$

Причем мембрана предполагается закрепленной на границе /C/, для которой существует функция Green'а. Начальные условия задаются в виде

$$u(x,y;t)\big|_{t=0} = \Phi(x,y),\quad \frac{\partial u}{\partial t}\Big|_{t=0} = \varphi(x,y)\ .$$

Как и в первой работе, дается доказательство единственности решения этой краевой задачи, в предположении возможности его продолжения за контур С и постепенного исчезновения скорости движения с течением времени. В третьей части нашей диссертации мы даем доказательство единственности для случая не только последействующей, но и релаксирующей мембраны, уже не пользуясь этим предположением.

Затем следует ряд положений, относящихся к свойствам собственных функций и собственных частот и, в общем, хорошо известных. Работа заканчивается доказательством того, что полученное в виде бесконечного ряда решение действительно удовлетворяет начальным и краевым условиям и основному уравнению (G).

Наконец, в третьей работе [32] «Основания теории деформаций упруго-вязких тел» автор распространяет выводы двух предыдущих статей на общий случай трехмерной задачи. Линейно сочетая законы Гука и Ньютона, тензор напряжений автор представляет в виде суммы тензоров: чисто упругого напряжения и напряжения вязкого. В результате получается векторное уравнение уравнения движения сплошной среды, с одной стороны содержащее в себе, как частный случай, Navier-Stokes'а гидродинамики вязких жидкостей, а с другой стороны – уравнения Lamé теории упругости.

Из системы уравнений, полученных автором в этой работе, те уравнения, которые рассматривались в первых двух статьях автора, получаются, как частные случаи.

## III. ВЫВОДЫ. ПРИМЕНЕНИЕ ТЕОРИИ К ВОПРОСУ О КОЛЕБАНИЯХ УПРУГО-ВЯЗКИХ МЕМБРАН, ПОСЛЕДЕЙСТВУЮЩИХ И РЕЛАКСИРУЮЩИХ ПО ЛИНЕЙНОМУ ЗАКОНУ.

Из того, что было сказано по поводу теорий Maxwell'а и Boltzmann'а в первой части этой статьи по линейному закону могут быть описаны при помощи , можно сделать следующие выводы.

1. Релаксационные явления по линейному закону могут быть описаны при помощи дифференциального соотношения

$$\dot{\varepsilon} = \nu\sigma + \frac{1}{E}\dot{\sigma}, \qquad (M)$$

где $\nu$ и $E$ – постоянные, зависящие от рода и свойств вещества, $\sigma$ - напряжение, точка указывает производные по времени, $\varepsilon$ – деформация.

2. Явления последействия могут быть представлены или интегральным линейным уравнением Boltzmann'а:

$$\varepsilon(t) = \frac{1}{E}\sigma(t) + \int_{-\infty}^{t} B(t-\tau)\varepsilon(\tau)d\tau, \qquad /B/$$

где $\varepsilon$ – деформация, $B(t-\tau)$ – «наследственная» (V.Volterra) функция, или при помощи соотношения Tompson'а :

$$\sigma = E\varepsilon + \mu\dot{\varepsilon}, \qquad /T/$$

где $\mu$ – физическая постоянная.

3. Оба явления, и релаксация, и последействие, могут быть одновременно описаны при помощи одного линейного соотношения вида

$$\dot{\sigma} + r\sigma = b\dot{\varepsilon} + nb\dot{\varepsilon}, \qquad /I/$$

где $\sigma$ и $\varepsilon$ – напряжение и соответствующая деформация, $r$ – скорость релаксации, $n$ – скорость последействия, $b$ – модуль упругости.

4. Следуя А.Ю. Ишлинскому, можно было бы рассмотреть еще более общий случай, кроме линейных законов релаксации и последействия включающий

явление упрочнения /линейного/ и не упускающий из поля зрения чисто упругого поведения тела по закону Гука до предела упругости. Это соотношение таково:

$$\dot{\sigma}+r\sigma=b\dot{\varepsilon}+nb\dot{\varepsilon}\pm(r-n)k \quad \text{при } |r\sigma-nb\varepsilon|>(r-n)\sigma_s, \qquad \dot{\sigma}=b\dot{\varepsilon}$$

при $|r\sigma-nb\varepsilon|<(r-n)\sigma_s$ ,

разности плюс или минус определяется знаком $r\sigma-nb\varepsilon$ . Здесь $k$ есть пластическая постоянная /П.М.Огибалов,l.c.,38/, $\sigma_s$ – некоторое постоянное значение напряжения, а разность $r-n$ скоростей релаксации и последействия считается положительной, как это имеет место в действительности.

Мы будем исходить из дифференциального соотношения /I/.

Так как

$$\dot{\sigma}+r\sigma=e^{-rt}\frac{d}{dt}(e^{rt}\sigma),$$

$$\dot{\varepsilon}+n\varepsilon=e^{-nt}\frac{d}{dt}(e^{nt}\varepsilon),$$

то /I/ можно представить в виде

$$e^{-(r-n)t}\frac{d}{dt}(e^{rt}\sigma)=b\frac{d}{dt}(e^{nt}\varepsilon)$$

Если ввести новые переменные $\sum$ и $E$ посредством соотношений:

$$e^{rt}\sigma=\sum \ , \ e^{nt}\varepsilon=E, \qquad /2/$$

то предыдущее принимает вид:

$$e^{-(r-nt)}d\sum=bdE,$$

откуда интегрированием получаем:

$$E(t)-E(-\infty)=\frac{1}{b}\int_{-\infty}^{t}e^{-(r-n)\tau}d\sum(\tau).$$

Здесь интеграл правой части следует понимать в смысле Stieltjes’a. Можно считать, далее, что для естественного состояния тела при имеем $E(-\infty)=0$ и тогда

$$E(t)=\int_{-\infty}^{t}\frac{1}{b}e^{-(r-n)\tau}d\sum(\tau) \ . \qquad /3/$$

Мы назовем $\sum$ и $E$ соответственно «приведенными» напряжением и деформацией. Соотношение /3/ позволяет определить приведенную деформацию $E$ в любой момент времени, если известно в функции времени приведенное напряжение.

Наоборот, из $d\sum = be^{(r-n)t}dE$ следует:

$$\sum(t) = \int_{-\infty}^{t} be^{(r-n)\tau} dE(\tau) \ , \qquad /4/$$

где интеграл опять понимается в Stieltjes'овском смысле. При помощи /4/ находится приведенное напряжение, когда известна в функции времени приведенная деформация. Постоянство во времени одной из двух приведенных величин $\sum$ или $E$ влечет за собой, на основании, /3/ или /4/, постоянство другой.

На формулы /3/ и /4/ можно смотреть, как на своего рода обобщение закона Гука. Они устанавливают взаимно-однозначное соответствие между процессами изменения деформации и напряжения в последующем и релаксирующем теле.

Дифференцируя /4/ по $x$, считая $b, n$ и $r$ постоянными, получаем

$$\frac{\partial \sum(t)}{\partial x} = \int_{-\infty}^{t} be^{(r-n)\tau} d\left(\frac{\partial E(\tau)}{\partial x}\right),$$

или, на основании /2/:

$$e^{rt} \frac{\partial \sigma}{\partial x} = \int_{-\infty}^{t} be^{(r-n)\tau} d\left(e^{n\tau} \frac{\partial \varepsilon}{\partial x}\right).$$

В случае одномерной задачи о колебаниях стержня с учетом релаксации и последействия:

$$\frac{\partial \sigma}{\partial x} = \rho \frac{\partial^2 u}{\partial t^2}, \quad \varepsilon = \frac{\partial u}{\partial x}$$

и предыдущее принимает вид:

$$\rho e^{rt} \frac{\partial^2 u}{\partial t^2} = \int_{-\infty}^{t} be^{(r-n)\tau} d\left(e^{n\tau} \frac{\partial^2 u}{\partial x^2}\right). \qquad /5/$$

Аналогично, дифференцированием по $x$ уравнения /3 / и используя те же зависимости для

$\sigma, u$ и $\varepsilon$ мы получили бы

$$e^{nt} \frac{\partial^2 u}{\partial x^2} = \int_{-\infty}^{t} \frac{1}{b} e^{-(r-n)\tau} d\left(\rho e^{r\tau} \frac{\partial^2 u}{\partial \tau^2}\right). \qquad /6/$$

Уравнения /5/ и / 6/, которым должно удовлетворять смещение $u$, вследствие того, что

$$d\left(e^{n\tau} \frac{\partial^2 u}{\partial x^2}\right) = e^{n\tau} \left\{ n \frac{\partial^2 u}{\partial x^2} + \frac{\partial^3 u}{\partial x^2 \partial t} \right\} d\tau$$

$$d(\rho e^{r\tau}\frac{\partial^2 u}{\partial \tau^2}) = \rho e^{r\tau}\left\{r\frac{\partial^2 u}{\partial \tau^2}+\frac{\partial^3 u}{\partial \tau^3}\right\}d\tau ,$$

можно представить в виде

$$\frac{\partial^2 u}{\partial t^2} = \int_{-\infty}^{t}\frac{b}{\rho}e^{-r(t-\tau)}\left\{n\frac{\partial^2 u}{\partial x^2}+\frac{\partial^3 u}{\partial x^2 \partial \tau}\right\}d\tau ,$$

$$\frac{\partial^2 u}{\partial x^2} = \int_{-\infty}^{t}\frac{\rho}{b}e^{-n(t-\tau)}\left\{r\frac{\partial^2 u}{\partial \tau^2}+\frac{\partial^3 u}{\partial \tau^3}\right\}d\tau \quad . \qquad /7/$$

Каждое из двух уравнений /7/ есть следствие одного и того же третьего, более общего. Действительно, дифференцируя /7/ по $t$, имеем как из одного, так и из другого:

$$\frac{\partial^3 u}{\partial x^2 \partial t} = \frac{\rho^r}{b}\frac{\partial^2 u}{\partial t^2}+\frac{\rho}{b}\frac{\partial^3 u}{\partial t^3}-n\frac{\partial^2 u}{\partial x^2} . \qquad /8/$$

Это есть именно то уравнение, описывающее колебания релаксирующего и последействующего стержня, которое изучал А.Ю.Ишлинский [37]. От того, которое изучал я в своей первой работе[30], оно отличается дополнительным членом $\frac{\rho}{b}\frac{\partial^3 u}{\partial t^3}$ ,

Учитывающим линейную релаксацию. То, что в моей работе [30] обозначено через $M$ и $\eta$ /упругость и вязкость/, здесь представлено как $\frac{bn}{r}$ и $\frac{b}{r}$, так что

$$\frac{b}{\eta} = r \quad \text{и} \quad \frac{M}{\eta} = n . \qquad /9/$$

Отсюда видно, что предположение

$$r > n$$

равнозначно с допущением

$$b > M .$$

На основании /9/, первое из двух уравнений /7/ напишем с обозначением констант, употребленным в моей статье [30]:

$$\rho\frac{\partial^2 u}{\partial t^2} = \int_{-\infty}^{t} re^{-r(t-\tau)}\left\{M\frac{\partial^2 u}{\partial x^2}+\eta\frac{\partial^3 u}{\partial x^2 \partial \tau}\right\}d\tau$$

В случае двухмерной задачи с мембраной мы имели бы:

$$\rho\frac{\partial^2 u}{\partial t^2}=\int_{-\infty}^{t} re^{-r(t-\tau)}\left\{M\left[\frac{\partial^2 u}{\partial x^2}+\frac{\partial^2 u}{\partial y^2}\right]+\eta\left[\frac{\partial^3 u}{\partial x^2\partial\tau}+\frac{\partial^3 u}{\partial y^2\partial\tau}\right]\right\}d\tau\ . \qquad /10/$$

Второму уравнению из /7/ в задаче о мембране соответствует:

$$M\left[\frac{\partial^2 u}{\partial x^2}+\frac{\partial^2 u}{\partial y^2}\right]=\int_{-\infty}^{t}\frac{\rho}{\eta}e^{-n(t-\tau)}\left\{\frac{\partial^2 u}{\partial \tau^2}+\frac{1}{r}\frac{\partial^3 u}{\partial \tau^3}\right\}d\tau\,. \qquad /11/$$

Вместо /8/ получилось бы:

$$\rho\frac{\partial^3 u}{\partial t^3}=Mr\left[\frac{\partial^2 u}{\partial x^2}+\frac{\partial^2 u}{\partial y^2}\right]+\eta r\left[\frac{\partial^3 u}{\partial x^2\partial t}+\frac{\partial^3 u}{\partial y^2\partial t}\right]-r\rho\frac{\partial^2 u}{\partial t^2}\,,$$

или, после деления на $r$:

$$\eta\left[\frac{\partial^3 u}{\partial x^2\partial t}+\frac{\partial^3 u}{\partial y^2\partial t}\right]+M\left[\frac{\partial^2 u}{\partial x^2}+\frac{\partial^2 u}{\partial y^2}\right]=\rho\left[\frac{\partial^2 u}{\partial t^2}+\frac{1}{r}\frac{\partial^3 u}{\partial t^3}\right]\ . \qquad /П/$$

Таково уравнение поперечных колебаний упруго-вязкой мембраны с учетом линейного закона не только упругого последействия, как в моей второй работе /31/, но и релаксации.

Мы рассмотрим более общее уравнение вынужденных колебаний:

$$\eta\left[\frac{\partial^3 u}{\partial x^2\partial t}+\frac{\partial^3 u}{\partial y^2\partial t}\right]+M\left[\frac{\partial^2 u}{\partial x^2}+\frac{\partial^2 u}{\partial y^2}\right]+f(x,y;t)=\rho\left[\frac{\partial^2 u}{\partial t^2}+\frac{1}{r}\frac{\partial^3 u}{\partial t^3}\right]\ . \qquad /Ш/$$

Прежде всего, считая $r$ достаточно большим, *покажем*, что *оно не может иметь более одного решения*, которое удовлетворяло бы начальным и краевым условиям задачи: $u|_{t=0}=\Phi(x,y)$; $\frac{\partial u}{\partial t}|_{t=0}=\varphi(x,y)$; $u|_{C}=0$ при $t\geq 0$. /Ш$^{\mathrm{I}}$/

<u>*В самом деле*</u>, если бы это было не так и $u_1$ и $u_2$ были два различных решения задачи /Ш/,/Ш$^{\mathrm{I}}$/ , то однородное уравнение /П/, очевидно, имело бы решение $u=u_1-u_2$ , удовлетворяющее условиям:

$$u|_{t=0}=0;\qquad \frac{\partial u}{\partial t}|_{t=0}=0;\ u|_{C}=0 \text{ при } t\geq 0\,, \qquad /П^{\mathrm{I}}/$$

и отличное от тождественного нуля. Покажем, что такого решения $u$ уравнения /П/ при условиях /П$^{\mathrm{I}}$/ не существует.

В формуле Green'а

$$\iiint\left\{\frac{\partial P}{\partial x}+\frac{\partial Q}{\partial y}+\frac{\partial R}{\partial t}\right\}dxdydt=\iint Pdydt+Qdtdx+Rdxdy\,,$$

где за основную область в системе декартовых прямоугольных координат $OXYT$ принята часть пространства, ограниченная плоскостями $t=0$ и $t=T_0(T>0$ и произвольно выбрано) и цилиндрической поверхностью с контуром С мембраны в качестве направляющей и с образующей, параллельной оси $OT$, положим

$$P=\eta\frac{\partial^2 u}{\partial x\partial t}\frac{\partial u}{\partial t},\quad Q=\eta\frac{\partial^2 u}{\partial y\partial t}\frac{\partial u}{\partial t},\quad R=\eta\frac{\partial^2 u}{\partial t\partial t}\frac{\partial u}{\partial t}.$$

Так как

$$\frac{\partial}{\partial x}\left(\eta\frac{\partial u}{\partial t}\frac{\partial^2 u}{\partial x\partial t}\right)=\eta\frac{\partial u}{\partial t}\frac{\partial^3 u}{\partial x^2\partial t}+\eta\left(\frac{\partial^2 u}{\partial x\partial t}\right)^2$$

и аналогично для двух других координат, получим

$$\iiint\eta\frac{\partial u}{\partial t}\left[\frac{\partial^3 u}{\partial x^2\partial t}+\frac{\partial^3 u}{\partial y^2\partial t}+\frac{\partial^3 u}{\partial t^3}\right]dxdydt+\iiint\eta\frac{\partial u}{\partial t}\left[\left(\frac{\partial^2 u}{\partial x\partial t}\right)^2+\left(\frac{\partial^2 u}{\partial y\partial t}\right)^2+\left(\frac{\partial^2 u}{\partial t^2}\right)^2\right]dxdydt=$$

$$=\iint\eta\frac{\partial u}{\partial t}\left[\frac{\partial^2 u}{\partial x\partial t}dydt+\frac{\partial^2 u}{\partial y\partial t}dtdx+\frac{\partial^2 u}{\partial t^2}dxdy\right].\qquad /12/$$

Из $u|_C=0$ при всех $t\geq 0$ следует, что $\frac{\partial u}{\partial t}|_C=0$. Поэтому слагаемые правой части, соответствующие боковой поверхности цилиндра, равны нулю. Кроме того, на торцах $t=0$ и $t=T_0$ цилиндра $dt=0$, а на торце $t=0$, в дополнение к этому, еще и $\frac{\partial u}{\partial t}=0$. Таким образом, от поверхностного интеграла правой части остается

$$\iint\eta\left[\frac{\partial u}{\partial t}\frac{\partial^2 u}{\partial t^2}\right]_{t=T}dxdy=\iint\frac{\eta}{2}\left[\frac{\partial u}{\partial t}\left(\frac{\partial u}{\partial t}\right)^2\right]_{t=T}dxdy\ .$$

Используя исходное уравнение /П/, перепишем /12/ в виде:

$$\iiint\eta\left[\left(\frac{\partial^2 u}{\partial x\partial t}\right)^2+\left(\frac{\partial^2 u}{\partial y\partial t}\right)^2+\left(\frac{\partial^2 u}{\partial t^2}\right)^2\right]dxdydt+\iiint\frac{\rho}{2}\frac{\partial}{\partial t}\left(\frac{\partial u}{\partial t}\right)^2dxdydt+\iiint\frac{\rho}{r}\frac{\partial u}{\partial t}\frac{\partial^3 u}{\partial t^3}dxdydt+$$

$$+\iiint\eta\frac{\partial u}{\partial t}\frac{\partial^3 u}{\partial t^3}dxdydt-\iiint M\frac{\partial u}{\partial t}\left[\frac{\partial^2 u}{\partial x^2}+\frac{\partial^2 u}{\partial y^2}\right]dxdydt=$$

$$=\iint\frac{\eta}{2}\left[\frac{\partial}{\partial t}\left(\frac{\partial u}{\partial t}\right)^2\right]_{t=T}dxdy\ .$$

Интегрируя второй интервал левой части по t, получаем:

$$\iiint\frac{\rho}{2}\frac{\partial}{\partial t}\left(\frac{\partial u}{\partial t}\right)^2dxdydt=\iint\frac{\rho}{2}\left[\left(\frac{\partial u}{\partial t}\right)^2_{t=T}-\left(\frac{\partial u}{\partial t}\right)^2_{t=0}\right]dxdy\ .$$

Имея в виду, что $\frac{\partial u}{\partial t}\Big|_{t=0}=0$ согласно начальному условию, перепишем предыдущее соотношение так:

$$\iiint\eta\left[\left(\frac{\partial^2u}{\partial x\partial t}\right)^2+\left(\frac{\partial^2u}{\partial y\partial t}\right)^2+\left(\frac{\partial^2u}{\partial t^2}\right)^2\right]dxdydt-\iiint M\frac{\partial u}{\partial t}\left[\frac{\partial^2u}{\partial x^2}+\frac{\partial^2u}{\partial y^2}\right]dxdydt+$$

$$+\iint\frac{\rho}{2}\frac{\partial}{\partial t}\left(\frac{\partial u}{\partial t}\right)^2_{t=T}dxdy+\iiint(\eta+\frac{\rho}{r})\frac{\partial u}{\partial t}\frac{\partial^3u}{\partial t^3}dxdydt=$$

$$=\iint\frac{\eta}{2}\left[\frac{\partial}{\partial t}\left(\frac{\partial u}{\partial t}\right)^2\right]_{t=T}dxdy,$$

или

$$\iiint\eta\left[\left(\frac{\partial^2u}{\partial x\partial t}\right)^2+\left(\frac{\partial^2u}{\partial y\partial t}\right)^2+\left(\frac{\partial^2u}{\partial t^2}\right)^2\right]dxdydt-\iiint M\frac{\partial u}{\partial t}\left[\frac{\partial^2u}{\partial x^2}+\frac{\partial^2u}{\partial y^2}\right]dxdydt=$$

$$=\frac{1}{2}\iint\left[\eta\frac{\partial}{\partial t}\left(\frac{\partial u}{\partial t}\right)^2-\rho\left(\frac{\partial u}{\partial t}\right)^2\right]_{t=T}dxdy-\iiint\left(\eta+\frac{\rho}{r}\right)\frac{\partial u}{\partial t}\frac{\partial^3u}{\partial t^3}dxdydt.$$

Выполняя в последнем интеграле интегрирование по частям, после приведения подобных членов получим

$$\iiint\eta\left[\left(\frac{\partial^2u}{\partial x\partial t}\right)^2+\left(\frac{\partial^2u}{\partial y\partial t}\right)^2\right]dxdydt-\iiint M\frac{\partial u}{\partial t}\left[\frac{\partial^2u}{\partial x^2}+\frac{\partial^2u}{\partial y^2}\right]dxdydt=$$

$$=-\iint\frac{\rho}{2}\left(\frac{\partial u}{\partial t}\right)^2_{t=T}dxdy-\frac{\rho}{r}\iiint\frac{\partial u}{\partial t}\frac{\partial^3u}{\partial t^3}dxdydt. \qquad /13/$$

Преобразуем второй интеграл левой части:

$$\iiint M\frac{\partial u}{\partial t}\frac{\partial^2u}{\partial x^2}dxdydt=\iiint M\left[\frac{\partial}{\partial x}\left(\frac{\partial u}{\partial t}\frac{\partial u}{\partial x}\right)-\frac{1}{2}\frac{\partial}{\partial t}\left(\frac{\partial u}{\partial x}\right)^2\right]dxdydt,$$

$$\iiint M\frac{\partial u}{\partial t}\frac{\partial^2u}{\partial y^2}dxdydt=\iiint M\left[\frac{\partial}{\partial y}\left(\frac{\partial u}{\partial t}\frac{\partial u}{\partial y}\right)-\frac{1}{2}\frac{\partial}{\partial t}\left(\frac{\partial u}{\partial y}\right)^2\right]dxdydt$$

Складывая все почленно, имеем:

$$\iiint M\frac{\partial u}{\partial t}\left[\frac{\partial^2u}{\partial x^2}+\frac{\partial^2u}{\partial y^2}\right]dxdydt=\iiint M\left[\left[\frac{\partial}{\partial x}\left(\frac{\partial u}{\partial t}\frac{\partial u}{\partial x}\right)+\frac{\partial}{\partial y}\left(\frac{\partial u}{\partial t}\frac{\partial u}{\partial y}\right)\right]-\frac{1}{2}\frac{\partial}{\partial t}\left[\left(\frac{\partial u}{\partial x}\right)^2+\left(\frac{\partial u}{\partial y}\right)^2\right]\right]dxdydt$$

или, согласно формуле Green'а,

здесь первые два элемента интеграла правой части дают нуль, так как в плоскостях $t = 0$ и $t = T_0$ имеем $dt = 0$, а на боковой поверхности $\frac{\partial u}{\partial t} = 0.$

Поэтому

$$\iiint M \frac{\partial u}{\partial t}\left[\frac{\partial^2 u}{\partial x^2} + \frac{\partial^2 u}{\partial y^2}\right] dxdydt = -\iint \frac{M}{2}\left[\left(\frac{\partial u}{\partial x}\right)^2 + \left(\frac{\partial u}{\partial y}\right)^2\right]_{t=T} dxdy.$$

Соотношение /13/ принимает теперь вид:

$$\iiint \eta \left[\left(\frac{\partial^2 u}{\partial x \partial t}\right)^2 + \left(\frac{\partial^2 u}{\partial y \partial t}\right)^2\right] dxdydt = -\iint \frac{M}{2}\left[\left(\frac{\partial^2 u}{\partial x^2}\right)^2 + \left(\frac{\partial^2 u}{\partial y^2}\right)^2\right]_{t=T} dxdy -$$

$$-\iint \frac{\rho}{2}\left(\frac{\partial u}{\partial t}\right)^2_{t=T} dxdy - \frac{\rho}{r}\iiint \frac{\partial u}{\partial t}\frac{\partial^3 u}{\partial t^3} dxdydt \quad . \qquad /14/$$

Допустим, что скорость $r$ релаксации настолько велика, что последний член правой части может быть опущен без значительной погрешности. Тогда будем иметь:

$$\iiint \eta \left[\left(\frac{\partial^2 u}{\partial x \partial t}\right)^2 + \left(\frac{\partial^2 u}{\partial y \partial t}\right)^2\right] dxdydt = -\iint \left\{\frac{M}{2}\left[\left(\frac{\partial^2 u}{\partial x^2}\right)^2 + \left(\frac{\partial^2 u}{\partial y^2}\right)^2\right] + \frac{\rho}{2}\left(\frac{\partial u}{\partial t}\right)^2\right\}_{t=T} dxdy \ .$$

Вследствие не отрицательности слагаемых и противоположности знаков слева и справа должно быть

$$\left.\frac{\partial u}{\partial x}\right|_{t=T} = 0, \quad \left.\frac{\partial u}{\partial y}\right|_{t=T} = 0, \qquad \left.\frac{\partial u}{\partial t}\right|_{t=T} = 0$$

при любом $T \geq 0$.

Отсюда следует, что $u = const$. А так как $u = 0$ при $t = 0$, то $u \equiv 0$, *что и доказывает единственность решения уравнения /Ш/ при условиях /Ш$^I$/.*

При этом существенным оказалось предположение, что скорость релаксации $r$ настолько велика, что член

$$-\frac{\rho}{r}\iiint \frac{\partial u}{\partial t}\frac{\partial^3 u}{\partial t^3} dxdydt$$

не влияет на знак правой части в равенстве /14/.

А.Ю.Ишлинский в своей работе [27] о колебаниях стержня дает несколько иное доказательство единственности решения, предполагая, что скорость релаксации больше скорости последействия $n$.

Вернемся к однородному уравнению /П/ свободных колебаний мембраны:

$$\eta\left[\frac{\partial^3 u}{\partial x^2 \partial t}+\frac{\partial^3 u}{\partial y^2 \partial t}\right]+M\left[\frac{\partial^2 u}{\partial x^2}+\frac{\partial^2 u}{\partial y^2}\right]=\rho\left[\frac{\partial^2 u}{\partial t^2}+T\frac{\partial^3 u}{\partial t^3}\right], \qquad /П/$$

где $T=\frac{1}{r}$ есть период релаксации, предполагаемый очень малым. Будем искать частное решение для /П/ в виде произведения из функции координат $Z(x,y)$ на функцию от времени $U(t)$. После подстановки такого вида решения $u$ в /П/ мы получим:

$$\eta\left[\frac{\partial^2 Z}{\partial x^2}+\frac{\partial^2 Z}{\partial y^2}\right]U'+M\left[\frac{\partial^2 Z}{\partial x^2}+\frac{\partial^2 Z}{\partial y^2}\right]U=\rho\left[ZU''+\frac{1}{r}ZU'''\right],$$

или, разделяя переменные:

$$\frac{\frac{\partial^2 Z}{\partial x^2}+\frac{\partial^2 Z}{\partial y^2}}{Z}=\frac{\rho\left[U''+\frac{1}{r}U'''\right]}{MU+\eta U'}.$$

Поскольку левая часть зависит только от координат, а правая–только от времени и поскольку между координатами $x,y$ и временем $t$ не может быть никакой зависимости, каждое из этих отношений есть постоянная величина, которую мы обозначим $-\lambda$._

Таким образом, приходим к следующим двум уравнениям нашей задачи:

$$\frac{\partial^2 Z}{\partial x^2}+\frac{\partial^2 Z}{\partial y^2}=-\lambda Z, \qquad /IY/$$

$$\frac{\rho}{r}U'''+\rho U''+\lambda\eta U'+\lambda MU=0. \qquad /Y/$$

Уравнение /IY/ было подробно исследовано в нашей работе: «К вопросу о малых колебаниях упруго-вязких мембран» /см.[31], стр.473-480/. Оно определяет последовательность собственных значений $\lambda_k$ и собственных функций $Z_k(x,y)$, которые можно считать нормированными и попарно ортогональными. Я позволю себе не останавливаться здесь на нем.

Уравнению /Y/ соответствует характеристическое уравнение:

$$\frac{\rho}{r}x^3+\rho x^2+\lambda\eta x+\lambda M=0, \qquad /15/$$

которое не имеет положительных корней, так как все его коэффициенты положительны. Легко установить границы, в которых лежат его отрицательные корни. Именно, переписав это уравнение в виде

$$x^2 = -\frac{\lambda}{\rho}\frac{\eta x + M}{\frac{1}{r}x + 1},$$

видим, что действительный корень $x$ этого уравнения должен сообщить противоположные знаки линейным функциям

$$\eta x + M \text{ и } \frac{1}{r}x + 1.$$

Следовательно, действительный корень $x$ /15/ должен лежать между точками пересечения с осью $OX$ прямых

$$y_1 = \eta x + M, \quad y_2 = \frac{1}{r}x + 1.$$

А так как эти прямые пересекают $OX$ в точках с абсциссами

$$x' = -\frac{M}{\eta} = -n, \quad x'' = -r,$$

То каждый действительный корень $\alpha$ уравнения /15/ удовлетворяет неравенствам

$$-r \le \alpha \le -n. \qquad /16/$$

А.Ю.Ишлинский [37] доказывает это неравенство из соотношений Виетта. Далее, по способу, применяемому Ишлинским, можно доказать, что *если /15/ имеет два сопряженных комплексных корня и один отрицательный, то вещественная часть комплексного корня отрицательна.*

<u>В самом деле</u>, пусть

$\alpha + \beta i, \alpha - \beta i, \gamma;$ ($\alpha, \beta$-действительные)

три корня уравнения /15/. Тогда в силу соотношений Виетта:

$2\alpha + \gamma = -r.$

Так как самое меньшее возможное значение для $\gamma$ есть $-r$ по /16/, то самое большое значение для $2\alpha$ есть 0, так что $\alpha \le 0$, что и утверждалось.

Если назвать $\alpha_{1,m}\alpha_{2,m}, \alpha_{3,m}$ корни характеристического уравнения /15/, соответствующие собственному значению $\lambda = \lambda_m (m = 1,2,...)$, то общее решение краевой задачи о последействующей и релаксирующей мембране представится в виде:

$$u(x,y;t) = \sum_{m=1}^{\infty}\left\{C_{1,m}e^{\alpha_{1,m}t} + C_{2,m}e^{\alpha_{2,m}t} + C_{3,m}e^{\alpha_{3,m}t}\right\}Z_m(x,y),$$

Причем постоянные С с индексами должны быть подобраны так, чтобы удовлетворялись начальные условия:

$$u\big|_{t=0} = \Phi(x,y), \quad \dot{u}\big|_{t=0} = \varphi(x,y), \quad \ddot{u}\big|_{t=0} = \psi(x,y).$$

Предполагая, что $\Phi(x,y), \varphi(x,y)$ и $\psi(x,y)$ могут быть разложены в ряд по собственным функциям $Z_m(x,y)$ задачи:

$$u(x,0)=\sum_{m=1}^{\infty}a_m Z_m(x,y),\ \ \dot{u}(x,y;0)=\sum_{m=1}^{\infty}b_m Z_m(x,y),\ \ \ddot{u}(x,y;0)=\sum_{m=1}^{\infty}c_m Z_m(x,y),$$

получим, в качестве уравнений для определения C:

$$\begin{aligned} &C_{1,m}+C_{2,m}+C_{3,m}=a_m\\ &\alpha_{1,m}C_{1,m}+\alpha_{2,m}C_{2,m}+\alpha_{3,m}C_{3,m}=b_m\\ &\alpha_{1,m}^2C_{1,m}+\alpha_{2,m}^2C_{2,m}+\alpha_{3,m}^2C_{3,m}=c_m \end{aligned}$$

Наконец, решение задачи о вынужденных колебаниях упруго-вязкой мембраны с реаксацией и последействием сообразно линейному закону осуществляется по способу почти не отличающемуся от того, который был применен мной во второй из предложенных работ [31]. При этом все ряды, которые получаются дифференцированием полученного решения столько раз, сколько это требуется уравнением /Ш/, сходятся равномерно и абсолютно.

Москва, октябрь 1942 г. *А.ГЕРАСИМОВ*

ЦИТИРОВАННАЯ ЛИТЕРАТУРА

# **8.**ОБОБЩЕНИЕ ЛИНЕЙНЫХ ЗАКОНОВ ДЕФОРМИРОВАНИЯ И ЕГО ПРИМЕНЕНИЕ К ЗАДАЧАМ ВНУТРЕННЕГО ТРЕНИЯ*

Явления упругого последействия в свете теории наследственности описываются линейным интегральным уравнением Больцмана [1]

$$\sigma(t) = E\varepsilon(t) + \int_0^\infty G(\tau)\varepsilon(t-\tau)d\tau ,$$

где $E$ – упругая постоянная и $G(\tau)$ – наследственная функция, определяемая из опыта. Экспериментальные исследования показывают, что особого внимания заслуживает частный вид этого уравнения, соответствующий только наследственной части напряжения,

$$\sigma(t) = \int_0^\infty G(\tau)\varepsilon(t-\tau)d\tau .$$

Не меньший интерес представляет также тот случай, когда напряжение $\sigma(t)$ зависит от всех предшествующих, надлежащим образом взвешенных значений скоростей деформации, но не деформации.

Для таких процессов деформирования зависимость между $\sigma$ и $\varepsilon$

$$\sigma(t) = \int_0^\infty K(\tau)\dot{\varepsilon}(t-\tau)d\tau \quad ,$$

если ограничиваться только наследственной частью напряжения.

Наследственная функция для некоторых материалов (волокнистой структуры) должна иметь вид

$$K(\tau) = \frac{A}{\tau^\alpha} ,$$

где постоянная $A > 0$ и $\alpha$ лежит между нулем и единицей. Положив еще

$$A = \frac{\kappa}{\Gamma(1-\alpha)} ,$$

где постоянная $\kappa > 0$ зависит от свойств вещества и $\Gamma$ есть эйлеров интеграл второго рода, получим

$$\sigma(t) = \kappa \frac{1}{\Gamma(1-\alpha)} \int_0^\infty \frac{\dot{\varepsilon}(t-\tau)d\tau}{\tau^\alpha} = \kappa \frac{\partial^\alpha \varepsilon(t)}{\partial t^\alpha} \qquad (0 < \alpha < 1) \qquad (1)$$

Так как производная от $\varepsilon(t)$ по $t$ порядка $\alpha$ будет

$$\frac{1}{\Gamma(1-\alpha)} \int_0^\infty \frac{\dot{\varepsilon}(t-\tau)d\tau}{\tau^\alpha} .$$

---

* *Доложено в Институте механики АН СССР 29 мая 1947 г.*

Это линейное соотношение между $\sigma$ и $\varepsilon$ при $\alpha = 0$ обращается в закон Гука, при $\alpha = 1$ – в закон Ньютона для внутреннего трения. Мы будем исходить из зависимости (1) для любого $\alpha$ между 0 и 1.

*<u>Пример 1.</u>* Рассмотрим задачу о движении жидкости между параллельными плоскостями, из которых одна неподвижна, а другая движется прямолинейно параллельно первой по заданному закону. Примем зависимость между $\sigma$ и $\varepsilon$ для жидкости в виде (1), где $\alpha$ имеет определенное значение $(0 < \alpha < 1)$.

Пусть (фиг.1) плоскость $x = 0$ неподвижна, плоскость $x = l$ движется в самой себе направлении оси $y$ по закону

$$y(x,t)\big|_{x=l} = \varphi(t) \qquad (\varphi(0) = 0, \varphi'(0) = 0)$$

где $\varphi(t)$ – данная функция времени.

Начальные и граничные условия предположим соответственно в виде

$$y(x,t)\big|_{t=0} = 0, \quad \frac{\partial y}{\partial t}\Big|_{t=0} = 0, \quad y(x,t)\big|_{x=0}, y(x,t)\big|_{x=l} = \varphi(t) \qquad (2)$$

Рассмотрим элемент объема жидкости $AB$, ограниченный двумя гранями с площадями, равными единице в плоскостях $x = const$ и $x + dx = const$. Уравнение движения получим в виде

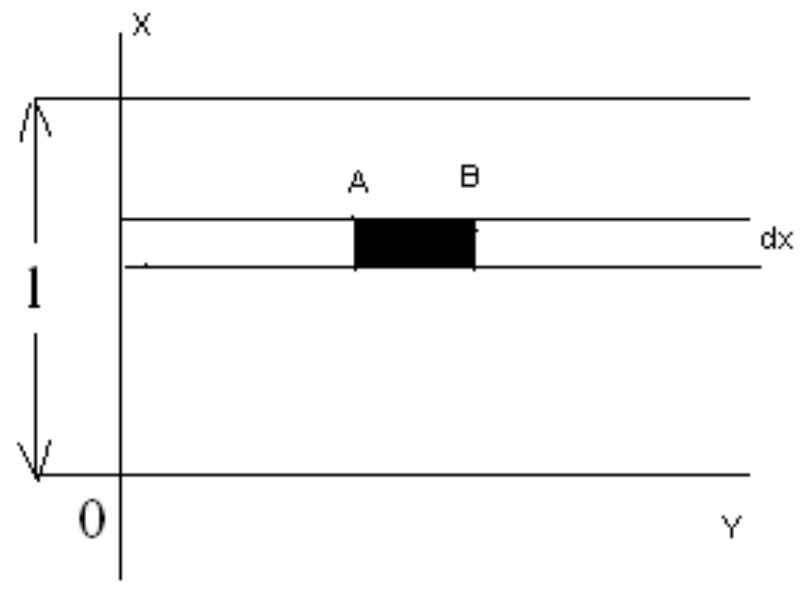

Фиг.1

$$\frac{\partial \sigma}{\partial x} = \rho \frac{\partial^2 y}{\partial t^2} \qquad (3)$$

где $\rho = const$ есть плотность жидкости.

Дифференцируя (1) по $x$ и помня, что $\varepsilon = \partial y / \partial x$, после подстановки в (3) получим

$$\kappa \frac{\partial^\alpha}{\partial t^\alpha}\left(\frac{\partial^2 y}{\partial x^2}\right) = \rho \frac{\partial^2 y}{\partial t^2} \qquad (4)$$

Будем искать решение $y(x,t)$ уравнения (4) при условиях (2) по способу Хевисайда. Пусть $Y(x, p)$ –изображение для $y(x,t)$.

В области изображений уравнению (4) при условиях (2) будет отвечать уравнение

$$\frac{d^2 Y}{dx^2} = c^2 p^{2m} Y \qquad (c^2 = \frac{\rho}{\kappa}, 2m = 2 - \alpha, \frac{1}{2} < m < 1) \qquad (5)$$

при условиях

$$Y\big|_{x=0} = 0, \quad Y\big|_{x=l} = \Phi(p) \stackrel{\bullet}{\underset{\bullet}{\leftarrow}} \varphi(t) \qquad (6)$$

Решение уравнения (5) при условиях (6) имеет вид

$$Y(x,p)=\Phi(p)\frac{shcp^m x}{shcp^m l} \qquad (7)$$

По теореме Бореля получим для оригинала

$$y(x,t)=\frac{d}{dt}\int_0^t K(t-\theta)\varphi(\theta)d\theta \qquad (8)$$

где

$$K(\vartheta)\xrightarrow{\bullet}\frac{shcp^m x}{shcp^m l} \quad \text{или} \quad K(\vartheta)=\frac{1}{2\pi i}\int_L \frac{e^{p\vartheta}}{p}\frac{shcp^m x}{shcp^m l}dp \qquad (9)$$

Здесь $L$ есть, например, прямая, проходящая параллельно мнимой оси $p$ справа от всех особенностей подинтегральной функции. Применяя известное разложение, имеем

$$\frac{shcp^m x}{shcp^m l}=\frac{x}{l}+\frac{2}{\pi}\sum_{k=1}^{\infty}\frac{(-1)^k}{k}\sin\frac{k\pi x}{l}\sum_{j=0}^{\infty}(-1)^j\left(\frac{k\pi}{cl}\right)^{2j}p^{-2jm} \qquad (10)$$

Так как при $\tau \geq 0$

$$p^{-2jm}\xleftarrow{\bullet}\frac{\tau^{2jm}}{\Gamma(2jm+1)} \qquad (j=0,1,2,...)$$

то

$$K(t-\theta)=\frac{x}{l}+\frac{2}{\pi}\sum_{k=1}^{\infty}\frac{(-1)^k}{k}\sin\frac{k\pi x}{l}\sum_{j=0}^{\infty}(-1)^j\left(\frac{k\pi}{cl}\right)^{2j}\frac{(t-\theta)^{2jm}}{\Gamma(2jm+1)} \qquad (11)$$

Действительная часть функций

$$\sum_{j=0}^{\infty}(-1)^j\left(\frac{k\pi}{cl}\right)^{2j}\frac{\tau^{2jm}}{\Gamma(2jm+1)} \qquad (k=1,2,...) \qquad (12)$$

для $\frac{1}{2}<m<1$ вообще многозначна; но при $m=\frac{1}{2}$ или $m=1$ имеет лишь одну действительную ветвью. Рассмотрим частные случаи.

При $m=\frac{1}{2}$ , или $\alpha=1$,т.е. для жидкости, следующей закону Ньютона $\sigma=\kappa\, d\varepsilon(t)/dt$, из (11) получается

$$K(t-\theta)=\frac{x}{l}+\frac{2}{\pi}\sum_{k=1}^{\infty}\frac{(-1)^k}{k}\sin\frac{k\pi x}{l}\exp[-\frac{k^2\pi^2(t-\theta)}{c^2l^2}]$$

Для распределения смещения имеем

$$y(x,t)=\frac{d}{dt}\int_0^t\left\{\frac{x}{l}+\frac{2}{\pi}\sum_{k=1}^{\infty}\frac{(-1)^k}{k}\sin\frac{k\pi x}{l}\exp[-\frac{k^2\pi^2(t-\theta)}{c^2l^2}]\right\}\varphi(\theta)d\theta$$

или

$$y(x,t)=\left\{\frac{x}{l}+\frac{2}{\pi}\sum_{k=1}^{\infty}\sin\frac{k\pi x}{l}\right\}\varphi(t)- \qquad (13)$$

$$-\frac{2\pi}{c^2l^2}\sum_{k=1}^{\infty}(-1)^k k\sin\frac{k\pi x}{l}\int_0^t\varphi(\theta)\exp[-\frac{k^2\pi^2(t-\theta)}{c^2l^2}]d\theta$$

Как известно,

$$\frac{x}{l}+\frac{2}{\pi}\sum_{k=1}^{\infty}\frac{(-1)^k}{k}\sin\frac{k\pi x}{l}=\begin{cases}0(x<l)\\1(x=l)\end{cases} \qquad (14)$$

Поэтому из (13) получаем

$$y(x,t)=-\frac{2}{\pi^2c^2}\sum_{k=1}^{\infty}(-1)^k k\sin\frac{k\pi x}{l}\int_0^t\varphi(\theta)\exp[-\frac{k^2\pi^2(t-\theta)}{c^2l^2}]d\theta \quad (0\le x<l)$$

$$y(x,t)\big|_{x=l}=\varphi(t)$$

В таком виде решение для этого случая дано А. И. Лурье [2].

При $\alpha=0$, или $m=1$, т.е. когда между плоскостями находится упругая среда, деформируемая на сдвиг, имеем

$$K(t-\theta)=\frac{x}{2}+\frac{2}{\pi}\sum_{k=1}^{\infty}\frac{(-1)^k}{k}\sin\frac{k\pi x}{l}\sin\frac{k\pi(t-\theta)}{cl}$$

Для распределения смещений получаем

$$y(x,l)=\left\{\frac{x}{l}+\frac{2}{\pi}\sum_{k=1}^{\infty}\frac{(-1)^k}{k}\sin\frac{k\pi x}{l}\right\}\varphi(t)-$$

$$-\frac{2}{cl}\sum_{k=1}^{\infty}(-1)^k k\sin\frac{k\pi x}{l}\int_0^t\varphi(\theta)\sin\frac{k\pi(t-\theta)}{cl}d\theta$$

На основании (14) найдем

$$y(x,t)=-\frac{2}{cl}\sum_{k=1}^{\infty}(-1)^k k\sin\frac{k\pi x}{l}\int_0^t\varphi(\theta)\sin\frac{k\pi(t-\theta)}{cl}d \quad (0\le x\le l)$$

$$y(l,t)=\varphi(t)$$

Это, как известно, уравнение вынужденных поперечных волн в упругой среде.

Таким образом, предельные случаи $\alpha = 1$и $\alpha = 0$не выпадают из найденной выше общей формы решения для любого $\alpha$ между 0 и 1.Поэтому в зависимости (1) можно считать $0 \le \alpha \le 1$.

Вернемся опять к общему выражению (11) для $K(\tau)$. При $\alpha = \frac{1}{2}$ $(m = \frac{3}{4})$ оказалось бы, что каждому значению $\tau \ge 0$ отвечают два действительных значения каждой из функций (12), входящих в выражение (11) для $K(\tau)$. Эти значения будут

$$\sum_{k=0}^{\infty} (-1)^k \frac{\left(a_m \tau^{3/4}\right)^{2k}}{\Gamma(3k/2+1)}, \qquad \sum_{k=0}^{\infty} \frac{\left(a_m \tau^{3/4}\right)^{2k}}{\Gamma(3k/2+1)}, \qquad \left(a_m = \frac{m\pi}{cl}, m = 1,2,\ldots\right) \tag{15}$$

где радикалы понимаются в арифметическом смысле.

Разложения (15) сходятся равномерно для всех $\tau$ на любом интервале конечной длины.

Но, в то время как функции, определяемые первым рядом (15), остаются ограниченными для $\tau \ge 0$, функции, определенные вторым рядом, неограниченно растут с ростом $\tau$. Сохраняя лишь одно первое выражение и отбрасывая второе, как не обеспечивающее должной сходимости ядра $K(\tau)$ при $\tau \to \infty$, введем обозначение

$$C_n(\tau) = \sum_{k=0}^{\infty} a_n^{2k} \frac{(-1)^{k} \tau^{3k/2}}{\Gamma(3k/2+1)} \qquad \left(a_n = \frac{m\pi}{cl}, n = 1,2,\ldots\right) \tag{16}$$

Тогда

$$y(x,t) = \frac{2}{\pi} \sum_{n=1}^{\infty} \frac{(-1)^n}{n} \sin\frac{m\pi x}{l} \int_0^t C_n'(t-\theta)\varphi(\theta)d\theta \qquad (0 \le x \le l) \tag{17}$$

$$y(l,t) = \varphi(t)$$

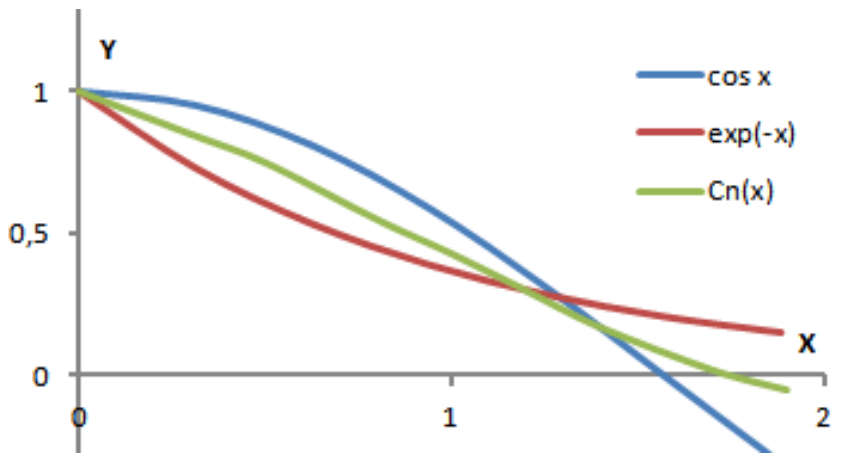

Фиг.2 *Примерные графики функций* $\cos x$ , $\exp(-x)$ *и функции* $C_n(x)$ *при* $a_n = 1$.

Зная распределение смещений (17) для случая $a_n = 1/2$, можно найти величину сдвига $\varepsilon = \partial y/\partial x$, а затем, беря от этого выражения производную по $t$ половинного порядка, получить при помощи (1) распределение напряжений. Однако способ Хевисайда позволяет найти напряжение в любом слое, в частности, на границах без знания распределения смещений. Покажем это для случая равномерного движения верхней плоскости $x = l$

$$\varphi(\theta) = \upsilon\theta \qquad (\upsilon = const > 0)$$

Из решения (7) в области изображений дифференцированием по $x$ и подстановкой в результате $x = l$ имеем

$$\varepsilon\Big|_{x=l} = \frac{\partial y}{\partial x}\Big|_{x=l} \doteq \frac{dY}{dx}\Big|_{x=l} = = cp^m ctghcp^m l\Phi(p)$$

Чтобы найти изображение для $\sigma\big|_{x=l}$, нужно это выражение умножить на $\kappa p^{\alpha} = \kappa p^{2-2m}$. Имеем

$$\sigma\big|_{x=l} \doteq \kappa cp^{5/4}\Phi(p) ctghclp^{3/4}$$

Пользуясь известными рядами и имея в виду, что $\upsilon t \doteq \upsilon/p$, легко найдем

$$\kappa cp^{2-m}\frac{\upsilon}{p} ctghclp^m = \kappa c\upsilon\left\{p^{1-m} + 2\sum_{k=0}^{\infty}\sum_{j=0}^{\infty}(-1)^j \frac{[(2k+2)cl]^j p^{1+(j-1)m}}{j!}\right\}$$

Помня, что $m-$ не целое, следовательно, допускающее возможность применения известной операционной формулы для перехода к оригиналу, получим

$$\sigma\big|_{x=l} = \kappa c\upsilon\left\{\frac{1}{\Gamma(m)t^{1-m}} + 2\sum_{k=0}^{\infty}\sum_{j=0}^{\infty}(-1)^j \frac{[(2k+2)cl]^j}{j!\Gamma(-jm+m)t^{1+(j-1)m}}\right\} \qquad (18)$$

При $t=0$ напряжение $\sigma\big|_{x=l}$ обращается в бесконечность порядка $(1-m)$. Это и естественно, так как в вязкой жидкости нельзя мгновенно сообщить граничной плоскости $x=l$ конечную скорость, применяя усилие конечной величины (Лурье [2]).

*Пример 2.* Рассмотрим движение жидкости между коаксильными цилиндрическими поверхностями, вращающимися по данному закону. Мы опять будем исходить из зависимости (1). Пусть $r_1, r_2-$ радиусы поверхностей внутреннего и внешнего, а $r-$ промежуточного цилиндров, $L-$ длина цилиндров, $\rho-$ плотность жидкости. Далее, пусть $\varphi(r,t)-$ подлежащее определению угловое смещение жидкого слоя на радиусе $r$. Для деформации, очевидно, имеем

$$\varepsilon = r\frac{\partial\varphi}{\partial r} \qquad (19)$$

Нетрудно убедиться, что результирующий момент сил инерции для слоя между поверхностями $r$ и $r+dr$ есть

$$-\rho 2\pi L dr\frac{\partial^2\varphi}{\partial t^2}r^3 \qquad (20)$$

Момент сил напряжения при сдвиге в этом слое равен

$$2\pi L dr\frac{\partial}{\partial r}(r^2\sigma) \qquad (21)$$

Согласно принципу Даламбера

$$\frac{\partial}{\partial r}(r^2\sigma) = \rho r^2 \frac{\partial^2 \varphi}{\partial t^2} \tag{22}$$

Умножив выражение (1) для $\sigma$ на $r^2$ и продифференцировав по $r$, а затем учитывая (20), имеем

$$\frac{\partial}{\partial r}(r^2\sigma) = \kappa \frac{\partial}{\partial r}[r^2 \varepsilon^{(\alpha)}] \tag{23}$$

Сравнивая (22) и (23), получим уравнение движения

$$\rho r^2 \frac{\partial^2 \varphi}{\partial t^2} = \kappa \frac{\partial}{\partial r}[r^3 \frac{\partial}{\partial r} \frac{\partial^\alpha \varphi}{\partial t^\alpha}] \tag{24}$$

Его решение должно удовлетворять начальным условиям, которые предположим в виде

$$\varphi(r,t)\big|_{t=0} = 0, \quad \frac{\partial \varphi}{\partial t}\Big|_{t=0} = 0 \qquad (r_1 \le r \le r_2) \tag{25}$$

=и краевым условиям, выражающим прилипание частиц на стенках:

$$\varphi(r,t)\big|_{r=r_1} = \varphi_1(t), \quad \varphi(r,t)\big|_{r=r_2} = \varphi_2(t) \tag{26}$$

Здесь $\varphi_1$ и $\varphi_2$ – данные функции, обращающиеся в нуль вместе с их производными первого порядка при $t = 0$.

Применим опять способ Хевисайда. Пусть $\Phi(r,p) \xleftarrow{\bullet} \varphi(r,t)$. Уравнению (24) при условиях (25) и (26) в области изображений соответствует

$$c^2 r^3 p^{2m} \Phi = \frac{\partial}{\partial r}[r^3 \frac{\partial \Phi}{\partial r}] \qquad (c^2 = \frac{\rho}{\kappa}, m = \frac{2-\alpha}{2}, \frac{1}{2} < m < 1) \tag{27}$$

Уравнение (27) сопровождается условиями

$$\Phi(r_1,p) = \Phi_1(p), \quad \Phi(r_2,p) = \Phi_2(p), \quad \Phi_1(p) \xleftarrow{\bullet} \varphi_1(t), \ \Phi_2(p) \xleftarrow{\bullet} \varphi_2(t) \tag{28}$$

Вводя дуговое смещение $F(r,p)$ в области изображений под условием

$$F(r,p) = r\Phi(r,p) \tag{29}$$

вместо (27) и (28) получим

$$\frac{d^2F}{dr^2} + \frac{1}{r}\frac{dF}{dr} - \left(c^2 p^{2m} + \frac{1}{r^2}\right)F = 0$$

$$(F(r_1,p) = r_1\Phi_1(p), \quad F(r_2,p) = r_2\Phi_2(p)) \tag{30}$$

Общий интеграл бесселева уравнения (30) имеет вид

$$F(r,p) = AJ_1(\lambda i r) + BN_1(\lambda i r) \qquad (\lambda = cp^m, i^2 = -1) \tag{31}$$

Постоянные $A$ и $B$ надо выбрать так, чтобы было

$$J_1(\lambda i r_1)A + N_1(\lambda i r_1)B = r_1\Phi_1(p), \quad J_1(\lambda i r_2)A + N_1(\lambda i r_2)B = r_2\Phi_2(p) \tag{32}$$

Из (32) можно будет определить $A$ и $B$ в функциях от $p$, если только $p$ не является корнем уравнения

$$\Delta = J_1(\lambda i r_1)N_1(\lambda i r_2) - J_1(\lambda i r_2)N_1(\lambda i r_1) = 0 \tag{33}$$

Тогда получим

$$\Phi(r,p)=\frac{1}{\Delta}\{[r_1\Phi_1(p)J(\lambda ir)N(\lambda ir_2)-r_2\Phi_2(p)J(\lambda ir)N(\lambda ir_1)]+$$
$$+[r_2\Phi_2(p)J(\lambda ir_1)N(\lambda ir)-r_1\Phi_1(p)J(\lambda ir_2)N(\lambda ir)]\} \qquad (34)$$

Здесь и в дальнейшем индекс порядка функций Бесселя опущен. Пользуясь затем теоремой Бореля, получим

$$\varphi(r,t)=\frac{r_1}{r}\frac{d}{dt}\int_0^t\varphi_1(\vartheta)K_1(t-\vartheta,r)d\vartheta+\frac{r_2}{r}\frac{d}{dt}\int_0^t\varphi_2(\vartheta)K_2(t-\vartheta,r)d\vartheta \qquad (35)$$

где $K_1(\theta)$ и $K_2(\theta)$ суть оригиналы, соответствующие функциям

$$\kappa_1(p,r)=\frac{J(\lambda ir)N(\lambda ir_2)-J(\lambda ir_2)N(\lambda ir)}{\Delta} \qquad (36)$$

$$\kappa_2(p,r)=\frac{J(\lambda ir_1)N(\lambda ir)-J(\lambda ir)N(\lambda ir_1)}{\Delta}$$

т.е.

$$K_1(\theta,r)=\frac{1}{2\pi i}\int_L\frac{e^{p\theta}}{p}\kappa_1(p,r)dp,$$

При условии, что кривая $L$ в плоскости $p$, вообще произвольная, лежит правее всех особенностей подынтегральных выражений, в частности, правее всех нулей левой части (33), так что условия возможности решения (32) относительно $A$ и $B$ соблюдаются.

Ограничимся частным случаем, когда $r_2-r_1=\delta_0$ мало по сравнению с $r_1$ (подшипник скольжения). Введем обозначения

$$r_1\Phi_1(p)=F_1,\quad r_2\Phi_2(p)=F_2,\qquad r\Phi(r,p)=F,\ J(\lambda ir_1)=J$$
$$N(\lambda ir_1)=N,\quad \lambda ir_1=x,\qquad \lambda i(r_2-r_1)=\delta,\ \lambda i(r-r_1)=\xi$$

Тогда

$$\lambda ir_2=x+\delta,\quad J(\lambda ir_2)=J+J'\delta,\quad J(\lambda ir)=J+J'\xi$$
$$\lambda ir=x+\xi,\quad N(\lambda ir_2)=N+N'\delta,\ N(\lambda ir)=N+N'\xi$$

После вычислений имеем

$$\Phi(r,p)=\frac{r_1}{r}\frac{r_2-r}{\delta_0}\Phi_1(p)+\frac{r_2}{r}\frac{r-r_1}{\delta_0}\Phi_2(p) \qquad (37)$$

Применения теоремы Бореля не требуется, и сразу находим

$$\varphi(r,t)=\frac{r_1}{r}\frac{r_2-r}{\delta_0}\varphi_1(t)+\frac{r_2}{r}\frac{r-r_1}{\delta_0}\varphi_2(t) \qquad (38)$$

Отсюда для деформации $\varepsilon$ получаем следующее выражение

$$\varepsilon=r\frac{\partial\varphi}{\partial r}=\frac{r_1r_2}{\delta_0 r}\psi(t)\quad(\psi(t)=\varphi_2(t)-\varphi_1(t)) \qquad (39)$$

Наконец, пользуясь (1), получим выражение для напряжения

$$\sigma=\frac{\kappa r_1r_2}{\delta_0 r}\frac{\partial^\alpha\psi}{\partial t^\alpha} \qquad (40)$$

В частности напряжения на поверхностях будут

$$\sigma\Big|_{r=r_1}=\frac{\kappa r_2}{\delta_0}\frac{\partial^\alpha\psi}{\partial t^\alpha},\qquad\qquad \sigma\Big|_{r=r_2}=\frac{\kappa r_1}{\delta_0}\frac{\partial^\alpha\psi}{\partial t^\alpha} \tag{41}$$

Из (41) между прочим видно, что второе меньше первого в $r_2/r_1$ раз; основание этого обстоятельства чисто геометрическое: закон действия и противодействия требует равенства обеих сил, из которых одна распределена по меньшей, другая по большей площади; тогда концентрация силы, т.е. напряжение, должна быть больше на той поверхности, которая меньше. При $\alpha=1$ имеем вместо (41) закон Кулона: напряжение оказывается пропорциональным относительной скорости.

В более общем случае допустим, что относительное угловое смещение $\psi(t)=\varphi_2(t)-\varphi_1(t)$ представляется рядом синусом

$$\psi(t)=\sum_k a_k\sin\omega_k t\qquad (a_k=const) \tag{42}$$

Тогда из (41) найдем

$$\sigma\Big|_{r=r_1}=\frac{\kappa r_2}{\delta_0}\sum_k a_k\omega_k^\alpha\sin(\omega_k t+\alpha\frac{\pi}{2}) \tag{43}$$

Выражение (43) дает ход изменения во времени для напряжений на внутренней поверхности при любом $\alpha$ в интервале $0\le\alpha\le1$ в зависимости (1). Следует, однако, иметь в виду, что движение смазки в подшипнике должно быть ламинарным и, следовательно, относительная угловая скорость $d\psi(t)/dt$ должна быть достаточно малой.

Известные случаи линейных зависимостей с постоянными коэффициентами между напряжением $\sigma$ и деформацией $\varepsilon$ могут быть представлены в единой форме при помощи интегралов Стильтьеса по непрерывно (или скачкообразно) меняющемуся порядку производных по времени.

В самом деле, каждому деформируемому телу по отношению к данной деформации поставим в соответствие две функции $e(\alpha)$ и $s(\alpha)$ действительного переменного $\alpha$, определенные для всех действительных значений $\alpha$. Предполагается, что обе они имеют ограниченное изменение на любом конечном интервале значений $\alpha$. Будем называть $e(\alpha)$ инерцией данного тела относительно $\sigma$ и $s(\alpha)$ инерцией по отношению к деформации. Тогда всякий линейный закон, связывающий $\sigma$ с $\varepsilon$ посредством интегральной, или дифференциальной, или дискретной зависимости с не меняющимися по времени коэффициентами, является лишь частным случаем общей зависимости

$$\int_{-\infty}^{+\infty}\frac{\partial^2\sigma}{\partial t^\alpha}de(\alpha)=\int_{-\infty}^{+\infty}\frac{\partial^2\varepsilon}{\partial t^\alpha}ds(\alpha) \tag{44}$$

Здесь предполагается, что оба интеграла Стильтьеса существуют; там, где $e(\alpha)$ и $s(\alpha)$ не определены, их нужно считать равными нулю. В случае абсолютно упругого тела (закон Гука) в (1) надо положить

$$de(0)=1,\quad ds(0)=E\qquad (de(\alpha)=0, ds(\alpha)=0, \alpha\neq0)$$

В случае, когда тело, не будучи вполне упруго, релаксирует со скоростью $q$ и последействует со скоростью $n$, т.е. когда оно ведет себя по закону (А.Ю.Ишлинский [3])

$$\sigma + \frac{1}{q}\dot{\sigma} = E\left(\varepsilon + \frac{1}{n}\dot{\varepsilon}\right) \qquad (q, n - const > 0, q > n) \tag{45}$$

то чтобы получить (45), в уравнении (1) достаточно положить

$$de(0) = 1, \qquad de(1) = \frac{1}{q}, \qquad ds(0) = E, \qquad ds(1) = \frac{E}{n}$$

$$(de(\alpha) = 0, \quad ds(\alpha) = 0, \quad \alpha \neq 0, \quad \alpha \neq 1)$$

Ясно, что и вторичные релаксационно-последейственные процессы деформирования, описываемые соотношением

$$\omega + \frac{1}{q}\dot{\sigma} + \lambda\ddot{\sigma} = E\left(\varepsilon + \frac{1}{n}\dot{\varepsilon} + \mu\ddot{\varepsilon}\right) \tag{46}$$

также укладываются в рамки зависимости (1), если положить $de(\alpha) = ds(\alpha) = 0$ для всех $\alpha$, кроме 0,1 и 2, а для этих последних считать

$$de(0) = 1, \qquad de(1) = \frac{1}{q}, \qquad de(2) = \lambda,$$

$$ds(0) = E, \quad ds(1) = \frac{E}{n}, \quad ds(2) = E\mu$$

В заключение этой статьи мы даем решение уже рассмотренной в первом примере задачи о жидкости между параллельными плоскостями, когда жидкость деформируется на чистый сдвиг движением одной плоскости относительно другой. Зависимость между $\sigma$ и $\varepsilon$ предположим в виде (44) при некоторых определенных функциях инерции $e(\alpha)$ и $s(\alpha)$, заданных на некотором (ограниченном или неограниченном) множестве. После того, что было сказано выше по поводу определения этих функций, не нарушая общности рассуждений, можно продолжить область определения $e(\alpha)$ и $s(\alpha)$ на всю действительную ось $\alpha$.

Пусть плоскость $x = 0$ (фиг.1) неподвижна, плоскость $x = l$ движется в самой себе в направлении оси $y$ по заданному закону; пусть $\rho$ есть постоянная плотность жидкости.

Краевые условия задачи будут те же, что и раньше:

$$y(x,t)\big|_{x=0} = 0, \quad y(x,t)\big|_{x=l} = \varphi(t) \tag{47}$$

Начальные условия таковы

$$\frac{\partial^{\omega} y}{\partial t^{\omega}}\Big|_{t=0} = 0 \qquad (\omega \in M) \tag{48}$$

где $M$ – конечное или бесконечное, дискретное или континуальное, ограниченное или нет множество значений $\alpha$, определяемое свойствами $e(\alpha)$ и $s(\alpha)$, т.е. свойствами жидкости.

Функция $\varphi(t)$ должна быть задана в согласии с (48).

Как и в первом примере, уравнения движения

$$\frac{\partial \sigma}{\partial x} = \rho \frac{\partial^2 y}{\partial t^2} \qquad (49)$$

мы комбинируем с тем, которое получается из (44) дифференцированием по $x$. Имеем

$$\int_{-\infty}^{+\infty} \frac{\partial^\alpha}{\partial t^\alpha}\left(\frac{\partial \sigma}{\partial x}\right) de(\alpha) = \int_{-\infty}^{+\infty} \frac{\partial^\alpha}{\partial t^\alpha}\left(\frac{\partial \varepsilon}{\partial x}\right) ds(\alpha) \qquad (50)$$

Принимая во внимание, что $\partial\varepsilon/\partial x = \partial^2 y/\partial x^2$, получим из (49) уравнение движения в виде

$$\rho \int_{-\infty}^{+\infty} \frac{\partial^{2+\alpha} y}{\partial x^{2+\alpha}} ds(\alpha) = \int_{-\infty}^{+\infty} \frac{\partial^\alpha}{\partial t^\alpha}\left(\frac{\partial^2 y}{\partial x^2}\right) ds(\alpha) \qquad (51)$$

Если $Y(x,p) \overset{\bullet}{\underset{\bullet}{\leftarrow}} y(x,t)$, то в области изображений уравнению (51) соответствует

$$\frac{d^2 Y}{dx^2} \int_{-\infty}^{+\infty} p^\alpha ds(\alpha) = Y\rho p^2 \int_{-\infty}^{+\infty} p^\alpha de(\alpha) \qquad (52)$$

При этом функция $Y(x,p)$ должна удовлетворять условиям

$$Y(x,p)\big|_{x=0} = 0, \quad Y(x,p)\big|_{x=l} = \Phi(p) \quad (\Phi(p) \overset{\bullet}{\underset{\bullet}{\leftarrow}} \varphi(t)) \qquad (53)$$

Решение уравнения (52) при условиях (53) будет

$$Y(x,p) = \Phi(p)\frac{shApx}{shApl}$$

$$\left( A = \sqrt{\frac{J_e}{J_s}}, J_e = \int_{-\infty}^{+\infty} p^\alpha de(\alpha), J_s = \int_{-\infty}^{+\infty} p^\alpha ds(\alpha) \right) \qquad (54)$$

Определив теперь $K(\theta)$ как такую функцию, изображением которой является множитель при $\Phi(p)$ в (54) и пользуясь теоремой Бореля, найдем окончательно

$$y(x,t) = \frac{d}{dt}\int_0^t K(t-\theta)\varphi(\theta)d\theta$$

и задача будет решена.

*Поступила в редакцию 5 VI 1947*

# **9.**КИНЕТИКА ПРОЦЕССА ВЫТЯГИВАНИЯ.

## 1. Стационарный процесс

В основу этой теоретической работы берем факт, установленный известными опытами Л.Н. Гинзбурга, М.М. Слатинцева и О.М. Спесивцевой [1-3]. При стационарном процессе индивидуальные скорости отдельных волокон в границах вытяжного поля представлены почти исключительно только двумя значениями, причем каждое волокно меняет свою скорость мгновенно (и не больше одного раза) с меньшего значения на большее. Это положение мы будем называть в дальнейшем *принципом отсутствия промежуточных скоростей.*

Мы учитываем, что, строго говоря, ни один из двух зажимов аппарата не является совершенным, т.е. в заднем зажиме имеются волокна, уже получившие большую из двух возможных скоростей. Также в переднем зажиме (по движению ленты) имеются волокна, еще сохранившие меньшую скорость.

Кроме того, будет учтена способность волокон растягиваться. Мы показываем, как из этих положений, без привлечения каких-либо иных гипотез, методами классической механики подвижных сред строится вся кинетика процесса. Существенно отметить, что в развиваемой теории понятие длины волокна не играет той роли, которая ему обычно приписывается.

**1.***Общие замечания. Обозначения. Основные формулы.* Речь идет о простейшем вытяжном аппарате с двумя парами валиков. Растягиваемый волокнистый объект мы называем лентой. Ось $x$ направим по движению ленты. Начало координат О помещаем в осевой плоскости заднего зажима. Длину поля (разводку) обозначим $l$ и для простоты примем за единицу длины.

Обозначим через $v_0, v_1 (v_1 > v_0)$ – наблюдаемые значения скоростей отдельных волокон (считаем, что никаких иных скоростей волокна не имеют); через $v(0), v(1)$ – средние арифметические (по совокупности) скорости в плоскостях заднего и соответственно переднего зажимов.

Из-за наличия медленных волокон в плоскости переднего зажима

$$v(1) < v_1 \tag{1.1}$$

Из-за наличия быстрых волокон в плоскости заднего зажима:

$$v(0) > v_0 \tag{1.2}$$

Вместе с тем, по смыслу работы аппарата

$$v(1) > v(0) \tag{1.3}$$

Пусть $v(x)$ – средняя арифметическая скорость волокон, пронизывающих сечение поля с абсциссой $x, 0 \le x \le 1$.

Обозначим $q-$ секундный расход массы (масса волокнистого материала, проходящего через сечение поля в единицу времени – секунду). В стационарном процессе $q-$ одно и то же для всех сечений поля. Вместе с тем $q$ является линейной плотностью количества движения.

Обозначим через $\psi(s)ds-$ вероятность того, что в пучке волокон выбранное наудачу волокно имеет длину между $s$ и $s+ds$, при этом

$$\int_0^\infty \psi(s)ds = 1; \qquad (1.4)$$

$\Lambda-$ средняя линейная плотность массы волокна, обратно пропорциональная номеру волокна; $\lambda(x)$ – линейная плотность массы ленты в сечении поля; величина, обратно пропорциональная номеру ленты в этом сечении.

По определению линейной плотности $q$ количества движения

$$\lambda(x)v(x) = q \qquad (1.5)$$

Отношение

$$\frac{\lambda(x)}{\Lambda} = n(x) \qquad (1.6)$$

дает число волокон, пронизывающих сечение $x$.

На основании (1.5) –(1.6) имеем

$$n(x) = \frac{q}{\Lambda v(x)} \qquad (1.7)$$

Отношение

$$\frac{n(0)}{n(1)} = \frac{v(1)}{v(0)} = b > 1 \qquad (1.8)$$

есть вытяжка (фактическая). Ее нужно отличать от предельной вытяжки $B$, определяемой посредством предельных скоростей $v_0$ и $v_1$ отношением

$$B = \frac{v_1}{v_0} = \frac{n_0}{n_1} \qquad (n_0 = \frac{q}{\Lambda v_0}, n_1 = \frac{q}{\Lambda v_1}) \quad . \qquad (1.9)$$

Так как $v(1) < v_0$ и $v(0) > v_0$, то

$$b < B \qquad (1.10)$$

Фактическая вытяжка всегда меньше предельной вследствие наличия некоторого количества быстрых волокон в заднем и медленных – в переднем зажимах. Если бы $v_0$ была общей скоростью входа для всех волокон и $v_1-$ общей скоростью их выхода, то фактическая вытяжка совпала бы с предельной. В каждом сечении поля пронизывающие его волокна имеют одну из двух возможных скоростей: $v_0$ или $v_1$.

Пусть $n_0(x)$– число волокон, пронизывающих сечение $x$ и имеющих скорость $v_0$, а $n_1(x)$– число волокон, пронизывающих сечение $x$ и имеющих скорость $v_1$. Согласно (1.7)

$$n_0(x)+n_1(x)=n(x)=\frac{q}{\Lambda v(x)} \qquad (1.11)$$

С другой стороны, можно показать, что *средняя арифметическая скорость* $v(x)$ волокон, пронизывающих сечение $x$, есть не что иное, как *скорость центра массы этих волокон.* В самом деле, если $\psi(s)ds$ есть относительное число волокон с длиной между $s$ и $s+ds$ (оно не должно зависеть от скорости), то сечение $x$ пронизывается $n_0(x)\psi(s)ds$ волокнами с длиной $s$ и скоростью $v_0$. Их масса есть $\Lambda s n_0(x)\psi(s)ds$. То же сечение пронизывается $n_1(x)\psi(s)ds$ волокнами с длиной $s$ и скоростью $v_1$. Их масса есть $\Lambda s n_1(x)\psi(s)ds$.

Скорость центра массы волокон всевозможных длин, пронизывающих сечение $x$, определится как

$$\frac{v_0\Lambda n_0(x)\int_0^\infty s\psi(s)ds+v_1\Lambda n_1(x)\int_0^\infty s\psi(s)ds}{\Lambda n(x)\int_0^\infty s\psi(s)ds} \qquad (1.12)$$

что совпадает с

$$\frac{v_0 n_0(x)+v_1 n_1(x)}{n_0(x)+n_1(x)}=v(x)$$

Последнее выражение есть средняя арифметическая скорость по совокупности всех волокон, пронизывающих сечение $x$.

Из последнего выражения следует, что

$$v_0 n_0(x)+v_1 n_1(x)=v(x)n(x)=\frac{q}{\Lambda} \qquad (1.13)$$

(1.11) и (1.13) дают:

$$n_0(x)=\frac{q}{\Lambda(v_1-v_0)}\frac{v_1-v(x)}{v(x)}, \quad n_1(x)=\frac{q}{\Lambda(v_1-v_0)}\frac{v(x)-v_0}{v(x)} \qquad (1.14)$$

В частности, на входе в поле $x=0$ имеется $n_0(0)$ волокон со скоростью $v_0$ и $n_1(0)$ волокон со скоростью $v_1$

$$n_0(0)=\frac{q}{\Lambda(v_1-v_0)}\frac{v_1-v(0)}{v(0)}, \quad n_1(0)=\frac{q}{\Lambda(v_1-v_0)}\frac{v(0)-v_0}{v(0)} \qquad (1.15)$$

Аналогично, в переднем зажиме $x=1$ имеется $n_0(1)$ волокон со скоростью $v_0$ и $n_1(1)$ волокон со скоростью $v_1$

$$n_0(1)=\frac{q}{\Lambda(v_1-v_0)}\frac{v_1-v(1)}{v(1)},\ n_1(1)=\frac{q}{\Lambda(v_1-v_0)}\frac{v(1)-v_0}{v(1)} \qquad (1.16)$$

Отсюда вытекает, что фактическая вытяжка $b$ связана с предельной $B$ соотношением

$$b=\frac{n(0)}{n(1)}=\frac{n_0(0)+n_1(0)}{n_0(1)+n_1(1)}=\frac{n_0}{1+(B-1)\alpha}:\frac{n_0}{B-(B-1)\beta}=\frac{B-(B-1)\beta}{1+(B-1)\alpha},$$

где

$$\alpha=\frac{n_1(0)}{n(0)}=\frac{v(0)-v_0}{v_1-v_0} \qquad (0<\alpha<1) \qquad (1.17)$$

является мерой *проскальзывания волокон в заднем зажиме* и

$$\beta=\frac{n_0(1)}{n(1)}=\frac{v_1-v(1)}{v_1-v_0} \qquad (0<\beta<1) \qquad (1.18)$$

есть мера *проскальзывания волокон в переднем зажиме*.

Числа $\alpha$ и $\beta$ зависят от свойств обрабатываемого материала и от условий процесса. В большинстве случаев они лишь немного отличаются от нуля (особенно $\beta$). Вообще имеет место соотношение $\alpha+\beta<1$.

Подобно тому, как это было сделано для средней арифметической скорости, средняя квадратичная скорость $u(x)$ может быть представлена в виде:

$$[u(x)]^2=\frac{\frac{1}{2}v_0^2\Lambda n_0(x)\int s\psi(s)ds+\frac{1}{2}v_1^2\Lambda n_1(x)\int s\psi(s)ds}{\frac{1}{2}\Lambda n(x)\int s\psi(s)ds}$$

не отличающемся от

$$\frac{v_0^2 n_0(x)+v_1^2 n_1(x)}{n_0(x)+n_1(x)}=[u(x)]^2 \qquad (1.19)$$

Подставляя сюда $n_0(x)$ и $n_1(x)$ из (1.14), получим

$$[u(x)]^2=(v_0+v_1)v(x)-v_0v_1 \qquad (1.20)$$

Откуда следует, что

$$[u(x)]^2-[v(x)]^2=[v_1-v(x)][v(x)-v_0]>0 \qquad (1.21)$$

Для отношения квадратов обеих средних скоростей находим

$$\frac{[u(x)]^2}{[v(x)]^2}=\sigma(x)=\frac{v_0+v_1}{v(x)}-\frac{v_0v_1}{[v(x)]^2}$$

Оно достигает максимума в той точке $x^*$ поля, для которой средняя арифметическая скорость $v^*$ равна средней гармонической из предельных:

$$v = \frac{2v_0 v_1}{v_0 + v_1} \tag{1.22}$$

или, для которой средняя квадратичная скорость $u^*$ равна средней геометрической из предельных:

$$u = \sqrt{v_0 v_1}$$

Это максимальное отношение $\sigma_m$ получается равным

$$\sigma_m = \frac{(v_0 + v_1)^2}{4v_0 v_1} = \frac{(B+1)^2}{4B} \tag{1.23}$$

и зависит только от вытяжки. Например, для $B = 3$ имеем $\sigma_m = 1.7$; для $B = 100$ получается $\sigma_m = 25.6$.

Наконец, заметим, что из (1.14) получается

$$\frac{dn_0(x)}{dx} = \frac{v_1}{v_1 - v_0}\frac{dn}{dx}, \quad \frac{dn_1(x)}{dx} = -\frac{v_0}{v_1 - v_0}\frac{dn}{dx}$$

так что

$$-\frac{dn_0(x)}{dx} = +\frac{dn_1(x)}{dx}B \tag{1.24}$$

Градиент падения числа медленных волокон в $B$ раз больше градиента возрастания числа быстрых. Интегрируя (1.24), имеем

$$n_0(x) + Bn_1(x) = const$$

Но на входе в поле

$$n_0(0) + Bn_1(0) = \frac{q}{\Lambda v_0} = n_0$$

Поэтому во всех сечениях поля должно быть

$$n_0(x) + Bn_1(x) = n_0 \tag{1.25}$$

**2.***Диссипативная функция.* Легко видеть, что выражение

$$E(x) = \frac{\lambda(x)[u(x)]^2}{2} = \frac{q}{2v}[(v_0 + v_1)v - v_0 v_1] = \frac{qv_0}{2} + \frac{qv_1}{2} - \frac{qv_0 v_1}{2v(x)}$$

представляет линейную плотность кинетической энергии в сечении $x$ поля.

С другой стороны, выражение

$$E'(x) = \frac{\lambda(x)[v(x)]^2}{2} = \frac{qv(x)}{2}$$

дает ту проблематическую линейную плотность кинетической энергии в рассматриваемом сечении $x$, которую следовало бы приписать, если бы все волокна, пронизывающие это сечение, имели бы единую скорость, равную $v(x)$. Различие между обоими выражениями обусловлено наличием разброса собственных скоростей отдельных волокон по обоим предельным значениям $v_0$ и $v_1$. Величина

$$E(x) - E'(x) = D(x), \tag{2.1}$$

равная

$$\frac{\lambda\left(u^2 - v^2\right)}{2} = \frac{qv_0}{2} + \frac{qv_1}{2} - \frac{qv_0v_1}{2v} - \frac{qv}{2},$$

на основании (1.20) может быть представлена в виде

$$D(v) = \frac{q}{2}\frac{(v_1 - v)(v - v_0)}{v} \tag{2.2}$$

Разность $D(v)$ есть доля кинетической энергии в единице длины ленты, выбранной около сечения $x$ поля и обусловленной исключительно тем обстоятельством, что часть $n_0(x)$ волокон проходит через это сечение со скоростью $v_0$, тогда как другая их часть $n_1(x)$ имеет в этом сечении скорость $v_1$.

У входа в поле, где почти все волокна имеют скорость $v_0$, $D(v)$, будучи положительной величиной, близка к нулю. То же имеет место и вблизи выхода из поля, где почти нет медленных волокон. Но в некотором сечении $x^*$, лежащем где-то посредине поля, $D(v)$ достигает единственного максимума, обусловленного наибольшим разбросом скоростей волокон по обоим предельным значениям $v_0$ и $v_1$. Так как

$$\frac{dD(v)}{dx} = \frac{q}{2}\frac{v_0v_1 - v^2}{v^2}\frac{dv}{dx}$$

то этот максимум достигается в той точке поля $x^*$, где

$$v^2 = v_0v_1$$

т.е. там, где средняя арифметическая скорость $v(x)$ становится равной средней геометрической из скоростей $v_0$ и $v_1$ отдельных волокон. Функцию

$$D(v) = \frac{q}{2}\frac{(v_1 - v)(v - v_0)}{v}$$

мы называем *диссипативной функцией.*

Она может быть связана с вероятностями наличия в сечении $x$ как медленных, так и быстрых волокон. В самом деле, формулы (1.14):

$$n_0(x) = \frac{q}{\Lambda(v_1 - v_0)}\frac{v_1 - v(x)}{v(x)}, \quad n_1(x) = \frac{q}{\Lambda(v_1 - v_0)}\frac{v(x) - v_0}{v(x)} \tag{2.3}$$

дающие числа волокон, пронизывающих сечение $x$ и идущих со скоростями, соответственно равными $v_0$ и $v_1$, могут быть переписаны в виде:

$$\frac{n_0(x)}{n(x)} = \frac{v_1 - v(x)}{v_1 - v_0} = \frac{n_0}{n_0 - n_1}\frac{n(x) - n_1}{n(x)} = \frac{B}{B-1}\frac{n(x) - n_1}{n(x)}$$

$$\frac{n_1(x)}{n(x)} = \frac{v(x) - v_0}{v_1 - v_0} = \frac{n_1}{n_0 - n_1}\frac{n_0 - n(x)}{n(x)} = \frac{1}{B-1}\frac{n_0 - n(x)}{n(x)} \tag{2.4}$$

Левые части этих равенств дают условные вероятности того, что через сечение $x$ проходит $n_0(x)$ медленных волокон и соответственно $n_1(x)$ быстрых в предположении, что полное число волокон в сечении $x$ есть $n(x)$. Обозначая эти вероятности соответственно $p_0(x)$ и $p_1(x)$, имеем

$$\frac{n(x)-n_1}{n(x)} = \frac{B-1}{B} p_0(x), \qquad \frac{n_0 - n(x)}{n(x)} = (B-1) p_1(x) \tag{2.5}$$

С другой стороны, диссипативная функция $D(v)$ может быть выражена не через скорости, а через числа волокон. Именно, так как

$$v = \frac{q}{\Lambda n},$$

то

$$D(v) = \frac{q}{2}\frac{(v_1 - v)(v - v_0)}{v} = \frac{q^2}{2\Lambda n_0 n_1} n \frac{(n-n_1)}{n}\frac{(n_0 - n)}{n} \tag{2.6}$$

Перемножая почленно оба равенства (2.5), получим:

$$\frac{(n-n_1)(n_0-n)}{n^2} = \frac{(B-1)^2}{B} p_0(x) p_1(x) \quad ,$$

вследствие чего (2.6) принимает вид

$$D(v) = \frac{q^2}{2\Lambda n_0 n_1} n \frac{(n-n_1)}{n}\frac{(n_0-n)}{n} = \frac{(B-1)^2}{B}\frac{q^2}{2\Lambda n_0 n_1} n(x) p_0(x) p_1(x) \tag{2.7}$$

Таким образом, функции $D(v)$ при прочих равных условиях пропорциональна полному числу волокон в данном сечении и обеим вероятностям наличия среди них как медленных, так и быстрых.

*3.Растягивающая сила.* Как говорилось выше, $v(x)$ есть скорость центра массы волокон, пронизывающих сечение $x$ поля и $q = \lambda(x)v(x)-$ постоянная (в стационарном процессе) линейная плотность количества движения или секундный расход массы.

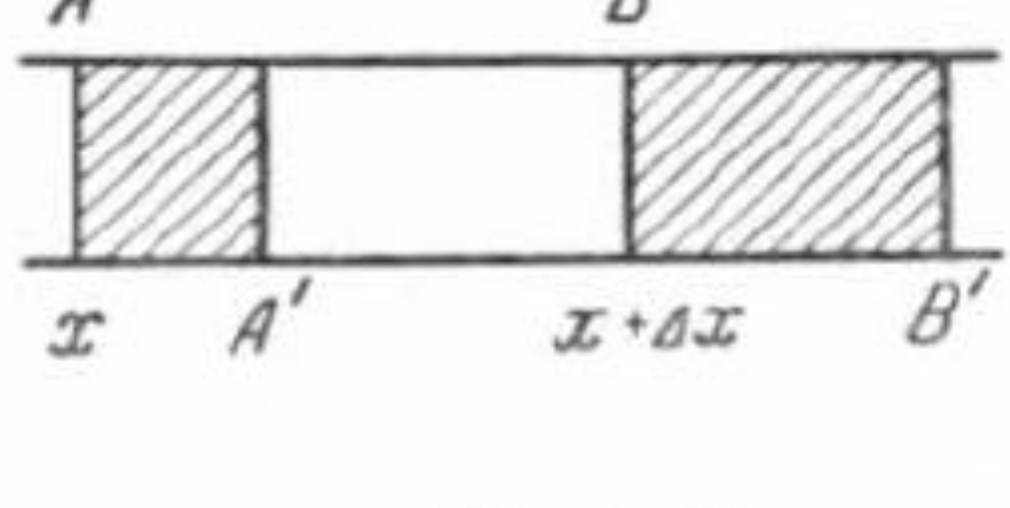

Фиг. 1

Рассмотрим маленький отрезок $AB$ ленты между абсциссами $x$ и $x+dx$ (фиг.1).

Если за время $\Delta t$ центр массы, лежавший в $A$, переместится в сечение $A'$, то центр массы, который находился в начале этого промежутка времени в сечении $B$, переместится в сечение $B'$. При поступательном движении рассматриваемого элемента со скоростью $v(x)$ оба перемещения $AA'$ и $BB'$ центров массы были бы одинаковы. Но вследствие перераспределения скоростей волокон по значениям $v_0$ и $v_1$ этого не будет и $BB'$ окажется больше $AA'$. В самом деле,

$$AA' = v(x)\Delta t, \qquad BB' = v(x+\Delta x)\Delta t$$

и, следовательно, рассматриваемый отрезок длины ленты, в начале промежутка времени равный $\Delta x$, по прошествии времени $\Delta t$ станет равным

$$AA' = AB - AA' + BB' = \Delta x - v(x)\Delta t + v(x+\Delta x)\Delta t$$

и получит прирост длины, равный $v(x+\Delta x)\Delta t - v(x)\Delta t$, что с точностью до малых первого порядка относительно $\Delta x$ равно $(dv/dx)\Delta x\Delta t$.

Так как каждой единице длины ленты соответствует количество движения $q$, то этому приросту длины рассматриваемого элемента соответствует прирост в количестве движения, равный $q(dv/dx)\Delta x\Delta t$.

Прирост количества движения за единицу времени получается равным $q(dv/dx)\Delta x$, а будучи отнесен к единице длины ленты, составляет $q(dv/dx)$.

Но прирост количества движения за единицу времени равен движущей силе. Поэтому выражение

$$f = q\frac{dv}{dx} \qquad (3.1)$$

представляет линейную плотность силы, уравновешивающей в единице длины ленты сопротивление последней к растягиванию и работающей в процессе вытягивания против сил взаимного сцепления между волокнами.

Выражение (3.1) будем называть *линейной плотностью растягивающей силы* или просто *растягивающей силой.* Легко видеть, что растягивающая сила на всей длине поля от зажима до зажима равна

$$F = f_1 - f_0 = \int_0^1 q\frac{dv}{dx}dx = qv_1 - qv_0 = qv_0(B-1)$$

Она пропорциональна избытку вытяжки над единицей и секундному расходу массы, как и должно быть.

При значении $B=1$ растягивающая сила была бы равна нулю и аппарат работал бы лишь на вовлечение новых масс в поступательное движение без растягивания.

*4. Уравнение процесса.* В нормально протекающем стационарном процессе вытягивания между растягивающей силой и диссипативной функцией должна быть функциональная зависимость. В наиболее общем случае плотность $f$ растягивающей силы, определенная посредством (3.1), представляется в виде степенного ряда по диссипативной функции, выражение для которой дано в (2.2):

$$f = a_0 + a_1 D + a_2 D^2 + \ldots \qquad (4.1)$$

Здесь коэффициенты $a_0, a_1, \ldots$, вообще говоря, зависят от абсциссы $x$. Однако, как показывает исследование этого вопроса, достаточно близкая к действительности картина получается уже в том простом частном случае, когда, во-первых, ограничиваются лишь двумя первыми членами вышенаписанного ряда и, во-вторых, считают коэффициенты положительными постоянными. Мы полагаем

$$f = r + \delta D \tag{4.2}$$

где $r$ и $\delta$ – неотрицательные постоянные. Они зависят от условий процесса и от внутренних свойств обрабатываемого материала, в частности от модальной длины волокна.

Подставляя в (4.2) вместо $f$ и $D$ их выражения (3.1) и (2.2), получим дифференциальное уравнение процесса в виде:

$$q\frac{dv}{dx} = r + \frac{q\delta}{2}\frac{(v_1 - v)(v - v_0)}{v}$$

или

$$\frac{dv}{dx} = -\frac{\delta}{2}\frac{1}{v}\left[v^2 - 2\left(\frac{r}{\delta q} + \frac{v_0 + v_1}{2}\right)v + v_0 v_1\right]$$

Здесь переменные легко разделяются:

$$\frac{\delta}{2}dx = -\frac{vdv}{v^2 - 2\left(\frac{r}{\delta q} + \frac{v_0 + v_1}{2}\right)v + v_0 v_1} \tag{4.3}$$

Корни трехчлена в знаменателе $w_0$ и $w_1$

$$w_0 = \frac{r}{\delta q} + \frac{v_0 + v_1}{2} - \sqrt{\left(\frac{r}{\delta q} + \frac{v_0 + v_1}{2}\right)^2 - v_0 v_1}$$

(4.4)

$$w_1 = \frac{r}{\delta q} + \frac{v_0 + v_1}{2} + \sqrt{\left(\frac{r}{\delta q} + \frac{v_0 + v_1}{2}\right)^2 - v_0 v_1}$$

оказываются положительными. В частности, при $r = 0$ мы имели бы $w_0 = v_0$, $w_1 = v_1$. С ростом $r$ от нуля $w_1$ монотонно возрастает, так как

$$\frac{dw_1}{dr} = \frac{1}{\delta q} + \frac{\left(\frac{r}{\delta q} + \frac{v_0 + v_1}{2}\right)\frac{1}{\delta q}}{\sqrt{\left(\frac{r}{\delta q} + \frac{v_0 + v_1}{2}\right)^2 - v_0 v_1}} > 0$$

А так как $w_0 w_1 = v_0 v_1$, то $w_0$ монотонно убывает с ростом $r$. Уравнение процесса (4.3) можно теперь представить в виде:

$$\frac{\delta}{2}dx = -\frac{vdv}{(v - w_1)(v - w_0)} \quad , \tag{4.5}$$

или, развертывая в элементарные дроби:

$$\frac{\delta}{2}dx = -\frac{w_1}{w_1 - w_0}\frac{dv}{v - w_1} + \frac{w_0}{w_1 - w_0}\frac{dv}{v - w_0} \qquad (4.6)$$

Интегрируя, находим:

$$\frac{\delta}{2}x + C = -\frac{w_1}{w_1 - w_0}\ln(v - w_1) + \frac{w_0}{w_1 - w_0}\ln(v - w_0) \qquad (4.7)$$

где $C$ – произвольная постоянная.

При $x = 0$ должно быть $v = v(0)$. Поэтому

$$C = -\frac{w_1}{w_1 - w_0}\ln[v(0) - w_1] + \frac{w_0}{w_1 - w_0}\ln[v(0) - w_0]$$

Вычитая из (4.7), найдем

$$\frac{\delta}{2}x = \frac{w_1}{w_1 - w_0}\ln\frac{v(0) - w_1}{v - w_1} + \frac{w_0}{w_1 - w_0}\ln\frac{v - w_0}{v(0) - w_0} \qquad ([\delta] = \frac{1}{длина})$$

Так как $v(0)$, как и $v(x)$, меньше $v_1$ и, следовательно, меньше $w_1$, то этот результат удобнее писать в виде:

$$\frac{\delta}{2}x = \frac{w_1}{w_1 - w_0}\ln\frac{w_1 - v(0)}{w_1 - v} + \frac{w_0}{w_1 - w_0}\ln\frac{v - w_0}{v(0) - w_0} \qquad (4.8)$$

На выходе из поля ($x = 1$) имеем $v = v(1)$. Поэтому должно быть

$$\frac{\delta}{2} = \frac{w_1}{w_1 - w_0}\ln\frac{w_1 - v(0)}{w_1 - v(1)} + \frac{w_0}{w_1 - w_0}\ln\frac{v(1) - w_0}{v(0) - w_0} \qquad (4.9)$$

Разделив (4.8) на (4.9), а затем разделив числитель и знаменатель результата на $w_0/(w - w_0)$, получим

$$x = \frac{\ln\left[\left(\frac{w_1 - v(0)}{w_1 - v}\right)^{\varepsilon}\frac{v - w_0}{v(0) - w_0}\right]}{\ln\left[\left(\frac{w_1 - v(0)}{w_1 - v(1)}\right)^{\varepsilon}\frac{v(1) - w_0}{v(0) - w_0}\right]} \qquad (\frac{w_1}{w_0} = \varepsilon > 1) \qquad (4.10)$$

или, отсюда, потенцируя, получим

$$\frac{v - w_0}{(w_1 - v)^{\varepsilon}} = PQ^{x} \qquad (4.11)$$

где положительные постоянные $P$ и $Q$, имеют значения:

$$P = \frac{v(0) - w_0}{[w_1 - v(0)]^{\varepsilon}}, \quad Q = \frac{v(1) - w_0}{v(0) - w_0}\left(\frac{w_1 - v(0)}{w_1 - v(1)}\right)^{\varepsilon} \qquad (4.12)$$

Это и есть уравнение распределения средних арифметических скоростей по сечениям поля. В частности, если $r = 0$, так что

$$w_0 = v_0, \quad w_1 = v_1, \quad \varepsilon = \frac{v_1}{v_0} = B \qquad (4.13)$$

то

$$P = \frac{v(0) - v_0}{[v_1 - v(0)]^{B}}, \quad Q = \frac{v(1) - v_0}{v(0) - v_0}\left(\frac{v_1 - v(0)}{v_1 - v(1)}\right)^{B} \qquad (4.14)$$

Как увидим дальше, в этом случае растягивание волокон предполагается отсутствующим ($r = 0$), но тогда приходится признать наличие проскальзывания в обоих зажимах.

Для коэффициентов проскальзывания мы имели:

$$\alpha = \frac{n_1(0)}{n(0)} = \frac{v(0) - v_0}{v_1 - v_0}, \qquad \beta = \frac{n_0(1)}{n(1)} = \frac{v_1 - v(1)}{v_1 - v_0} \tag{4.15}$$

Отсюда следует

$$v(0) - v_0 = (v_1 - v_0)\alpha, \qquad v(1) - v_0 = (v - v_0)(1 - \beta)$$

$$v_1 - v(0) = (v_1 - v_0)(1 - \alpha), \qquad v_1 - v(1) = (v_1 - v_0)\beta$$

Поэтому

$$P = \frac{(v_1 - v_0)\alpha}{(v_1 - v_0)^B (1-\alpha)^B} = \frac{\alpha}{v_0^{B-1}(B-1)^{B-1}(1-\alpha)^B},\ Q = \frac{(1-\alpha)^B (1-\beta)}{\alpha\beta^B} \tag{4.16}$$

и уравнение распределения скоростей принимает вид:

$$v_0^{B-1}(B-1)^{B-1}\frac{v - v_0}{(v_1 - v)^B} = \frac{\alpha}{(1-\alpha)^B}\left[\frac{(1-\alpha)^B(1-\beta)}{\alpha\beta^B}\right]^x \tag{4.17}$$

Таким образом, имеющее место растягивание отдельных волокон формально можно свести к проскальзыванию волокон в зажимах.

*5.Обсуждение результатов. Частные случаи.* Прежде всего заметим, что развиваемая теория процесса содержит следующий частный случай. Представим себе, что вытягиваемая лента такова, что в каждом сечении поля $x$ все волокна имеют единую скорость $v(x)$, меняющуюся вдоль поля от значения $v_0$ до $v_1$. Предполагается, следовательно, что никакого разброса скоростей по двум значениям $v_0$ и $v_1$ нет. Это имеет место в том случае, когда лента ведет себя в поле подобно растягиваемому резиновому шнуру, поступающему в поле со скоростью $v_0$ и выходящему с большей скоростью $v_1$. Тогда профиль скорости в каждом сечении поля есть плоскость. Волокна растягиваются, но не проскальзывают одно по другому. Этот род вытяжки Васильев назвал «вытяжкой первого рода».

Тогда растягивающая сила на единицу длины, или плотность растягивающей силы постоянна во всех сечениях поля и мы имеем дифференциальное уравнение процесса в виде:

$$q\frac{dv}{dx} = r \tag{5.1}$$

Отсюда получается:

$$v(x) = v_0\left(B-1\right)x + v_0\big) \qquad , \qquad (5.\ 2)$$

причем

$$r = q\left(v_1 - v_0\right) = qv_0\left(B-1\right) \qquad (5.\ 3)$$

Последнее является не чем иным, как выражением закона Гука: деформирующая сила $r$ прямо пропорциональна деформации $B-1$. Во-вторых, заметим, что принятое общее уравнение процесса:

$$f = r + \delta D \qquad (5.4)$$

может быть представлено в каждом из двух следующих видов:

(1) либо $$\frac{dv}{dx} = \frac{r}{q} + \frac{\delta}{2}\frac{\left(v_1 - v\right)\left(v - v_0\right)}{v} \qquad (5.5)$$

где $v_0$ и $v_1$ – наблюдаемые скорости волокон

(2) либо

$$\frac{dv}{dx} = \frac{\delta}{2}\frac{\left(w_1 - v\right)\left(v - w_0\right)}{v} \qquad (5.6)$$

где «условные» скорости $w_0$ и $w_1$ связаны с наблюдаемыми $v_0$ и $v_1$ так:

$$\frac{w_0 + w_1}{2} = \frac{r}{\delta q} + \frac{v_0 + v_1}{2},\ w_0 w_1 = v_0 v_1,\ (w_0 < v_0, w_1 > v_1) \qquad (5.7)$$

При помощи последних соотношений всегда можно перейти от двучленного уравнения (5.6) к трехчленному (5.5) и обратно.

При $r = 0$ из (5.7) следует

$$w_0 = v_0, \qquad w_1 = v_1 \qquad , \qquad (5.\ 8)$$

так что оба вида уравнения процесса перестают отличаться друг от друга. Это имеет место тогда, когда волокна считаются нерастяжимыми и вытягивание ленты происходит исключительно за счет проскальзывания волокон друг по другу. Существенно отметить, что появление в этом случае логарифмических особенностей на обоих концах поля заставляет считать тогда, что $v(0) > v_0$ и соответственно $v(1) < v_1$.

Наконец, в случае растягивающихся волокон ($r > 0$) уравнение процесса (5.4), конечно, удобнее писать в виде (5.6):

$$\frac{dv}{dx} = \frac{\delta}{2}\frac{\left(w_1 - v\right)\left(v - w_0\right)}{v}$$

отличающемуся от

$$\frac{dv}{dx} = \frac{\delta}{2}\frac{\left(v_1 - v\right)\left(v - v_0\right)}{v} \qquad (5.9)$$

только тем, что в предыдущем условии вместо наблюдаемых скоростей $v_0$ и $v_1$ фигурируют на тех же правах «условные» скорости $w_0$ и $w_1$, вычисляемые посредством (4.4) по $r, \delta, v_0$ и $v_1$. При этом

$$w_0 < v_0, \qquad w_1 > v_1 \qquad (5.10)$$

Вместо растяжения, обусловленного $r > 0$, формально вводится проскальзывание в зажимах при неизменной длине волокон, и коэффициенты $\alpha$ и $\beta$ этого фиктивного проскальзывания вычисляются по формулам:

$$\alpha = \frac{v_0 - w_0}{w_1 - w_0}, \qquad \beta = \frac{w_1 - v_1}{w_1 - w_0} \tag{5.11}$$

Таким образом, дело сводится к тому же дифференциальному уравнению процесса для вытягивания ленты из неизменяемых волокон

$$\frac{dv}{dx} = \frac{\delta}{2}\frac{(v_1 - v)(v - v_0)}{v} \tag{5.12}$$

при наличии проскальзывания в обоих зажимах с известными коэффициентами $\alpha$ и $\beta$, определяемыми посредством (5.11), где $w_0$ и $w_1$ играют роль наблюдаемых скоростей (они вычисляются по (4.4)), а $v_0$ и $v_1$ являются скоростями $v(0)$ и $v(1)$ на входе и на выходе.

Решение этого уравнения при данных граничных условиях легко выполняется по способу, описанному в § 4. В дальнейшем мы будем считать процесс вытягивания сопровождающимся проскальзыванием волокон в обоих зажимах и сохраним обозначения $v_0$ и $v_1$ для наблюдаемых скоростей $v(0)$ и $v(1)$ – для средних арифметических скоростей на входе и соответственно на выходе. Коэффициенты проскальзывания $\alpha$ и $\beta$ будем считать известными. Уравнение кривой распределения средних скоростей $v(x)$ по сечениям $x$ поля в конечном виде будем писать так:

$$[v_0(B-1)]^{B-1}\frac{v - v_0}{(v_1 - v)^B} = \frac{\alpha}{(1-\alpha)^B}\left[\frac{(1-\alpha)^B(1-\beta)}{\alpha\beta^B}\right]^x \tag{5.13}$$

Величина, обратная $\delta$, имеет размерность длины. Легко видеть, что $\delta$ зависит только от вытяжки (предельной) $B$ и от обоих коэффициентов проскальзывания $\alpha$ и $\beta$; а именно:

$$\frac{1}{\delta} = S = \frac{B-1}{2\ln\left[(1-\alpha)^B(1-\beta)/(\alpha\beta^B)\right]} \tag{5.14}$$

Дифференциальное уравнение процесса (5.13) выглядит при этом так:

$$\frac{B-1}{2\ln\left[(1-\alpha)^B(1-\beta)/(\alpha\beta^B)\right]}\frac{dv}{dx} = \frac{(v_1 - v)(v - v_0)}{v} \tag{5.15}$$

Из него видно, что $dv/dx > 0$; следовательно, $v$ монотонно растет в поле от значения $v(0)$ на входе до $v(1)$ на выходе. Градиент скорости $dv/dx$, будучи положительным, сначала растет, достигает наибольшей величины в середине поля и затем уменьшается у выхода из поля, если только единственная точка перегиба кривой распределения скоростей, имеющая абсциссу $x^*$, для которой

$$v^{*2} = v_0 v_1 \tag{5.16}$$

попадает между обоими зажимами $x = 0$ и $x = 1$. Далее мы рассматриваем различные частные случаи видов кривой распределения скоростей.

*6. Кривая утонения.* Подставляя в (5.15) вместо $v$ его выражение через полное число $n$ волокон в данном сечении поля $v = q/(\Lambda n)$ и делая аналогичную замену для $v_0$ и $v_1$, получим дифференциальное уравнение

$$-\frac{(B-1)n_0^2}{B\ln\left[(1-\alpha)^B(1-\beta)\alpha\beta^{-B}\right]}\frac{dn}{dx} = \left[(n_0+n_1)n^2 - n^3 - n_0 n_1 n\right] \qquad (6.1)$$

распределения числа волокон $n(x)$ по сечениям поля или «кривой утонения». Правая часть в (6.1) неотрицательна, как и множитель при $dn/dx$ в левой части. Следовательно, $n(x)$ монотонно убывает с ростом $x$.

Дифференцируя (6.1) и пользуясь (5.16), найдем

$$\frac{d^2n}{dx^2} = \left[\frac{B}{(n-1)n_0^2}\ln\frac{(1-\alpha)^B(1-\beta)}{\alpha\beta^B}\right]^2 (n_0-n)(n-n_1)n\left[2(n_0+n_1)n - 3n^2 - n_0n_1\right] \qquad (6.2)$$

Корни трехчлена $-n^2 + \frac{2}{3}(n_0+n_1)n - \frac{1}{3}n_0n_1$ таковы

$$n_0^* = \tfrac{1}{3}\left[n_0 + n_1 + \sqrt{n_0^2+n_1^2-n_0n_1}\right],\ n_1^* = \tfrac{1}{3}\left[n_0 + n_1 - \sqrt{n_0^2+n_1^2-n_0n_1}\right] \qquad (6.3)$$

или

$$n_0^* = \tfrac{1}{3}n_1\left[B+1+\sqrt{B^2-B+1}\right], \qquad n_1^* = \tfrac{1}{3}n_1\left[B+1-\sqrt{B^2-B+1}\right]$$

Замечая, что

$$B-1 < \sqrt{B^2-B+1} < B$$

можем положить

$$\sqrt{B^2-B+1} = B - z, \qquad z = B - \sqrt{B^2-B+1} \quad (0<z<1) \qquad (6.4)$$

Тогда эти корни представятся в виде

$$n_0^* = \tfrac{1}{3}n_1[2B+1-z], \qquad n_1^* = \tfrac{1}{3}n_1[1+z] \qquad (6.5)$$

Так как коэффициент при $n^2$ написанного выше трехчлена есть -1, то изображающая его парабола обращена вершиной вправо от вертикальной, смотрящей вверх оси $n$. При $n > n_0^*$ или $n < n_1^*$ этот трехчлен, а значит, и $d^2n/dx^2$ отрицателен, при $n_0^* > n > n_1^*$ – положителен. Следовательно, кривая утонения монотонно спускаясь к оси $x$ абсцисс, обращена выпуклостью вправо (по направлению оси $x$), пока $n$ остается больше $n_0^*$ или становится меньше $n_1^*$. Между корнями $n_0^*$ и $n_1^*$ она выпукла влево. Вообще говоря, она может иметь две точки перегиба, ординаты которых равны $n_0^*$ и $n_1^*$. Однако, как сейчас увидим, в пределах вытяжного поля может оказаться не более какой-либо одной из них. В самом деле, теоретически возможны следующие комбинации:

1) $n_0^* < n(1)$, так что и подавно $n_1^* < n(1)$.–Обе точки перегиба лежат справа от выхода из поля, $n > n_0^*$ над всем полем, кривая утонения выпукла вправо.

2) $n(0) > n_0^* > n(1)$, $n_1^* < n(1)$.– Лишь $n_0^*$ лежит в поле, над первой частью поля кривая утонения выпукла вправо, у выхода– влево.

3) $n(0) > n_0^* > n_1^* > n(1)$.– Обе точки перегиба лежат в поле.

4) $n_0^* > n(0), n(1) < n_1^* < n(0)$.– В поле лежит лишь $n_1^*$.

5) $n_1^* > n(0)$, так что и подавно $n_1^* > n(1)$.–Обе точки перегиба лежат вне поля, перед входом.

Из этих пяти комбинаций реальны только две первые.

Действительно, если $n_1^* > n(1)$, то, принимая во внимание, что

$$n_1^* = \frac{n_1}{3}(1+z),\ n(1) = \frac{n_0}{B-(B-1)\beta} = \frac{B}{B-(B-1)\beta} n_1 \quad ,$$

мы имеем

$$\frac{1+z}{3} > \frac{B}{B-(B-1)\beta}$$

Так как $z$ –правильная дробь, то $\frac{1}{3}(1+z) \le \frac{2}{3}$ при всяком $z$. В то же время наименьшее значение для

$$\frac{B}{B-(B-1)\beta}$$

есть единица.

Поэтому из перечисленных комбинаций нужно отбросить все те, для которых $n_1^* > n(1)$. Остаются следующие две:

1) $n_0^* < n(1)$.– Обе точки перегиба– вне поля, за выходом из него.
$d^2n/dx^2 < 0$. Кривая утонения выпукла вправо над всем полем (фиг.2а).

2) $n_1^* < n(1)$, $n(0) > n_0^* > n(1)$. Верхняя точка перегиба $(n_0^*)$ лежит в поле, нижняя – вне поля справа. Первая часть кривой утонения выпукла вправо, вторая – влево (фиг. 2,б).

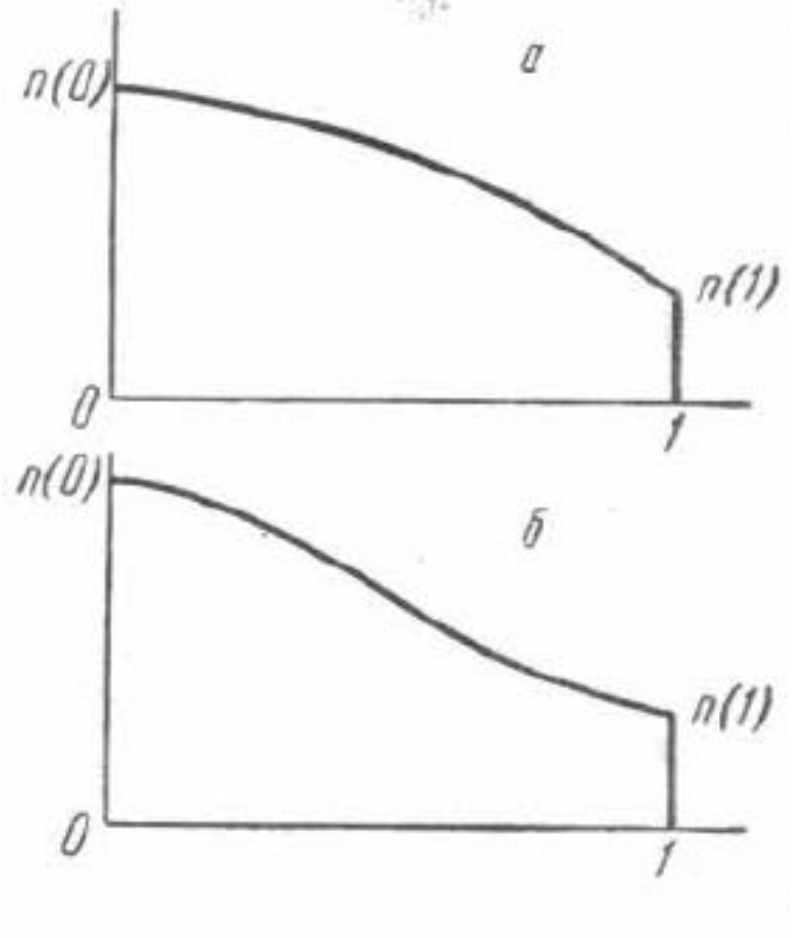

Фиг. 2

Первый случай реален только при достаточно большом коэффициенте проскальзывания $\beta$ на выходе из поля.

В самом деле, условие этого случая: $n_0^* < n(1)$ равноценно следующему:

$$n_0^* = \frac{n_1}{3}(2B+1-z) < n(1) = n_1 \frac{B}{B-(B-1)\beta}$$

или

$$2B+1-z < \frac{3B}{B-(B-1)\beta}$$

Отсюда получается, что кривая утонения без точки перегиба может иметь место только в случае, если

$$\beta > \frac{B}{B-1}\frac{B-2+\sqrt{B^2-B+1}}{B+1+\sqrt{B^2-B+1}} \qquad (6.6)$$

т.е. при достаточно больших коэффициентах проскальзывания $\beta$ на выходе, или, что то же, при достаточно малых коэффициентах проскальзывания $\alpha$ на входе, потому что всегда должно иметь место неравенство $\alpha + \beta < 1$. Например, при $B = 2$ кривая утонения не имеет точки перегиба, если $\beta > 0.36$; для $B = 6$ должно быть $\beta > 0.93$, так что на долю $\alpha$ приходится во втором случае меньше чем 0.07.

Кривая утонения с точкой перегиба получается, как легко видеть, при одновременном соблюдении неравенств

$$1 + (B-1)\alpha < \frac{3B}{B+1+\sqrt{B^2-B+1}}$$

$$B - (B-1)\beta > \frac{3B}{B+1+\sqrt{B^2-B+1}}$$

или

$$\alpha < \frac{2B-1-\sqrt{B^2-B+1}}{(B-1)[B+1+\sqrt{B^2-B+1}]}, \beta < \frac{B}{B-1}\frac{B-2+\sqrt{B^2-B+1}}{B+1+\sqrt{B^2-B+1}} \qquad (6.7)$$

т.е. при достаточно малом проскальзывании как на входе, так и на выходе. Поэтому вероятность возникновения случая перегиба кривой утонения больше, чем противоположного.

При наличии лишь двух значений скорости у отдельных волокон вопрос о проскальзывании в зажимах тесно связан с вопросом о том, в какой именно момент волокно меняет свою скорость и какой из двух зажимов оказывается «сильнее», выражаясь словами Зотикова. Коэффициенты проскальзывания могут, на наш взгляд, служить количественной характеристикой этих двух обстоятельств. Ниже мы приводим пять примеров расчета кривых распределения скорости и кривых утонения.

*7. Примеры.*

*Пример 1*. $v_0 = 1, B = 6, \alpha = 0.01, \beta = 0.3$

Вычисления по соответствующим формулам дают

$$v(0) = 1.05,\quad v(1) = 4.50,\ b = \frac{v(1)}{v(0)} \approx 4.5,\qquad S = 0.2\ ,\ n_0 = 10000.$$

Кривая распределения скоростей (фиг.3)

$$x = 1.103 + \frac{\lg(v-1) - 6\lg(6-v)}{4956}$$

| | | | | | | | | | |
|---|---|---|---|---|---|---|---|---|---|
| $v$ | 1.05 | 1.20 | 1.50 | 2.00 | 2.50 | 3.00 | 3.50 | 4.00 | 4.50 |
| $x$ | 0.000 | 0.139 | 0.253 | 0.374 | 0.480 | 0.586 | 0.703 | 0.835 | 1.000 |
| $n$ | 9524 | 8000 | 6667 | 5000 | 4000 | 3333 | 2857 | 2500 | 2222 |

Ордината точки перегиба кривой распределения скоростей равна $v^* \approx 2.4$

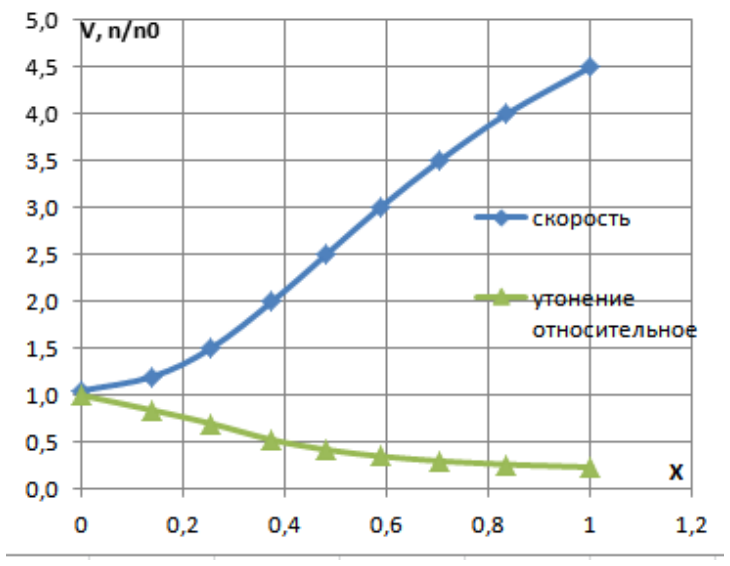

Фиг. 3

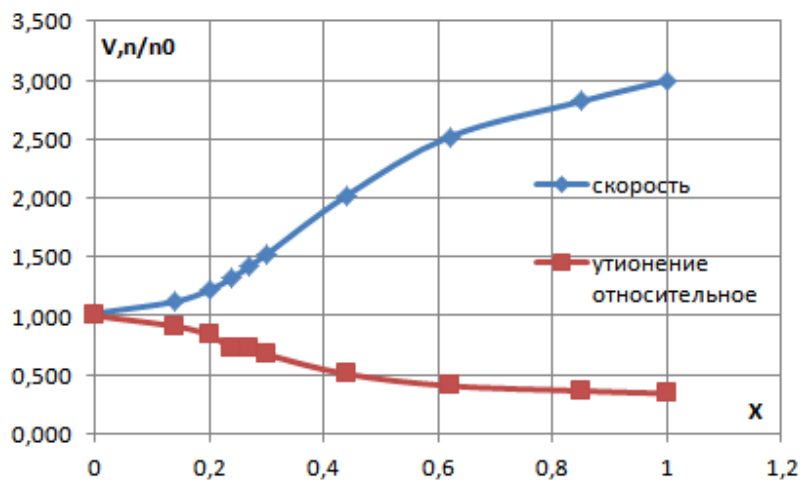

Фиг.4

*<u>Пример 2</u>*. $v_0 = 1, B = 3, \alpha = 0.1, \beta = 0.01$. Вычисления дают

$$v(0) = 1.02, \quad v(1) = 2.98, \quad b = 2.96, \quad S = 0.055 \quad , \quad n_0 = 10000.$$

Кривая распределения скоростей (фиг.4)

$$x = \frac{\lg 400}{\lg 970300} + \frac{\lg(v-1) - 3\lg(3-v)}{\lg 970300}$$

| | | | | | | | | | | |
|---|---|---|---|---|---|---|---|---|---|---|
| $v$ | 1.02 | 1.12 | 1.22 | 1.32 | 1.42 | 1.52 | 2.02 | 2.52 | 2.82 | 3.00 |
| $x$ | 0.000 | 0.14 | 0.20 | 0.24 | 0.27 | 0.30 | 0.44 | 0.62 | 0.85 | 1.000 |
| $n$ | 9804 | 8928 | 8197 | 7076 | 7042 | 6579 | 4950 | 3968 | 3546 | 3333 |

*<u>Пример 3</u>*. $v_0 = 1, B = 6, \alpha = 0.01, \beta = 0.01$ Вычисления дают

$$v(0) = 1.05, \quad v(1) = 5.95, \quad b = 5.9, \quad S = 0.071 \quad , \quad n_0 = 10000.$$

Кривая распределения скоростей (фиг.5)

$$x = 0.391 + \frac{\lg(v-1) - 6\lg(6-v)}{14}$$

| | | | | | | | | | | | | |
|---|---|---|---|---|---|---|---|---|---|---|---|---|
| $v$ | 1.05 | 1.20 | 1.50 | 2.00 | 2.50 | 3.00 | 3.50 | 4.00 | 4.50 | 5.00 | 5.50 | 5.95 |
| $x$ | 0.00 | 0.05 | 0.09 | 0.13 | 0.17 | 0.21 | 0.25 | 0.29 | 0.35 | 0.43 | 0.57 | 1.00 |
| $n$ | 9524 | 8333 | 6667 | 5000 | 4000 | 3000 | 2857 | 2500 | 2222 | 2000 | 1818 | 1680 |

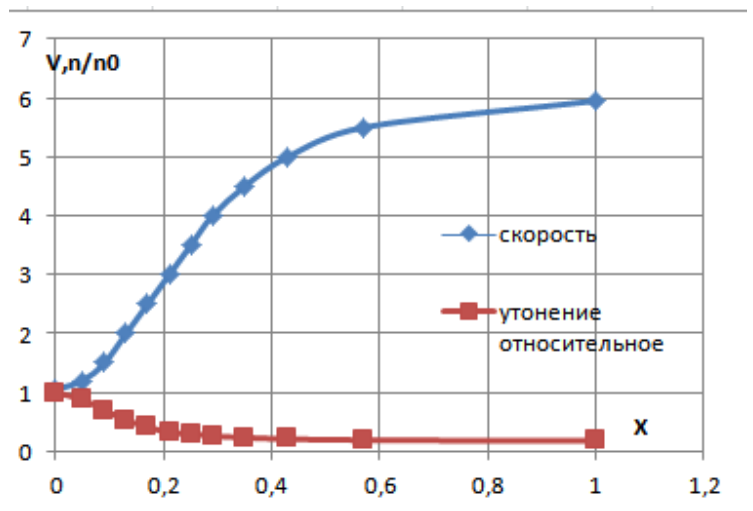

Фиг.5

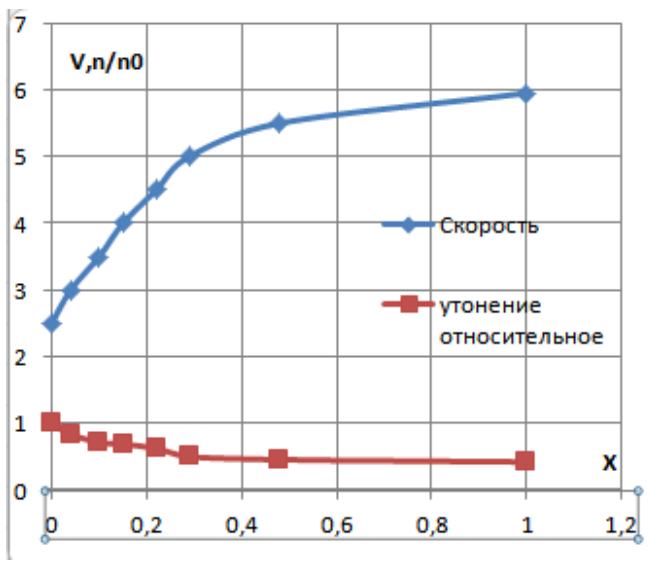

Фиг.6

*<u>Пример 4</u>*. $v_0 = 1, B = 6, \alpha = 0.3, \beta = 0.01$ Вычисления дают

$$v(0) = 2.5, \quad v(1) = 5.95, \quad b = \frac{v(1)}{v(0)} = 2.38, \quad S = \frac{1}{\delta} = 0.086 \quad , \quad n_0 = 10000.$$

Кривая распределения скоростей (фиг.6) : $x = \frac{\lg 1226}{\lg(3.882 \cdot 10^{11})} + \frac{\lg(v-1) - 6\lg(6-v)}{\lg(3.882 \cdot 10^{11})}$

| $v$ | 2.50 | 3.00 | 3.50 | 4.00 | 4.50 | 5.00 | 5.50 | 5.95 |
|---|---|---|---|---|---|---|---|---|
| $x$ | 0.00 | 0.04 | 0.10 | 0.15 | 0.22 | 0.29 | 0.48 | 1.00 |
| $n$ | 4000 | 3333 | 2857 | 2750 | 2500 | 2000 | 1818 | 1680 |

Точка перегиба на кривой скоростей, перед входом в поле, $v^* = 2.38$.

*Пример 5*. $v_0 = 1, B = 6, \alpha = 0.5, \beta = 0.8$ Вычисления дают

$$v(0) = 1.25, \quad v(1) = 2.00, \quad b = \frac{v(1)}{v(0)} = 1.60, \quad S = \frac{1}{\delta} = 0.90 \quad , \quad n_0 = 10000.$$

Кривая распределения скоростей (фиг.7) $x = 4.329 + \frac{\lg(v-1) - 6\lg(6-v)}{1.077}$

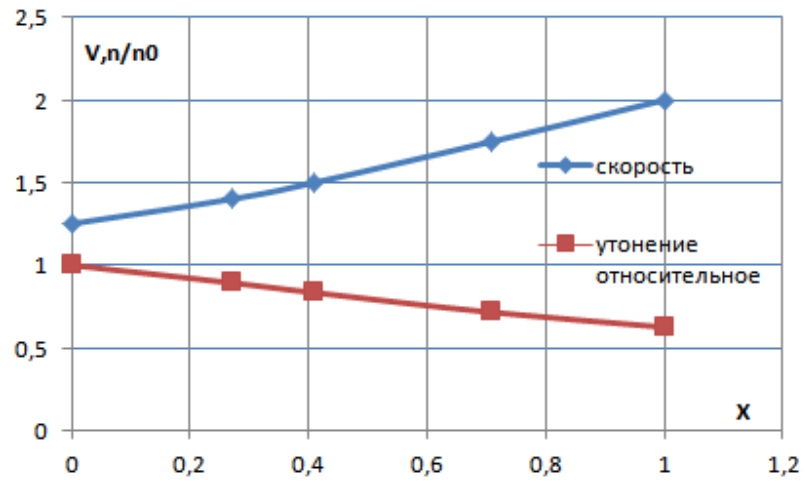

Фиг.7

| $v$ | 1.25 | 1.40 | 1.50 | 1.75 | 2.00 |
|---|---|---|---|---|---|
| $x$ | 0.00 | 0.27 | 0.41 | 0.71 | 1.00 |
| $n$ | 8000 | 7143 | 6667 | 5714 | 5000 |

Точка перегиба кривой распределения – за выходом из поля.

*Поступила 15 V 1956*

ЛИТЕРАТУРА

# **10.**КИНЕТИКА ПРОЦЕССА ВЫТЯГИВАНИЯ.

## II. Нестационарный процесс

В ранее опубликованной статье [1] изучался стационарный процесс вытягивания на основе опытно установленного принципа отсутствия промежуточных скоростей у отдельных волокон. Такой процесс никогда, строго говоря, не имеет места: неровнота ленты, органически связанная со структурой вытягивания объекта, является сопровождающим обстоятельством, в такой мере неизбежным, как и нежелательным. Для технологии процесса очень важно поэтому изучить явление с точки зрения возникающих систематических или случайных отклонений от стационарности.

Ниже рассматривается случай нестационарного процесса.

Как и раньше, принимается простейшая схема вытяжного аппарата из двух пар валиков, между которыми пропускается непрерывно бегущая лента. Начало координат O поместим в зажиме питающей (задней) пары, ось $Ox$ направим по движению ленты; длину поля, а значит, абсциссу зажима передней пары примем за единицу, которую будем обозначать в случае надобности буквой $l$.Пусть $\lambda$ обозначает среднюю линейную плотность массы одиночного волокна, являющуюся величиной, обратной номеру волокна, $v(0,t)$ – среднюю арифметическую по совокупности скорость ленты на входе в поле, $v(1,t)$ – аналогичную величину на выходе из поля. При этом $v(0,t)$ и $v(1,t)$ обозначают заданные условиями работы аппарата функции времени $t$,удовлетворяющие неравенствам

$$0 \le v(0,t) \le v(1,t), \qquad (t \ge 0)$$

$v(x,t)$ – среднюю (по совокупности) скорость в абсциссе $x$ в момент $t$; $\lambda(x,t)$ – среднюю линейную плотность массы ленты в сечении $x$ в момент $t$, являющуюся величиной, обратной номеру ленты в этом сечении в этот момент, $q(x,t)$ – линейную плотность количества движения в сечении $x$ в момент $t$, которая является одновременно мгновенным значением секундного расхода массы ленты в этом сечении в этот момент.

**1.***Уравнение неразрывности.* За время $dt$ в элемент поля $dx$ между сечениями $x$ и $x+dx$ сквозь сечение $x$ входит масса, равная $q(x,t)dt$. За то же время сквозь сечение $x+dx$ из элемента выходит масса $q(x+dx,t)dt$. Прирост массы в элементе за время $dt$ есть

$$q(x,t)dt - q(x+dx,t)dt = -\frac{\partial q}{\partial x}dxdt \qquad (q = \lambda v)\,.$$

С другой стороны, масса элемента $dx$, в момент $t$ представляемая в виде $\lambda(x,t)dx$, в момент $t+dt$ представится как $\lambda(x,t+dt)dx$. Ее прирост за время $dt$ есть

$$\lambda(x,t+dt)dx-\lambda(x,t)dx=\frac{\partial\lambda}{\partial t}dxdt$$

Сравнивая оба выражения, получим уравнение неразрывности

$$\frac{\partial(\lambda v)}{\partial x}=-\frac{\partial\lambda}{\partial t} \qquad (1.1)$$

**2.***Уравнение движения.* Будем исходить из общих соображений, разработанных И.В. Мещерским [2-3].

Пусть масса движущегося тела меняется со временем. Если в момент $t$ масса тела есть $M$, скорость его в этот момент есть $v$ и если за время $dt$ тело присоединяет массу $dM$, которая до присоединения имела скорость $u$, то, называя буквой $F$ движущую силу в момент $t$, имеем

$$F=\frac{d(Mv)}{dt}-u\frac{dM}{dt} \qquad (2.1)$$

Это соотношение выражает второй закон динамики для поступательного движения с учетом реакции вовлекаемых в движение масс. Рассмотрим отрезок ленты, в момент $t$ заключенный между абсциссами 0 и $x$. Его количество движения в этот момент есть

$$\int_0^x q(\xi,t)d\xi .$$

В момент $t+dt$ линейная плотность количества движения в сечении $\xi$ станет равной $q(\xi,t+dt)$, а передний конец рассматриваемого отрезка переместится за время $dt$ из сечения $x$ в сечение $x+dx=x+v(x,t)dt$ и новое значение количества движения станет

$$\int_0^{x+vdt} q(\xi,t+dt)d\xi$$

Разность между этим интегралом и предыдущим, написанная в несколько измененном виде

$$\int_0^{x+vdt} q(\xi,t+dt)d\xi-\int_0^{x} q(\xi,t+dt)d\xi+\int_0^{x} q(\xi,t+dt)d\xi-\int_0^{x} q(\xi,t)d\xi$$

по теореме о среднем равна

$$q(x,t)v(x,t)dt+dt\int_0^x\frac{\partial q(\xi,t)}{\partial t}d\xi . \qquad (2.2)$$

Она дает прирост количества движения рассматриваемого отрезка ленты за время $dt$. Соответствующий прирост за единицу времени есть

$$q(x,t)v(x,t)+\int_0^x \frac{\partial q(\xi,t)}{\partial t}d\xi$$

Выражение (2.2) соответствует члену $d(Mv)/dt$ в общем уравнении (2.1) Мещерского.

Второй же член $udM/dt$ в правой части (2.1) равен секундному приросту массы, умноженному на скорость этой массы до ее присоединения. В рассматриваемом случае секундный прирост массы на входе в поле есть $q(0,t)$, скорость его в момент $t$ вступления в поле есть $v(0,t)$. Значит, вместо $udM/dt$ в (2.1) нужно подставить величину $q(0,t)v(0,t)$. Это выражение получается прямо из (2.2) подстановкой $0$ вместо $x$.

Следовательно, движущая сила $F(x,t)$, действующая в момент $t$ на отрезок $(0,x)$ ленты, определяется равенством

$$F(x,t)=q(x,t)v(x,t)-q(0,t)v(0,t)+\int_0^x \frac{\partial q(\xi,t)}{\partial t}d\xi \qquad (2.3)$$

или, после дифференцирования по $x$

$$\frac{\partial[F-qv]}{\partial x}=\frac{\partial q}{\partial t}$$

Так как

$$q=\lambda v \qquad (2.4)$$

то вместо предыдущего можно написать также

$$\frac{\partial[F-\lambda v^2]}{\partial x}=\frac{\partial(\lambda v)}{\partial t} \qquad (2.5)$$

Уравнение (2.3) или (2.5) и выражает общий закон движения деформируемой ленты.

**3.** *Обсуждение полученных уравнений.* Вопрос о четвертом уравнении.
Система трех независимых уравнений (1.1), (2.3) и (2.4) связывает четыре функции: $F,q,\lambda,v$ переменных $x$ и $t$. Эта система должна быть дополнена четвертым уравнением, не зависящим от трех предыдущих и связывающим, вообще говоря, все четыре функции. Это четвертое уравнение можно было бы назвать в известном смысле уравнением состояния ленты. Оно должно отражать *связь между структурными свойствами объекта и параметрами процесса.*

Например, в стационарном случае вытягивания ровничной ленты можно было бы составить очень простое выражение для диссипативной функции, от которой зависит градиент силы $dF/dx$. В нестационарном процессе дело обстоит сложнее, хотя и здесь в принципе возможен путь через выражение для диссипативной функции.

В предлагаемой статье мы не делаем никаких гипотез о физических свойствах ленты и ограничиваемся рассмотрением лишь тех результатов, которые получаются из трех уравнений: (1.1), (2.3) и (2.4). Некоторые из этих результатов интересны. Например, в стационарном режиме $q$ не зависит не только от $t$, но и от $x$. Тогда (2.3) дает

$$F(x) = q[v(x) - v(0)] \quad .$$

Так как правая часть представляет величину, пропорциональную скорости деформации предельного растяжения ленты на отрезке $(0, x)$, то предыдущее соотношение является своеобразным выражением одного из известных законов деформирования.

Широко практикуемый в настоящее время метод радиоактивных изображений (например, изотопов фосфора) позволяет экспериментально определить каждую из трех функций: $\lambda(x,t), q(x,t), v(x,t)$ в отдельности. Таков, например, аппарат Письманика для определения неровноты пряжи, действующий экземпляр которого можно видеть в павильоне атомной энергии на Всесоюзной промышленной выставке в Москве. Другое средство для решения той же экспериментальной задачи – так называемая «установка Б», функционирующая на кафедре физики Московского текстильного института.

Цель предлагаемой статьи – показать, что и каким образом может дать система общих уравнений (1.1), (2.3) и (2.4), когда опыт даст сведения о пространственно-временном распределении или линейной плотности массы $\lambda(x,t)$, или линейной плотности количества движения $q(x,t)$, или, наконец, распределения скоростей $v(x,t)$ по сечениям поля в различные моменты времени.

**4.***Определение $q$, $v$ и $F$ по известной $\lambda(x,t)$ и секундному расходу массы $q(0,t)$ на входе в поле.* Если $\lambda(x,t)$ известна для всех абсцисс $x$ в поле $(0 \geq x \geq l)$ и для всех моментов времени $t$ $(0 \leq t)$, то, прежде всего, из уравнения неразрывности (1.1), учитывая (2.4), находим

$$q(x,t) = q(0,t) - \int_0^x \frac{\partial \lambda(\xi,t)}{\partial t} d\xi \qquad (4.1)$$

Тогда из $q = \lambda v$ получаем

$$v(x,t)=\frac{1}{\lambda(x,t)}(q(0,t)-\int_0^x\frac{\partial\lambda(\xi,t)}{\partial t}d\xi) \tag{4.2}$$

При $x=0$ и $x=l$ это дает

$$v(0,t)=\frac{q(0,t)}{\lambda(0,t)}, \qquad v(l,t)=\frac{1}{\lambda(l,t)}(q(0,t)-\int_0^l\frac{\partial\lambda(\xi,t)}{\partial t}d\xi) \tag{4.3}$$

Для окончательной вытяжки $B(l,t)$ в момент $t$ получаем

$$B(l,t)=\frac{v(l,t)}{v(0,t)}=\frac{\lambda(0,t)}{\lambda(l,t)}\left[1-\frac{1}{q(0,t)}\int_0^l\frac{\partial\lambda(\xi,t)}{\partial t}d\xi\right] \tag{4.4}$$

В обычно употребляемых вытяжных аппаратах с двумя парами валиков движения обеих пар управляются единым механизмом. Вследствие этого оба движения координированы так, что с точностью до малых случайных отклонений отношение $B(l,t)$ скоростей ленты на выходе и на входе не должно зависеть от времени. Вообще же вытяжка $B(l,t)$ является функцией времени, колеблющейся по законам случая около некоторого постоянного значения в данном режиме. Выражение (4.4) показывает, что, вообще говоря, лишь в стационарном процессе отношение скоростей на выходе и на входе равно обратному отношению линейных плотностей массы. Наконец, зная $\lambda(x,t)$, $q(x,t)$ и $v(x,t)$, найдем распределение движущих сил $F(x,t)$ по уравнению

$$F(x,t)=q(x,t)v(x,t)-q(0,t)v(0,t)+\int_0^x\frac{\partial q(\xi,t)}{\partial t}d\xi \tag{4.5}$$

и задача будет решена до конца. Отметим, что найденное решение в такой постановке вопроса оказывается единственным.

**5.** *Определение $\lambda$, $v$ и $F$ по известному секундному расходу массы $q(x,t)$ и распределению линейной плотности массы $\lambda(x,0)$ в начальный момент времени.*

Эта задача решается по способу, почти не отличающемуся от предыдущего; из (1.1), учитывая $q=\lambda v$, имеем

$$\lambda(x,t)=\lambda(x,0)-\int_0^t\frac{\partial q(x,\tau)}{\partial x}d\tau \tag{5.1}$$

Здесь $\lambda(x,0)$ играет роль постоянной интеграции по $t$, подобно тому, как в предыдущей задаче $q(0,t)$ играло роль постоянной интеграции по $x$. Зная $q(x,t)$ и $\lambda(x,t)$, найдем

$$v(x,t)=\frac{q(x,t)}{\lambda(x,t)} \qquad (5.2)$$

и, в частности, на входе и выходе из поля

$$v(0,t)=\frac{q(0,t)}{\lambda(0,t)}, \quad v(l,t)=\frac{q(l,t)}{\lambda(l,t)}=q(l,t)\left[\lambda(l,0)-\int_0^t \frac{\partial q(x,\tau)}{\partial x}\Big|_{x=l}d\tau\right]^{-1} \qquad (5.3)$$

Отсюда получаем для окончательной вытяжки $B(l,t)$ в момент времени $t$

$$B(l,t)=\frac{v(l,t)}{v(0,t)}=\frac{q(l,t)}{q(0,t)}\lambda(0,t)\left[\lambda(l,0)-\int_0^t \frac{\partial q}{\partial x}\Big|_{x=l}d\tau\right]^{-1} \qquad (5.4)$$

Опять видим, что лишь в стационарном режиме, когда $q=const$ для всех сечений поля, имеет место равенство

$$B=\frac{v(l)}{v(0)}=\frac{q\lambda(0)}{q\lambda(l)}=\frac{\lambda(0)}{\lambda(l)}$$

Наконец, по распределению расходов массы и скоростей находим распределение движущих сил по уравнению (4.5)

$$F(x,t)=q(x,t)v(x,t)-q(0,t)v(0,t)+\int_0^x \frac{\partial q(\xi,t)}{\partial t}d\xi \qquad (5.5)$$

Решение единственно и в этом случае.

**6.** *Определение* $\lambda(x,t)$ *по* $v(x,t)$. Теперь представим себе, что из опыта стало известно распределение скоростей по сечениям поля в различные моменты времени. Таким образом, будет известна функция $v(x,t)$ для области $0\le x\le l$ и $t\ge 0$.

Посмотрим, как в таком случае может быть определена для тех же значений обоих аргументов линейная плотность массы $\lambda(x,t)$.

Уравнение неразрывности может быть записано в виде

$$\frac{\partial q}{\partial x}+\frac{\partial \lambda}{\partial t}=0 \qquad (6.1)$$

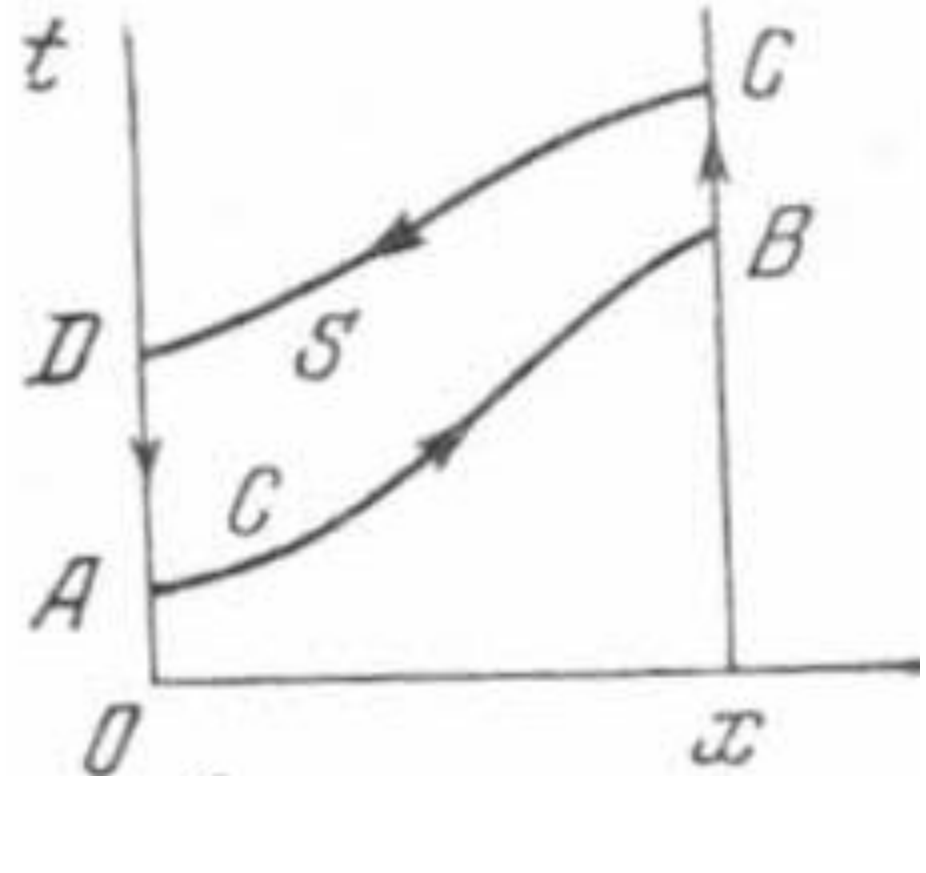

Контур $ABCD$ состоит из отрезков $BC$ и $DA$, прямых $x=0$, $x=const$ и двух кривых: $AB$ и $CD$, на каждой из которых переменные $x$ и $t$ связаны условием

$$dx=vdt \qquad (6.2)$$

Применив к $C$ и $S$ теорему Грина, имеем

$$\iint_S \left[\frac{\partial q}{\partial x}+\frac{\partial \lambda}{\partial t}\right]dxdt=\int_C qdt-\lambda dx=0=\int_C \lambda(vdt-dx) \qquad (6.3)$$

Докажем теперь, что в области $S$ (фиг.1) может существовать не больше одной функции $\lambda(x,t)$, удовлетворяющей уравнению (6.1) при данной $v(x,t)$ и принимающей на отрезке $DA$ (где $x=0$) границы $C$ данные значения $\varphi(t)$.

В самом деле, если бы существовали две такие функции: $\lambda_1$ и $\lambda_2$, то их разница $\lambda=\lambda_1-\lambda_2$ была бы отлична от тождественного нуля в $S$ и на $C$, удовлетворяла бы уравнению $F(\lambda)=0$ и обращалась бы в нуль на $DA$. Между тем, подставив эту разницу $\lambda$ в (6.3), мы имели бы

$$\left[\int_{AB}+\int_{BC}+\int_{CD}+\int_{DA}\right]\lambda(vdt-dx)=0 \qquad (6.4)$$

Но на $AB$ и $CD$ интегралы исчезают вследствие (6.2). Интеграл вдоль $DA$ обращается в нуль вследствие $\lambda=0$ на $DA$. Остается

$$\int_{BC}\lambda(vdt-dx)=0$$

или, так как на $BC$ имеем $x=const$,

$$\int_{BC}\lambda vdt=0 \qquad (6.5)$$

Функции $\lambda$ и $v$ не отрицательны, причем $v>0$. Поэтому непременно $\lambda\equiv 0$ на $BC$ и, следовательно, $\lambda_1\equiv\lambda_2$ всюду в полу-полосе $0\le x\le l, t\ge 0$ на основании произвола в выборе границы $ABCD$ области $S$.

Тем самым показано, что существует не больше одного распределения линейных плотностей массы $\lambda(x,t)$ в вытяжном поле, если это распределение задано распределением скоростей и условием

$$\lambda(x,t)\big|_{x=0}=\varphi(t)=\lambda(0,t) \qquad (6.6)$$

на входе в поле.

Если удастся найти решение уравнения (6.1) при условии (6.6) с заданной функцией времени, то, как и в рассмотренных ранее случаях, определим затем $q=\lambda v$ и, наконец,
$F(x,t)$. В такой постановке задачи найденное ее решение будет единственно возможным.

Итак, *ставим задачу.*

*Линейная плотность поступающего в вытяжной аппарат материала является заданной функцией времени* $t \geq 0$. *Распределение скоростей* $v(x,t)$ *по сечениям* $x$ *поля* $0 \leq x \leq l$ *и времени* $t \geq 0$ *известно. Найти распределение линейной плотности* $\lambda(x,t)$ *для* $0 \leq x \leq l$ *и для* $t \geq 0$.

Рассмотрим несколько более общее уравнение вместо (6.1)

$$F_\alpha(\lambda) \equiv v\frac{\partial \lambda}{\partial x} + \frac{\partial v}{\partial x}\lambda + \alpha\frac{\partial \lambda}{\partial t} = 0 \qquad (6.7)$$

где $\alpha$ – безразмерный параметр. При $\alpha = 1$ это уравнение совпадает с уравнением (6.1).

Теорема о единственности решения без труда распространяется и на обобщенное уравнение (6.7).

Будем искать решение для задачи

$$F_\alpha(\lambda) = 0\,,\; \lambda(x,t)\big|_{x=0} = \lambda(0,t) = \varphi(t) \qquad (6.8)$$

в виде ряда по степеням параметра $\alpha$

$$\lambda(x,t) = \lambda_0(x,t) + \alpha\lambda_1(x,t) + \alpha^2\lambda_2(x,t) + \ldots \qquad (6.9)$$

Составляя формально ряды для производных, подставляя в (6.7) и сравнивая множители при одинаковых степенях параметра $\alpha$, имеем последовательность рекуррентных уравнений

$$\frac{\partial(v\lambda_0)}{\partial x} = 0\,,\quad \frac{\partial(v\lambda_k)}{\partial x} = -\frac{\partial \lambda_{k-1}}{\partial t} \quad (k = 1,2,\ldots) \qquad (6.10)$$

Из этой системы определим сначала $\lambda_0(x,t)$, от которой потребуем, чтобы она обращалась в заданную функцию $\varphi(t)$ при $x = 0$. Затем шаг за шагом найдем все $\lambda_k(x,t)$, $k = 1,2,\ldots$, требуя, чтобы каждая из них обращалась в нуль при $x = 0$. Тогда ряд (6.9) формально даст решение задачи (6.8). Чтобы решение было не только формальным, но и фактическим, нужно иметь уверенность в равномерной и абсолютной сходимости его во всей полуполосе $0 \leq x \leq l\,, t \geq 0$. Если это действительно так при $\alpha = 1$, то ряд

$$\lambda(x,t) = \sum_{k=0}^{\infty}\frac{(-1)^k}{v(x,t)}\int_0 \frac{\partial}{\partial t}\ldots\frac{1}{v(x,t)}\int_0 \frac{\partial}{\partial t}\left[\varphi(t)\frac{v(0,t)}{v(x,t)}\right]\underbrace{dx\ldots dx}_{k\cdot раз} \qquad (6.11)$$

является решением задачи

$$\frac{\partial \lambda}{\partial t} + v\frac{\partial \lambda}{\partial x} + \lambda\frac{\partial v}{\partial x} = 0\,,\quad \lambda(x,t)\big|_{x=0} = \varphi(t)$$

Мы видели уже, что это решение будет единственно возможным. Впрочем, непосредственно из физических соображений ясно, что знание желаемого распределения скоростей и условия на входе еще недостаточно для того, чтобы процесс мог фактически осуществляться. Этому отвечает то обстоятельство, что многократно применяемый оператор

$$\frac{-1}{v}\int_0^x \frac{\partial}{\partial t}\ldots dx$$

в ряде, дающем $\lambda(x,t)$, не ограничен, и надлежащая сходимость ряда, вообще говоря, ничем не обеспечена.

В каждом конкретном случае сходимость ряда должна проверяться применительно к этому случаю.

ЛИТЕРАТУРА

## **11.**О СКОРОСТЯХ ВОЛОКОН В ПРОЦЕССЕ ВЫТЯГИВАНИЯ

В ранее сообщенной работе [2] были даны три динамических уравнения, которые должны совместно удовлетворяться четырьмя функциями от координаты $x$ и времени $t$ не зависимо от того, что вытягивается. Недостающее четвертое уравнение должно существенно отразить физическую природу вытягиваемого объекта. В названной работе[2] говорилось, что это уравнение можно было бы получить экспериментально, путем промеров секундного расхода массы и линейной плотности ленты в различных сечениях поля и в различные моменты времени. В предлагаемой теперь статье мы приводим недостающее уравнение. Обозначения сохраняем введенные раньше [1-2].

Пусть $v(0,t)$ есть скорость входа ленты в поле в момент $t$ $(t \geq 0)$, $v(l,t)$ – скорость ее выхода из поля в тот же момент, причем

$$v(l,t) \geq v(0,t) \geq v(0,\tau) > 0 \qquad (0 \leq \tau \leq t) \tag{1}$$

Называя $v(x,t)$ $(0 \leq x \leq l)$ среднюю (по совокупности волокон) скорость ленты в сечении $x$ в момент $t$, определяем вытяжку $B$ в этот момент в этом сечении отношением

$$B(x,t) = \frac{v(x,t)}{v(0,t)} \tag{2}$$

Пусть $n(x,t)$ – число волокон, пронизывающих сечения $x$ в момент $t$. Если плотность вероятности для распределения длин волокон в сечении $x$ в момент $t$ есть $\psi(x,t,s)$, то

$$n(x,t)\psi(x,t;s)ds \tag{3}$$

является числом тех волокон из $n(x,t)$, которые имеют длину, бесконечно мало отличающуюся от $s$.

Пронизывая сечение $x$, они по-разному расположены в момент $t$ относительно этого сечения. Обозначив $\alpha$ длину переднего конца волокна $(0 \leq \alpha \leq s)$ и вводя плотность вероятности для распределения длин передних концов волокон с длиной $s$ в сечении $x$ в момент $t$; $\omega(x,t,s;\alpha)$, получим для числа волокон из группы (4), имеющих при общей длине $s$ длину переднего конца, практически не отличающуюся от $\alpha$, выражение

$$n(x,t)\psi(x,t;s)ds\omega(x,t,s;\alpha)d\alpha \tag{4}$$

При этом для всех возможных $x$ и $t$ должно быть

$$\int_0^S \omega(x,t,s;\alpha)d\alpha = 1 \quad , \quad \int_0^S \psi(x,t;s)ds = 1, \tag{5}$$

где $S$ – наибольшая длина волокна.

Мы считаем, *во-первых*, что процесс, вообще нестационарный, протекает *без последействия*, т.е. распределение скоростей волокон зависит лишь от мгновенных значений, определяющих процесс обстоятельств, а не от тех, которые эти обстоятельства имели в предшествующие моменты времени. Такое допущение можно сделать потому, что каждое отдельное волокно практически лишено инерции и тотчас воспринимает тот режим, который ему в данный продиктован. *Во-вторых*, мы принимаем, что два любых волокна с одинаковой длиной, занимающие в данный момент времени один и тот же участок поля, находятся в одинаковых кинематических условиях, т.е. в точках с одинаковыми абсциссами имеют одинаковые мгновенные скорости. *В-третьих*, в отличие от принятой в ранней работе [1] гипотезы об отсутствии промежуточных скоростей у отдельных волокон, мы допускаем теперь наличие переходов от одной скорости к другой по ряду промежуточных значений. Наконец, *в четвертых*, опять-таки в отличие от [1], считаем, что проскальзывания волокон нет ни на входе, ни на выходе.

Волокна, в момент $t$ не покинувшие входного зажима, имеют в этот момент скорость входа $v(0,t)$ (контролируются входной парой валиков). Волокна, которые в момент $t$ уже зажаты между парой передних валиков, имеют скорость выхода $v(l,t)$. Прочие волокна той части ленты, которая в момент $t$ находится в вытяжном поле, имеют собственные скорости, которые являются предметом специальных исследований (работа В.А. Протасовой в радиоизотопной лаборатории МТИ и более ранние: Гинзбурга, Слатинцева, Спесивцевой и др.). Такие волокна находятся в момент $t$ в состоянии «плавания». Из числа $n(x,t)$ волокон, пронизывающих в момент $t$ сечение $x$, найдутся контролируемые задней парой, если только длина заднего конца $\beta = s - \alpha$ окажется больше $x$. Если*) найдутся такие волокна, для которых длина переднего конца $\alpha$ будет в момент $t$ больше, чем $l - x$, то они в этот момент будут контролироваться передней парой валиков. Наконец, те волокна из числа $n(x,t)$, для которых $s - x < \alpha < l - x$, окажутся «плавающими» в момент $t$. Впрочем, может случиться, что одно или даже два из этих трех множеств волокон окажется пустым.

Из сказанного следует, что, обозначая неизвестную пока скорость «плавания» одиночного волокна просто буквой $w$, для средней (по совокупности) скорости волокон в сечении $x$ в момент $t$ можем написать

$$v(x,t) = v(0,t) \iint\limits_{0<\alpha<s-x} \psi(x,t;s)\omega(x,t,s;\alpha)dsd\alpha + \tag{6}$$

$$+ \iint\limits_{s-x<\alpha<l-x} w\psi(x,t;s)\omega(x,t,s;\alpha)dsd\alpha + v(l,t) \iint\limits_{l-x<\alpha<S} \psi(x,t;s)\omega(x,t,s;\alpha)dsd\alpha$$

Здесь $w$ – некоторая неизвестная функция переменных $x, t, s$ и $\alpha$.

Может случиться, что один или даже два из интегралов правой части в выражении (6) отсутствуют. Уравнение (6) в совокупности с тремя уравнениями, упомянутыми в начале этой статьи[1], с исчерпывающей полнотой решало бы

**Примечание. Нарисованная картина существенно предполагает неизменяемость движущихся волокон по длине, в частности, их нерастяжимость; в дальнейшем мы освободимся от этого допущения.*

общую задачу о вытягивании ровничной ленты, если бы функции $\psi, \omega$ и $w$ были в точности известны. Между тем даже в стационарном случае мы не имеем в настоящее время достаточных данных по этому поводу. Окончательное слово будет сказано результатами хорошо поставленного опыта. Пока опытных данных совсем или почти совсем нет.

Мы ставим вопрос, есть ли надобность делать различие в кинематическом поведении плавающих волокон и волокон, контролируемых на входе или на выходе из прибора. Рассмотрим сначала волокно, плавающее в момент $t$. Пусть оно пронизывает в этот момент сечение $x$, имеет длину $s$ и длину переднего конца $\alpha$. Если разводка (длина поля) есть $l$, то условие плавания обозначат, что $s \le x + \alpha \le l$.

Назовем $w$ мгновенную скорость точки волокна, находящейся в момент $t$ точно в сечении $x$. По смыслу процесса

1) скорость $w$ тем больше отличается от мгновенной скорости $v(0,t)$ на входе, чем больше расстояние $x+\alpha-s$ заднего конца плавающего волокна от входа;
2) $w$ тем меньше отличается от мгновенной скорости $v(l,t)$ на выходе из поля, чем меньше расстояние $l-x-\alpha$ переднего конца волокна от выхода;
3) из двух упомянутых в 1) и 2) разностей абсцисс $x+\alpha-s$ и $l-x-\alpha$ вторая влияет на $v(l,t)-w$ в большей степени, чем первая на $w-v(0,t)$;
4) Высказанное в 1) , 2) и 3) остается в силе, каковы бы ни были возможные значения $x$ и $t$.

Все сказанное можно изобразить уравнением

$$\frac{[v(l,t)-w]^{\mu}}{w-v(0,t)} = [v(l,t)-v(0,t)]^{\mu-1}\left[\frac{l-(x+\alpha)}{(x+\alpha)-s}\right]^{\nu} \qquad (7)$$

где $\mu$ и $\nu$ – положительные постоянные, причем $\mu > 1$. Не зависящий от $x$ множитель

$$[v(l,t)-v(0,t)]^{\mu-1}$$

уравнивает размерности справа и слева. Условие отсутствия проскальзывания в зажимах соблюдено, так как из (7) следует:

$$w = v(l,t) \quad (x+\alpha = l), \qquad w = v(0,t) \qquad (x+\alpha = s) \qquad (8)$$

т.е. когда либо передний, либо соответственно задний конец плавающего волокна приходится на плоскость либо переднего, либо заднего зажима.

Кроме того, по условию плавания при $x = 0$ имеем $\alpha \ge s$, значит

$$\alpha\big|_{x=0} = s \qquad (9)$$

При $x = l$ имеем $\alpha \le 0$, так что

$$\alpha|_{x=l} = 0 \qquad (10)$$

Таким образом, *на входе в поле* ($x = 0$) *распределение длин* $\alpha$ *передних концов плавающих волокон совпадает с распределением их полных длин.*

*На выходе из поля* $x = l$ *длины задних концов* $s - \alpha$ плавающих волокон распределены так же, *как их полные длины*.

Далее, из (7) следует, что в любой данный момент $t$ все частицы данного плавающего волокна имеют одну и ту же мгновенную скорость $w$, которую можно поэтому называть мгновенной скоростью данного плавающего волокна. В самом деле, рассмотрим плавающее волокно, занимающее положение между абсциссами $x + \alpha$ и $x + \alpha - s$ в момент $t$. Скорость $w$ его точки, лежащей в абсциссе $x$, определяется из (7).

Рассмотрим другую его точку $y$ в том же положении волокна. Для сечения $y$ длина $\alpha'$ переднего конца рассматриваемого волокна в сумме с $y$ дает абсциссу переднего конца, которая есть $x + \alpha$. Итак

$$y + \alpha' = x + \alpha$$

От подстановки (10) в уравнение (7) мы получим то же значение мгновенной скорости $w$, что и раньше. Между тем, по смыслу уравнения (7), эта $w$ была бы мгновенной скоростью частицы волокна, находящейся в сечении $y$.

Из всего сказанного следует, что всякое плавающее волокно, либо начинающее процесс плавания, либо заканчивающее его, ведет себя как контролируемое либо на входе, либо соответственно на выходе. Говоря иначе, контролирование на входе в кинематическом смысле непрерывно переходит в плавание, которое, в свою очередь, непрерывно переходит в контролирование на выходе.

Если приписывать всякому волокну длину, равную длине той его части, которая в данный момент находится внутри вытяжного поля, то справедливым окажется и обратное суждение: каждое волокно, контролируемое в данный момент на входе в поле, ведет себя как плавающее со скоростью $v(0,t)$ и имеющее длину $s = \alpha$, растущую в данный момент пропорционально $v(0,t)$. Аналогично, каждое волокно, контролируемое в данный момент на выходе из поля, тоже ведет себя как плавающее, но со скоростью $v(l,t)$ и имеющее длину $s - \alpha$ в этот момент и уменьшающее ее пропорционально $v(l,t)$.

Однако такая *концепция* требует радикального пересмотра законов распределения как длин $s$ волокон, так и длин их передних концов $\alpha$. Мы полагаем, что связанные с этим вопросом трудности можно обойти, пользуясь следующими соображениями:

1) Распределение средних скоростей волокон по сечениям поля и во времени описывать, исходя из уравнения (7).

2) Отказаться от понятий «длины волокна» и «длины переднего конца» его при расчетах средних скоростей в вытяжном поле.
3) Отказаться при изучении этого вопроса и от понятия «контролирования».
4) Ленту в процессе ее обработки в вытяжном поле рассматривать, как агрегат материальных точек (частиц), размеры волокон считать, следовательно, равными нулю по сравнению с длиной поля (разводкой).
5) Эти частицы считать связанными в направлении оси поля в бесконечно протяженные в оба конца цепочки. Связывающие силы могут быть бесконечно разнообразны по величине, и их распределение вдоль поля и во времени подчинено законам счисления зависимых случайных величин.
6) Все частицы, упомянутые в 4), только «плавают», характер этого движения описан уравнением (7), в котором надо положить $\alpha = s = 0$.
7) Специфику вытягиваемого объекта – ленты – характеризовать не путем попыток установления законов распределения длин волокон и длин их передних концов, а надлежащим подбором двух отвлеченных положительных постоянных $\mu$ и $\nu$, о которых будет идти речь дальше.
8) Полученные из тщательно обставленного опыта (см. выше) сведения о распределении скоростей волокон в вытяжном поле стараться уложить в уравнения типа (7) с выбором $\mu$ и $\nu$,что, по-видимому, можно сделать и что сделать проще, чем определять для каждого отдельного случая законы распределения $s$ и $\alpha$, фигурирующие в общем выражении (6) для $v(x,t)$.

Принимая во внимание сказанное и полагая $\alpha \equiv 0$, $s \equiv 0$ в уравнении (7) имеем

$$\frac{[v(l,t)-w]^{\mu}}{w-v(0,t)} = [v(l,t)-v(0,t)]^{\mu-1}\left[\frac{l-x}{x}\right]^{\nu} \qquad (11)$$

Отсюда можно, говоря в принципе, найти мгновенную скорость $w$ для всех тех частиц ленты, которые в момент $t$ находятся в сечении поля с абсциссой $x$.

Значение $w$ отождествляем с искомой средней скоростью ленты в момент $t$ в сечении $x$. Законность этого очевидна: при $\alpha \equiv 0$, $S = 0$ обращается в $v(x,t) = w$. В частности, для стационарного процесса без проскальзывания в зажимах вместо (11) имеем

$$\frac{[v_1-w]^{\mu}}{w-v_0} = v_0^{B-1}[B-1]^{\mu-1}\left[\frac{l-x}{x}\right]^{\nu}$$

где $v_0$ и $v_1 > v_0$ – постоянные скорости на входе и на выходе, $B$ – вытяжка, $w$– искомая скорость в абсциссе $x$, $0 \le x \le l$. Дифференцируя равенство (7), находим

$$\frac{\partial w}{\partial x} = \frac{\nu l}{x(l-x)}\,\frac{[w-v(0,t)][v(l,t)-w]}{[(\mu-1)w+\{B(t)-\mu\}v(0,t)]} \qquad (12)$$

Заметим, что при $\mu > 1$ для всех возможных $B(t)$

$$0 < (\mu - 1)w + [B(t) - \mu]v(0,t)$$

Следовательно,

$$\frac{\partial w}{\partial x} > 0 \qquad (t \geq 0, 0 < x < l)$$

и $w(x,t)$ в каждый момент времени монотонно возрастает внутри поля.

При $x = 0$, как и при $x = l$, производная (12) приобретает неопределенный вид. Однако, если $v > \mu$, то из(11) следует, что при $x \to l$ разность $v(l,t) - w$ является малой более высокого порядка, чем $l - x$. Также при $v > \mu$ разность $w - v(0,t)$ стремится обратиться в ноль быстрее, чем $x$. Это значит, что при $v > \mu$ касательная к кривой скоростей в каждый момент времени параллельна оси $x$ как на входе, так и на выходе. Будучи выпуклой вниз у входа в поле, кривая во всякий момент становится выпуклой вверх у выхода.

Можно несколько *обобщить предыдущее*, достигая этим и несколько большей простоты: мы допустим, что $\mu > 1$*есть не постоянная, а функция времени*, равная вытяжке. Кроме того, пусть и $v$ *есть функция времени,* равная

$$v = \kappa B(t),$$

где $\kappa > 1$ не зависит ни от $x$, ни от $t$, причем $[\kappa] = 1$. Исходное уравнение (11) перепишется в виде

$$\frac{[v(l,t) - w]^{B(t)}}{w - v(0,t)} = [v(l,t) - v(0,t)]^{B(t)-1}\left[\frac{l-x}{x}\right]^{\kappa B(t)} \qquad (13)$$

а вместо (12) будем иметь тогда

$$\frac{\partial w}{\partial x} = \frac{\kappa B(t) l}{[B(t) - 1]x(l-x)} \frac{[w - v(0,t)][v(l,t) - w]}{w} \qquad (14)$$

В обоих последних соотношениях фигурирует теперь, кроме мгновенной вытяжки $B(t)$, единственная абсолютная отвлеченная постоянная $\kappa > 1$, которую надлежит рассматривать как функцию свойств волокнистого вещества ленты. Проектируемые нами радиоизотопные промеры будут иметь целью выбрать $\kappa$ так, чтобы под опытные данные можно было подвести зависимость (13). Из (14) можно получить соответствующее выражение для случая стационарного процесса. Не противоречит ли оно тому, что мы писали раньше [1]?

В названной работе [1] мы рассматривали такой стационарный процесс, при котором не исключалось проскальзывание волокон в зажимах аппарата. Собственно говоря, вытяжное поле при таких условиях не было ограничено отрезком $(0,l)$ оси $x$, а простиралось над всей осью $x$ от $-\infty$ до $+\infty$. В абсциссах $x=0$ и $x=l$ скорости $v(0)$ и $v(l)$ должны были отличаться от «наблюдаемых» $v_0$ и соответственно $v_1$, которые играли роль «достигаемых в бесконечности»:

$$v_0 = v(-\infty), \quad v_1 = v(+\infty)$$

В случае, рассматриваемом теперь для стационарного процесса, проскальзывание волокон в зажимах исключено – средние скорости на концах поля $v(0)$ и $v(l)$ в точности совпадают с $v_0$ и соответственно с $v_1$. Заменив в (14) обозначения $x$ на $x'$ и переходя затем от координаты $x'$ к новой переменной, названной опять $x$ по уравнению

$$x = -l\ln\frac{l-x'}{x'}, \qquad \frac{ldx'}{x'(l-x')} = \frac{dx}{l}$$

мы *отображаем конечный отрезок* $(0,l)$ оси $x'$ *на всю бесконечно простирающуюся в оба конца ось* $x$. При $x'=0$ имеем $x=-\infty$, при $x'=l$ получаем $x=+\infty$.

После подстановки вместо (14) получим

$$\frac{dw}{dx} = \frac{2\kappa B}{(B-1)l}\frac{(w-v_0)(v_1-w)}{2w}$$

т.е. то самое уравнение процесса, которое мы получили (из других соображений) в цитированной статье. Между прочим обнаруживается и смысл введенной там «динамической длины» волокна

$$S = \frac{B-1}{2\kappa B}l$$

Из сказанного следует, что вопрос о *наличии или отсутствии промежуточных скоростей* у волокон в стационарном процессе не столь существенен, как думают обычно.

Уравнение (13) рекомендуем в качестве недостававшего четвертого в дополнение к тем трем, которые были нами предложены в предыдущем сообщении [2].

Из четырех независимых уравнений:

$$q = \lambda v, \quad \frac{\partial q}{\partial x} = -\frac{\partial \lambda}{\partial t}, \quad \frac{\partial (F-qv)}{\partial x} = -\frac{\partial q}{\partial t}.$$

$$\frac{[v(l,t)-v]^{B(t)}}{v-v(0,t)} = [v(l,t)-v(0,t)]^{B(t)-1}\left[\frac{l-x}{x}\right]^{\kappa B(t)}$$

мы имеем возможность (правда, лишь в принципе) определить все четыре неизвестные функции $v, \lambda, q$ и $F$ независимых переменных $x$ и $t$ с надлежащими краевыми и начальными условиями, как это намечено было в цитированной статье [2].

## **12.**КВАЗИСТАЦИОНАРНЫЙ ПРОЦЕСС РАБОТЫ ВЫТЯЖНОГО ПРИБОРА

В предыдущей статье мы дали систему уравнений:

$$q=\lambda v\,,\quad \frac{\partial q}{\partial x}=-\frac{\partial \lambda}{\partial t}\,,\quad \frac{\partial (F-qv)}{\partial x}=-\frac{\partial q}{\partial t}\,.$$

$$\frac{[v(l,t)-v]^{B(t)}}{v-v(0,t)}=[v(l,t)-v(0,t)]^{B(t)-1}\left[\frac{l-x}{x}\right]^{\kappa B(t)} \qquad (1)$$

где $q-$ секундный расход массы ленты, $\lambda-$линейная плотность, $v-$скорость, $F$ – сила, $B(t)-$ полная вытяжка в момент $t$, равная

$$B(t)=\frac{v(l,t)}{v(0,t)} \qquad (2)$$

положительная постоянная $\kappa$ подбирается из опыта и является динамической характеристикой вещества.

Система (1) определяет все четыре неизвестные функции $v,\lambda,q$ и $F$ независимых переменных, как бы не менялись с течением времени скорости $v(0,t)$ и $v(l,t)$ на входе и на выходе. Но конструкция прибора такова, что за исключением случайных малых отклонений полная вытяжка остается постоянной для любого $t$, потому что скорость $v(l,t)$ на выходе во всякий момент пропорциональна скорости $v(0,t)$ на входе. Имея это в виду, полагаем

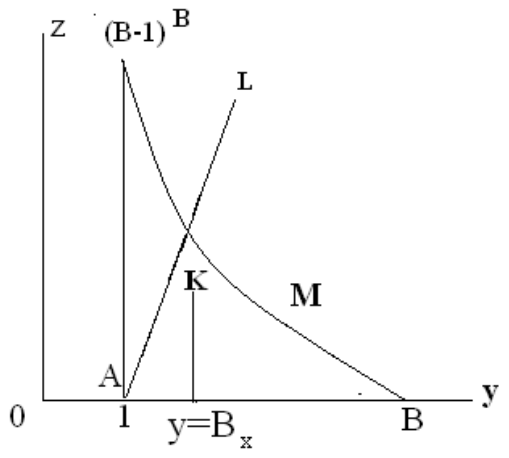

$$v(0,t)=\sigma(t)v(0)\ ,\ v(l,t)=\sigma(t)v(l) \qquad (3)$$

где $\sigma(t)$ – данная функция $t$, причем $\sigma(t)>0$ и $\sigma(0)=1$. С другой стороны, *отсутствие последействия в волокнах* позволяет считать, что и вообще

$$v(x,t)=\sigma(t)v(x)\ , \qquad (4)$$

так что переменные $x,t$ в выражении для скорости разделены естественными условиями процесса.

Вследствие однородности последнее уравнение (1) переписывается

$$\frac{(B-B_x)^B}{B_x-1}=(B-1)^{B-1}\left(\frac{l-x}{x}\right)^{\kappa B} \qquad (5)$$

Для любого $B\geq 1$ и $x\ (0\leq x\leq l)$ графическим приемом можно определить $B_x$. Положим

$$(B-1)^{B-1}\left(\frac{l-x}{x}\right)^{\kappa B}=C\,,\quad B_x=y\,,\ (B-y)^B=z_1,\ G(y-1)=z_2$$

Построим кривые $M(z_1,y),L(z_2,y)$. На отрезке $(1,B)$ они имеют вид, схематически показанный на фигуре. Парабола порядка $B$ зависит только от $B$.

С ростом $x$ прямая $L$ поворачивается от вертикального положения для $x=0$ к горизонтальному при $x=l$ по часовой стрелке. Точка пересечения $K$ обеих кривых имеет абсциссу $B_x$. Итак, для всякого $t \geq 0$ и $x$ $(0 \leq x \leq l)$ можно считать известной

$$v(x,t)=\sigma(t)v(x)=v(0,t)B_x \qquad (6)$$

Легко понять, что при сделанных оговорках секундный расход массы, один и тот же во всех сечениях в данный момент, может быть функцией только времени $t$. Из второго уравнения (1) видим, что тогда $\lambda$ есть функция только $x$. На основании первого уравнения (1) это позволяет заключить, что

$$q=q_0\sigma(t),\ \lambda=\frac{q}{v(x,t)}=\frac{q_0}{v(x)} \qquad (7)$$

где $q_0$ – секундарный расход массы при $t=0$. Тогда

$$\frac{\partial q}{\partial t}=q_0\sigma'(t)$$

и третье уравнение (1) сводится к

$$F-qv=q_0\sigma'(t)x+A(t),$$

где $A(t)$ – произвольная функция времени. Но при $x=0$ имеем

$$F(0,t)=q_0\sigma(t)v(0,t)+A(t)$$

Отсюда

$$F(x,t)-F(0,t)=q(t)\left[v(x,t)-v(0,t)\right]+\frac{\partial q}{\partial t}x \qquad (8)$$

и вместе с тем решение задачи доведено до конца.

*Остановимся еще на одном вопросе.* Первые три уравнения написанной выше системы (1) верны для процесса вытягивания любого линейно-протяженного тела. Вводя для краткости обозначение

$$qv-F=\Phi \qquad (9)$$

и дифференцируя второе из уравнений (1) по $t$, а третье – по $x$, получим

$$\frac{\partial^2\Phi}{\partial x^2}=\frac{\partial^2\lambda}{\partial t^2}. \qquad (10)$$

Это указывает на специфический характер процесса. Например, если бы процесс был таков, что каждая из двух функций $\Phi$ и $\lambda$ зависела бы от одной и той же величины $u$, в свою очередь являющейся функцией только двух независимых переменных $x$ и $t$, так что

$$\Phi=\Phi(u),\quad \lambda=\lambda(u) \qquad (11)$$

то (10) приняло бы вид хорошо изученного волнового одномерного уравнения

$$\frac{\partial}{\partial x}\left[\Phi'(u)\frac{\partial u}{\partial x}\right]=\frac{\partial}{\partial t}\left[\lambda'(u)\frac{\partial u}{\partial t}\right], \qquad (12)$$

решаемого известными методами.

СОДЕРЖАНИЕ

9 785001 220886

**Новоженова Ольга Георгиевна**

# БИОГРАФИЯ И НАУЧНЫЕ ТРУДЫ АЛЕКСЕЯ НИКИФОРОВИЧА ГЕРАСИМОВА. О ЛИНЕЙНЫХ ОПЕРАТОРАХ, УПРУГО-ВЯЗКОСТИ, ЭЛЕВТЕРОЗЕ И ДРОБНЫХ ПРОИЗВОДНЫХ

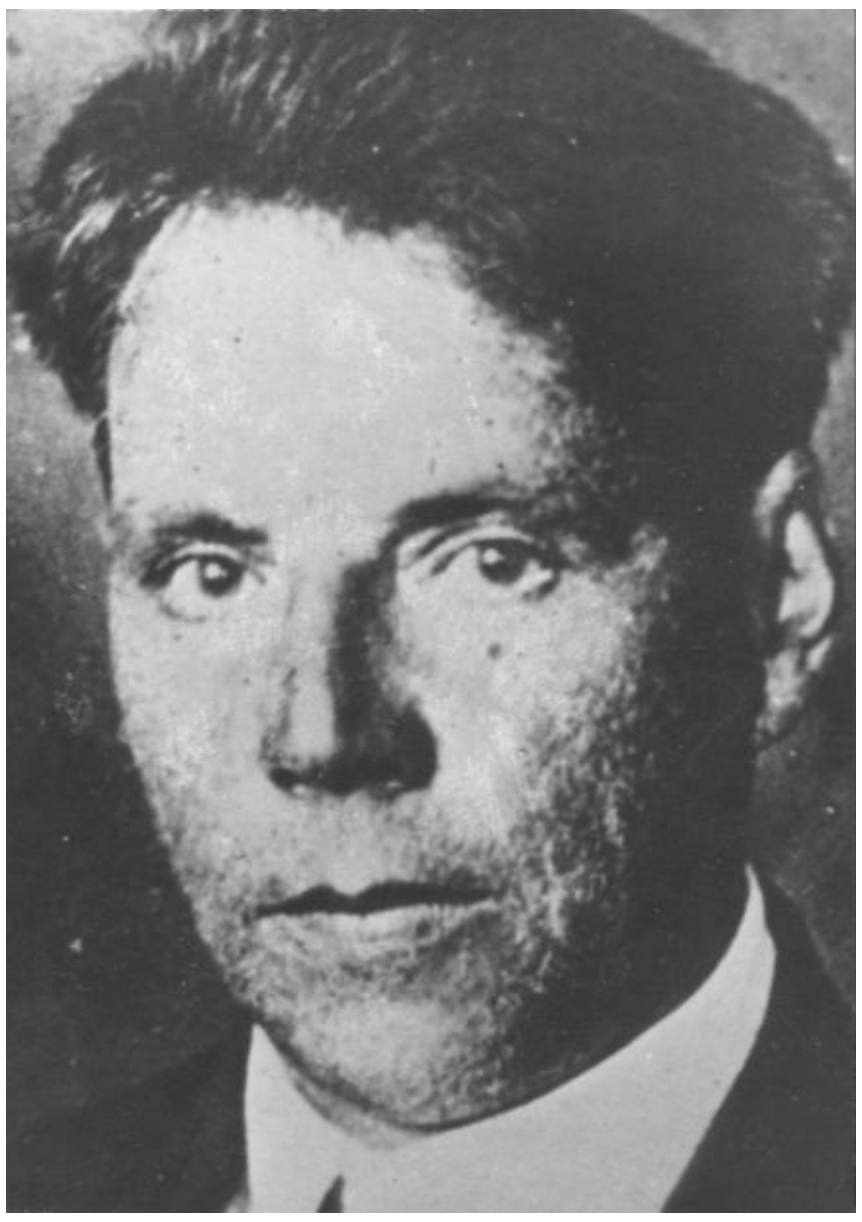

**АЛЕКСЕЙ НИКИФОРОВИЧ ГЕРАСИМОВ**
**(24.03.1897–14.03.1968)**

окончил в 1936 г. экстерном физико-математический факультет I-го Московского Государственного Университета с дипломом I-й степени, защитив под руководством проф. И.И. Привалова диплом на тему:
«Принцип соответствия в теории линейных операторов».
В 1943 г. также экстерном написал и защитил под руководством проф. В.В.Степанова кандидатскую диссертацию на тему: «Некоторые задачи теории упругости с учетом последействия и релаксации по линейному закону». Писал докторскую диссертацию, из которой была опубликована только одна статья «Обобщение линейных законов деформирования и его применение к задачам внутреннего трения» в журнале «Прикладная математика и механика»,1948, №3, 251-260.
В этой статье А.Н.Герасимов впервые в мировой литературе использовал дробную производную для решения задач вязко – упругости.

***Килбас А.А.***

На наш взгляд, это не совсем верно. Правильнее их называть дробными производными Герасимова-Капуто, так как в 1948 году советский механик А.Н. Герасимов ввел частную производную вида (1.9′) относительно $t$ на всей оси:

$$({}^{C}D^{\alpha}_{-,t}u)(x,t) = \frac{1}{\Gamma(\alpha)}\int_{-\infty}^{t}\frac{u_y(x,y)dy}{(t-y)^{\alpha}} \qquad (1.9''),$$

$$(t>0,\ x\in\mathbf{R};\ 0<\alpha<1).$$

в своей работе

*Теория и приложения дифференциальных уравнений дробного порядка (Курс лекций). Методологическая школа-конференция «Математическая физика и нанотехнологии», посвященная 40-летию возрождения Самарского государственного университета (Самара, 4-9 октября 2009 г.)*

**Kiryakova V.**

We mark now **50 years of the paper**:
M. Caputo, Linear model of dissipation whose Q is almost frequency independent-II. *Geophysical J. Royal Astronomic Soc.* **13** (1967), 529–539 (Reproduced in: *Fract. Calc. Appl. Anal.* **11**, No 1 (2008), 3–14),
in which the author **introduced a variant of the fractional derivative**, called usually in literature after his name, as **Caputo derivative**.

• It is **120th anniversary of the birth of the Russian (Soviet) mechanician Alexey N. Gerasimov (born 24 March 1897)**, a pioneer of using Fractional Calculus in solid mechanics, and in differential equations with fractional order partial derivatives. See the historical survey by O. Novozhenova, dedicated to the anniversary, in this journal's issue, 790–809. In a report on May 29, 1947 at the Institute of Mechanics of the USSR Academy of Sciences, **70 years ago, Gerasimov introduced the same kind of fractional derivative**. On the base of this report, the paper:
A.N. Gerasimov, Generalization of laws of the linear deformation and their application to problems of the internal friction (In Russian). *Prikladnaya Matematika i Mekhanika* **12**, No 3 (1948), 251–260,
was published, **but remained not popular** for long time.

**EDITORIAL**
**FCCA RELATED NEWS, EVENTS AND BOOKS**
**(FCAA–VOLUME 20–3–2017)**/*Fract.Calc. and Appl. Anal., 2017, vol.20,N 3, p.567*